

\input amstex
\documentstyle{amams} 

\document

\catcode`\@=11
\font\twelvemsb=msbm10 scaled 1100

\font\ninemsb=msbm10 scaled 800
\newfam\msbfam
\textfont\msbfam=\twelvemsb  \scriptfont\msbfam=\ninemsb
  \scriptscriptfont\msbfam=\ninemsb
\def\msb@{\hexnumber@\msbfam}
\def\Bbb{\relax\ifmmode\let\next\Bbb@\else
 \def\next{\errmessage{Use \string\Bbb\space only in math
mode}}\fi\next}
\def\Bbb@#1{{\Bbb@@{#1}}}
\def\Bbb@@#1{\fam\msbfam#1}
\catcode`\@=12

 \catcode`\@=11
\font\twelveeuf=eufm10 scaled 1100
\font\teneuf=eufm10
\font\nineeuf=eufm7 scaled 1100
\newfam\euffam
\textfont\euffam=\twelveeuf  \scriptfont\euffam=\teneuf
  \scriptscriptfont\euffam=\nineeuf
\def\euf@{\hexnumber@\euffam}
\def\frak{\relax\ifmmode\let\next\frak@\else
 \def\next{\errmessage{Use \string\frak\space only in math
mode}}\fi\next}
\def\frak@#1{{\frak@@{#1}}}
\def\frak@@#1{\fam\euffam#1}
\catcode`\@=12



\annalsline{154}{2001}

\received{November 24, 1999}

\revised{May 16, 2000}

\startingpage{1}


\title{Deviation of ergodic averages\\ for area-preserving flows\\
       on surfaces of higher genus} 

\shorttitle{Deviation of ergodic averages} 


\acknowledgements{ This paper rests on the work of several mathematicians, H. Masur,
J. Smillie, W. Veech and A. Zorich among them, and it was strongly inspired by the work of A. Zorich
and M. Kontsevich-A.~Zorich. I wish to thank particularly A. Eskin, J. Smillie and
A. Zorich for their
interest during the slow progress of this work and for discussing with me their ideas on the subject.
I am grateful to R.~Gunning for his encouragement and his enthusiasm
during several discussions on the content of
Section 4, which improved very much my understanding of moduli spaces.
I am also very grateful to J. Mather for his patience in listening to a number of sometimes very
confused and tentative presentations of these results and for his help in clarifying my ideas with
numerous questions and suggestions. This research was supported by NSF grant DMS-9704791. }


\author{Giovanni Forni}

\institutions{Princeton University, Princeton,  NJ\\
{\eightpoint {\it E-mail address\/}: gforni@math.princeton.edu}\\
{\eightpoint {\it Current address\/}}:   Northwestern University\\
{\eightpoint {\it E-mail address\/}: gforni@math.northwestern.edu}}


\centerline{\bf Table of Contents}
\def\sni#1{\smallbreak\noindent{#1}. }
\def\ssni#1{\vglue-1pt\noindent\hskip18pt {#1}.}
\medbreak
\noindent {Introduction}
\sni{1} {The Kontsevich-Zorich cocycle}
\sni{2} {Variational formulas}
\sni{3} {A lower bound for the second Lyapunov exponent}
\sni{4} {The determinant locus}
\sni{5} {The Kontsevich-Zorich formula revisited and other formulas for the\hfill\break 
\phantom{aaa} Lyapunov exponents}
\sni{6} {Basic currents for measured foliations}
\ssni{6.1} {Basic currents and invariant distributions}
\ssni{6.2} {Weighted Sobolev spaces of currents }
\ssni{6.3} {A geometric estimate of the Poincar\'e constant}

\sni{7} {The structure of the space of basic currents of finite order}
\ssni{\hskip1pt 7.1} {Basic currents with non-vanishing cohomology class}
\ssni{\hskip1pt 7.2} {Basic currents with vanishing cohomology class}
\sni{8} {The non-uniform hyperbolicity of the Kontsevich-Zorich cocycle and an\hfill\break
\phantom{aaa}\ application to currents}
\ssni{\hskip3pt 8.1} {The stable and unstable sub-bundles of the Kontsevich-Zorich cocycle
\hglue.62in   as
bundles of basic currents}
\ssni{\hskip3pt 8.2} {The Kontsevich-Zorich cocycle is non-uniformly hyperbolic}
\ssni{\hskip3pt 8.3} {The Oseledec's theorem for the bundle of closed currents of order $1$}
\sni{9} {The deviation of ergodic averages and open questions}
\ssni{\hskip1pt 9.1} {The $L^2$ mean deviation and the cohomological equation}
\ssni{\hskip1pt 9.2} {Sobolev estimates for the (first) return orbits}
\ssni{\hskip1pt 9.3} {Special sequences of `close' return times}
\ssni{\hskip1pt 9.4} {The main theorem on the deviation of ergodic averages}
\ssni{\hskip1pt 9.5} {The deviation for area-preserving vector fields and open questions}

\bigbreak
\intro

Quantities of physical interest in a mechanical system are often described 
as time averages along the trajectories of the system. Since such averages are in general difficult 
to compute, a fundamental idea which goes back to Boltzmann is to replace them with averages 
over the phase space. In the mathematical formulation, a flow $(M,\Phi_X^{\tau},\mu)$ on a manifold
$M$, generated by a smooth vector field $X$, is {\it ergodic }if the time averages of any integrable
function converge with probability one to the average with respect to the invariant probability 
measure $\mu$ as `time' goes to infinity. It is therefore a relevant problem to establish 
quantitative estimates on the speed of convergence and on the asymptotic behaviour of time averages 
of functions in ergodic systems. However, a basic result in ergodic theory [50, \S 3.2 B, Th.\ 2.3] 
states that for a general integrable function the speed of convergence can be arbitrarily slow. It 
is then reasonable to restrict the class of functions under consideration to functions with a certain
degree of smoothness.

 For many systems with strongly chaotic, {\it hyperbolic }behaviour 
the
{\it Central Limit Theorem }(CLT) is known to hold. Examples include 
smooth 
uniformly hyperbolic systems, like geodesic flows on compact manifolds 
of negative curvature (Sinai, Ratner), discontinuous uniformly
hyperbolic systems, like billiards with convex scatterers (Sinai's billiards)
and non-uniformly hyperbolic systems, like Bunimovich's billiards (the reader
can consult the survey by M. Denker \ref\Dk or the more recent papers by 
L.-S. Young \ref\Ygtwo, \ref\Ygthree). The CLT states that, for any H\"older 
continuous function $f$ on $M$,
$${1\over{ \sqrt{\Cal T}}}\,\left\{\,\int_0^{\Cal T} f\bigl(\Phi_X^{\tau}(p)\bigr)
\,d\tau\,\,-\,\, {\Cal T}\int_M f\,d\mu\,\right\} \,\,\to \,\, {\Cal N}(0,\sigma_f)
\,\,, \eqnu $$ 
in the sense of probability distributions, where ${\Cal N}(0,\sigma_f)$ is the 
normal (Gaussian) distribution  with zero mean and variance $\sigma_f^2$. A 
related result, also known for some examples of hyperbolic systems \ref\Dk, 
is the so-called {\it Law of Iterated Logarithms }(LIL). The LIL provides the 
following almost everywhere upper bound on the {\it deviation of ergodic 
averages:}
$$ \int_0^{\Cal T} f\bigl(\Phi_X^{\tau}(p)\bigr)
\,d\tau\,\,-\,\, {\Cal T}\int_M f\,d\mu\,=\,{\Cal O}\bigl(\sqrt{ {\Cal T} 
\log\log {\Cal T}}\bigr)\,\,. \speqnu {0.1'}$$
If CLT and LIL hold, then the stochastic process $\{f \circ \Phi^{\tau}_X\}_{\tau\in
{\Bbb R}}$ on the probability space $(M,\mu)$ behaves to some extent as if it were 
truely random (i.e. such as outcomes from flipping a coin or the Brownian motion). 

In the non-hyperbolic case less is known. If the system is {\it uniquely ergodic}, the
convergence of ergodic averages is uniform, with respect to $p\in M$, for all continuous functions. 
It is therefore natural to ask for {\it pointwise }estimates on the deviation of ergodic averages. 
A fundamental and well understood example, for its relevance in Hamiltonian mechanics, is given 
by regular quasi-periodic motions on tori, which have zero topological entropy, are never mixing but 
generically are uniquely ergodic. In particular, for a generic smooth aperiodic conservative flow on 
the $2$-dimensional torus $T^2$, the Denjoy-Koksma inequality [25, Chap.\ VI, \S 3] and the 
metric (measure) theory of continued fractions [36, Chap.\ III, Ths. 30 \& 31] imply that, if $f$ 
is a function of bounded variation, then for all $p\in T^2$ and any $\alpha>0$,
$$\left|\,\int_0^{\Cal T} f\bigl(\Phi_X^{\tau}(p)\bigr)\,d\tau\,\,-\,\, 
{\Cal T}\int_M f\,d\mu\,\right|\leq C_{\alpha} \,\hbox{Var}(f) \,\log{\Cal T}
\bigl(\log\log {\Cal T}\bigr)^{1+\alpha}\,\,. \eqnu$$
In addition, if $f$ is sufficiently smooth, a Fourier series argument shows that  
$$\left|\,\int_0^{\Cal T} f\bigl(\Phi_X^{\tau}(p)\bigr)\,d\tau\,\,-\,\, 
{\Cal T}\int_M f\,d\mu\,\right|\leq C_f\,\,, \speqnu{0.2'}$$
as a consequence of the existence of smooth solutions of the cohomological equation $Xu=
f-\int_M f\,d\mu$. Conservative aperiodic flows on tori are examples of {\it elliptic }dynamics,
characterized by very slow or no divergence of nearby orbits. Hyperbolic systems display 
exponential divergence of orbits, which is the fundamental source of their strongly chaotic 
properties.

The above results leave open all intermediate cases. A lemma proved by M. Ratner [51, 
Lemma 3.1], based on the Chebychev inequality, states that polynomial decay of correlations implies
polynomial speed of convergence of ergodic averages. According to Ratner, if a function $f$ has
polynomial decay of correlation with exponent $0<\alpha_f<1$, then for any $\varepsilon>0$ there are
$P:=P_{\varepsilon}\subset M$ with $\mu(P)\geq 1-\varepsilon$ and ${\Cal T}_P>0$ such that, if $p\in P$, 
$$
\left|\,\int_0^{\Cal T} f\bigl(\Phi_X^{\tau}(p)\bigr) \,d\tau\,\,-\,\, {\Cal T}\int_M f\,d\mu\,\right| \leq 
{\Cal T}^{1-\alpha'}\,\,, \eqnu $$
for all ${\Cal T}\geq {\Cal T}_P$, where $\alpha':=\alpha'_f=\alpha_f/8$. Ratner's result was proved with
the horocycle flow on a compact surface of constant negative 
curvature in mind  and for which the required
polynomial decay of correlations was known \ref\Rttwo. In the horocycle case
Ratner's result was improved upon by M. Burger. He proved [7, Th.\ 2 (C)] 
that, if $M$ is the unit cotangent bundle of a compact Riemannian surface
$R$ of constant negative curvature, then $(0.3)$ holds uniformly for all 
$p\in M$ with exponent $0<\alpha'<1$ determined by the smallest non-zero 
eigenvalue of the Laplace-Beltrami operator on $R$. The results of Ratner 
[52] and Burger [7] are based on the Fourier analysis for the Lie group 
${\rm PSL}(2,{\Bbb R})$. The horocycle flow, which is uniquely ergodic and 
mixing, is an example of {\it parabolic }dynamics, characterized by 
polynomial divergence of nearby orbits.

In this paper we establish results on the asymptotic behaviour 
of ergodic 
averages for a class of conservative flows with isolated canonical saddle singularities 
on higher genus surfaces. Such flows were originally introduced and studied mainly by 
the Russian school as a simple model of the dynamics of Hamiltonian systems in the 
case of a `multivalued' Hamiltonian, as in the presence of an electromagnetic field 
[48, \S 6], or other first integral, as in the case of billiards in rational polygons and related 
systems \ref\ZK. L.-S. Young \ref\Ygone proved that all continuous flows on surfaces have 
zero topological entropy. A.V. Kocergin [37, \S 5, Th.\ 3] established the mixing property for 
ergodic flows in case all singularities are (degenerate) {\it mixing saddle points }[37,  
Def.\ 5.3]. To the author's best knowledge, no satisfactory estimates for the decay of 
correlations are known and mixing is an open question in case one of the saddles is of 
Morse type (a recent result of B. Fayad [16, I, Th.\ 3.1] for flows on the torus suggests 
that flows with mixing saddle points should have polynomial decay of correlations).  In 
contrast to the asymptotic behaviour of ergodic averages, which essentially depends only on the orbit
foliation, the mixing property for flows with singularities is extremely 
sensitive to time changes. In fact, if the return time function (to a transverse interval) 
is of bounded variation, then the flow is not mixing \ref\Kttwo. In this case the flow cannot 
however be continuous. The (unique) ergodicity on the complement of the singular set
for generic flows, known in the related case of interval exchange transformations as the
{\it Keane conjecture}, was finally proved  by H. Masur \ref\Mstwo and W. Veech \ref\Vcone.  

 A. Zorich, motivated by a question of S. P. Novikov [48, \S 6] on the semiclassical 
trajectories of an electron on the Fermi surface of a metal [82, \S\S 1.3 \& 1.4], has discovered 
in numerical experiments a striking non-standard behaviour of ergodic averages. In his papers 
[78]--[82] and his joint paper with M. Kontsevich \ref\KZone he has proved several partial results 
and formulated conjectures that draw a fairly complete picture of the new phenomenon. In
\ref\Ftwo the author has proved the existence of distributional obstructions, not given by 
measures, to solving the cohomological equation $Xu=f$. As a consequence, the estimate 
$(0.2')$ can be proved only for a finite codimensional subspace of the space of smooth 
functions (described by the vanishing of a finite number of invariant distributions), in 
agreement with Zorich's conjectures. This paper essentially completes Zorich's program 
on the deviation of ergodic averages. The approach followed is Zorich's and, especially,
Kontsevich-Zorich's \ref\KZone, except for the relevant role played by invariant distributions 
and by the methods developed in \ref\Ftwo. The results obtained below, together with those 
of \ref\Ftwo, show that the dynamics of generic conservative flows on higher genus surfaces 
differs significantly from the genus one case. It appears in fact to be closer to the behaviour 
expected in the parabolic case, represented by the horocycle flow. This phenomenon can 
heuristically be explained by the divergence of nearby orbits produced by saddle-like 
singularities that are unavoidable in the higher genus case [33, Chap.\ 8, \S 4].

 Let $M$ be a compact orientable surface of genus $g\geq 2$, let $(\Sigma,\imath)$ be a 
divisor on $M$, in the sense that $\Sigma:=\{p_1, \ldots,p_{\sigma}\}$ is a finite subset of $M$ and 
$\imath:=(\imath_1, \ldots,\imath_{\sigma})\in {\Bbb Z}^{\sigma}$, with the properties that $\imath_k 
<0$ for all $k\in\{1, \ldots,\sigma\}$ and $\sum_k\imath_k=2-2g$. Let ${\Cal F}^{\imath}(M,\Sigma)$
be the space of orientable measured foliations, in the sense of W.Thurston \ref\Th, with
(possibly degenerate) saddle-like singularities of canonical type at $\Sigma$, of index $\imath_k$
at $p_k$ for each $k\in\{1, \ldots,\sigma\}$, and no regular leaves homologous to zero. Let $\omega$ 
be any smooth non-negative \hbox{$2$-form} on~$M$, positive on $M\setminus\Sigma$, allowed to have
zeroes at
$\Sigma$. Let ${\Cal E}^{\imath}_{\omega}(M,\Sigma)$ be the space of all smooth $\omega$-preserving
vector fields on $M\setminus\Sigma$ such that the form $\eta_X:=\imath_X
\omega$ is smooth (and closed) on $M$ and the orbit foliation ${\Cal F}_X:=\{\eta_X=0\}\in {\Cal F}
^{\imath}(M,\Sigma)$. The main goal of the paper is to prove the following result, conjectured by
M.~Kontsevich and A. Zorich \ref\KZone on the basis of numerical experiments:

\proclaimtitle{Theorem 9.6} 
\proclaim{ Theorem} For {\rm `}\/almost all\/{\rm '} $\,X\in {\Cal E}^{\imath}_{\omega}(M,\Sigma)${\rm ,}
the flow $\Phi_X^{\tau}$ has a deviation spectrum in the following sense{\rm .} There exists a finite set 
of exponents
$$\lambda'_1(X)=1> \lambda'_2(X)>\cdots >\lambda'_s(X)>0\,\,\eqnu$$
and a splitting of the space ${\Cal I}_X^1(M)$ of $X$\/{\rm -}\/invariant distributions of order $1$ 
{\rm (}\/i.e.\ solutions
${\Cal D}\in H^{-1}(M)$ of the equation $X{\Cal D}=0$)
$${\Cal I}_X^1(M)={\Cal I}_X^1(\lambda'_1)\oplus {\Cal I}_X^1(\lambda'_2)\oplus\cdots \oplus 
{\Cal I}_X^1(\lambda'_s) \speqnu {0.4'}$$
such that the following holds{\rm .} Let $f\in H^1_0(M\setminus\Sigma)$ be a weakly differentiable 
function
with $L^2$ partial derivatives{\rm ,} supported in $M\setminus\Sigma${\rm ,} such that
$${\Cal D}^X(f)=0\,\,,\,\,\,\,\hbox{ for all }\,{\Cal D}^X\in {\Cal I}_X(\lambda'_1)\oplus\cdots\oplus
{\Cal I}_X(\lambda'_i)\,\,.\eqnu $$
If $i<s$ in $(0.5)${\rm ,} then  for all $p\in M$ with regular forward trajectory{\rm ,}
$$\limsup_{{\Cal T}\to +\infty}{{\log |\int_0^{\Cal T} f(\Phi_X(p,\tau))\,d\tau|}\over
             {\log {\Cal T}}}\,\leq \, \lambda'_{i+1}(X)\eqnu $$
and there exists ${\Cal D}^X_{i+1}\in {\Cal I}_X(\lambda'_{i+1})\setminus\{0\}$ such that{\rm ,} if 
$\,{\Cal D}^X_{i+1}(f)\not= 0${\rm ,} then equality holds in $(0.6)$ for almost all $p\in M${\rm .}

If $i=s$ in $(0.5)${\rm ,} then{\rm ,}  for all $p\in M$ with regular forward trajectory{\rm ,}
$$\limsup_{{\Cal T}\to +\infty}{{\log |\int_0^{\Cal T} f(\Phi_X(p,\tau))\,d\tau|}\over
              {\log {\Cal T}}}\,= \, 0\,\,.\speqnu{0.6'}$$
The multiplicities $\,m_i(X):=\hbox{{\rm dim}}\,{\Cal I}_X(\lambda'_i)$ satisfy the conditions $\sum_i 
m_i(X)=g$ and $\,m_1(X)=1${\rm ,} since ${\Cal I}_X(\lambda'_1)={\Bbb R}\cdot\omega\,${\rm }. The deviation 
spectrum with multiplicities $\{\bigl(\lambda'_i(X),m_i(X)\bigr)\,|\,i=1, \ldots,s\}$ is locally constant
on ${\Cal E}^{\imath}_{\omega}(M,\Sigma)${\rm .}
\endproclaim

 The {\it measure class }on ${\Cal E}^{\imath}_{\omega}(M,\Sigma)$ is the (standard) one
obtained by considering the Lebesgue measure on the A. Katok's fundamental classes of area-preserving
vector fields. The {\it fundamental class }of a vector field $X\in {\Cal E}_{\omega}(M,\Sigma)$ is, 
following \ref\Ktone, the cohomology class of the closed $1$-form $\eta_X:=\imath_X\omega$ in
$H^1(M,\Sigma;{\Bbb R})$. The pull-back of the Lebesque measure class on the cohomology 
vector space under the map ${\Cal E}^{\imath}_{\omega}(M,\Sigma)\to H^1(M,\Sigma;{\Bbb R})$, 
given by $X\to [\eta_X]$, yields the measure class considered in the statement of Theorem 0.1. 
The symbols $H^1(M)$, $H^{-1}(M)$ denote the standard $L^2$ {\it Sobolev spaces }on the compact 
manifold $M$ \ref\Ad. The space $H^1(M)$ is defined as the space of $L^2$ weakly differentiable
functions with $L^2$ partial derivatives and the space $H^{-1}(M)$ as the dual of the Banach space
$H^1(M)$. The space $H^1_0(M\setminus\Sigma)\subset H^1(M)$ is the subspace of functions with
compact support in $M\setminus\Sigma$. 

In the case $i=0$, the statement of Theorem 0.1 follows immediately from the 
(unique) ergodicity of `almost all' vector fields (the {\it Keane Conjecture}, proved by \ref\Mstwo,
\ref\Vcone and later on by several other authors by different methods) and from Birkhoff's ergodic
theorem. In fact, the continuous embedding $H^1(M) \subset L^1(M,\omega)$ holds by the H\"older
inequality. In case the function $f$ has zero average, by the theory of sign changes of ergodic integrals
[50, Chap.\ 3(C)], the lower limit of the quantity in $(0.6)$ is $-\infty$ and the limit does not exist for
almost all $p\in M$. In the related case of interval exchange transformations, A. Zorich \ref\Zrthree,
\ref\Zrfive proved a statement similar to Theorem 0.1 for linear combinations of characteristic 
functions of the sub-intervals. He conjectured that all the exponents in $(0.4)$ are non-zero and 
that the spectrum is {\it simple}, in the sense that all the multiplicities $m_1=\cdots =m_s=1$ (hence 
$s=g$), as suggested by the numerical evidence and the results he had presented in \ref\Zrone,
\ref\Zrtwo, \ref\Zrthree. We were unable to prove the part of the conjecture concerning the simplicity 
of the spectrum, except in the case $g=2$. In \ref\KZone M. Kontsevich and A. Zorich conjectured that, 
in the case of area-preserving flows on surfaces, Zorich's theorem could be generalized to smooth
 functions. The role played by invariant distributions was not part of the Kontsevich-Zorich picture. 
The author regards it as his most important contribution to the subject.

The proof of Theorem 0.1 is based on several results, concerning the Teichm\"uller geodesic
flow, which, we believe, are of independent interest. Following the conjectural approach outlined in 
\ref\KZone, we have proved the non-uniform hyperbolicity of a `renormalization' cocycle, introduced 
in \ref\KZone, and earlier in a different form in \ref\Zrtwo, \ref\Zrthree. The {\it Kontsevich-Zorich 
cocycle }is a flow on a symplectic vector bundle over the moduli space ${\Cal M}_g$ of holomorphic
quadratic differentials on Riemann surfaces of genus $g\geq 2$. The fiber of the vector bundle at a
quadratic differential $q$ is the cohomology vector space with real coefficients $H^1(M_q,{\Bbb R})$ 
of the Riemann surface $M_q$ carrying $q$. The cocycle is defined as a lift of the Teichm\"uller
geodesic flow, obtained by parallel transport of cohomology classes with respect the canonical flat
(Gauss-Manin) connection of the bundle. Since the Kontsevich-Zorich cocycle is defined on a 
symplectic vector bundle (of dimension $2g$), it has a symmetric Lyapunov spectrum:
$$\lambda_1=1\geq \lambda_2\geq\cdots\geq \lambda_g\geq 0 \geq -\lambda_g\geq\cdots
-\lambda_2\geq  -\lambda_1=-1\,\,. \eqnu$$
The numbers $\lambda'_1(X)>\cdots>\lambda'_s(X)>0$ appearing in the statement of Theorem 0.1 are the 
{\it distinct }positive Lyapunov exponents of the Kontsevich-Zorich cocycle with respect to 
the canonical absolutely continuous invariant measure of the Teichm\"uller flow on the appropriate
connected component of a stratum of the moduli space determined by $X\in
{\Cal E}^{\imath}_{\omega}(M,\Sigma)$. We prove
the following result:

\proclaimtitle{Corollary 2.2 \& Th.\ 8.5}
\proclaim {Theorem}
{\rm (i)} The first Lyapunov exponent $\lambda^{\mu}_1=1$ of the Kontsevich\/{\rm -}\/Zorich cocycle{\rm ,}
 with
respect to any ergodic probability measure $\mu$ of the Teichm{\rm \"{\it u}}ller flow on the stratum of 
squares
of holomorphic differentials{\rm ,} is simple{\rm .} In fact{\rm ,} the second Lyapunov exponent
 $\lambda^{\mu}_2<1${\rm .}

\smallbreak
{\rm (ii)} The Kontsevich\/{\rm -}\/Zorich cocycle is non\/{\rm -}\/uniformly hyperbolic with respect to the 
canonical absolutely continuous invariant probability measure of the Teichm{\rm \"{\it u}}ller
 flow{\rm ,} in the sense that{\rm ,} 
for {\rm (}\/Lebesgue\/{\rm )} almost all quadratic differentials on each connected component of a
stratum of squares{\rm ,}
$$\lambda_1=1>\lambda_2\geq \cdots \geq \lambda_g >0\,\,. \speqnu {0.7'}$$
\endproclaim

 The Lyapunov spectrum of the Kontsevich-Zorich cocycle contains all the non-trivial 
information on the Lyapunov spectrum of the Teichm\"uller geodesic flow, which was proved by
W. Veech \ref\Vctwo to be non-uniformly hyperbolic on each connected component of a stratum of
quadratic differentials, with respect to the canonical absolutely continuous invariant probability
measure. In fact, the explicit formulas $(1.2')$ hold and Theorem 0.2 is therefore equivalent to the 
following:

\specialnumber{ 0.2'}
\proclaim{Theorem} 
{\rm (i)} All probability invariant measures of the Teichm{\rm \"{\it u}}ller geodesic flow{\rm ,} in particular 
all periodic trajectories{\rm ,} on any stratum of squares{\rm ,} are non\/{\rm -}\/uniformly hyperbolic{\rm ,}
 in the sense that
the Lyapunov exponent $0$\break {\rm (}\/corresponding to the direction of the flow\/{\rm )} is simple{\rm
.}
\smallbreak
{\rm (ii)} The multiplicity of the Lyapunov exponent $1${\rm ,} on every connected component of a 
stratum of squares having $\sigma\in \{1, \ldots,2g-2\}$ distinct zeroes{\rm ,} is almost everywhere exactly 
equal to $\sigma-1${\rm ,} with respect to the canonical absolutely continuous invariant
 probability measure{\rm .}
\endproclaim

The proof of Theorem 0.2 (and of Theorem 0.2$'$), carried out in Sections 2--5 and Section 8, is based 
on explicit equations for the Kontsevich-Zorich cocycle in appropriate coordinates, established in 
Section 2. By applying such equations, we are able to obtain a remarkably simple proof of Theorem 0.2
(i) (Corollary 2.2) and to compute formulas for the hyperbolic Laplacian of the logarithm
of the volume on any isotropic $k$-plane in $H^1(M,{\Bbb R})$ over Teichm\"uller disks (in \S 3 
and \S 5). In the particular case $k=g$, we obtain a new version of the formula obtained in \ref
\KZone. The argument is thus reduced to a geometric problem concerning the determinant of the 
directional derivatives of the period matrix, which is solved in Section~4 by studying quadratic 
differentials near special boundary points of the moduli space. The complete proof of the
non-uniform hyperbolicity of the Kontsevich-Zorich cocycle (Theorem 8.5) depends on additional
information on the invariant sub-bundles (which could be in fact deduced, as remarked by A. Eskin, 
from \ref\Vctwo), studied in Section~6 and Section~8.       

 Any holomorphic quadratic differential $q$ on a Riemann surface $M_q$ can be identified
with a pair $({\Cal F}_q,{\Cal F}_{-q})$ of transverse {\it measured foliations }(in the sense of 
\ref\Th), the {\it horizontal foliation }${\Cal F}_q:=\{\Im(q^{1/2})=0\}$, the {\it vertical foliation
 }${\Cal F}_{-q}:=\{\Re(q^{1/2})=0\}$. Let $q$ be the square of a holomorphic (abelian) differential 
which is a {\it regular }point (in the sense of the Oseledec's theorem \ref\Os, [32, Th.\ S.2.9]) 
of the Kontsevich-Zorich cocycle. Let $E^+_q\,[E^-_q]\subset H^1(M_q,{\Bbb R})$ the unstable 
[the stable] subspace at $q$ of the Kontsevich-Zorich cocycle. Such spaces have a dynamical
interpretation, discovered by A. Zorich \ref\Zrone, in terms of the foliations ${\Cal F}_{\pm q}$. Let
${\widehat\gamma}_{\pm q}^{\Cal T}$ be the closed curve on $M$ obtained as the union of a leaf 
of length ${\Cal T}>0$ of the foliation ${\Cal F}_{\pm q}$ and of a `short' segment of length
 $\leq\hbox{diam}\,(M_q)$. Then the distance of the homology class $[{\widehat \gamma}_{\pm q}
^{\Cal T}]\in H_1(M,{\Bbb R})$ from the Poincar\'e dual of the subspace $E^{\pm}_q$ is 
(uniformly) bounded for all ${\Cal T}>0$. The cycles contained in the Poincar\'e dual of the 
subspace $E^{\pm}_q$ have been called [47, \S 7.9.3] the {\it Zorich cycles }of the measured 
foliation ${\Cal F}_{\pm q}$. The Schwartzman \ref\Sn {\it asymptotic cycle }of the foliation 
${\Cal F}_{\pm q}$ is a particular Zorich cycle. It follows from Theorem 0.2(ii) that the dimension 
of the space of Zorich cycles is `generically' equal to the genus $g\geq 2$. 

\pagegoal=50pc
We have proved that all cohomology classes in the $E^{\pm}_q$, i.e. all the
Poincar\'e duals of Zorich cycles of ${\Cal F}_{\pm q}$, can be represented by (closed) currents 
with special properties. Let $\Sigma_q\subset M$ be the finite set of the zeroes of the quadratic
differential $q$. A current of dimension (and degree) equal to $1$ on the $2$-dimensional manifold 
$M\setminus\Sigma_q$ is a continuous linear functional on the space of smooth $1$-forms with 
compact support \ref\dR, [56, Chap.\ IX]. A {\it basic current }for the measured foliation 
${\Cal F}_{\pm q}$ is defined as a distributional generalization of the notion of basic form, well 
known in the geometric theory of foliations since the work of\break B. L. Reinhart \ref\Reone, \ref\Retwo 
(see also [2, \S 1.5 \& \S 7], [63, Chap.\ 4]). A current $C$ of dimension (and degree) equal to $1$ 
is basic for the measured foliation ${\Cal F}_{\pm q}|_{M\setminus\Sigma_q}$
 
\smallbreak\noindent if, for all vector
fields
$X$ tangent to ${\Cal F}_{\pm q}$ with compact support in $M\setminus\Sigma_q$,
\vglue-16pt
$$ \imath_X C={\Cal L}_X C=0\,\,. \eqnu$$ 
\vglue-4pt
The operation of contraction $\imath_X$ and Lie derivative ${\Cal L}_X$ are extended to currents in 
the standard distributional sense [56, Chap.\ IX, \S 3]. If a current $C$ on $M\setminus\Sigma_q$
(which is $2$-dimensional) is basic then it is closed, hence it represents, by the {\it generalized 
de Rham theorem }[11, Th.\ 12], [56, Chap.\ IX, \S 3, Th.\ I], a cohomology class in $H^1(M
\setminus\Sigma_q,{\Bbb R})$. Let $X^{\pm}\in {\Cal E}^{\imath}_{\omega}(M,\Sigma_q)$ be vector fields
respectively tangent to the leaves of ${\Cal F}^{\pm}_q$. There is a
one-to-one correspondence between 
the space of $X^{\pm}$-invariant distributions on $M\setminus\Sigma_q$ and the space of basic currents for 
${\Cal F}_{\pm q}|_{M\setminus\Sigma_q}$. In fact, if ${\Cal D}^{\pm}$ is any $X^{\pm}$-invariant
distribution, then its contraction $\imath_{X^{\pm}} {\Cal D}^{\pm}$ is a current of dimension (and degree)
equal to $1$ which is basic for ${\Cal F}_{\pm q}$. The $X^{\pm}$-invariant measures are a particular case
of 
$X^{\pm}$-invariant distributions and the related basic currents where implicitly considered by A. Katok 
\ref\Ktone.

\pagegoal=48pc
 The space of currents on $M\setminus\Sigma_q$ is filtered by a sequence of {\it weighted
Sobolev spaces }${\Cal H}^{-s}_q(M)$, $s\geq 0$. A basic current for ${\Cal F}_{\pm q}$, of {\it finite 
order }$s\geq 0$, is a current $C\in {\Cal H}^{-s}_q(M)$ of dimension (and degree) equal to $1$ such 
that the identities $(0.8)$ hold for all vector fields $X$ tangent to ${\Cal F}_{\pm q}$ for which the 
operators $\imath_X:{\Cal H}^{-s}_q(M)\to {\Cal H}^{-s}_q(M)$ and ${\Cal L}_X:{\Cal H}^{-s}_q(M) \to
{\Cal H}^{-s-1}_q(M)$ are well defined. Basic currents of order $0$ correspond to {\it absolutely continuous}
invariant measures. It can be proved (see Lemma 6.2) that the cohomology class of a basic current for
${\Cal F}_{\pm q}$, of order $s\geq 0$, belongs to the cohomology vector space $H^1(M,{\Bbb R})\subset 
H^1(M\setminus\Sigma_q,{\Bbb R})$. The following result describes the invariant sub-bundles of the 
Kontsevich-Zorich cocycle in terms of basic currents:

\vglue-9pt
\proclaimtitle {Theorem 8.3}
\proclaim {Theorem} For almost all quadratic differential $q$ in every stratum of squares of 
holomorphic differentials{\rm ,} the invariant unstable {\rm [}\/stable\/{\rm ]} subspace at $q$ of the
Kontsevich\/{\rm -}\/Zorich  cocycle is the subspace $E^+_q\,[E^-_q]\subset H^1(M_q,{\Bbb R})$ given by the
cohomology classes of  the basic currents for ${\Cal F}_q$ $[{\Cal F}_{-q}]$ which are of order $1${\rm .}
\endproclaim

Theorem 0.3 is in fact valid with respect to any ${\rm SL}(2,{\Bbb R})$ invariant probability 
measure on the moduli space and it requires no information on the dimension of the stable or
unstable sub-bundle. The proof is based on a Cheeger-type estimate  of the first non-trivial
eigenvalue of the Dirichlet form of the flat metric induced by a quadratic differential, proved
in Section 6.3, and on the {\it logarithmic law }for geodesics in the moduli space, proved in \ref\Msthree.

 Theorem 0.2 and Theorem 0.3 imply that the dimension of the space of basic current of order
$1$ for a `generic' measured foliation on a compact orientable surface $M$ is equal to the genus
$g\geq 2$ of $M$. It is interesting to compare such a result with the results we have proved in Section 7 
on the space of all basic currents of finite order.

\proclaimtitle{Theorems 7.1 \& 7.7}  
\proclaim {Theorem} Let $\eta_{\Cal F}$ denote the smooth closed $1$\/{\rm -}\/form
such that ${\Cal F}:=\{\eta_{\Cal F}=0\}$\/{\rm . }\/

 {\rm (i)}  There exists $r>1$ such that\/{\rm ,}\/ for almost all orientable measured foliations ${\Cal
F}\in  {\Cal F}^{\imath}_{\omega}(M,\Sigma)$\/{\rm ,}\/ the space $H^{1,s}_{\Cal F}(M,{\Bbb R})$ of
cohomology  classes of ${\Cal F}$\/{\rm -}\/basic currents of any order $s\geq r$ {\rm (}\/called the
$H^{-s}$ basic cohomology  of ${\Cal F}${\rm )} has codimension $1$ in the cohomology $H^1(M,{\Bbb
R})${\rm .} In fact{\rm ,} 
$$H^{1,s}_{\Cal F}(M,{\Bbb R}):=\{c\in H^1(M,{\Bbb R})\,|\, c\wedge [\eta_{\Cal F}]=0\}\,\,. 
\eqnu$$ 

{\rm (ii)} Let $X$ be any smooth vector field on $M\setminus\Sigma$ such that $\imath_X
\eta_{\Cal F}\equiv 1${\rm .}  Let $C$ be an ${\Cal F}$\/{\rm -}\/basic current of order $s\geq 0$ such that 
the cohomology class $[C]=0\in H^1(M,{\Bbb R})${\rm .} Then there exists an ${\Cal F}$\/{\rm -}\/basic current
$\widehat C$ of order $s-1$ such that $C=d\,\imath_X {\widehat C}${\rm . (}\/The restriction of the
operator 
$d\circ \imath_X$ to ${\Cal F}$\/{\rm -}\/basic currents does not depend on the choice of the vector
field $X${\rm ).}
\endproclaim

By Theorem 0.4(i), the regularity restriction on basic current in the statement of 
Theorem 0.3 is crucial. Theorem 0.4(ii) implies that there are only a finite number of `primitive' 
basic currents for any measured foliation. However, since $2g-1>g$, by Theorem 0.4(i), not all of 
them have a dynamical significance, according to Theorem 0.1 or Theorem 0.3. The question whether all 
basic currents (with non-vanishing cohomology class) have a dynamical interpretation is left 
unanswered (see Question 9.10). The proof of Theorem 0.4 is based on the results of \ref\Ftwo.

 In order to prove Theorem 0.1, we introduce a cocycle $G^c_t$ over the Teich\"muller flow 
$G_t$, defined on the infinite dimensional bundle $Z^{-1}(M)$ over the moduli space with fiber given 
by the space of all closed currents of order $1$. The projection of the cocycle $G^c_t$ under the
map $Z^{-1}_q(M)\to H^1(M_q,{\Bbb R})$, given by the generalized de Rham theorem, coincides with the
Kontsevich-Zorich cocycle.  There is no general Oseledec's theorem for infinite dimensional (Hilbert)
bundles and the available ones (see the survey \ref\Sl) do no seem to apply to the case at hand. 
However, we derive from Theorem 0.2 and Theorem 0.3 the following Oseledec-type splitting:

\proclaimtitle{Theorem 8.7}
\proclaim {Theorem} There is a $G^c_t$\/{\rm -}\/invariant splitting
$$ Z^{-1}(M)\equiv {\Cal B}^1_+(M)\oplus {\Cal B}^1_-(M)\oplus {\Cal E}^{-1}(M)\,\,.\eqnu$$
\medbreak
\item{\rm (i)} ${\Cal B}^1_{\pm}(M)$ is the bundle of ${\Cal F}_{\pm q}$-basic 
currents of order $1$ over the moduli space of squares of holomorphic differentials{\rm .} 

\smallbreak
\item{\rm (ii)} The restriction of the cocycle $G^c_t$ to the bundle ${\Cal B}^1_+(M)\oplus 
{\Cal B}^1_-(M)$ is measurably isomorphic to the Kontsevich\/{\rm -}\/Zorich cocycle on
 the cohomology bundle{\rm .}

\smallbreak
\item{\rm (iii)} The Lyapunov spectrum of the cocycle $G^c_t$ on ${\Cal B}^1_{\pm}(M)$ 
consists of the set of Lyapunov exponents $\pm\{\lambda_1, \ldots,\lambda_g\}$ given by $(0.7')${\rm .}

\smallbreak
\item{\rm (iv)} The fiber of the bundle ${\Cal E}^{-1}(M)$ at $q$ consists of all exact currents of\break 
order~$1$. The Lyapunov spectrum of the cocycle $G^c_t$ on ${\Cal E}^{-1}(M)$ is reduced to the single
Lyapunov exponent $0${\rm .}

\endproclaim

 The strategy of the proof of Theorem 0.1 is then clarified by the remark that the return 
orbits of a flow $\Phi^{\tau}_X$, $X\in {\Cal E}^{\imath}_{\omega}(M,\Sigma)$, closed by a (short) 
transverse segment, are closed currents of order $1$, by the Sobolev embedding theorem [1, Th.\
5.4]. Hence, their evolution under the action of the cocycle $G^c_t$ can be studied by applying 
Theorem 0.5. Since $G^c_t$ is a cocycle over the Teichm\"uller flow, which plays the role of a
{\it renormalization dynamics }for the flows of the vector fields $X\in {\Cal E}^{\imath}_{\omega}
(M,\Sigma)$, Theorem 0.1 can be derived from Theorem 0.5 by standard techniques. The $X$-invariant 
distributions appearing in the statement of Theorem 0.1 are in one-to-one correspondence with the 
basic currents of order $1$ for the orbit foliation ${\Cal F}_X$ given by 
Theorem 0.2(ii) and Theorem~0.3.

\vglue6pt
\section{The Kontsevich-Zorich cocycle}
\vglue6pt

Let $T_g$ be the Teichm\"uller space of all marked closed Riemann surfaces of genus
$g\geq 2$. It can be viewed as the quotient of the space of all complex structures on a compact
orientable topological surface of genus $g$, with respect to the equivalence relation induced by
the natural action of the group $\hbox{Diff}^+_0(M)$ of orientation-preserving diffeomorphisms
{\it homotopic to the identity}. For a general introduction to Teichm\"uller theory and
references to the relevant literature we refer to the survey article by L. Bers [6]. In 
particular, $T_g$ has a natural complex structure (Alhfors-Bers), isomorphic to that of a bounded 
domain of holomorphy of the Euclidean space ${\Bbb C}^{3g-3}$, and it can be equipped with a natural
Finsler metric (the Teichm\"uller metric). The quotient of the bundle of non-zero holomorphic
quadratic differentials on marked closed Riemann surfaces of genus $g$, with respect to the
natural action of $\hbox{Diff}^+_0(M)$ by pull-back, is a complex manifold $Q_g$, which can be
identified with the cotangent bundle of $T_g$ (minus the zero section). We will be concerned with
subvarieties of $Q_g$ given by holomorphic quadratic differentials with a prescribed pattern of
zeroes.

 Let $\kappa=(k_1, \ldots ,k_{\sigma})$ satisfy the following properties: 
each $k_i$ is
a positive even integer and $\sum k_i=4g-4$. We define $Q_{\kappa}$ to be the subset of $Q_g$
consisting of quadratic differentials $q$ such that $(1)$ $q$ is the square of a holomorphic
(abelian) differential; $(2)$ the distinct zeroes of $q$ have orders $(k_1, \ldots ,k_{\sigma})$.
The set $Q_{\kappa}$ is a complex analytic subvariety of $Q_g$, which is called a {\it stratum}
of~$Q_g$. Since $\sum k_i=4g-4$, the number of strata in $Q_g$ is finite. The strata of $Q_g$
which correspond to non-orientable quadratic differentials will not be considered here. They
can be studied by passing to suitable branched double covers over $M$. If $q\in Q_{\kappa}$,
we denote by $M_q$ be the Riemann surface carrying~$q$ and by $\Sigma_q$ the finite set of
distinct zeroes of $q$.

Let $\Gamma_g:=\hbox{Diff}^+(M)/\hbox{Diff}_0^+(M)$ be the {\it pure mapping class
group }of the compact orientable surface $M$ of genus $g\geq 2$. It is defined as the group of
isotopy equivalence classes of orientation-preserving diffeomorphisms of $M$. The mapping class
group acts by pull-back on $T_g$, $Q_g$, $Q_{\kappa}$. Such actions are
properly discontinuous
and given by biholomorphic isometries (with respect to the Teichm\"uller metric). The quotient
spaces $R_g:=T_g/\Gamma_g$, ${\Cal M}_g:=Q_g/\Gamma_g$  are, respectively, the {\it moduli spaces}
of complex structures and of quadratic differentials on surfaces of genus $g$. Since there are
Riemann surfaces with non-trivial (but finite) automorphism group, $R_g$, $Q_g$ are complex
orbifolds (hence in particular complex analytic spaces) of dimension $3g-3$, i.e.\ for each point
$p\in R_g$  $[p\in {\Cal M}_g]$  there is a neighbourhood ${\Cal U}_p$, a finite subgroup
$G\subset\hbox{GL}(3g-3,{\Bbb C})$  $[\hbox{GL}(6g-6,{\Bbb C})]$  and a neighbourhood
$\tilde{\Cal U}_p$ of the origin in ${\Bbb C}^{3g-3}$ $[{\Bbb C}^{6g-6}]$, such that
${\Cal U}_p\equiv \tilde{\Cal U}_p/G$. For each symbol $\kappa$, the quotient ${\Cal M}_{\kappa}
:=Q_{\kappa}/\Gamma_g$ is called a stratum of the moduli space ${\Cal M}_g$. From their
definition, it follows that strata are complex analytic subvarieties of ${\Cal M}_g$. Moreover,
each stratum ${\Cal M}_{\kappa}$ can be equipped with the structures listed below [66], 
[68], [38]:
\medbreak
\item{(1)} {\it a natural complex affine orbifold structure modeled on the vector space}
$H^1(M,\Sigma_{\kappa};{\Bbb C})$, {\it where }$\Sigma_{\kappa}$ {\it is a finite subset of}
$M$ {\it with }card($\Sigma_{\kappa})=\sigma_{\kappa}$;
\item{(2)} {\it an absolutely continuous measure }$\mu_{\kappa}$;
\item{(3)} {\it a locally quadratic non-holomorphic function }$A:{\Cal M}_{\kappa}\to{\Bbb R}^+$;
\item{(4)} {\it a non-holomorphic action of the group }$GL_+(2,{\Bbb R})$.
\medbreak
\noindent The structure $(1)$ is given by the {\it period map}. There exists a locally
defined continuous square root on $Q_{\kappa}$, such that if $q\in Q_{\kappa}$, $q^{1/2}$ is a
holomorphic (abelian) differential with zeroes at $\Sigma_q$. The cohomology class of $q^{1/2}$
in $H^1(M_q,\Sigma_q;{\Bbb C})$ is defined by integration of $q^{1/2}$ along relative cycles in
$(M_q,\Sigma_q)$. The period map gives a natural complex affine structure on $Q_{\kappa}$, modeled
on the vector space $H^1(M,\Sigma_{\kappa};{\Bbb C})$, where $\Sigma_{\kappa}$ is a finite set with
$\hbox{card}(\Sigma_{\kappa})=\sigma$. In fact, the cohomology $H^1(M_{\tilde q},\Sigma_{\tilde q};
{\Bbb C})$ can be canonically identified with $H^1(M_q,\Sigma_q;{\Bbb C})$, by parallel transport
with respect to the so-called {\it Gauss-Manin connection }[23, \S 2], on a neighbourhood
${\Cal U}_q$ of $q\in Q_{\kappa}$. The Gauss-Manin connection in this context is simply given by the
property that parallel sections are the holomorphic sections of the holomorphic vector bundle
over $Q_{\kappa}$ with fiber $H^1(M_q,\Sigma_q;{\Bbb C})$ which are {\it locally constant}. Hence
the period map $\Pi_q$ is defined on ${\Cal U}_q$ into $H^1(M_q,\Sigma_q;{\Bbb C})$ and is a
(local) biholomorphism of ${\Cal U}_q$ onto an open subset of $H^1(M_q,\Sigma_q;{\Bbb C})$. This
result can be proved by studying holomorphic deformations of quadratic differentials in terms
analogous to the Kodaira-Spencer theory of deformations of complex structures [27, Chap.
IV, \S 2]. The complex affine structure on $Q_{\kappa}$, just described, induces, by quotient
with respect to the action of the mapping class group $\Gamma_g$, a complex affine orbifold structure
on the stratum ${\Cal M}_{\kappa}$ in the moduli space.

The construction of the measure $\mu_{\kappa}$ in $(2)$ is also based on the period map.
In fact, the Lebesgue measure on $H^1(M_q,\Sigma_q;{\Bbb C})$, uniquely normalized by
the condition that the volume of the $2g-1+\sigma$ complex torus obtained as a quotient over
the integer lattice ${\Bbb C}\otimes_{\Bbb Z}H^1(M_q,\Sigma_q;{\Bbb Z})$ is equal to $1$, gives a
smooth measure on $Q_{\kappa}$, which projects by push-forward to a measure $\mu_{\kappa}$ on~${\Cal M}_{\kappa}$.

The total area function $A:Q_g\to {\Bbb R}^+$ is defined for any quadratic
differential $q$ as the total area of the surface $M_q$ with respect to the area form
$\omega_q$ of the metric $R_q$ induced by $q$ [58, 5.3]. The function $A$ is invariant
under the action of $\Gamma_g$, hence it can be viewed as a function on ${\Cal M}_g$.

Finally, the group ${\rm GL}_+(2,{\Bbb R})$ of $2\times 2$ matrices with positive
determinant acts on $Q_{\kappa}$ through the linear transformations on the pairs of real-valued\break
$1$-forms $\bigl(\Re(q^{1/2})\,,\,\,\Im(q^{1/2})\bigr)$. In the affine coordinate
system described by $(1)$, it is the action of the group ${\rm GL}_+(2,{\Bbb R})$ on the vector space
$$H^1(M_q,\Sigma_q;{\Bbb C})\equiv {\Bbb C}\otimes_{\Bbb R} H^1(M_q,\Sigma_q;{\Bbb R})
\equiv {\Bbb R}^2\otimes_{\Bbb R} H^1(M_q,\Sigma_q;{\Bbb R}) \eqnu$$
through the first factor in the tensor product. Such an action commutes with the action of the
mapping class group $\Gamma_g$ and it induces therefore an action on~${\Cal M}_{\kappa}$.

 It follows immediately from the definitions that the subgroup ${\rm SL}(2,{\Bbb R})$ of
$2\times 2$ matrices with determinant equal to $1$ preserves the measure $\mu_{\kappa}$ and the area
function $A$. Consequently, on the hypersurface $Q^{(1)}_{\kappa}:=A^{-1}(1)\cap Q_{\kappa}$ we can define
an induced smooth measure $\mu_{\kappa}^{(1)}:= \mu_{\kappa}/dA$ which is preserved by the action
of the group ${\rm SL}(2,{\Bbb R})$. The action of the $1$-parameter subgroup of diagonal matrices
$G_t:=\hbox{diag}(e^t, e^{-t})$ gives a measure-preserving flow on $Q_g$, which preserves each stratum
$Q_{\kappa}$ and each hypersurface $Q^{(1)}_{\kappa}$. Since the action of ${\rm GL}_+(2,{\Bbb R})$
commutes with the action of the mapping class group $\Gamma_g$, $G_t$ induces a measure preserving
flow on the moduli space ${\Cal M}_g$, which preserves each stratum ${\Cal M}_{\kappa}$ and each
hypersurface ${\Cal M}^{(1)}_g:= A^{-1}(1)\cap {\Cal M}_g$, ${\Cal M}^{(1)}_{\kappa}:=A^{-1}(1)\cap
{\Cal M}_{\kappa}$. Such a flow is known as the {\it Teichm\"uller geodesic flow}. As discovered by
W.Veech [68] the strata are not always connected. A complete description of all connected
components of all strata has been given conjecturally in [38] and has been recently proved
in [39] for all the strata of squares. In this paper we will not rely on such a description.
By a theorem of H. Masur [43] and W. Veech [66] the Teichm\"uller flow is ergodic in the
following sense:

\proclaimtitle{[43], [66]}
\proclaim {Theorem}  The total volume of the measure $\mu_{\kappa}^{(1)}$
on ${\Cal M}^{(1)}_{\kappa}$ is finite and the Teichm{\rm \"{\it u}}ller geodesic flow $G_t=\hbox{diag}(e^t, e^{-t})$
is ergodic on each connected component of ${\Cal M}^{(1)}_{\kappa}$ with respect to the invariant 
measure $\mu_{\kappa}^{(1)}${\rm .} \endproclaim

In [38] M. Kontsevich and A. Zorich introduced a natural cocycle over the Teichm\"uller
geodesic flow defined on a bundle over ${\Cal M}^{(1)}_{\kappa}$ with fiber $H^1(M_q,{\Bbb R})$. Such 
a cocycle plays the role of a {\it renormalization group }dynamics for flows along measured foliations 
and for interval exchange transformations, as it will appear from the content of Section~9. Let $H^1_g(M,
{\Bbb C})$ be the holomorphic vector bundle with fiber $H^1(M_q,{\Bbb C})$ over the space $Q^{(1)}_g$ of 
holomorphic quadratic with unit total area. Since $Q^{(1)}_g$ is simply connected, the vector bundle 
$H^1_g(M,{\Bbb C})$ can be trivialized by parallel transport with respect to the Gauss-Manin 
connection, i.e.\ $H^1_g(M,{\Bbb C})\equiv Q^{(1)}_g\times H^1(M,{\Bbb C})$ (Riemann surfaces in the 
Teichm\"uller space are marked). The Teichm\"uller geodesic flow $G_t$ on $Q^{(1)}_g$ lifts trivially 
to the product bundle by the identity action on the second factor. Hence, it lifts to a flow on 
$H^1_g(M,{\Bbb C})$ (in other terms, cohomology classes are moved by parallel transport along the 
Teichm\"uller flow). Such a flow induces, by quotient with respect to the mapping class group 
$\Gamma_g$, a flow on the complex affine (orbifold) vector bundle ${\Cal H}^1_g(M,{\Bbb C}):=H^1_g(M,
{\Bbb C})/\Gamma_g$ which preserves real and imaginary parts of the bundle. The restriction of the 
flow to the real sub-bundle ${\Cal H}^1_g(M,{\Bbb R}):=H^1_g(M,{\Bbb R})/\Gamma_g$, introduced in 
[38], will be called the {\it Kontsevich-Zorich cocycle }and denoted by $G^{KZ}_t$. Such a cocycle
leaves invariant each stratum ${\Cal H}^1_{\kappa}(M,{\Bbb R})$, defined as the restriction to ${\Cal 
M}_{\kappa}^{(1)}$ of the (orbifold) vector bundle ${\Cal H}^1_g(M,{\Bbb R})$. 

Since $H^1(M,{\Bbb R})$ has a natural symplectic structure, given by the wedge product
on the de Rham cohomology, the structure group of the bundles ${\Cal H}^1_{\kappa}(M,{\Bbb R})$
can be reduced to $Sp(2g,{\Bbb R})\subset {\rm GL}(2g,{\Bbb R})$. In this case Lyapunov exponents form
a symmetric subset of the real line. Kontsevich and Zorich conjectured in [38], on the basis
of numerical computations (see also [79]), that the {\it Lyapunov spectrum is simple}, that is:

\nonumproclaim{Kontsevich-Zorich Conjecture} The $2g$ Lyapunov exponents of the Kontsevich\/{\rm -}\/Zorich cocycle 
have the following properties\/{\rm :}\/
$$1=\lambda_1>\lambda_2>\cdots >\lambda_g>0>\lambda_{g+1}=-\lambda_g>\cdots 
>\lambda_{2g}=-\lambda_1=-1\,\,. \eqnu$$
\endproclaim

\numbereddemo{Remark} In [66] W. Veech proved that the Teichm\"uller geodesic flow is
non-uniformly hyperbolic on each stratum ${\Cal M}^{(1)}_{\kappa}$. Since the bundle
${\Cal H}^1_{\kappa}(M,{\Bbb R})$ does not coincide with the whole tangent bundle of
${\Cal M}^{(1)}_{\kappa}$, Veech's theorem does not imply the Kontsevich-Zorich Conjecture.
By $(1.1)$ the Lyapunov exponents of the Teichm\"uller geodesic flow can be written as follows, 
in terms of the exponents $(1.2)$ of the Kontsevich-Zorich cocycle [78, \S 5], [38, \S 7]:
$$ \eqalign{ 2&\geq 1+\lambda_2\geq 1+\lambda_3\geq \cdots\geq 1+\lambda_g \geq \underbrace{1=\cdots= 1}_{\sigma_{\kappa}-1}   \cr
&\geq 1-\lambda_g\geq \cdots \geq 1-\lambda_2\geq 0\geq -1+\lambda_2\geq \cdots \geq -1+\lambda_g  \cr
\cr 
&\geq\underbrace{-1=\cdots= -1}_{\sigma_{\kappa}-1}\geq -1-\lambda_g\geq -1-\lambda_{g-1}\geq \cdots
\geq -1-\lambda_2\geq  -2 \,\,. \cr} \speqnu{1.2'}$$
The non-uniform hyperbolicity of the {\it Teichm\"uller flow} is therefore equivalent to the 
statement that $\lambda_2<\lambda_1=1$. \enddemo

In this paper we prove that the Kontsevich-Zorich cocycle is {\it non-uniformly hyperbolic}, that is
$\lambda_g>0$. Our method also yields a new proof of Veech's theorem $\lambda_2<1$ (Corollary 2.2).

  \vglue-4pt  \section{Variational formulas}
\vglue-4pt

  Let $q\in {\Cal M}^{(1)}_{\kappa}$. Any cohomology class $c\in H^1(M_q,{\Bbb R})$ can be represented
by a harmonic differential. In fact $c=[\Re(h^+)]$, where $h^+$ is a holomorphic differential on $M_q$. The vector
space $H^1(M_q,{\Bbb R})$ can be endowed with the Hodge norm induced by the Hermitian structure on $M_q$. The
Kontsevich-Zorich cocycle along the Teichm\"uller orbit of $q$ can be written in the Hodge representation in the 
form of an O.D.E. on a Hilbert space of square integrable functions. We recall the basic setting established
in [18].

 Let $R_q:=|q|^{1/2}$ be the {\it flat {\rm (}\/smooth\/{\rm )} metric }induced by $q$ on $M$, degenerate at $\Sigma_q$,
and let $\omega_q$ be the corresponding area form [58, \S 5.3]. At any point $p\in M_q\setminus\Sigma_q$,
there exists a canonical holomorphic coordinate $z:=x+iy$ such that $q=dz^2$ [58, \S 5.1].
Consequently,
$$R_q=(dx^2+dy^2)^{1/2}\,\,\,\,\hbox{ and }\,\,\,\, \omega_q=
dx\wedge dy\,\, .\eqnu$$
At a zero $p\in\Sigma_q$ of order $k$, $q=z^kdz^2$ with respect to a canonical coordinate. Hence,
$$R_q=|z|^{k/2}\,(dx^2+dy^2)^{1/2}\,\,\,\,\hbox{ and }\,\,\,\,
\omega_q=|z|^k\,dx\wedge dy\,\,. \speqnu {2.1'}$$
Let ${\Cal F}_{\pm q}$ be the {\it horizontal and the vertical foliation }of $q$ respectively, defined
by ${\Cal F}_q:=\{\Im(q^{1/2})=0\}$, ${\Cal F}_{-q}:=\{\Re(q^{1/2})=0\}$. Such foliations are naturally 
measured foliations in the sense of Thurston \ref\Th with the transverse measures $|\Im(q^{1/2})|$, 
$|\Re(q^{1/2})|$ respectively. There exists a unique {\it positively oriented frame }$\{S,T\}$ of the 
tangent bundle $TM$ on $M\setminus\Sigma_q$, orthonormal with respect to $R_q$, such that $S$, $T$ are
tangent to the horizontal, respectively the vertical, foliation of $q$. The vector fields $S$, $T$ are
smooth on $M\setminus\Sigma_q$ and blow up at $\Sigma_q$. Let 
$$\eqalign{ \eta_T:&=\Re(q^{1/2})=-\imath_T\omega_q\,\, , \cr
            \eta_S:& =\Im(q^{1/2})=\imath_S\omega_q\,\,. \cr} \eqnu$$ 
Then $\omega_q=\eta_T\wedge\eta_S$ and, in addition, since $\eta_S$ and $\eta_T$ are closed $1$-forms,
the vector fields $S$, $T$ preserve the area form $\omega_q$. 

Let $L^2_q(M):=L^2(M,\omega_q)$ be the Hilbert space of square-integrable functions with respect to
the area form. The {\it Cauchy-Riemann operators} $\,\partial^{\pm}_q:=S\pm i\,T\,$ are closed on the domain
$H^1(M)\subset L^2_q(M)$ (the Sobolev space of weakly differentiable functions with $L^2$ derivatives on the
compact manifold $M$) and have closed ranges $R^{\pm}_q$ of finite codimension. In fact, the adjoints of the
Cauchy-Riemann operators satisfy the identity $(\partial^{\pm}_q)^{\ast}=-\partial^{\mp}_q$, hence,  by
Hilbert space theory, $R^{\pm}_q$ are the orthogonal complements of the kernels ${\Cal M}^{\mp}_q\subset
L^2_q(M)$, which consist of anti-meromorphic, respectively meromorphic, $L^2_q$ functions (with poles at
$\Sigma_q$) [18, Prop. 3.2]. There exist therefore two orthogonal decompositions:
$$ L^2_q(M)=R^-_q \oplus^{\perp} {\Cal M}^+_q =R^+_q \oplus^{\perp} {\Cal M}^-_q \eqnu $$
and two associated (orthogonal) projections $\pi^{\pm}_q:L^2_q(M)\to {\Cal M}^{\pm}_q$.
The spaces ${\Cal M}^{\pm}_q$, endowed with the real Euclidean structure induced by the Hilbert structure on
$L^2_q(M)$, are {\it isomorphic }to the cohomology $H^1(M_q,{\Bbb R})$, endowed with the {\it Hodge inner
product }determined by the metric $R_q$ (by the area form $\omega_q$). In fact, the ${\Bbb R}$-linear maps 
$c^{\pm}_q:{\Cal M}^{\pm}_q\to H^1(M_q, {\Bbb R})$ defined by 
$$\eqalign{ & c^+_q(m^+):= [\Re(m^+q^{1/2})]\,\,, \cr
            & c^-_q(m^-):= [\Re(m^-{\overline q}^{1/2})] \,\,\cr} \eqnu $$
are isomorphisms of vector spaces. In addition
$$|\!|c^{\pm}_q(m^{\pm})|\!|_q^2:=\int_M c^{\pm}_q(m^{\pm})\,\wedge\,\ast c^{\pm}_q(m^{\pm}) \,=\,
|m^{\pm}|_0^2 \,\,, \eqnu$$    
where $|\!|\cdot|\!|_q$ is the {\it Hodge norm }on $H^1(M_q,{\Bbb R})$ (the Hodge $\ast$-operator on $M_q$ is given 
by ${\ast}\eta_T=\eta_S$, ${\ast}\eta_S=-\eta_T$) and $|\cdot|_0$ the norm on $L^2_q(M)$. These assertions follow
from the remark that, if $h^{\pm}$ is any holomorphic, respectively anti-holomorphic, differential on
$M_q$, then $h^+/q^{1/2}$, $h^-/{\bar q}^{1/2}$ are respectively a meromorphic, an anti-meromorphic
function with poles at $\Sigma_q$ which belong to $L^2_q(M)$. The fact that all real cohomology classes can be
represented as the real part of a holomorphic (or anti-holomorphic) differential is a classical property of
Riemann surfaces [13, III.2]. A computation similar to $(2.5)$ shows that the {\it symplectic form }on
$H^1(M_q,{\Bbb R})$, given by the wedge product, can be written in ${\Cal M}^{\pm}_q$ as follows:
$$\eqalignno{& c^{\pm}_q(m_1^{\pm}) \wedge c^{\pm}_q(m_2^{\pm})= \Im
(m^{\pm}_1,m^{\pm}_2)_q \,\,, &\speqnu{2.5'}\cr
\noalign{\noindent 
where $(\cdot,\cdot)_q$ is the inner product in $L^2_q(M)$.}}
$$ 

By $(2.3)$ any function $u\in L^2_q(M)$ has a 
decomposition $u=\partial^+_q v + \pi^-_q(u)$ with $v\in H^1(M)$. 
The operator $U_q:L^2_q(M)\to L^2_q(M)$ given by
$$\eqalign{U_q(u):&=\partial^-_q v - {\overline{\pi^-_q(u)}}\,\,, \cr
           \hbox{ if }\,\,\,\,u&= \partial^+_q v + \pi^-_q(u) \cr} \eqnu $$
is a ${\Bbb R}$-linear isometry. In fact, by [18, Prop. 3.2] and the
orthogonality of the decomposition $(2.3)$, we have
$$|U_q(u)|_0^2=|\partial^-_q v|_0^2+|\pi^-_q(u)|_0^2 = |\partial^+_q
v|_0^2+|
\pi^-_q(u)|_0^2=|u|_0^2\,\,.
\speqnu {2.6'}$$

 Let $q_t:=G_t(q)$ denote the orbit of the quadratic differential 
$q\in Q^{(1)}_{\kappa}$ under the Teichm\"uller geodesic flow $G_t$. By the 
definition of the Teichm\"uller flow, $q_t$ is given by the equations:
$$\eqalign{ &\eta_T(t):=\Re(q_t^{1/2})=e^t\,\Re(q^{1/2})=e^t\eta_T\,\, , \cr
&\eta_S(t):=\Im(q_t^{1/2})=e^{-t}\,\Im(q^{1/2})=e^{-t}\eta_S\,\,. \cr}
\eqnu $$
In all the arguments which follow it is crucial that the area form
$\omega_q$ of the metric $R_q$ and, consequently, the inner product
of the Hilbert space $L^2_q(M)$, are invariant under the Teichm\"uller
geodesic flow. Let $\{S_t,T_t\}$ and $\partial^{\pm}_t:=S_t\pm i\,T_t$ be 
the orthonormal system and the Cauchy-Riemann operators determined by $q_t$. We have:
$$\partial^{\pm}_t=S_t\pm i\,T_t = e^{-t}S \pm i\,e^t T \,\,.
\speqnu{2.7'}$$
Let ${\Cal M}^{\pm}_t:=N(\partial^{\pm}_t) \subset L^2_q(M)$ denote 
the subspaces of meromorphic, respectively anti-meromorphic, $L^2_q$
functions on the Riemann surface $M_t$ (endowed with the complex 
structure determined by $q_t$). The Kontsevich-Zorich cocycle
$G_t^{KZ}$ is described, along the orbit $q_t:=G_t(q)$ of the 
Teichm\"uller flow, by the following variational formula:

\proclaim{Lemma} The ordinary differential equation 
$$u'=U_{q_t}(u) \eqnu $$ 
is well\/{\rm -}\/defined in $L^2_q(M)$ and satifies the following properties{\rm .}
\medbreak
\item {\rm (1)} Solutions of the Cauchy problem for  {\rm (2.8)}  exist for all 
times and are uniquely determined by the initial condition{\rm .}
\item {\rm (2)} If $u_t\in L^2_q(M)$ is any solution of  {\rm (2.8)}  with initial
condition $u_0\in {\Cal M}^+_q$, then $u_t \in {\Cal M}^+_t${\rm ,} for all 
$t\in {\Bbb R}${\rm .}
\item {\rm (3)} Let $m^+_t\in {\Cal M}^+_t$ be the unique solution of  {\rm (2.8)} 
with initial condition $m^+_0=m^+\in {\Cal M}^+_q${\rm .} For all $t\in {\Bbb
R}${\rm ,} we have

$$ G_t^{KZ}\bigl( c^+_q(m^+)\bigr) \,=\, c^+_{q_t}(m^+_t)\,\,. \speqnu
{2.8'}$$ 
\endproclaim

\demo{{P}roof of {\rm (1)}} The O.D.E. $(2.8)$ is well defined since the Hilbert 
space $L^2_q(M)$ is invariant along orbits of the Teichm\"uller flow.
The function $F(t,u):=U_{q_t}(u)$ is uniformly Lipschitz in the second
coordinate, since $U_q$ is an isometry, and smooth on ${\Bbb R}\times
L^2_q(M)$. In fact, the Cauchy-Riemann operators
$\partial^{\pm}_t$, densely defined on the Hilbert space $L^2_q(M)$
with common domain $H^1(M)$, depend smoothly on $t\in {\Bbb R}$, as 
it can be seen by explicitly writing them in terms of the fixed orthonormal 
frame $\{S,T\}$. The standard local existence and uniqueness results for 
O.D.E. on Banach spaces therefore apply [40, Chap. IV]. Since for any 
solution $|u'_t|_0=|u_t|_0$, it follows that 
$$|u_t|_0 \leq |u_0|_0 + \int_0^t |u_s|_0\,ds \,\,, \eqnu $$
hence by Gronwall's lemma $|u_t|_0$ does not blow up in finite time.
Local solutions can be therefore continued for all times. 
 \enddemo

\demo{{P}roof of {\rm (2)}}
A computation shows that the flow of $(2.8)$ preserves meromorphic
functions. In fact, since by $(2.7)$ and $(2.7')$
$${d\over{dt}}\,q_t^{1/2} = {\overline{q_t}}^{1/2}\,\,,\,\,\,\, 
{d\over{dt}}\,\partial^{\pm}_t =
-\partial^{\mp}_t \,\,, \eqnu $$
it follows that, if $u_t$ is a solution of the O.D.E. $(2.8)$,
we have:
$$ {d\over{dt}} (\partial^+_t u_t)= -\partial^-_t(u_t)+
\partial^+_t(u'_t)=0\,\,, 
\eqnu $$
in the weak sense in $L^2_q(M)$. 

\demo{Proof of {\rm (3)}} Let $\pi^{\pm}_t:L^2_q(M) \to {\Cal M}^{\pm}_t
\subset L^2_q(M)$ be the orthogonal projections. If $m^+_t\in
{\Cal M}^+_t$ is a solution of $(2.8)$, there exists a function
$v_t\in H^1(M)$ (unique up to additive constants) such that:
$$\eqalign{   m^+_t &= \partial^+_t v_t \,+\, \pi^-_t(m^+_t) \,\,,\cr
{d\over{dt}} m^+_t &=\partial^-_t v_t \,-\,{ \overline{\pi^-_t(m^+_t)} }
\,\,. \cr} \eqnu$$
If $v_t$ is chosen with zero average for all $t\in {\Bbb R}$ then the
map $t\to v_t\in H^1(M)$ is smooth. By $(2.12)$  
$${d\over{dt}}\Re(m^+_t q_t^{1/2})= \Re\bigl(({d\over{dt}}m^+_t +
{\overline{m^+_t}})q_t^{1/2}\bigr)=2\Re(dv_t)\,\equiv 0\,\in
H^1(M,{\Bbb R})\,.\eqnu$$
Hence $(2.8')$ follows by the definition of the Kontsevich-Zorich
cocycle $G^{KZ}_t$. In fact, over the Teichm\"uller space, the
cocycle acts as the identity on cohomology classes.
\enddemo
 \pagebreak

  Let $c\in H^1(M_q,{\Bbb R})$ and $c_t:=G^{KZ}_t(c)$. Let
$c_t=[\Re(m^+_t q_t^{1/2})]$, $m^+_t \in {\Cal M}^+_t$. By Lemma
2.1, $m^+_t$ is a solution of $(2.8)$ or, equivalently, of $(2.12)$.
 
\specialnumber{2.1'}
\proclaim{Lemma} The evolution of the Hodge norm under the Konsevich\/{\rm -}\/Zorich 
cocycle is described by the following formulas\/{\rm :}\/
$$\eqalign{ {d\over{dt}}|m^+_t|_0^2 &= -2\Re B_q(m^+_t):= 
-2\Re \int_M (m^+_t)^2\omega_q \,\, , \cr
             {d^2\over{dt^2}}|m^+_t|_0^2 &= 4\{|\pi^-_t(m^+_t)|_0^2 
-\Re \int_M (\partial^+_t v_t)
               (\partial^-_tv_t)\omega_q \}\,\,.\cr} \eqnu $$
\endproclaim

\demo{Proof} By $(2.5)$ and the invariance of the $L^2_q$ inner
product under the Teichm\"uller flow, $|\!|c_t|\!|=|m^+_t|_0$. A
computation, based on $(2.12)$, shows that
$$\eqalign{ {d\over{dt}}|m^+_t|_0^2&=2\Re (m^+_t, {d\over{dt}}m^+_t)_q=
-\Re(m^+_t,{\overline{\pi^-_t(m^+_t)}})_q  \cr 
&=-2\Re (m^+_t,\overline{m^+_t})_q=-2\Re\int_M (m^+_t)^2\omega_q \,\,,\cr} 
\eqnu$$
$$\align {d\over{dt}}B_q(m^+_t)&=2\int_M m^+_t\,{d\over{dt}}m^+_t\,\omega_q \cr
&=2\int_M \bigl(\partial^+_tv_t+\pi^-_t(m^+_t)\bigr)\bigl(\partial^-_tv_t-{\overline{\pi^-_t(m^+_t)}}
\bigr)\omega_q\speqnu{2.15'} \cr &= -2|\pi^-_t(m^+_t)|_0^2+2\int_M (\partial^+_t v_t)(\partial^-_tv_t)\omega_q\,\,. 
\cr
\noalign{\vskip-24pt}
\endalign $$
 \enddemo 
\vglue12pt
Calculations of first and second variation 
formulas for various functionals on the Teichm\"uller 
space are widespread in Teichm\"uller theory (see for instance [64], [72], [60], [61]). 
In particular, M.Taniguchi [60] proved variational formulas for the extremal length of the homology class of 
a given curve (equivalently, for the Hodge norm of the 
dual cohomology class), with respect to fairly general 
quasi-conformal deformations. In the special case of Teichm\"uller deformations, not explicitly considered in [60], [61], Taniguchi's variational formulas agree, by Lemma 2.1 in the case of the second variation, 
with $(2.14)$ (and $(3.12)$ below). We thank the referee 
and M.Wolf for bringing the above mentioned references 
to our attention.

 Lemma 2.1$'$ immediately implies Veech's result, recalled in
Remark 1.2. In fact, we will obtain an upper bound for the second Lyapunov 
exponent of the Kontsevich-Zorich cocycle in terms of the function 
$\Lambda^+: {\Cal M}^{(1)}_{\kappa}/SO(2,{\Bbb R}) \to {\Bbb R}$, defined as follows:
$$\Lambda^+(q):= \max\left\{ {{|B_q(m^+)|}\over {|m^+|_0^2}} \,\,\Big|\,\,m^+\in {\Cal M}^+_q\setminus\{0\}\,,
\,\,\int_M m^+\,\omega_q=0\,\right\}\,\,. \eqnu $$
\specialnumber{2.2}
\proclaim {Corollary} Let $\mu$ be any $G_t$\/{\rm -}\/ergodic probability measure on ${\Cal M}^{(1)}_{\kappa}${\rm .} The
second Lyapunov exponent of the Kontsevich\/{\rm -}\/Zorich cocycle with respect to the measure $\mu$ satisfies the
following inequality\/{\rm :}\/
$$\lambda_2^{\mu}\leq \int_{{\Cal M}^{(1)}_{\kappa}} \Lambda^+ \,d\mu \,< \,\lambda_1^{\mu}= 1\,\,. 
\eqnu $$
\endproclaim

\demo{Proof} The Kontsevich-Zorich cocycle has two $1$-dimensional invariant continuous sub-bundles 
$E_{\pm 1}\subset {\Cal H}_{\kappa}^1(M,{\Bbb R})$ with Lyapunov exponents $\pm 1$. At any $q\in {\Cal M}^
{(1)}_{\kappa}$
$$\eqalign{  E_1(q) & :={\Bbb R}\cdot [\Im(q^{1/2})]\,\,,\cr
             E_{-1}(q) & :={\Bbb R}\cdot [\Re(q^{1/2})]\,\,. \cr} \eqnu $$
It follows in particular that, for all $G_t$-ergodic probability
measures $\mu$, $\lambda_1^{\mu}=1$, $\lambda_{2g}^{\mu}=-\lambda_1^{\mu}
=-1$. Let $E_0\subset {\Cal H}_{\kappa}^1(M,{\Bbb R})$ be the continuous 
invariant sub-bundle defined as follows. At any $q\in {\Cal M}^{(1)}
_{\kappa}$,
$$E_0(q):=\{ c\in H^1(M_q,{\Bbb R})\,|\,c\wedge [\Re(q^{1/2})]=c\wedge [\Im(q^{1/2})]=0\} \,\,. 
\speqnu {2.18'}$$
A cohomology class $c=[\Re(m^+q^{1/2})]\in E_0(q)$ if and only if $\int_M m^+\omega_q=0$. There is an
invariant continuous splitting ${\Cal H}_{\kappa}^1(M,{\Bbb R})= E_1\oplus E_{-1} \oplus E_0$. In 
addition, the bundle $E_0$ is invariant under the action of the circle group  $SO(2,{\Bbb R})$ on 
${\Cal M}^{(1)}_{\kappa}$.  Let $c\in E_0(q)$. Let $q_t:=G_t(q)$ and $c_t:=G^{KZ}_t(c)$. Then $c_t=
[\Re(m^+_tq_t)^{1/2}]$, where $m^+_t\in {\Cal M}^+_t$ has zero average. Hence, by Lemma $2.1'$,
$${d\over{dt}}\log |m^+_t|_0= -{{\Re B_q(m^+_t)}\over {|m^+_t|_0^2}}\leq \Lambda^+(q_t)\,\,. 
\eqnu $$ 
By taking time averages, it follows that
$${1\over {\Cal T}} \log {{|m^+_{\Cal T}|_0}\over {|m^+|_0}} \leq {1\over {\Cal T}}\int_0^{\Cal T}
\Lambda^+(q_t)\,dt \,\,. \speqnu{2.19'}$$
By Oseledec's theorem \ref\Os, [32, Th. S.2.9], and Birkhoff's ergodic theorem (for the ergodic measure $\mu$), by taking the limit in
$(2.19')$ as ${\Cal T}\to +\infty$, we obtain the first of the two
inequalities in $(2.17)$. In order to
prove the remaining inequality we will show that $\Lambda^+(q)<1$ for {\it all }$q\in {\Cal M}_{\kappa}
^{(1)}$. By Schwartz inequality, we have:
$$|B_q(m^+)|=|(m^+,{\overline {m^+}})_q|\leq |m^+|_0^2\,\,, \eqnu$$
where equality holds if and only if there exists $\lambda\in {\Bbb C}$ ($|\lambda|=1$) such that 
${\overline{m^+}}=\lambda\,m^+$. However, in that case $m^+$ would be meromorphic and anti-meromorphic, 
hence constant, and, having zero average, $m^+=0$. It follows that, if $m^+\not=0$, the inequality in 
$(2.20)$ is strict, hence, by compactness of the unit sphere in ${\Cal M}^+_q\subset L^2_q(M)$, 
$\Lambda^+(q)<1$. 
\enddemo

 \section{A lower bound for the second Lyapunov exponent}   

 The Teichm\"uller space $Q^{(1)}_g$ is foliated by the orbits of the action of ${\rm SL}(2,{\Bbb 
R})\subset {\rm GL}_+(2,{\Bbb R})$ described in Section~1.  The leaves are $3$-dimensional and are isometric to
the unit cotangent bundle $S^{\ast}(D)\equiv {\rm SL}(2,{\Bbb R})$ of the Poincar\'e disk $D:=D(0,1)$ 
(endowed with the hyperbolic metric $|dz|/(1-|z|^2)$ of constant curvature $-4$). The circle group
${\rm SO}(2,{\Bbb R})$ is a subgroup of ${\rm SL}(2,{\Bbb R})$ and acts therefore on $Q^{(1)}_g$. The
quotient $Q^{(1)}_g/{\rm SO}(2,{\Bbb R})$ is then foliated by $2$-dimensional leaves isometric to the
Poincar\'e disk.  In fact, given any $q\in Q^{(1)}_g$, there exists an isometric embedding
$j_q:D\to Q^{(1)}_g/{\rm SO}(2,{\Bbb R})$ such that $j_q(0)= [q]\in Q^{(1)}_g/{\rm SO}(2,{\Bbb R})$. Let
$q_{\theta}:=e^{\imath\theta}q$ be the rotation of the quadratic differential $q$ by an angle
$\theta\in [0,2\pi]$ and, if $z=re^{\imath\theta}$, let 
$$q_z := G_t(q_{\theta})\,\,,\,\,\,\,t:={1\over 2}\log{{1+r}\over{1-r}}\,\,.\eqnu$$ 
The embedding $j_q$ is given by
$$j_q(z):= [q_z]\in Q^{(1)}_g/{\rm SO}(2,{\Bbb R})\,\,, \,\,\,\,z\in D\,\,.  \speqnu{3.1'}$$ 
Let $\pi_g:Q^{(1)}_g\to T_g$ be the canonical projection onto the Teichm\"uller space $T_g$ of
(marked) Riemann surfaces, defined by $\pi_g(q)=M_q$, the Riemann surface carrying the holomorphic
quadratic differential $q\in Q^{(1)}_g$. The composed map $\pi_g\circ j_q$ is an isometric
embedding and its image is called a {\it Teichm\"uller disk} [46, \S 2.6.5]. Since the
action of ${\rm GL}_+(2,{\Bbb R})$ commutes with the action of the mapping class group $\Gamma_g$, the
foliation of the space $Q^{(1)}_g/{\rm SO}(2,{\Bbb R})$ by Teichm\"uller disks projects to a foliation
${\Cal T}_g$ of ${\Cal M}^{(1)}_g/{\rm SO}(2,{\Bbb R})$. Any leaf ${\Cal T}$ of ${\Cal T}_g$ is a
complex curve endowed with a hyperbolic metric of constant curvature $-4$.  It is in fact the
quotient of the Poincar\'e disk by a discrete group of isometries. The generic leaf is isometric
to the Poincar\'e disk [38, \S 7]. The foliated unit cotangent bundle $S^{\ast}({\Cal
T}_g)$ gives a foliation of the moduli space ${\Cal M}^{(1)}_g$, which can be naturally identified
with the orbit foliation of the action of ${\rm SL}(2,{\Bbb R})$. The Teichm\"uller geodesic flow
preserves the foliation $S^{\ast}({\Cal T}_g)$ (since the latter is the orbit foliation of the
full group ${\rm SL}(2,{\Bbb R})$) and its restriction to each leaf $S^{\ast}({\Cal T})$ coincides with
the geodesic flow of the Poincar\'e metric on ${\Cal T}$. The foliations ${\Cal T}_g$ and
$S^{\ast}({\Cal T}_g)$ preserve each stratum ${\Cal M}^{(1)}_{\kappa}/{\rm SO}(2,{\Bbb R})$ and ${\Cal
M}^{(1)}_{\kappa}$ respectively.  

In [38] M. Kontsevich and A. Zorich outlined the
proof of a formula for the sum $\lambda_1+...+\lambda_g$ of the non-negative Lyapunov exponents of
the Kontsevich-Zorich cocycle, based on the computation on a Teichm\"uller disk of the (Poincar\'e)
Laplacian of the logarithm of the volume on a generic Lagrangian subspace of the fiber
$H^1(M,{\Bbb R})$. A new version of their formula will be proved in Section 5. Following [38],
the basic tool in proving bounds for the Lyapunov exponents of the Kontsevich-Zorich cocycle is given
by the elementary lemma below, which does not rely on Sullivan's results on the Brownian motion (on
the Poincar\'e disk) [59] exploited in [38]. We remark that a similar idea of using
sub-harmonic functions to prove lower bounds for Lyapunov exponents was introduced, in a different
context, by M. Herman [26].  

\proclaim{Lemma} Let $\triangle_h$ be the {\rm (}\/hyperbolic\/{\rm )}
Laplacian on the Poincar{\rm \'{\it e}} disk $D$ and let $L:D\to {\Bbb R}$ be a smooth solution of the Poisson
equation $\triangle_h L=\Lambda${\rm ,} where $\Lambda$ is a smooth bounded function on $D${\rm .} Then 
$${1\over {2\pi}} {{\partial} \over {\partial t}}\int_0^{2\pi} L(t,\theta)\,d\theta = {1\over
2}\tanh(t)\, {1\over {|D_t|}} \int_{D_t}\Lambda\,\omega_P \,\,,\eqnu $$ 
where $(t,\theta)$ are geodesic polar coordinates on the Poincar{\rm \'{\it e}} disk{\rm ,} $|D_t|$ denotes the 
Poincar{\rm \'{\it e}} area of the disk $D_t$ of geodesic radius $t>0${\rm ,} centered at the origin{\rm ,} and 
$\omega_P$ is the Poincar{\rm \'{\it e}} area form{\rm .}  \endproclaim 

\demo{Proof}  In geodesics polar coordinates centered at the origin, the Poincar\'e Laplacian 
can be written as follows [24, Intr. \S 4.2]:  
$$ \triangle_h={{\partial^2} \over {\partial t^2}} + 2 \coth(2t) {{\partial}\over {\partial t}} 
+ {4\over{\sinh^2(2t)}} {{\partial^2} \over {\partial \theta^2}}\,\,.  \eqnu $$ 
Let $L_r$, $\Lambda_r$ be the circular averages of the functions $L$, $\Lambda$ respectively:  
$$ L_r(t):={1\over {2\pi}}\int_0^{2\pi} L(t,\theta)\,d\theta\,\,\,\,,\,\,\,\,\,\, \,\,
\Lambda_r(t):= {1\over{2\pi}}\int_0^{2\pi}\Lambda(t,\theta)\,d\theta\,\,.  \eqnu $$ 
By averaging the equation $\triangle_hL =\Lambda$, we obtain the following O.D.E.:  
$$ {{\partial^2} \over {\partial t^2}}L_r + 2 \coth(2t) {{\partial} \over {\partial t}}L_r= 
\Lambda_r \,\,.  \speqnu{3.4'}$$ 
The equation $(3.4')$ can be explicitly solved by setting $M_r:=\partial L_r/\partial t$.  Since
$M_r(0)=0$ the solution has the following expression:  
$$M_r(t)= {1\over {\sinh(2t)}} \int_0^t\Lambda_r(s)\sinh(2s)ds\,\,.  \eqnu $$ 
We recall the following elementary formulas from hyperbolic geometry:  
$$\omega_P= {1\over 2} \sinh(2t) dt\,d\theta \,\,\,\,,\,\,\,\,\,\,\,\,
|D_t|= \pi {{\cosh(2t)-1}\over 2} \,\,.  \eqnu $$ 
Hence, by $(3.5)$ and $(3.6)$,
\medbreak
\noindent (3.7) \hfill${\displaystyle {{\partial} \over {\partial t}}L_r = {{\cosh(2t)-1} \over {2\sinh(2t)}}\,{1\over {|D_t|}}
\int_{D_t}\Lambda\,\omega_P = {1\over 2}\tanh(t)\,{1\over {|D_t|}} \int_{D_t}\Lambda\,\omega_P
\,\,.  }$
\enddemo
 \advance\eqcount by 1

Let $q\in Q^{(1)}_{\kappa}$, $M_z$ be the Riemann surface carrying the quadratic 
differential $q_z$, given by $(3.1)$, and let ${\Cal M}^{\pm}_z$ be the spaces of meromorphic, 
respectively anti-meromorphic, functions on $M_z$ which belong to $L^2_q(M)$.  We remark that the 
$L^2$ structure induced by $q_z$ is independent of $z\in D$, since the area form of a quadratic 
differential is invariant under the action of ${\rm SL}(2,{\Bbb R})$.  Let $\pi^{\pm}_z:L^2_q(M)\to 
{\Cal M}^{\pm}_z$ denote the orthogonal projections.  Let $c\in H^1(M_q,{\Bbb R})$. There exists 
a uniquely determined meromorphic function $m^+_z\in {\Cal M}^+_z$ such that $c=[\Re(m^+_zq_z^{1/2})]$.  
The function $z\to m^+_z\in L^2_q(M)$ is
holomorphic on $D$, but it depends on the choice of the base point $q\in Q^{(1)}_{\kappa}$. The
norm $|m^+_z|_0$ is invariant under the action of the circle group on $Q^{(1)}_{\kappa}$, hence
it gives a well defined real analytic function on the Teichm\"uller disk $j_q(D)\subset
Q^{(1)}_g/{\rm SO}(2,{\Bbb R})$, independent of the choice of the base point.  

\proclaim{Lemma} The following formula holds\/{\rm : }\/ 
$$ \triangle_h \log |m^+_z|_0= 4 {{|\pi^-_z(m^+_z)|_0^2} \over {|m^+_z|_0^2}}-2{{|B_q(m^+_z)|^2} 
\over {|m^+_z|_0^4}} \geq 2 {{|\pi^-_z(m^+_z)|_0^2} \over {|m^+_z|_0^2}}\,\,.  \eqnu $$ 
\endproclaim 

\demo{Proof} If $[q']\in j_q(D)$, then $j_{q'}(D)=j_q(D)$.  Hence it suffices to prove the formula
at the origin $z=0$.  Let $c=[\Re(m^+q^{1/2})] \in H^1(M,{\Bbb R})$, where $m^+\in {\Cal M}^+_q$. 
Let $q_{\theta}:= e^{\imath\theta}q$ and $q_{(t,\theta)}:=G_t (q_{\theta})$. Let $X_{\theta}$ be the 
normalized directional derivative at $z=0$ defined, for any smooth function $\phi:D\to {\Bbb R}$, 
by 
$$X_{\theta}\phi:= {d\over {dr}}\phi (re^{\imath\theta})\,|_{r=0}\,\,.  \eqnu $$ 
By Lemma $2.1'$ we have:
$$\eqalign{X_{\theta}|m^+_z|_0^2&=-2\Re \left[e^{-\imath\theta}B_q(m^+)\right]=-2\Re \left[e^{-\imath\theta}
\int_M(m^+)^2\omega_q\,\right]\,\,,\cr X_{\theta}^2 |m^+_z|_0^2 &= 4\left\{|\pi^-_q(m^+)|_0^2 
-\Re \left[e^{-2\imath\theta}\int_M (\partial^+_q v) (\partial^-_qv)\omega_q\,\right]\right\}\,\,.\cr} \eqnu $$ 
In fact, let ${\Cal M}^{\pm}_{\theta}\subset L^2_q(M)$ be the kernels of the Cauchy-Riemann 
operators $\partial^{\pm}_{\theta}$ determined by the quadratic differential $q_{\theta}$. Let
$\pi^{\pm}_{\theta}: L^2_q(M)\to {\Cal M}^{\pm}_{\theta}$ be the orthogonal projections.  Any
cohomology class $c\in H^1(M,{\Bbb R})$ can be represented as $c=[\Re(m^+_{\theta}q^{1/2}_{\theta})]$, 
where $m^+_{\theta}\in {\Cal M}^+_{\theta}$.  We have:
$$\eqalign{ \partial^{\pm}_{\theta} &= e^{\pm \imath{{\theta} \over 2}} \partial^{\pm}_q \,\,, \cr
\pi^{\pm}_{\theta}&= \pi^{\pm}_q\,\,, \cr m^+_{\theta} &= e^{-\imath{{\theta} \over 2}} m^+\,\,,\cr}
\eqnu $$ 
hence, in particular, ${\Cal M}^{\pm}_{\theta}={\Cal M}^{\pm}_q$ and, if $m^+=\partial^+_q v + 
\pi^-_q(m^+)$, $$ \eqalign { m^+_{\theta}&= \partial^+_{\theta} v_{\theta}+\pi^-_{\theta}(m^+_{\theta})
\,\,, \cr v_{\theta}&= e^{-\imath\theta} v\,\,.\cr} \speqnu{3.11'}$$ 
The formulas $(3.10)$ follow therefore from the formulas $(2.14)$, taken at $t=0$, for the quadratic
differential $q_{\theta}$.  

 The hyperbolic gradient and the hyperbolic Laplacian can be written in terms of the 
directional derivatives $X_{\theta}$. Hence, by $(3.10)$, at $z=0$, 
$$\eqalign{\nabla_h|m^+_z|_0^2&=(X_0|m^+_z|_0^2,X_{\pi\over 2}|m^+_z|_0^2)=
-2\bigl(\Re B_q(m^+),\Im B_q(m^+))\,\,,\cr \triangle_h|m^+_z|_0^2&={1\over {\pi}}\int_0^{2\pi}
X_{\theta}^2|m^+_z|_0^2\,d\theta= 8|\pi^-_q(m^+)|_0^2 \,\,.  \cr} \eqnu $$ 
Since, given any positive smooth function $\phi$ on $D$, 
$$\triangle_h\log\phi={1\over 2}\left\{{{\triangle_h(\phi^2)}\over{\phi^2}}- {{|\nabla_h\phi^2|^2}
\over{\phi^4}}\right\}\,\,, \eqnu $$ 
the identity in $(3.8)$ at $z=0$ follows from $(3.12)$ by taking $\phi(z):=|m^+_z|_0$ in
$(3.13)$.  In addition, $|B_q(m^+_z)|\leq |\pi^-_z(m^+_z)|_0 |m^+_z|_0$, since
$$|B_q(m^+_z)|=|(m^+_z,{\overline{m^+_z}})_q|=|(\pi^-_z(m^+_z),{\overline{m^+_z}})_q|\leq
|\pi^-_z(m^+_z)|_0 |m^+_z|_0 \,\,.  \eqnu $$ 
The inequality in $(3.8)$ follows from $(3.14)$ and the proof is therefore complete.  
\enddemo  

The second Lyapunov exponent of the Kontsevich-Zorich cocycle, with respect to any ${\rm SL}(2,{\Bbb R})$-invariant ergodic
measure, can be bounded from  below in terms of the non-negative function
$\Lambda^-:{\Cal M}^{(1)}_{\kappa}/{\rm SO}(2,{\Bbb R})\to {\Bbb R}$, defined as follows:  
$$\Lambda^-(q):=\min\left\{ {{|\pi^-_q(m^+)|_0^2}\over {|m^+|_0^2}}\,\,\big|\,\, m^+\in {\Cal M}^+_q\setminus
\{0\} \right\}\,\,.  \eqnu $$ 

\vglue8pt

\proclaim {Theorem} Let $\mu$ be any $G_t$\/{\rm -}\/ergodic probability measure on ${\Cal M}^{(1)}_{\kappa}$ such that for
$\mu$\/{\rm -}\/almost all orbits of the circle group the Haar {\rm (}\/Lebesgue\/{\rm )} measure is absolutely continuous
with respect to the conditional measure induced by $\mu${\rm .} The second Lyapunov exponent 
$\lambda_2^{\mu}$ of the Kontsevich\/{\rm -}\/Zorich cocycle{\rm ,} with respect to the measure $\mu${\rm ,} satisfies the 
following lower bound\/{\rm :}  
$$\lambda_2^{\mu}\,\,\geq\,\, \int_{{\Cal M}^{(1)}_{\kappa}} \Lambda^-\,d\mu\,\,.\eqnu $$
\endproclaim

\vglue8pt

\demo{Proof} Since the Lebesgue measures on the orbits of the circle action are absolutely continuous with 
respect to the conditional measures induced by $\mu$, by Fubini's theorem, Birkhoff's ergodic theorem and Oseledec's theorem, for
$\mu$-almost all $q\in {\Cal M}^{(1)}_{\kappa}$, $q_{\theta}$ is a Birkhoff generic point and a regular point (in the sense of
Oseledec's theorem) for almost all $\theta\in [0,2\pi]$. The Lyapunov exponent $\lambda_2^{\mu}$ is the top Lyapunov exponent
of the restriction of the Kontsevich-Zorich cocycle to the invariant sub-bundle $E_0$, introduced in the proof of Corollary 2.2. 
We recall that $E_0$ is invariant under the action of the circle group ${\rm SO}(2,{\Bbb R})$  and a cohomology class
$c=[\Re(m^+q^{1/2})]\in E_0(q)$ if and only if $\int_M m^+\omega_q=0$. Let
$c\in E_0(q)$. Let $q_z$ be given by $(3.1)$ and $m^+_z\in {\Cal M}^+_z$ be a (zero average)
meromorphic function in $L^2_q(M)$ such that $c=[\Re(m^+_z q_z^{1/2})]$. Let $G^{KZ}_z(c)$ be the
cohomology class $c$ at the point $q_z \in Q^{(1)}_{\kappa}$ (the parallel transport of cohomology
classes is trivial over the Teichm\"uller space $Q^{(1)}_{\kappa}$) and let $L:D\to{\Bbb R}$ be
the smooth function $L(z):=\log|\!|G_z^{KZ}(c)|\!|=\log|m^+_z|_0$.  Let $z\equiv (t,\theta)$ with
respect to geodesic polar coordinates, that is $z=re^{\imath\theta}$, $r=\tanh(t)$. By Lemmas 3.1, 3.2
and definition $(3.15)$, we have:  
$${1\over {2\pi}} {{\partial}\over {\partial t}}\int_0^{2\pi}\log|\!|G_{(t,\theta)}^{KZ}(c)|\!| 
\,d\theta \geq \tanh(t)\,{1\over{|D_t|}}\int_{D_t}\Lambda^-\bigl(G_s(q_{\theta})\bigr)\,
\omega_P(s,\theta)\,\,. \eqnu $$ 
Let $\sigma_q$ be the normalized canonical (Haar) measure on the unit sphere $E_0^{(1)}(q)$ of 
the Euclidean space $E_0(q)\subset H^1(M_q,{\Bbb R})$, endowed with the Hodge inner product 
induced by the metric $R_q$. The measure $\sigma_q$ is invariant under the action of the circle 
group ${\rm SO}(2,{\Bbb R})$ on the space $Q_g$. By averaging the inequality in $(3.17)$ over $c\in 
E_0^{(1)}(q)$ with respect to $\sigma_q$, then over the moduli space ${\Cal M}_{\kappa}^{(1)}$ 
with respect to the $G_t$-invariant measure $\mu$, we obtain:  
$$\eqalign{ {1\over {2\pi}}{{\partial}\over {\partial t}}\int_{{\Cal M}_{\kappa}^{(1)}}\int_{E_0^{(1)}(q)}
\int_0^{2\pi}\log|\!|G_{(t,\theta)}^{KZ}(c)|\!|\,&d\theta d\sigma_q d\mu \cr \geq  
&\tanh(t)\,\int_{{\Cal M}_{\kappa}^{(1)}}\Lambda^-\,d\mu\,\,.\cr} \speqnu{3.17'}$$ 
Finally, by taking the time average of $(3.17')$ over the interval $[0,{\Cal T}]$, we have 
$$\eqalign{ {1\over {2\pi{\Cal T}}} \int_{{\Cal M}_{\kappa}^{(1)}}\int_{E_0^{(1)}(q)}\int_0^{2\pi} 
\log|\!|G_{({\Cal T},\theta)}^{KZ}(c)|\!|&\,d\theta d\sigma_q d\mu  \cr \geq &{{\log\cosh{\Cal T}}
\over {\Cal T}} \int_{{\Cal M}_{\kappa}^{(1)}}\Lambda^-\,d\mu\,\,.\cr} \eqnu $$ 
Hence, by passing to the limit as ${\Cal T}\to +\infty$ in $(3.18)$, we obtain $(3.16)$. In fact, 
since $c\in E_0$, by Oseledec's theorem, 
\medbreak
\noindent ($3.18'$) \hfill ${\displaystyle \lim_{{\Cal T}\to +\infty}{1\over {2\pi{\Cal T}}}\int_{{\Cal
M}_{\kappa}^{(1)}}\int_{E_0^{(1)}(q)}
\int_0^{2\pi}\log|\!|G_{({\Cal T},\theta)}^{KZ}(c)|\!|\,d\theta d\sigma_q d\mu\,\,\leq\,\, 
\lambda_2^{\mu} \,\,.  }$ 
\enddemo  
\vglue12pt
Proving that $\lambda_2^{\mu}>0$ requires a deeper analysis, carried out
in the next section, of the projection operator $\pi^-_q:{\Cal M}^+_q\to {\Cal M}^-_q$, which
enters crucially in the lower bound given by Theorem 3.3, as the quadratic differential varies
over the moduli space.\pagebreak

 \section{The determinant locus}

 Let $M$ be a marked Riemann surface of genus $g\geq 2$ and let $\{a_1,\ldots ,
a_g,\break b_1,\ldots ,b_g\}$ be a {\it canonical homology
basis} [13, III.1] giving the marking. Let $\{\theta_1,\ldots ,\theta_g\}$ be the dual basis of holomorphic differentials on $M$, that is the
basis characterized by the property that $\theta_j(a_k)=\delta_{jk}$. The $g\times g$ complex matrix $\Pi_{ij}:=\theta_i(b_j)$ is called
the {\it period matrix }of the marked Riemann surface $M$. The period matrix defines a holomorphic map $\Pi:T_g\to \Sigma_g$ on the
Teichm\"uller space 
$T_g$ with values in the Siegel space $\Sigma_g$ of symmetric complex matrices with positive definite imaginary part. Let $q\in Q_g$ be a 
holomorphic quadratic differential on $M$. Let $\mu_q:=|q|/q$ be the canonical Beltrami differential 
on $M$ associated to $q$, which represents the deformation of the complex structure of the Riemann
surface $M_q$ in the direction determined by the Teichm\"uller flow at $(M_q,q)$. The equation
$$ \det \bigl({{d\Pi}\over {d\mu_q}})=0 \eqnu $$
defines a real analytic hypersurface $D_g$, of real codimension equal to $2$, of the Teichm\"uller 
space $Q_g$. In addition, $D_g$, is invariant under the action of the mapping class group $\Gamma_g$. 
The quotient ${\Cal D}_g:=D_g/\Gamma_g$ is therefore a real analytic hypersurface of codimension $2$ 
of the moduli space ${\Cal M}_g$, which will be called the {\it determinant locus}. We remark that, 
since $(4.1)$ is invariant under the multiplicative (holomorphic) action of the group ${\Bbb C}^{\ast}
:={\Bbb C}\setminus\{0\}$ on the space $Q_g$, the quotients $D_g/{\Bbb C}^{\ast}$, ${\Cal D}_g/
{\Bbb C}^{\ast}$ are real analytic hypersurfaces of $Q^{(1)}_g/{\rm SO}(2,{\Bbb R})$, ${\Cal M}^{(1)}_g/
{\rm SO}(2,{\Bbb R})$ respectively. Let $D_{\kappa}:=D_g\cap Q_{\kappa}$, $D^{(1)}_{\kappa}:=D_g\cap 
Q^{(1)}_{\kappa}$, ${\Cal D}_{\kappa}:=D_{\kappa}/\Gamma_g$ and ${\Cal D}^{(1)}_{\kappa}:=
D^{(1)}_{\kappa}/\Gamma_g$.

The answer to the question whether the lower bound established in Theorem 3.3 is {\it strictly }positive depends on the geometry of
the determinant locus. In fact, the answer will be affirmative if the support of the probability 
measure $\mu$, ergodic with respect to the Teichm\"uller flow on ${\Cal M}^{(1)}_{\kappa}$, is not 
contained in ${\Cal D}^{(1)}_{\kappa}$. Let $q\in {\Cal M}^{(1)}_{\kappa}$ and $\pi^-_q:L^2_q(M)\to 
{\Cal M}_q^-$ be the orthogonal projection onto the subspace of anti-meromorphic functions. Let
$$ 1\equiv \Lambda_1(q)\geq \Lambda_2(q) \geq \cdots  \geq \Lambda_g(q) \geq 0 \eqnu $$ 
be the eigenvalues of the non-negative hermitian form $H_q:{\Cal M}_q^+\times {\Cal M}_q^+ \to {\Bbb R}$
defined by 
$$H_q(m^+_1,m^+_2):=\bigl(\pi^-_q(m^+_1), \pi^-_q(m^+_2) \bigr)_q\,\,. \speqnu{4.2'}$$
The functions $\Lambda_1,\ldots ,\Lambda_g$ are continuous on ${\Cal M}^{(1)}_{\kappa}$ and invariant under the
action of ${\rm SO}(2,{\Bbb R})$.

\proclaim{Lemma} The zero set $\{\Lambda_g=0\}={\Cal D}^{(1)}_{\kappa}${\rm .} \endproclaim 

\demo{Proof} Let $\{m^+_1,\ldots ,m^+_g\}$ be an orthonormal basis of ${\Cal M}_q^+$. The\break (symmetric) matrix 
$B$ of the projection operator $\pi^-_q$, with respect to the basis $\{m^+_1,\ldots ,m^+_g\}\subset {\Cal M}_q^+$ 
and $\{{\overline {m^+}}_1,\ldots ,{\overline {m^+}}_g\}\subset {\Cal M}_q^-$, and the matrix $H$ of the hermitian
form $H_q$, with respect to the basis $\{m^+_1,\ldots ,m^+_g\}$, are given by the following formulas:
$$B_{ij}=B_q(m^+_i,m^+_j)=(m^+_i,{\overline {m^+}}_j)_q\,\,\,\,,\,\,\,\,\,\,\,\, H=B^{\ast}B={\overline B}B
\,\,. \eqnu $$
Let $\{a_1,\ldots ,a_g,b_1,\ldots ,b_g\}$ be a canonical homology basis and $\{\theta_1,\ldots ,\theta_g\}$ be the dual 
basis of holomorphic differentials. The quotients $\phi^+_i:=\theta_i/q^{1/2}$ are meromorphic functions on
$M$ with poles at $\Sigma_q$, which belong to the space $L^2_q(M)$. The system $\{\phi^+_1,\ldots ,\phi^+_g\}$
is a basis of the space ${\Cal M}_q^+$. By Rauch's formula [28, Prop. A.3]:  
$${{d\Pi_{ij}}\over {d\mu_q}}=\int_M \theta_i\theta_j\,\mu_q=\int_M \phi^+_i\phi^+_j\,\omega_q=
B_q(\phi^+_i,\phi^+_j)\,\,. \eqnu $$
Since $\{\phi^+_1,\ldots ,\phi^+_g\}$ is a basis of ${\Cal M}_q^+$, there exists a non-singular $g\times g$
complex matrix $C=(c_{ij})$ such that
$$ \phi^+_i=\sum_{j=1}^g c_{ij} m^+_j\,\,\,\,,\,\,\,\,\,\,\,\, CC^{\ast}=\Im(\Pi)\,\,. \eqnu $$
In fact, by [13, III.2.3],
$$(\phi^+_i,\phi^+_j)_q={\imath\over 2}\,\int_M\theta_i\wedge {\overline {\theta_j}}={\imath\over 2}\,\sum_{k=1}^g \theta_i(a_k)
{\overline {\theta_j}}(b_k)-\theta_i(b_k){\overline {\theta_j}}(a_k)=\Im(\Pi_{ij})\,. 
\speqnu{4.5'}$$ 
By $(4.4)$ and $(4.5)$,
$$\eqalign{|\det\bigl({{d\Pi}\over {d\mu_q}}\bigr)|&= 
|\det(CBC^t)|= |\det C|^2\,|\det B| \cr
&=\det(\Im(\Pi))\,(\det H)^{1/2}=\det (\Im(\Pi))\,(\Lambda_1\cdots \Lambda_g)^{1/2}\,\,.\cr} 
\eqnu $$
Since $\Im(\Pi)$ is positive definite, $\Lambda_g(q)=0$ if and only if $q\in {\Cal D}^{(1)}_{\kappa}$.
\enddemo

\phantom{reason}

The basic idea in order to prove that the lowest eigenvalue $\Lambda_g$ in $(4.2)$ does not vanish 
identically is to study the determinant locus near the boundary of the moduli space. We briefly recall the 
methods introduced by Fay [15, III] and Masur \ref\Msone to analyze
the behaviour of holomorphic and quadratic differentials near the boundary of the moduli space. We warn the
reader that, in the author's opinion, all of the variational formulas contained in [15, III] are to some 
extent incorrect, as pointed out by A.Yamada in \ref\YaA, where the correct formulas are proved in great 
detail (compare, for instance, the formulas $(47)$, $(53)$ in [15] with, respectively, $(36)$, $(36)'$ 
and $(56)$ in \ref\YaA. However, our argument does not require {\it precise }variational formulas (only the 
dominant term in the appropriate expansions), hence it does not depend on which version, Fay's or Yamada's,
is considered to be the correct one. 

The moduli space $R_g=T_g/\Gamma_g$ can be compactified by adding Riemann surfaces with nodes
as a boundary. This compactification is explained in detail in \ref\Beone or \ref\Msone. A 
Riemann surface with nodes $M_0$ can be obtained by pinching a compact Riemann surface $M$ along disjoint
curves $\gamma_1 ,\ldots , \gamma_s$ with $1\leq s\leq\break 3g-3$. A connected component of the complement of the 
set of nodes is called a {\it part }of $M_0$; it is a finite non-singular Riemann surface with punctures. 
It is possible to assume that no part of $M_0$ is a sphere with just one or two punctures. The pinching
of each $\gamma_i$ produces two {\it paired }punctures on $M_0$, denoted by $p^{(1)}_i$, $p^{(2)}_i$.
Let $z^{(1)}_i$ and $z^{(2)}_i$ be holomorphic coordinates parametrizing disjoint punctured disks 
$U^{(1)}_i=\{0<|z^{(1)}_i|<1\}$ and $U^{(2)}_i=\{0<|z^{(2)}_i|<1\}$ in a neighbourhood of $p^{(1)}_i$ 
and $p^{(2)}_i$ respectively. There are local coordinates $(\tau_{s+1},\ldots ,\tau_{3g-3})$ in a neighbourhood 
of $0\in {\Bbb C}^{3g-3-s}$, parametrizing a neighbourhood of $M_0$ in its moduli space, such that, for each 
$M_{\tau}$ near $M_0$, the functions $z^{(1)}_i$, $z^{(2)}_i$ are still holomorphic coordinates in a
neighbourhood of the punctures. For each $(t_1,\ldots ,t_s)\in {\Bbb C}^s$ in a neighbourhood of zero,  
  consider the Riemann surface 
$M_{(t,\tau)}$ obtained by deleting all the neighbourhoods $\{0<|z^{(1)}_i|\leq |t_i|\}$ and 
$\{0<|z^{(2)}_i|\leq |t_i|\}$ from $M_{\tau}$ and gluing $z^{(1)}_i$ to $t_i/z^{(2)}_i$. If 
all $t_i\not=0$, the surface $M_{(t,\tau)}$ is a compact Riemann surface of genus $g$. The parameters
$(t,\tau)$ parametrize all Riemann surfaces in the compactified moduli space in a neighbourhood 
of $M_0$. Let $A_{t_i}:=\{|t_i|<|z^{(1)}_i|<1\}=\{|t_i|<|z^{(2)}_i|<1\}\subset M_{(t,\tau)}$ and let 
$\gamma_{t_i}$ be the curve defined by $|z^{(1)}_i|=|z^{(2)}_i|=|t_i|^{1/2}$. If $K$ is a compact 
subset of the complement of the set of the nodes in $M_0$, then for sufficiently small $(t,\tau)$ 
it is possible to regard $K$ as a compact subset of $M_{(t,\tau)}$ via the quasi-conformal map $\phi_{\tau}
:M_0\to M_{\tau}$ and the natural inclusion of $M_{(t,\tau)}\setminus\cup\{\gamma_{t_i}\}$ in~$M_{\tau}$.

There is a complex manifold $X$ of dimension $3g-2$ and a proper holomorphic map $\pi:X\to
{\Bbb C}^s\times {\Bbb C}^{3g-3-s}$ such that the fiber over $(t,\tau)$ is $M_{(t,\tau)}$. Such a
construction is explained by [15, III] or \ref\Msone, we recall it here
for the convenience of the reader. Let ${\tilde W}:=\cup M_{\tau}$
and ${\tilde Y}:= {\tilde W}\setminus \cup\phi_{\tau}(U^{(1)}_i\cup U^{(2)}_i)$, where the union is taken over $i=1,\ldots ,s$, $\tau\in U$ and
$U\subset {\Bbb C}^{3g-3-s}$ is a neighbourhood of the origin. Let 
$$W:={\tilde Y}\times D^s\,\cup\,\bigcup \,\{(p,t)\in\phi_{\tau}(U^{(k)}_i)\times D^s
\,|\,|z^{(k)}_i(p)|>|t_i|\,\}\,, \eqnu $$
where $D\subset {\Bbb C}$ denotes the open unit disk. Let then $V_i\subset {\Bbb C}^3$ be the non-singular
surface given by the equation $X_iY_i=t_i$ and let
$$X:=W\cup (V_1\times D^{s-1}\times U)\cup\cdots \cup (V_s\times D^{s-1}\times U)\,\,, \speqnu{4.7'}$$
where in the overlap $(p,t)$ is identified with 
$$\eqalign{ &(z^{(1)}_i(p), t_i/z^{(1)}_i(p),t_i,t_1,\ldots ,t_{i-1},t_{i+1},\ldots ,t_s) \cr \hbox{or }
\,\,&(t_i/z^{(2)}_i(p), z^{(2)}_i(p),t_i,t_1,\ldots ,t_{i-1}, t_{i+1},\ldots ,t_s)\,\,, \cr} \eqnu $$
according to whether $(p,t)\in W\cap\phi_{\tau}(U^{(1)}_i)\times D^s$ or $(p,t)\in W\cap\phi_{\tau}(U^{(2)}_i)
\times D^s$. It can be proved that $X$ is a complex manifold. Local coordinates at a point $(\phi_{\tau}(p),
t)\in W$ are given by $(\phi_{\tau}(z),t,\tau)$, where $z$ is a holomorphic coordinate for $M_0$ in a 
neighbourhood of $p$. The natural projection $\pi:X\to D^s\times U$ has fiber $M_{(t,\tau)}$. It is a 
standard procedure to introduce the coordinates    
$$x_i:={{z^{(1)}_i+z^{(2)}_i}\over 2}\,\,,\,\,\,\,y_i:={{z^{(1)}_i-z^{(2)}_i}\over 2} \eqnu $$
on $V_i$. Then $y_i=(x_i^2-t_i)^{1/2}$ represents $A_{t_i}\cap V_i$ as a branched double covering of a
neighbourhood of $x_i=0$ with branch points at $x_i=\pm t_i^{1/2}$. At
$t_i=0$ the fiber crosses itself
at the node corresponding to the pinching of the curve $\gamma_i$, given by the equations $x_i=y_i=t_i=0$,
and the normalization of $A_0$ is $U^{(1)}_i\cup U^{(2)}_i$. The $x_i$ are called the {\it pinching
coordinates }for the corresponding nodes. Holomorphic (abelian) differentials on $M_{(t,\tau)}$ are
customarily written in terms of these coordinates. A {\it regular }$d$-{\it differential }on a 
finite Riemann surface $M$ with punctures is by definition a meromorphic form of type $(d,0)$ on each 
part of $M$, holomorphic in the complement of the punctures, with poles of order at most $d$ at
the punctures and equal (opposite) residues at paired punctures if $d$ is even (odd). The 
{\it residue }of a regular $d$-differential $q$ at a point $p\in M$, denoted by 
$\hbox{Res}_q(p)$, is the residue at $z=0$ (in the standard sense) of the 
abelian differential $z^{d-1} q(z) dz$, where $z$ is a holomorphic coordinate
such that $z(p)=0$ and $q=q(z)dz^2$ [42, \S 4]. Only the cases $d=1$ and 
$d=2$ of abelian and quadratic differentials are relevant in this paper. Any 
regular $1$-differential $u(t,\tau)$ on $M_{(t,\tau)}$, holomorphic on $A_{t_i}$
and as a function on the total space $X$, can be written in terms of the pinching 
coordinate $x_i$ as follows [15, III], [42, 4.1]: 
$$u(t,\tau)=[a(t,\tau,x_i) + b(t,\tau,x_i)/(x_i^2-t_i)^{1/2}]\,dx_i\,\,,\eqnu $$
where $a$, $b$ are holomorphic functions. A direct computations shows that, under the change of
variables given by $(4.9)$, we have, with respect to the coordinate $z^{(k)}$ on $A_{t_i}$, $k=1,2$:
$$u(t,\tau)= [f^{(k)}(t,\tau,z^{(1)}+z^{(2)})\,+\, {{g^{(k)}(t,\tau,z^{(1)
}+z^{(2)})}\over {z^{(k)}}}]\,dz^{(k)}\,\,, \speqnu{4.10'}$$
where $f^{(k)}$, $g^{(k)}$ are holomorphic functions.

There is a normalized basis $\{\theta_1(t,\tau),\ldots , \theta_g(t,\tau)\}$ of regular 
$1$-differentials on $M_{(t,\tau)}$, holomorphic on $X$ [15, III], [42, Prop. 4.1], 
[71, Corollaries 1,4,5], with respect to a canonical homology basis $$\{a_1(t,\tau),\ldots ,
a_g(t,\tau),b_1(t,\tau),\ldots ,b_g(t,\tau)\}.$$ Let $\Pi_{(t,\tau)}$  be the corresponding period
matrix, which is well-defined (finite) if all $t_i\not =0$. Variational formulas for the period
matrix as $t_i\to 0$ were obtained in [15, III] and (corrected) in \ref\YaA. Let 
$q(t,\tau)\in Q_g$ be a holomorphic quadratic differential on $M_{(t,\tau)}$, holomorphic as
a function on $X$, which converges in the projective bundle $Q_g/{\Bbb R}^\ast$ as $t_i\to 0$, 
$i\in J\subset \{1,\ldots ,s\}$, to a regular quadratic differential having poles of order two 
at all the punctures and different from zero on every part of the corresponding pinched Riemann 
surface. Let $\mu_{(t,\tau)}:=|q(t,\tau)|/q(t,\tau)$ be the related family of Beltrami differentials. 
We prove a result concerning the limit behaviour of the {\it derivative }of the period matrix 
$\Pi_{(t,\tau)}$ in the direction of the Beltrami differential $\mu_{(t,\tau)}$ at the boundary 
of the moduli space. Let $M_{\tau}(t)$ be the image of the natural inclusion
of $M_{(t,\tau)}\setminus\cup\{\gamma_{t_i}\}$ in $M_{\tau}$. The surface 
$M_{\tau}(t)$ is an open subset of $M_{\tau}$ which can be described as follows:
$$M_{\tau}(t)=M_{\tau}\setminus\bigcup_{i=1}^s\bigl(\{|z_i^{(1)}|\leq 
|t_i|^{1/2}\}\cup\{|z_i^{(2)}|\leq |t_i|^{1/2}\}\bigr)\,\,. \eqnu $$

\proclaim{Lemma} As $t_i\to 0${\rm ,}  for all $i\in J\subset \{1,\ldots ,s\}${\rm ,}
$$ {{d\Pi^{ij}_{(t,\tau)}}\over {d\mu_{(t,\tau)}}}-\int_{M_{\tau}(t)}
\theta_i(t',\tau)\theta_j(t',\tau)\,\mu_{(t',\tau)}\,\,\to 0\,\,,\eqnu $$
where $t'_i=0$ if $i\in J${\rm ,} $t'_i=t_i\not=0$ if $i\not\in J${\rm .}
\endproclaim

\demo{Proof} By Rauch's formula $(4.4)$, the derivative of the period matrix
can be written as an integral over $M_{\tau}(t)$, identified to the complement
of the finite set of curves $\gamma_{t_i}$ in $M_{(t,\tau)}$. In addition, since
the integrand is uniformly bounded and converges on compact sets contained in 
the complement of the punctures by the variational formulas of [15, III] 
or \ref\YaA, it suffices to prove convergence on each annulus 
$$R^{(k)}_{t_i}:=\{|t_i|^{1/2}<|z^{(k)}_i|<1\}\subset A_{t_i}\,\,,\eqnu $$
for all $i\in J$, $k=1,2$.

 Let $R_t:=\{|t|^{1/2}<|z|<1\}$ one of the annuli in $(4.13)$. We 
should estimate, as $t\to 0$ in $\Bbb C$, an expression of the form
$$\int_{R_t} \bigl(\,\theta_i(t,z)\theta_j(t,z){{|q_t(z)|}\over {q_t(z)}}-
\theta_i(0,z)\theta_j(0,z){{|q_0(z)|}\over {q_0(z)}}\,\bigr)\,dz\wedge d\, {\overline z}  
\,\,,\eqnu $$
where $\theta_i(t,z)$ are meromorphic functions with at most a simple pole at 
$z=0$ given by $(4.10')$, that is
$$\theta_i(t,z)= f_i(t,z+t/z)\,+\, {{g_i(t,z+t/z)}\over {z}}={{F_i(t,z+t/z)}
\over {z}} \,\,, 
\eqnu $$
with $f_i$, $g_i$, hence $F_i$, holomorphic functions, and $q_t(z)$ is a 
meromorphic function with a double pole at $z=0$, which can be written 
as [42, 5.2--5.3]:
$$q_t(z)={{a(t,z+t/z)}\over {z}} \,+\, {{b(t,z+t/z)}\over {z^2}}=
{{A(t,z+t/z)}\over {z^2}}\,\,, \speqnu{4.15'}$$ 
with $a$, $b$, hence $A$, holomorphic functions and $A(0,0)\not=0$. We split the
difference in $(4.14)$ into the following three terms:
$$\eqalign{ & {\rm (I)}\,\, \int_{R_t}\bigl(\theta_i(t,z)-\theta_i(0,z)\bigr)\theta_j(t,z) 
{{|q_t(z)|}\over {q_t(z)}} \,dz\wedge d\, {\overline z}\,\,, \cr
            & {\rm (II)}\,\, \int_{R_t}\bigl(\theta_j(t,z)-\theta_j(0,z)\bigr)\theta_i(0,z) 
{{|q_t(z)|}\over {q_t(z)}} \,dz\wedge d\, {\overline z}\,\,, \cr
            & {\rm (III)}\,\, \int_{R_t}\theta_i(0,z)\theta_j(0,z)\bigr({{|q_t(z)|}\over {q_t(z)}}
-{{|q_0(z)|}\over {q_0(z)}}\bigl)\,dz\wedge d\, {\overline z}\,\,. \cr} \eqnu $$
By $(4.15)$, it follows that, since $t/z\to 0$ if $z\in R_t$ as $t\to 0$, there is a constant 
$C_0>0$ such that
$$|\theta_i(t,z)-\theta_i(0,z)|\leq C_0 {{|t|}\over {|z|^2}}\,\,\,\,,\,\,\,\,\,\,\,\,
|\theta_i(t,z)|\leq C_0 {1\over {|z|}}\,\,. \eqnu $$
By   H\"older inequality and (4.17), we obtain:
$$ \eqalign{ &|{\rm (I)}| \leq |\!|\theta_i(t,z)-\theta_i(0,z)|\!|_{L^2(R_t)}\,|\!|\theta_j(t,z)|\!|_{L^2(R_t)} 
\leq C_1 (-|t|\log |t|)^{1/2} \,\,, \cr
             &|{\rm (II)}| \leq |\!|\theta_j(t,z)-\theta_j(0,z)|\!|_{L^2(R_t)}\,|\!|\theta_i(0,z)|\!|_{L^2(R_t)} 
\leq C_1 (-|t|\log |t|)^{1/2} \,\,. \cr } \eqnu $$ 
In order to estimate the third term in $(4.16)$, we proceed as follows. By $(4.15')$, since 
$A(0,0)\not=0$, there exists $r>0$ such that $A(t,z)\not=0$, hence $|A(t,z)|/
A(t,z)$ is a real analytic
function, in the polydisk $\{|t|<r,\break |z|<r\}\subset {\Bbb C}^2$. There exists therefore a constant
$C_2>0$ such that, in the annulus $\{|t|^{1/2}<|z|<r/2\}$, we have: 
$$\left|{{|q_t(z)|}\over {q_t(z)}}-{{|q_0(z)|}\over {q_0(z)}}\right|= \left|{{|A(t,z+t/z)|}\over {A(t,z+t/z)}}-
{{|A(0,z)|}\over {A(0,z)}}\right| \leq C_2 {{|t|}\over {|z|}} \leq C_2\, |t|^{1/2}\,\,.\eqnu $$
It follows that
$$\eqalign { |{\rm (III)}| \leq &- C_3 |t|^{1/2}\log |t|   \cr 
                    +&\, \left|\int_{|z|\geq r/2}\theta_i(0,z)\theta_j(0,z)
\bigr({{|q_t(z)|}\over {q_t(z)}}-{{|q_0(z)|}\over {q_0(z)}}\bigl)\,
dz\wedge d\, {\overline z}\,\right|\,\,,  \cr }
\speqnu{4.19'}$$
where the integral on the right-hand side converges to zero as $t\to 0$ by the dominated 
convergence theorem. By $(4.18)$ and $(4.19')$ the proof is concluded. 
\enddemo 

By Lemma 4.2 and the variational formulas for the period matrix proved
in [15, III], \ref\YaA, it is possible to estimate the determinant
of the derivative of the period matrix and, by the identity $(4.6)$,   the determinant
of the hermitian form $H_q$  close to the boundary of the moduli space, by computing
them at the boundary. Let in particular ${\Cal S}_g$ be the component of the boundary
consisting of {\it regular }quadratic differentials on the Riemann sphere with $2g$ paired 
punctures, having poles of order $2$ at all punctures.

\specialnumber{4.2'}  
\proclaim{Lemma} Let $q_0\in {\Cal S}_g$. Let $M_0$ be the corresponding Riemann sphere with 
the $2g$ paired punctures $p^{(1)}_i${\rm ,} $p^{(2)}_i${\rm ,} $i=1,\ldots ,g${\rm .} Let $t\in {\Bbb C}^g$ be such that 
$t_i\not=0$ for all $i=1,\ldots ,g$ and $M_0(t)\subset M_0$ be the open subset defined as in $(4.11)$.
Let $\{\theta_1(0),\ldots ,\theta_g(0)\}$ be the basis of the space of regular $1$\/{\rm -}\/differentials on 
$M_0${\rm ,} dual to the canonical homology basis $\{a_1,\ldots ,a_g,b_1,\ldots b_g\}${\rm ,} where $a_i$ is represented 
by the {\rm (}\/oriented\/{\rm )} boundary of a disk centered at $p^{(1)}_i$ {\rm (}\/not containing other punctures\/{\rm )} and 
$b_i$ by a path joining $p^{(1)}_i$ to $p^{(2)}_i${\rm .} Then   
$$\int_{M_0(t)} \theta_i(0)\theta_j(0) {{|q_0|}\over {q_0}}\, -\,
{1 \over {2\pi}}\, {{|\hbox{Res}_{q_0}(p^{(1)}_i)|} 
\over {\hbox{Res}_{q_0}(p^{(1)}_i)}}\,\,\delta_{ij}\log |t_i| \eqnu $$  
is bounded  as $t \to 0$, for all $i,j\in \{1,\ldots ,g\}${\rm .}
\endproclaim 

\demo{Proof} The normalized basis $\{\theta_1(0),\ldots ,\theta_g(0)\}$ on the punctured Riemann
sphere $M_0$ can be explicitly written:
$$ \theta_i(0)= {{p^{(1)}_i-p^{(2)}_i}\over {2\pi\imath}}\,\,{{dz}\over 
{(z- p^{(1)}_i)(z- p^{(2)}_i)}}\,\,. \eqnu $$
In fact, $(4.21)$ gives the only basis of regular $1$-differential $\{\theta_1,\ldots ,\theta_g\}$ 
satisfying $\theta_i(a_j)=\delta_{ij}$. Let $R^{(k)}_i(t):=\{ |t_i|^{1/2} < |z-p^{(k)}_i|< r\}$, 
where $r>0$ has been chosen so that all the annuli $R^{(k)}_i(t)$ are disjoint and $q_0$ has
no zeroes inside the closed disks $D_r(p^{(k)}_i)$ of radius $r$ centered at the punctures. Let 
$$q_0={{A^{(k)}_i(z)}\over {(z-p^{(k)}_i)^2}}\, dz^2 \,\, \eqnu $$ 
on the disk $D_r(p^{(k)}_i)$. Since $q_0$ is a regular quadratic differential with poles of order
$2$ at all punctures, it follows that 
$$A^{(1)}_i(p^{(1)}_i)=\hbox{Res}_{q_0}(p^{(1)}_i)=
A^{(2)}_i(p^{(2)}_i)=\hbox{Res}_{q_0}(p^{(2)}_i)\,\not = \,0 \,\,. \speqnu{4.22'}$$
We have 
$$\int_{R^{(k)}_i(t)}\theta_i(0)^2{{|q_0|}\over {q_0}}=\int_{R^{(k)}_i(t)} 
{{F^{(k)}_i(z)} \over {|z-p^{(k)}_i|^2}}\,dz\wedge d\, {\overline z} \,\,,
\eqnu $$
where the function $F^{(k)}_i$ is given on the disk $D_r(p^{(k)}_i)$ by the formula:
$$ F^{(k)}_i(z):= {{(p^{(1)}_i-p^{(2)}_i)^2}\over {8\pi^2 \imath\,(z-p^{(h)}_i)^2}}
\,{{|A^{(k)}_i(z)|}\over {A^{(k)}_i(z)}} \,\,, \speqnu{4.23'}$$ 
where $k$, $h\in \{1,2\}$ and $k\not= h$. Hence $F^{(k)}_i$ is real analytic on 
the disk $D_r(p^{(k)}_i)$ and, by Taylor series expansion at $p^{(k)}_i$, we obtain that
$$\int_{R^{(k)}_i(t)} {{F^{(k)}_i(z)}\over{|z-p^{(k)}_i|^2}}
\,dz\wedge d\, {\overline z}\,-\, F^{(k)}_i(p^{(k)}_i)\int_{R^{(k)}_i(t)} 
{{dz\wedge d\, {\overline z}}\over {|z-p^{(k)}_i|^2}} \eqnu $$
is bounded as $t \to 0$. Similarly, if $i\not = j$, we obtain
$$\int_{R^{(k)}_i(t)}\theta_i(0)\theta_j(0){{|q_0|}\over {q_0}}=\int_{R^{(k)}_i(t)} 
{{F^{(k)}_{ij}(z)} \over {\left(\overline{z-p^{(k)}_i}\right)}}\,dz\wedge d\, {\overline z} \,\,,
\eqnu $$
where the function $F^{(k)}_{ij}$ is given on the disk $D_r(p^{(k)}_i)$ by the formula:
$$ F^{(k)}_{ij}(z):= {{(p^{(1)}_i-p^{(2)}_i)(p^{(1)}_j-p^{(2)}_j)}\over 
{8\pi^2 \imath\,(z-p^{(h)}_i)(z-p^{(1)}_j)(z-p^{(2)}_j)}}{{|A^{(k)}_i(z)|}\over 
{A^{(k)}_i(z)}} \,\,, \speqnu{4.25'}$$ 
where, as above, $k$, $h\in \{1,2\}$ and $k\not= h$. Since $F^{(k)}_{ij}$ is real analytic, 
hence bounded, on the disk $D_r(p^{(k)}_i)$, the integral in $(4.25)$ is bounded as 
$t \to 0$.

 Let $R_{ij}(t):= R^{(1)}_i(t)\cup R^{(2)}_i(t)\cup R^{(1)}_j(t)\cup
R^{(2)}_j(t)$. The integral
$$\int_{M_0(t)\setminus R_{ij}(t)} \theta_i(0)\theta_j(0) {{|q_0|}\over {q_0}} 
\eqnu $$
is bounded as $t \to 0$, for all  $i,j\in  \{1,\ldots ,g\}$, since the integrand is 
uniformly bounded on $M_0(t)\setminus R_{ij}(t)$. Our statement follows, by 
$(4.23)$, by the boundedness of the difference in $(4.24)$ and of the integrals
in $(4.25)$, $(4.26)$, by explicit computation of the term on the right of 
$(4.24)$ and by the formulas $(4.22')$ on the residues of $q_0$. 
\enddemo

\numbereddemo{Definition} A measured foliation ${\Cal F}$ on a compact orientable surface $M$ of
genus $g\geq 2$ is said to be periodic if all its regular leaves are closed (compact) curves.  
A periodic measured foliation ${\Cal F}$ will be called Lagrangian if the subspace ${\Cal L}({\Cal F}) 
\subset H_1(M,{\Bbb R})$, generated by the homology classes of the regular leaves of ${\Cal F}$, is 
a Lagrangian subspace of the homology vector space $H_1(M,{\Bbb R})$, endowed with the symplectic 
structure given by the intersection form. A periodic measured foliation ${\Cal F}$ is Lagrangian
if and only if it has $g\geq 2$ distinct regular leaves $\gamma_1,\ldots ,\gamma_g$ such that ${\widehat M}
:=M\setminus\cup \{\gamma_1,\ldots ,\gamma_g\}$ is homeomorphic to a sphere minus $2g$ (paired) disjoint 
disks. 
\enddemo    

Let ${\Cal L}^{\pm}_{\kappa}$ be the set of holomorphic quadratic 
differentials $q\in Q_{\kappa}$ such that the foliation ${\Cal F}_{\pm q}$ is 
Lagrangian. Let $\Lambda(M,{\Bbb R})$ be the Grassmannian of\break Lagrangian 
subspaces of the symplectic vector space $H_1(M,{\Bbb R})$. Let $\Lambda\in
\Lambda(M,{\Bbb R})$ and let 
$${\Cal L}^{\pm}_{\kappa,\Lambda}:=\{q\in {\Cal L}^{\pm}_{\kappa}\,|\,{\Cal L}
({\Cal F}_{\pm q})\cap\Lambda=\{0\}\,\}\,\,. \eqnu$$  

The final key step consists of the following crucial result:

\proclaim{Lemma} The set ${\Cal L}^+_{\kappa,\Lambda}\subset {\Cal L}^+
_{\kappa}$ $[{\Cal L}^-_{\kappa,\Lambda}\subset{\Cal L}^-_{\kappa}]$ of 
holomorphic quadratic differentials $q\in Q_{\kappa}$ such that the horizontal 
foliation ${\Cal F}_q$ {\rm [}\/the vertical foliation ${\Cal F}_{-q}$\/{\rm ]} is Lagrangian
and ${\Cal L}({\Cal F}_q)\cap \Lambda=\{0\}$ $[{\Cal L}({\Cal F}_{-q})\cap 
\Lambda=\{0\}]$ is dense in $Q_{\kappa}${\rm .}
\endproclaim

\demo{Proof} Let $\kappa:=(k_1,\ldots ,k_{\sigma})$, where all $k_i$ are even and $\sum k_i=4g-4$. Let 
${\Cal F}_{\kappa}(M)$ be the set of all isotopy equivalence classes of orientable measured foliations 
${\Cal F}$ on $M$ with canonical saddle-like singularities of multiplicities $(k_1,\ldots ,k_{\sigma})$ at
a finite set $\Sigma_{\Cal F}:=\{p_1,\ldots ,p_{\sigma}\}$. Let ${\Cal F}'_{\kappa}(M)\subset {\Cal F}
_{\kappa}(M)$ be the subset given by the measured foliations ${\Cal F}$ such that there exists $q\in 
Q_{\kappa}$ with ${\Cal F}_q={\Cal F}$ (or ${\Cal F}_{-q}={\Cal F}$). The set ${\Cal F}'_{\kappa}(M)$ 
is open in ${\Cal F}_{\kappa}(M)$ and the map $Q_{\kappa}\to {\Cal F}'_{\kappa}(M)\times {\Cal F}'
_{\kappa}(M)$, given by $q\to ({\Cal F}_q,{\Cal F}_{-q})$, is open. By A.Katok's {\it local 
classification theorem}, announced in [29, Th.\ 3] and proved in [32, Th.\ 14.7.4] or [47, Th.\ 7.11.7], 
the space ${\Cal F}_{\kappa}(M)$ is locally modeled on $H^1(M,\Sigma_{\kappa};{\Bbb R})$, while the space
$Q_{\kappa}$ is locally modeled on $H^1(M,\Sigma_{\kappa};{\Bbb C})$ (see \S 1).
Let $\Lambda\in \Lambda(M,{\Bbb R})$.  We claim that the set of Lagrangian 
foliations ${\Cal F}\in {\Cal F}_{\kappa}(M)$ such that ${\Cal L}({\Cal F})\cap
\Lambda\break =\{0\}$ is dense in ${\Cal F}_{\kappa}(M)$. The proof of the claim will
conclude the argument. 

Let ${\Cal F}:=\{ \eta_{\Cal F}=0\}\in {\Cal F}_{\kappa}(M)$ and let $\Sigma_{\Cal F}
\subset M$ denote its (finite) singularity set. If the cohomology class $[\eta_{\Cal F}]\in {\Bbb R} 
\cdot H^1(M,\Sigma_{\Cal F};{\Bbb Z})$, hence, in particular, if $[\eta_{\Cal F}]\in H^1(M,\Sigma_{\Cal F};
{\Bbb Q})$, then ${\Cal F}$ is periodic. In fact, in that case, the set of 
${\Cal F}$-length of isotopy classes of simple closed curves is discrete. 
Hence, all the regular leaves are closed by 
the Poincar\'e recurrence theorem. It follows that the periodic measured foliations are dense in 
${\Cal F}_{\kappa}(M)$. If ${\Cal F}$ is periodic, the surface can be decomposed, by cutting along 
the singular leaves, into a finite union of cylindrical components whose number is at most $3g-3$. 
Let $d:=\hbox{dim}\,{\Cal L}({\Cal F})\in \{1,\ldots ,g\}$, where ${\Cal L}({\Cal F})\subset 
H_1(M,{\Bbb Z})$ is given by Definition~4.3. If $d=g$, then ${\Cal F}$ is a 
Lagrangian measured foliation.

Let $P:H_1(M,{\Bbb R})\to H^1(M,{\Bbb R})$ be the (symplectic) map
given by the Poincar\'e duality. We claim that, if ${\Cal F}$ is periodic, then
$P^{-1}[\eta_{\Cal F}]\in {\Cal L}({\Cal F})$. More precisely, if $\{a_1,\ldots ,
a_s\}$ are the oriented waist curves of the cylinders $\{A_1,\ldots ,A_s\}$ of 
${\Cal F}$, which are respectively of heights $\{h_1,\ldots ,h_s\}$, then
$$P^{-1}[\eta_{\Cal F}]\,=\, \sum _{i=1}^s h_i\,[a_i] \,\in\, H_1(M,{\Bbb R})
\,. \eqnu$$
In fact, if $\gamma\subset M$ is any simple oriented closed curve, then $\gamma 
\cap A_i$ is homologous to $([a_i]\cap [\gamma])\cdot v_i$  relative to $\partial A_i$,
where $v_i$ is a positively oriented vertical segment joining the ends of $A_i$. Hence 
$(4.28)$ follows.

Let $\Lambda\in \Lambda(M,{\Bbb R})$ and let $\{a_1,\ldots ,a_k\}$ be a
maximal system of regular leaves of ${\Cal F}$ such that the system of homology
classes $\{[a_1],\ldots ,[a_k]\}$ is linearly independent in $H_1(M,{\Bbb R})$ and
the isotropic subspace $I(a_1,\ldots ,a_k)\subset H_1(M,{\Bbb R})$, generated by
$\{[a_1],\ldots ,[a_k]\}$, is transverse to $\Lambda$. If $k=g$, then ${\Cal F}$
is Lagrangian and ${\Cal L}({\Cal F}) \cap \Lambda= \{0\}$. If $k<g$, then
there exists a primitive homology class $h\in H_1(M,{\Bbb Z})$ such that 
$h\cap [a_1]=\cdots =h\cap [a_k]=0$ and $h\not\in I(a_1,\ldots ,a_k) \oplus \Lambda$.
In fact, since $\Lambda$ is Lagrangian, the dimension of the subspace $\{h\in
I(a_1,\ldots ,a_k)\oplus\Lambda\,|\,h\cap [a_1]=\cdots =h\cap [a_k]=0\}$ is equal to 
$g$. On the other hand, the dimension over $\Bbb Q$ of the subspace $\{h\in
H_1(M,{\Bbb Q})\,|\,h\cap [a_1]=\cdots =h\cap [a_k]=0\}$ is equal to $2g-k$. Since 
$2g-k>g$, there exists ${\widehat h}\in H_1(M,{\Bbb Q})$ such that ${\widehat h}\not\in
I(a_1,\ldots ,a_k)\oplus\Lambda$ and ${\widehat h}\cap [a_1]=\cdots ={\widehat h}\cap [a_k]=0$.
Let then $h\in H_1(M,{\Bbb Z})$ be the unique integer multiple of $\widehat h$ which
is a primitive integer class. 

Let $z=\sum_i n_i \gamma_i$ be an integer cycle on $M$, supported on a finite 
union of disjoint smooth simple closed curves $\{\gamma_1,\ldots ,\gamma_s\}$ such that $[z]=h$,
$\hbox{supp}(z)\cap a_1=\cdots =\hbox{supp}(z)\cap a_k=\emptyset$ and $\hbox{supp}
(z) \cap \Sigma_{\Cal F}= \emptyset$. Let ${\Cal V}(\gamma_i)\subset\subset{\Cal U}
(\gamma_i)$ be open tubular neighbourhoods of $\gamma_i$ in $M$ such that, for all $i\not=j
\in \{1,\ldots ,s\}$,
$$\eqalign{&\overline{{\Cal U}(\gamma_i)}\cap \overline{{\Cal U}(\gamma_j)} =\emptyset\,\,, \cr
&\overline{{\Cal U}(\gamma_i)}\cap (a_1\cup\cdots \cup a_k)=\emptyset\,\,,\cr
&\overline{{\Cal U}(\gamma_i)}\subset M\setminus\Sigma_{\Cal F}\,\,. \cr} \eqnu$$
Let ${\Cal U}^{\pm}(\gamma_i)$ be the two connected components of the open set
${\Cal U}(\gamma_i)\setminus\gamma_i$ and let ${\Cal V}^{\pm}(\gamma_i):={\Cal V}
(\gamma_i)\cap {\Cal U}^{\pm}(\gamma_i)$. Let $\phi_i:M\to {\Bbb R}$ be a function, 
smooth on $M\setminus\gamma_i$, with the following properties:
$$\eqalign{ \phi_i(p)&=0 \,\,,\,\,\,\, p\in {\Cal U}^-(\gamma_i)\cup M\setminus {\Cal U}^+ 
(\gamma_i)\,\, , \cr  \phi_i(p)&=1\,\,,\,\,\,\,p\in\overline{ {\Cal V}^+(\gamma_i)}\,\,.\cr} 
\speqnu{4.29'}$$
Let $r\in {\Bbb Q}\setminus\{0\}$. The smooth $1$-form $\eta_r:=\eta_{\Cal F}+r 
\sum_i n_i d\phi_i$ is closed and $\eta_r\to \eta_{\Cal F}$, as $r\to 0$, in the
space of smooth $1$-forms on $M$. Since $d\phi_i=0$ on $M\setminus {\Cal U}
(\gamma_i)$, $\eta_r\equiv\eta_{\Cal F}$ in a neighbourhood of $\Sigma_{\Cal F}$,
hence, if $r\not=0$ is sufficiently small, $\eta_r(p)=0$ if and only if $p\in \Sigma_{\Cal F}$, 
and the foliation ${\Cal F}_r:=\{\eta_r=0\}\in {\Cal F}_{\kappa}(M)$. Since $r\in {\Bbb Q}$, 
the fundamental class $[\eta_r]\in H^1(M, \Sigma_{\Cal F};{\Bbb Q})$, hence ${\Cal F}_r$ 
is periodic. The simple closed curves $a_1,\ldots ,a_k$ are regular leaves of ${\Cal F}_r$. In 
fact, $d\phi_i=0$ on $M\setminus {\Cal U}(\gamma_i)$ and $\cup\{a_1,\ldots ,a_k\}\subset 
M\setminus \cup\,{\overline {\Cal U}(\gamma_i)}$, hence $\eta_r\equiv \eta_{\Cal F}$ in a 
neighbourhood of  $\cup\{a_1,\ldots ,a_k\}$. Let $I_r\subset H_1(M,{\Bbb R})$ be the 
maximal isotropic subspace, generated by the regular leaves of ${\Cal F}_r$, such that
$I(a_1,\ldots ,a_k)\subset I_r$ and $I_r \cap \Lambda=\{0\}$. We claim that $I_r \not = 
I(a_1,\ldots ,a_k)$. If the latter claim is false, ${\Cal L}({\Cal F}_r)\subset I(a_1,\ldots ,a_k)
\oplus\Lambda$, hence, by $(4.28)$, $P^{-1}[\eta_r]\in I(a_1,\ldots ,a_k)\oplus\Lambda$, 
but, by $(4.29')$ and the definition of $\eta_r$, $P^{-1}[\eta_r]=P^{-1}[\eta_{\Cal F}] + 
rh \not\in I(a_1,\ldots ,a_k)\oplus \Lambda$. In fact, by $(4.28)$ and the definition of the 
system $\{a_1,\ldots ,a_k\}$, $P^{-1}[\eta_{\Cal F}]\in I(a_1,\ldots ,a_k)\oplus\Lambda$ and, 
by construction, $h\not\in I(a_1,\ldots ,a_k)\oplus\Lambda$.

By a finite iteration of the previous construction, we can show
that the closure in ${\Cal F}_{\kappa}(M)$ of the subset of all Lagrangian 
measured foliations ${\Cal F}$, such that ${\Cal L}({\Cal F})\cap\Lambda=\{0\}$,
contains the subset of all periodic measured foliations. Hence it coincides 
with the entire space ${\Cal F}_{\kappa}(M)$. Our claim is therefore proved. 
\enddemo

\specialnumber{4.4'}     
\proclaim {Lemma} Every connected component ${\Cal C}_{\kappa}$ of any stratum 
${\Cal M}_{\kappa}$ of the moduli space of quadratic differentials has a boundary point 
in ${\Cal S}_g$ with $2g$ {\rm (}\/paired\/{\rm )} real strictly positive residues. 
\endproclaim

\demo{Proof} Let $\{\gamma_1,\ldots ,\gamma_g\}$ be a system disjoint simple closed curves on $M$ with 
the property that ${\widehat M}:=M\setminus\cup\{\gamma_1,\ldots ,\gamma_g\}$ is homeomorphic to a sphere 
minus $2g$ disjoint disks. Let ${\Cal V}_{\kappa}(\gamma_i)$ be the open subset of quadratic 
differentials $q\in Q_{\kappa}$ such that the vertical foliation ${\Cal F}_q=\{\Re(q^{1/2})=0\}$ 
has a closed regular leaf $\gamma_i(q)$ isotopic to $\gamma_i$ and let ${\Cal V}_{\kappa}(\gamma_1,\ldots ,\gamma_g)
:= \cap {\Cal V}_{\kappa}(\gamma_i)$, $i=1,\ldots ,g$. By Lemma 4.4, since ${\Cal C}_{\kappa}\subset 
{\Cal M}_{\kappa}$ is open, the system $\{\gamma_1,\ldots ,\gamma_g\}$ can be chosen with the additional 
property that the pull-back of ${\Cal C}_{\kappa}$ to the Teichm\"uller space $Q_{\kappa}$ has 
non-empty intersection with ${\Cal V}_{\kappa}(\gamma_1,\ldots ,\gamma_g)$. We then construct a 
continuous deformation $\Phi:(0,1]^g\times {\Cal V}_{\kappa}(\gamma_1,\ldots ,\gamma_g) \to 
{\Cal V}_{\kappa}(\gamma_1,\ldots ,\gamma_g)$ having the properties listed below. Such a construction 
will conclude the proof. Let $q_t:=\Phi_t(q)$, $t=(t_1,\ldots ,t_g)\in (0,1]^g$. Then  
\smallbreak
\item{(1)} If $t_i=1$ for all $i=1,\ldots ,g$, then $q_t=q$.  \smallbreak
\item{(2)} If $t_i\to 0$ for all $i=1,\ldots ,g$, then the Riemann surface carrying $q_t$ converges
to a Riemann surface with nodes, pinched along the curves $\gamma_1,\ldots ,\gamma_g$, hence $q_t$
converges in the moduli space ${\Cal M}_g$ to a quadratic differential $q_0 \in {\Cal S}_g$.
\smallbreak
\item{(3)} The zero set of $q_t$ coincides with the zero set $\Sigma_q$ of $q$ for all $t\in (0,1]^g$
and the cohomology class $[\Im(q_t^{1/2})]=[\Im(q^{1/2})] \in H^1(M,\Sigma_q,{\Bbb R})$.
\smallbreak

   We will construct $\Phi_t(q)$ as a composition of deformations $$\Phi^{(i)}:(0,1] \times
{\Cal V}_{\kappa}(\gamma_i)\to {\Cal V}_{\kappa}(\gamma_i)  , \qquad   i=1,\ldots ,g,
$$
pinching along $\gamma_i$.
Let $q\in {\Cal V}_{\kappa}(\gamma_i)$ and let $\gamma_i(q)$ be the closed regular vertical 
$q$-trajectory isotopic to $\gamma_i$. Let then ${\widehat M}_i:=M\setminus \{\gamma_i(q)\}$. 
Let $\gamma^{(1)}_i$, $\gamma^{(2)}_i$ be the boundary components of ${\widehat M}_i$ produced by 
cutting along the regular trajectory $\gamma_i(q)$, endowed with the orientation induced by 
${\widehat M}_i$, and let $\phi_i:\gamma^{(1)}_i\to\gamma^{(2)}_i$ be the corresponding (orientation 
reversing) attaching map. Let $D^{(k)}_i$ be copies of the closed unit disk in the complex plane 
and $h^{(k)}_i$ be the meromorphic differentials on $D^{(k)}_i$, $k=1,2$, given by 
$$h^{(k)}_i:= {{c^{(k)}_i}\over {2\pi}}\,{{dz^{(k)}_i}\over {z^{(k)}_i}}\,\,, \eqnu $$ 
where $c^{(k)}_i\not=0$ are the periods of $\Im(q^{1/2})$ along $\gamma^{(k)}_i$. Since 
$\gamma^{(k)}_i$, $k=1,2$, can be identified to $\gamma_i(q)$ taken with opposite orientations, 
we have that $c^{(1)}_i + c^{(2)}_i =0$. The boundary circle of $D^{(k)}_i$ is a regular 
vertical trajectory of $h^{(k)}_i$. Hence by gluing the pairs $(D^{(k)}_i\setminus\{0\},
h^{(k)}_i)$ to the pair $({\widehat M}_i,q^{1/2})$ along the curves $\gamma^{(k)}_i$, we 
obtain a pair $(M_i,q_i^{1/2})$, where $M_i$ is a surface of genus $g-1$ with two (paired) 
punctures and $q_i^{1/2}$ is a regular $1$-differential with the same zeroes as $q^{1/2}$. Let 
$\psi^{(k)}_i:\partial D^{(k)}_i \to \gamma^{(k)}_i$ be the corresponding attaching map. Let $t_i
\in (0,1]$, let $S^{(k)}_i(t_i)$ be the circle of radius $t_i$ in the disk $D^{(k)}_i$ and let 
$\psi_i(t_i):S^{(1)}_i(t_i)\to S^{(2)}_i(t_i)$ be the attaching map defined as follows:
$$\psi_i(t_i)(z^{(1)}_i):= t_i\,(\psi^{(2)}_i)^{-1}\circ\phi_i\circ\psi^{(1)}_i
(t_i^{-1}z^{(1)}_i) \,\,. \speqnu {4.30'}$$
Let $(M_i(t_i),q_i(t_i)^{1/2})$ be the  pair obtained by restriction of $(M_i,q_i^{1/2})$ to the 
complement of the open disk $\{|z^{(k)}_i|<t_i\}\subset D^{(k)}_i$ in $M_i$, then by gluing
along the boundary curves in the way prescribed by the attaching map $(4.30')$.
If\break $t_i\not=0$, 
$M_i(t_i)$ is a compact Riemann surface of genus $g\geq 2$ and, since the boundary circles 
$S^{(k)}_i(t_i)$ are regular vertical trajectories of the quadratic differential $q_i$, $q_i(t_i)$ 
has the same zeroes as $q$, hence $q_i(t_i)\in Q_{\kappa}$. In addition, by the choice of the 
holomorphic differential $(4.30)$ and of the attaching map $(4.30')$, the cohomology class 
$[\Im \bigl(q_i^{1/2}(t_i)\bigr)]\equiv [\Im (q^{1/2})]\in H^1(M,\Sigma_q;{\Bbb R})$, for all 
$t_i\in (0,1]$. By construction, $\bigl(M_i(t_i),q_i(t_i)\bigr)=(M,q)$, if $t_i=1$, and 
$\bigl(M_i(t_i),q_i(t_i)\bigr)=(M_i,q_i)$, if $t_i=0$, where $q_i$ is a regular quadratic 
differential with two poles of order $2$ at the two (paired) punctures of the surface $M_i$, 
with real strictly positive (equal) residues. We then define:
$$\eqalign{&{\rm (a)}\,\,\Phi^{(i)}_{t_i}(q):= q_i(t_i)\,\,,\,\,\,\,\,\,\,\,\,\,\,\,\,\,\,\,\,\,\,\,\,\,
                        \,\,\,\,\,\,\,\,\,\,(t_i,q)\in (0,1]\times {\Cal V}_{\kappa}(\gamma_i)\,;\cr 
           &{\rm (b)}\,\,\Phi_t(q):=\Phi^{(1)}_{t_1}\circ \cdots \circ\Phi^{(g)}_{t_g}(q)\,\,,\,\,\,\,
             (t,q) \in (0,1]^g\times {\Cal V}_{\kappa}(\gamma_1,\ldots ,\gamma_g)\,\,.\cr}
\eqnu $$
It can be checked that the pinching deformation $\Phi:(0,1]^g\times {\Cal V}_{\kappa}(\gamma_1,
\ldots,\gamma_g) \to {\Cal V}_{\kappa}(\gamma_1,\ldots ,\gamma_g)$, given by $(4.31)$, is well-defined 
and has the required properties. In fact, $\Phi^{(i)}_{t_i}\bigl({\Cal V}_{\kappa}(\gamma_1,\ldots ,
\gamma_g)\bigr) \subset {\Cal V}_{\kappa}(\gamma_1,\ldots ,\gamma_g)$, for all $i=1,\ldots ,g$ and all 
$t_i\in (0,1]$. Since, by definition, $\Phi^{(i)}_{t_i}(q)=q$, if $t_i=1$, the property $(1)$ 
holds. Let $(M_0,q_0)$ be the regular quadratic differential constructed by gluing {\it all }the 
pairs $(D^{(k)}_i\setminus\{0\},h^{(k)}_i)$ to the pair $({\widehat M},q^{1/2})$ along the curves 
$\gamma^{(k)}_i$. Then, since $\widehat M$ is homeomorphic to a sphere minus $2g$ disjoint disks, 
$M_0$ is a punctured Riemann sphere (by the uniformization theorem) with $2g$ paired punctures 
and $q_0$ is a regular quadratic differential with the same zeroes as $q$ and poles of order 
$2$ at all punctures. Hence $q_0\in {\Cal S}_g$ and in particular, by $(4.30)$, it has strictly 
positive residues. Since $q_t:=\Phi_t(q)\to q_0$ as $t_i\to 0$, for all $i=1,\ldots ,g$, property 
$(2)$ is also proved. Finally, property $(3)$, which will be relevant in a similar construction 
carried out in Section 8.2, follows, as we have remarked above, from the choice of the differentials
$(4.30)$ and of the attaching maps $(4.30')$.
\enddemo

\proclaim{Theorem} No connected component ${\Cal C}_{\kappa}$ of a stratum ${\Cal M}_{\kappa}$ 
of holomorphic 
quadratic differentials is contained in the determinant locus{\rm .} In fact{\rm ,} the following stronger 
result holds{\rm .} Let ${\Cal C}^{(1)}_{\kappa}:={\Cal C}_{\kappa}\cap {\Cal M}^{(1)}_{\kappa}$. We have\/{\rm :}
$$\sup_{q\in {\Cal C}^{(1)}_{\kappa}} \Lambda_i(q)\,=\,1\,\,,\,\,\,\, \hbox{for all } i\in\{1,\ldots ,g\}\,. 
\eqnu $$
\endproclaim  

\demo{Proof} By the definition $(4.2)$ it is sufficient to prove the statement for $i=g$. By 
Lemma 4.4$'$, ${\Cal C}_{\kappa}$ has an accumulation point $(M_0,q_0)\in {\Cal S}_g$ with real 
strictly positive residues at all $2g$ (paired) punctures of the punctured Riemann sphere $M_0$. 
Let $(M_{(t,\tau)},q_{(t,\tau)})$, $(t,\tau)\in {\Bbb C}^g\times {\Bbb C}^{2g-3}$ be the holomorphic 
family, whose construction has been outlined above, parametrizing a neighbourhood of $(M_0,q_0)$ 
in the compactified moduli space of quadratic differentials. Then $(M_{(t,\tau)},q_{(t,\tau)})=
(M_0,q_0)$, if $(t,\tau)=(0,0)$, and $(M_{(t,\tau)},q_{(t,\tau)})\in {\Cal S}_g$, if $t=0$ and 
$\tau\in {\Bbb C}^{2g-3}$ is close to the origin. Let $\Pi_{(t,\tau)}$ be the period matrix of
$M_{(t,\tau)}$ with respect to a normalized basis $\{\theta_1(t,\tau),\ldots ,\theta_g(t,\tau)\}$ on 
$M_{(t,\tau)}$, holomorphic in a neigbourhood of the origin in ${\Bbb C}^g\times {\Bbb C}^{2g-3}$. 
The period matrix is well-defined if $t_i\not=0$ for all $i=1,\ldots ,g$. By [15, III] or 
[71, \S 3, Corollary 6] (see for instance [15, p.\ 54]), the period matrix satisfies the 
following asymptotics (taking into account that the normalization condition adopted by [15]
and \ref\YaA differs from ours by a factor $2\pi\imath$):
$$\Pi_{(t,\tau)}^{ij}-{1\over{2\pi\imath}}\,\delta_{ij}\,\log t_i \eqnu $$
is bounded as $t \to 0$,  uniformly on $\tau\in {\Bbb C}^{2g-3}$ in a compact 
neighbourhood of the origin. Let $\Lambda_i(t,\tau)$, $i\in \{1,\ldots ,g\}$, be the eigenvalues $(4.2)$ 
of the hermitian form $(4.2')$ at $(M_{(t,\tau)}, q_{(t,\tau)})$. By formula $(4.6)$, Lemmas 4.2 and 
$4.2'$ and by the boundedness of $(4.33)$, we obtain:
$$\lim_{(t,\tau)\to (0,0)} \Lambda_1(t,\tau)\cdots\Lambda_g(t,\tau)=1 \,\,.
\eqnu $$
Since any neighbourhood of $(M_0,q_0)$ in the compactified moduli space contains points of the 
connected component ${\Cal C}_{\kappa}$ and since the hermitian form $(4.2')$ and all its eigenvalues 
$\Lambda_i$ are invariant under the multiplicative ${\Bbb R}$-action on quadratic differentials, the 
statement is proved. 
\enddemo

\specialnumber{4.5'}
\proclaim{Corollary} The second Lyapunov exponents of the Kontsevich\/{\rm -}\/Zorich cocycle with 
respect to the ergodic invariant measure given by the restriction to any connected component 
${\Cal C}^{(1)}_{\kappa}$ of a stratum ${\Cal M}^{(1)}_{\kappa}$ is strictly positive{\rm .} In fact{\rm ,} 
it satisfies the estimate
$$\lambda_2({\Cal C}^{(1)}_{\kappa}) \geq {1\over { \mu^{(1)}_{\kappa}({\Cal C}^{(1)}_{\kappa})}} 
\int_{{\Cal C}^{(1)}_{\kappa}}\Lambda_g \,d\mu^{(1)}_{\kappa}\,>\, 0\,\,. \eqnu $$
\endproclaim

\demo{Proof} The first inequality in $(4.35)$ follows from the ergodicity of the Teichm\"uller flow 
on ${\Cal C}^{(1)}_{\kappa}$ (see Theorem 1.1) and from the lower bound proved in Theorem 3.3. In 
fact, $\Lambda^-\equiv \Lambda_g$, by the definitions $(3.15)$, $(4.2)$, $(4.2')$ and standard 
properties of hermitian forms. Finally, the strict positivity of the integral in $(4.35)$ holds by 
Theorem 4.5, since the function $\Lambda_g$ is continuous and the measure $\mu^{(1)}_{\kappa}$ is 
positive on any open set. \enddemo

 \section{The Kontsevich-Zorich formula revisited\\ and other formulas for the 
Lyapunov exponents}

 In \ref\KZone M. Kontsevich and A. Zorich obtained a  formula for the sum\break
$\lambda_1+\cdots +\lambda_g$ of the non-negative Lyapunov exponents of the 
Kontsevich-Zorich cocycle. In this section we prove a different version of the 
formula along the lines of Sections 2 and  3. The same method yields formulas for 
all the partial sums $\lambda_1+\cdots +\lambda_k$, $1\leq k \leq g$.

Let $G_k(M,{\Bbb R})\subset H^1(M,{\Bbb R})$ be the 
Grassmannian of isotropic subspaces of dimension $k\in \{1,\ldots  ,g\}$, which is a
compact manifold of dimension\break $(2g-k)k-(k^2-k)/2$. Let $q\in Q^{(1)}_{\kappa}$  
be an orientable quadratic differential and $I_k\in G_k(M,{\Bbb R})$ an 
isotropic $k$-plane. Let ${\Cal S}_k:=\{c_1,\ldots  ,c_k\}$ be an ordered basis 
of $I_k$. By $(2.4)$ such a basis corresponds to an ordered
system $\{m^+_1,\ldots  ,m^+_k\}$, linearly independent over ${\Bbb R}$, of 
meromorphic functions in $L^2_q(M)$. Let $\{v_1,\ldots  ,v_k\}\subset H^1(M)$
be a system of the functions determined (up to additive constants) by the
orthogonal decompositions 
$$m^+_i=\partial^+_q v_i + \pi^-_q(m^+_i)\,\,, \eqnu $$
where $\pi^-_q: L^2_q(M)\to {\Cal M}^-_q$ is the orthogonal projection operator
given by the splitting $(2.3)$. Let $A^k$, $H^k$, $B^k$ and $V^k$ be the $k\times k$ 
matrices defined as follows: 
$$\eqalign{A^k_{ij} &:=(m_i,m_j)_q\,\,, \cr
H^k_{ij} &:=(\pi^-_q(m^+_i),\pi^-_q(m^+_j))_q \,\,, \cr
B^k_{ij} &:=B_q(m^+_i,m^+_j):=(m^+_i,{\overline {m^+_j}})_q\,\,,\cr
V^k_{ij} &:= (\partial^+_q v_i,\partial^+_q {\overline{v_j}})_q=
(\partial^-_q v_i,\partial^-_q {\overline{v_j}})_q\,\,. \cr} 
\eqnu $$ 
Such matrices depend on the pair $(q,{\Cal S}_k)$, where $q\in Q^{(1)}_{\kappa}$
and ${\Cal S}_k$ is an ordered linearly independent isotropic system of length
$k\in \{1,\ldots  ,g\}$ of real cohomology classes. However they are invariant under
the natural action of the pure mapping class group $\Gamma_g$. By the isotropy
property and $(2.5')$, the matrix $A^k$ is {\it real }and symmetric, while $H^k$
is hermitian and $B^k$, $V^k$ are complex and symmetric.

Our goal is to compute the evolution of the $k$-dimensional volume on
any isotropic plane $I_k$ under the action of the Kontsevich-Zorich cocycle. Let\break
$q\in Q^{(1)}_{\kappa}$ and let ${\Cal S}_k:=\{c_1,\ldots  ,c_k\}$ be any linearly 
independent ordered isotropic system of cohomology classes.  

\proclaim{Lemma} Let $q_t:=G_t(q)${\rm ,} ${\Cal S}_k(t):=G^{KZ}_t({\Cal S}_k)$
and $A^k_t${\rm ,} $H^k_t${\rm ,} $B^k_t$ and $V^k_t$ be  the matrices associated to the pair
$\bigl(q_t,{\Cal S}_k(t)\bigr)$ according to $(5.2)${\rm .} The following formulas 
hold\/{\rm :}
$$\eqalign{ {d\over{dt}} \det(A^k_t)=&\  -2 \det(A^k_t) 
\hbox{\rm tr}[(A^k_t)^{-1}\Re B^k_t]  \, , \cr
{d^2\over{dt^2}} \det(A^k_t)= &\ 4\det(A^k_t)\,\{\hbox{\rm tr}^2
[(A^k_t)^{-1}\Re B^k_t] \cr
&- \ \hbox{\rm tr}[(A^k_t)^{-1}\Re B^k_t]^2+
\hbox{\rm tr}[(A^k_t)^{-1}\Re(H^k_t - V^k_t)]\}\,\,. \cr } \eqnu $$
\endproclaim

\demo{Proof} By the standard formulas for the derivatives of a determinant:
$$\eqalign{ {d\over{dt}}\det(A^k_t)=&\ \det(A^k_t)\,\hbox{\rm tr}
[(A^k_t)^{-1}{d\over{dt}}A^k_t]\,\,, \cr
{d^2\over{dt^2}}\det(A^k_t)=&\ \det(A^k_t)\{\hbox{\rm tr}^2
[(A^k_t)^{-1}{d\over{dt}}A^k_t] \cr
&+ \ \hbox{\rm tr}[(A^k_t)^{-1}{d^2\over{dt^2}}A^k_t
-\bigl((A^k_t)^{-1}{d\over{dt}}A^k_t\bigr)^2]\}\,\,. \cr } \eqnu $$
In addition, by the variational formulas $(2.12)$ for the Kontsevich-Zorich cocycle, proved in Lemma 2.1, we can compute 
$$ \eqalign{ {d\over{dt}}A^k_t&= -2\Re B^k_t\,\,,\cr
    {d^2\over{dt^2}}A^k_t&= 4\Re (H^k_t-V^k_t)\,\,. \cr} \speqnu{5.4'}$$
The formulas in $(5.3)$ are then obtained by $(5.4)$ and $(5.4')$. 
\enddemo  

We can then compute the following generalization of the formula of
Lemma 3.2. Let $q\in Q^{(1)}_{\kappa}$ and $q_z$, $z\in D$, be given by 
$(3.1)$. Let $I_k\in G_k(M,{\Bbb R})$ be an isotropic $k$-plane and ${\Cal S}_k$
be any given ordered basis of $I_k$. Let $A^k_z$, $H^k_z$, $B^k_z$ and $V^k_z$
be the matrices associated to the pair $\bigl(q_z,{\Cal S}_k \bigr)$ according
to $(5.2)$.

\proclaim {Lemma} The following formulas hold\/{\rm :}
$$\triangle_h\log |\det(A^k_z)|^{1/2}= 4\hbox{\rm tr}
[(A^k_z)^{-1}H^k_z]-2\hbox{\rm tr}[(A^k_z)^{-1}B^k_z(A^k_z)^{-1}
{\overline{B^k_z}}]\,\,. \eqnu $$
\endproclaim

\demo{Proof} The argument follows the proof of Lemma 3.2. It suffices to
prove $(5.5)$ at the origin of the Poincar\'e disk $D$. Let $X_{\theta}$ be
the directional derivative at $z=0$, defined by $(3.9)$. By Lemma 5.1 and
the identities $(3.11)$ and $(3.11')$, we have:
$$\eqalign{ X_{\theta}\det(A^k_z)=&\  -2 \det(A^k_0)
\hbox{\rm tr}[(A^k_0)^{-1}\Re (e^{-i\theta} B^k_0)]\,\,, \cr
X_{\theta}^2\det(A^k_z)=&\  4\det(A^k_0)\,\{\hbox{\rm tr}^2
[(A^k_0)^{-1}\Re(e^{-i\theta} B^k_0)] \cr
&-\ \hbox{\rm tr}[(A^k_0)^{-1}\Re(e^{-i\theta}
 B^k_0)]^2+\hbox{\rm tr}[(A^k_0)^{-1}\Re(H^k_0-e^{-2i\theta}V^k_0)]\}\,\,. 
\cr } \eqnu $$
We can write the hyperbolic gradient and Laplacian in terms of the directional
derivatives as in $(3.12)$:
$$\eqalign{\nabla_h \det(A^k_z)=&\ \bigl(X_0\det(A^k_z),
X_{{\pi}\over 2}\det(A^k_z)\bigr) \cr
=&\ -2\det(A^k_0)(\hbox{\rm tr}[(A^k_0)^{-1}\Re B^k_0],
\hbox{\rm tr}[(A^k_0)^{-1}\Im B^k_0])\,\,, \cr
\triangle_h \det(A^k_z)=&\ {1\over{\pi}}\int_0^{2\pi} 
X_{\theta}^2\det(A^k_z)\,d\theta \cr
= &\ 4\det(A^k_0)\,\{\hbox{\rm tr}^2
[(A^k_0)^{-1}\Re B^k_0]+\hbox{\rm tr}^2[(A^k_0)^{-1}\Im B^k_0] \cr
&-\ \hbox{\rm tr}[(A^k_0)^{-1}\Re B^k_0]^2- \hbox{\rm tr}[(A^k_0)^{-1}\Im B^k_0]^2
+2\hbox{\rm tr}[(A^k_0)^{-1}\Re H^k_0]\}\,\,. \cr } \speqnu{5.6'}$$
The desired formula $(5.5)$ (at $z=0$) can be derived from $(5.6')$ by applying
the formula $(3.13)$ with $\phi(z):=|\det(A^k_z)|^{1/2}$. \enddemo

\specialnumber{5.2'}
\proclaim {Lemma} The formula
$$\Phi_k:= 2\hbox{\rm tr}[(A^k)^{-1}H^k]-\hbox{\rm tr}[(A^k)^{-1} B^k(A^k)^{-1}
{\overline{B^k}}]\,\,,\eqnu $$
where the matrices $A^k${\rm ,} $H^k$ and $B^k$ are given by $(5.2)${\rm ,} is independent of the choice
of a basis of the isotropic plane $I_k$ and therefore defines a non\/{\rm -}\/negative function 
$\Phi_k:Q^{(1)}_{\kappa}\times G_k(M,{\Bbb R})\to {\Bbb R}${\rm ,} invariant under the
natural action of the modular group $\Gamma_g${\rm .} In addition{\rm ,} the following identities hold{\rm .} 
Let $I_1\subset I_2\subset \cdots \subset I_g \subset H^1(M,{\Bbb R})$ be a finite sequence of 
isotropic subspaces such that $\hbox{dim}(I_k)=k\in \{1,\ldots  ,g\}${\rm .} Let $q\in Q^{(1)}_{\kappa}$ 
and $\{m^+_1,\ldots  ,m^+_g\}$ be an orthonormal basis of the space ${\Cal M}^+_q$ of meromorphic 
functions in $L^2_q(M)$ such that{\rm ,} for each $k\in \{1,\ldots  ,g\}${\rm ,} $\{m^+_1,\ldots  ,m^+_k\}$ represents 
an orthonormal basis of $I_k${\rm ,} in the sense that{\rm ,} if $c_i=[\Re(m^+_iq^{1/2})]\in H^1(M,{\Bbb R})${\rm ,} 
then $\{c_1,\ldots  ,c_k\}$ is an orthonormal basis of $I_k$ with respect to the Hodge norm induced by 
the metric $R_q${\rm .} Then    
\medbreak
\item{\rm (1)} $\Phi_1(q,I_1)=2|\pi^-_q(m^+_1)|^2-|B_q(m^+_1)|^2$\,{\rm ,} \smallbreak
\item{\rm (2)} $\Phi_k(q,I_k)= \Phi_g(q,I_g)- \sum_{i,j=k+1}^g |B_q(m^+_i,m^+_j)|^2$\,{\rm ,} \smallbreak
\item{\rm (3)} $\Phi_g(q,I_g)\equiv \Lambda_1(q)+\cdots +\Lambda_g(q)$\,{\rm .}
\endproclaim

\demo{Proof} The formula $(5.7)$ depends in principle on the choice of a
basis $\{c_1,\ldots  ,c_k\}$ of the isotropic subspace $I_k$. However, since the matrix of a 
change of base is a real matrix, a computation using the commutativity property of the
trace shows that in fact there is no dependence on the choice of the base. In particular,
if we choose an orthonormal basis, then  by the isotropy of $I_k$, $(m^+_i,m^+_j)_q=
\delta_{ij}$ and we obtain the following formula:
$$ \Phi_k(q,I_k)= 2 \sum_{i=1}^k |\pi^-_q(m^+_i)|^2 - \sum_{i,j=1}^k |B_q(m^+_i,m^+_j)|^2\,\,.
\speqnu{5.7'}$$
The identity $(1)$ is a particular case of $(5.7')$. The identity $(2)$ can be computed on 
the basis of $(5.7')$ taking into account the following identity, which holds since
$\{m^+_1,\ldots  ,m^+_g\}$, and hence $\{{\overline {m^+_1}},\ldots  ,{\overline {m^+_g}}\}$, are 
orthonormal bases for the spaces of meromorphic, respectively anti-meromorphic, functions
in $L^2_q(M)$:
$$\pi^-_q(m^+_i)= \sum_{j=1}^g B_q(m^+_i,m^+_j)\,{\overline {m^+_j}} \,\,. \eqnu $$
Finally, if $k=g$, by $(5.7')$ and $(5.8)$,
$$\Phi_g(q,I_g)=\sum_{i=1}^g |\pi^-_q(m^+_i)|^2 =\Lambda_1(q)+\cdots +\Lambda_g(q)\,\,,
\eqnu $$
which proves $(3)$. The invariance of $\Phi_k$ under the natural action of the modular group follows
from the invariance of the matrices $A^k$, $H^k$, $B^k$, which enter in the definition $(5.7)$. 
\enddemo

Let $q\in Q^{(1)}_{\kappa}$, $z\equiv (t,\theta)\in D$ in geodesics polar 
coordinates on the\break Teichm\"uller disk centered at $q$ and let $q_z$ be defined as in $(3.1)$.
Let $G^{KZ}_z$ be given by the (trivial) parallel transport of cohomology classes from the 
origin at $q$ to $q_z$ along the Teichm\"uller disk centered at $q$ in $Q^{(1)}_{\kappa}$.   
By Lemma 3.1 and Lemma 5.2, 
$${1\over {2\pi}} {{\partial}\over {\partial t}}\int_0^{2\pi} 
\log |\det A^k_z|^{1/2}\,d\theta= \tanh(t) {1\over {|D_t|}} \int_{D_t}
\Phi_k\circ G^{KZ}_{\xi}(q,I_k)\,\omega_P(\xi) \,\,. \eqnu $$

As a first consequence of $(5.10)$, we prove the following version of the 
formula for the sum of the first $g$ exponents, given by M. Kontsevich and\break A. Zorich \ref\KZone:

\specialnumber{5.3}
\proclaim {Corollary} The Lyapunov exponents of the Kontsevich\/{\rm -}\/Zorich cocycle with respect to the
ergodic invariant measure given by the restriction to any connected component ${\Cal C}^{(1)^{\phantom{|}}}
_{\kappa}$ of a stratum ${\Cal M}^{(1)}_{\kappa}$ satisfy the following identity\/{\rm :}
$$\lambda_1+\cdots +\lambda_g= {1\over {\mu^{(1)}_{\kappa}}({\Cal C}^{(1)}_{\kappa})}
\,\int_{{\Cal C}^{(1)}_{\kappa}} (\Lambda_1+\cdots +\Lambda_g)\,d\mu^{(1)}_{\kappa}\,\,.
\eqnu $$
\endproclaim

\demo{Proof} Let $\sigma^k_q$ be the normalized canonical (Haar) measure on the Grassmannian
$G_k(M,{\Bbb R})$ of isotropic $k$-dimensional subspaces $I_k\subset H^1(M_q,{\Bbb R})$, endowed
with the Euclidean structure given by the Hodge inner product induced by the metric $R_q$. The
probability measure $\sigma^k_q$ is invariant under the action of the circle group ${\rm SO}(2,{\Bbb R})$
on $Q_g$. In the case $k=g$, since the function $\Phi_g$ does not depend on the Lagrangian plane
$I_g$, by averaging $(5.10)$ over $G_k(M,{\Bbb R})$ with respect to $\sigma_q^g$, then over a stratum 
${\Cal M}^{(1)}_{\kappa}$ with respect to any ergodic $G_t$-invariant measure $\mu$, we obtain: 
$$\eqalign{ {1\over {2\pi}} {{\partial}\over {\partial t}}\int_{{\Cal M}^{(1)}_{\kappa}}\int_{G_k(M,{\Bbb R})}
\int_0^{2\pi}&\log |\det A^g_z|^{1/2}\,d\theta d\sigma^g_q d\mu  \cr
&=\tanh(t) \int_{{\Cal M}^{(1)}_{\kappa}}(\Lambda_1+\cdots +\Lambda_g)\,d\mu \,\,. \cr} \eqnu $$
If $\mu$ is the normalized restriction to any connected component of ${\Cal M}^{(1)}_{\kappa}$ of 
the invariant measure $\mu^{(1)}_{\kappa}$, which is absolutely continuous and invariant under the 
action of the circle group on $Q^{(1)}_{\kappa}$, by applying Fubini's theorem to the left-hand
side of $(5.12)$ and averaging over $[0,{\Cal T}]$ with respect to time, we have:
$$\eqalign{ {1\over {\Cal T}} \int_{{\Cal M}^{(1)}_{\kappa}}\int_{G_k(M,{\Bbb R})} \log |\det &
\bigl(A^g_{\Cal T}(A^g_0)^{-1}\bigr)|^{1/2}\,d\sigma^k_q d\mu  \cr
= & {{\log\cosh {\Cal T}}\over {\Cal T}}
\int_{{\Cal M}^{(1)}_{\kappa}}(\Lambda_1+\cdots +\Lambda_g)\,d\mu \,\,. \cr} \speqnu{5.12'}$$
The identity $(5.12')$ yields $(5.11)$ in the limit ${\Cal T}\to +\infty$. In fact, by Oseledec's
theorem, the limit of the left-hand side is equal to the sum $\lambda_1+\cdots +\lambda_g$.     
\enddemo

\specialnumber{5.3'}
\demo{Remark} Since $\Lambda_1\equiv \lambda_1=1$, the formula $(5.11)$ is equivalent
to the following:
$$\lambda_2+\cdots +\lambda_g= {1\over {\mu^{(1)}_{\kappa}}({\Cal C}^{(1)}_{\kappa})}
\,\int_{{\Cal C}^{(1)}_{\kappa}} (\Lambda_2+\cdots +\Lambda_g)\,d\mu^{(1)}_{\kappa}\,\,.
\eqnu $$ 
In particular, such formula implies that $\lambda_2>0$ (Corollary $4.5'$) and it gives an
exact formula for the second Lyapunov exponent in the case $g=2$.
\enddemo 

The formulas $(5.10)$ directly imply the non-vanishing of about a half of Lyapunov
exponents in case $g>2$.

\specialnumber{5.4} 
\proclaim{Corollary} Let us consider the Lyapunov spectrum of the
Kontsevich\/{\rm -}\/Zorich cocycle
with respect to the restriction of the measure $\mu^{(1)}_{\kappa}$ to any connected 
component of ${\Cal M}^{(1)}_{\kappa}${\rm .} If $k<g/2${\rm ,}
$$1\equiv\lambda_1>\lambda_2\geq\cdots\geq \lambda_k\geq \lambda_{k+1}>0\,\,. \eqnu $$
\endproclaim

\demo{Proof} Let $k<g/2$ and let
$$M_k(q):= \sup\{\Phi_k(q,I_k)\,|\, I_k\in G_k(M,{\Bbb R})\} \,\,. \eqnu $$
We claim that any connected component of the stratum ${\Cal M}^{(1)}_{\kappa}$ has an
open subset where the following inequality holds:
$$M_k(q)<\Lambda_1(q)+\cdots +\Lambda_g(q) \,\,. \speqnu{5.15'}$$
In fact, if $(5.15')$ fails, by Lemma $5.2'$ and the compactness of the Grassmannian, there
exists an orthonormal system $\{m^+_1,\ldots  ,m^+_g\}\subset {\Cal M}^+_q$ such that $B_q(m^+_i,m^+_j)=0$ for all
$i,j\in \{k+1,\ldots  ,g\}$. By $(5.8)$, if the projection $\pi^-_q$ is injective, then
$g-k\leq k$. Since $k<g/2$ and, by Theorem 4.5, $\pi^-_q$ is injective on an open subset
of any connected component of ${\Cal M}^{(1)}_{\kappa}$, the claim is proved. 
\smallskip
\noindent By formula $(5.10)$,
$${1\over {2\pi}} {{\partial}\over {\partial t}}\int_0^{2\pi} 
\log |\det A^k_z|^{1/2}\,d\theta\leq \tanh(t) {1\over {|D_t|}} \int_{D_t}
M_k(q_{\xi}) \omega_P(\xi) \,\,. \eqnu $$
Then, by averaging successively over $G_k(M,{\Bbb R})$ with respect to $\sigma^k_q$, over a connected 
component ${\Cal C}^{(1)}_{\kappa}\subset {\Cal M}^{(1)}_{\kappa} $ with respect to the normalized 
restriction $\mu$ of the measure $\mu^{(1)}_{\kappa}$ and over the interval $[0,{\Cal T}]$ with 
respect to time, we obtain:
$$\eqalign{ {1\over {\Cal T}}\int_{{\Cal M}^{(1)}_\kappa}\int_{G_k(M,{\Bbb R})}\log |\det
\bigl(A^g_{\Cal T}(A^g_0)^{-1}\bigr)|^{1/2}&\,d\sigma^k_q\,d\mu \cr
\leq & {{\log\cosh {\Cal T}}\over {\Cal T}} \int_{{\Cal M}^{(1)}_\kappa}M_k(q)\,d\mu \,.\cr}
\speqnu{5.16'}$$
Finally, by taking the limit as ${\Cal T}\to +\infty$, we have:
$$\eqalign{ \lambda_1+\cdots +\lambda_k &\leq {1\over {\mu^{(1)}_{\kappa}({\Cal C}^{(1)}_{\kappa})}}
   \int_{{\Cal C}^{(1)}_{\kappa}} M_k(q)\,d\mu^{(1)}_{\kappa}  \cr
   &<{1\over {\mu^{(1)}_{\kappa}({\Cal C}^{(1)}_{\kappa})}}\int_{{\Cal C}^{(1)}_{\kappa}}
   (\Lambda_1+\cdots+\Lambda_g)\,d\mu^{(1)}_{\kappa}=\lambda_1+\cdots +\lambda_g\,\,. \cr} \eqnu $$
Hence, $\lambda_{k+1}>0$. In fact, a strictly positive lower bound can be derived from $(5.17)$. 
\enddemo

The investigation of higher Lyapunov exponents requires some kind of control on the
invariant sub-bundles of the Kontsevich-Zorich cocycle. A result in that direction will be
proved in the next section. We conclude this section by deriving one more consequence of
formula $(5.10)$, an exact formula for the exponents in terms of the invariant sub-bundles.

\specialnumber{5.5}
\proclaim{Corollary} Assume $\lambda_k>\lambda_{k+1}\geq 0$, $k\in \{1,\ldots  ,g-1\}${\rm .} Let
$E^{\pm}_k\subset {\Cal H}^1_{\kappa}(M,{\Bbb R})$ be the $k$\/{\rm -}\/dimensional invariant sub\/{\rm -}\/bundles 
corresponding respectively to the sets of Lyapunov exponents $\pm \{\lambda_1,\ldots  ,\lambda_k\}$
in the Oseledec\/{\rm '}\/s splitting of the Kontsevich\/{\rm -}\/Zorich cocycle{\rm .} Then the following formula holds\/{\rm :}\/
$$\lambda_1+\cdots +\lambda_k={1\over {\mu^{(1)}_{\kappa}({\Cal C}^{(1)}_{\kappa})}}
\int_{{\Cal C}^{(1)}_{\kappa}} \Phi_k\bigl(q,E^+_k(q)\bigr)\,d\mu^{(1)}_{\kappa}\,\,. 
\eqnu $$
\endproclaim

\demo{Proof} By Oseledec's theorem, for almost all $q\in {\Cal M}^{(1)}_{\kappa}$, almost
all\break $\theta\in [0,2\pi]$ and $\sigma^k_q$-almost all $k$-dimensional isotropic subspaces 
$I_k\in G_k(M,{\Bbb R})$,
$$ \hbox{dist} \bigl(G^{KZ}_{(t,\theta)}(q,I_k),G^{KZ}_t(q_{\theta},E^+_k(q_{\theta}))\bigr) 
\,\,\to\,\, 0\,\,, \eqnu $$
as $t\to +\infty$. Consequently, since the function $\Phi_k$ is bounded and continuous,
by Fubini's theorem and the dominated convergence theorem, 
$$\eqalign{ {1\over {|D_t|}} \int_{G_k(M,{\Bbb R})}\int_{D_t} |\Phi_k\circ &G^{KZ}_{(s,\theta)}(q,I_k)
  \cr -\,&\Phi_k\circ G^{KZ}_s\bigl(q_{\theta},E^+_k(q_{\theta})\bigr)|\,\omega_P(s,\theta)\,
d\sigma^k_q\,\,\to \,\, 0\,,\cr} \speqnu{5.19'}$$
as $t\to +\infty$, for almost all $q\in {\Cal M}^{(1)}_{\kappa}$. Hence, by averaging $(5.10)$ 
over $G_k(M,{\Bbb R})$ with respect to the measure $\sigma^k_q$, then with respect to the 
normalized restriction $\mu$ of the invariant measure $\mu^{(1)}_{\kappa}$ to a connected component 
${\Cal C}^{(1)}_{\kappa}\subset{\Cal M}^{(1)}_{\kappa}$ and applying Fubini's theorem, we find 
that, since $E^+_k$ is an invariant sub-bundle of $G^{KZ}$ and $\mu$ is ${\rm SL}(2,{\Bbb R})$-invariant,  
$$\eqalign{ {\partial\over{\partial t}} \int_{{\Cal M}^{(1)}_{\kappa}}\int_{G_k(M,{\Bbb R})}\log 
&|\det(A^k_t)|^{1/2} \,d\sigma^k_q d\mu   \cr -\,&\tanh(t) \int_{{\Cal M}^{(1)}_{\kappa}} 
\Phi_k\bigl(q,E^+_k(q)\bigr)\,d\mu\,\,\to \,\, 0\,,\cr} \eqnu $$
as $t\to +\infty$. Finally, by averaging over $[0,{\Cal T}]$ with respect to time and by 
Oseledec's theorem, we obtain $(5.18)$. \enddemo

 \section{Basic currents for measured foliations}

  The invariant unstable and stable subspaces of the Kontsevich-Zorich cocycle at a quadratic
differential $q\in {\Cal M}^{(1)}_{\kappa}$ will be described in terms of distributional invariants 
of, respectively, its horizontal and vertical foliations, ${\Cal F}_q$ and ${\Cal F}_{-q}$. Such
invariants will be defined as a distributional generalization of the notion of {\it basic form}, 
well known in the geometric theory of foliations \ref\Reone, \ref\Retwo, [63, Chap.\ 4], [2, \S 1.5 \& 
\S 7] and will therefore be called {\it basic currents}. In this section we prove a few 
properties of basic currents and of the {\it weighted Sobolev spaces }introduced in [18, \S 2]
which will be relevant in the study of the Kontsevich-Zorich cocycle. The most crucial is a version 
of the Poincar\'e inequality with an explicit geometric estimate for the constant. Such an estimate, 
equivalent to a lower bound for the first non-trivial eigenvalue of the Dirichlet form of the metric 
$R_q$, is obtained following Cheeger's method \ref\Cg, \ref\BGM. 

\demo{{\rm 6.1.} Basic currents and invariant distributions}
  Let $\Cal F$ be a measured foliation on a compact orientable surface $M$ of genus $g\geq 
2$ in the sense of Thurston \ref\Th, [14]. Let $\Sigma \subset M$ be its (finite) set of 
singularities. 
\enddemo

\numbereddemo{Definition} A basic current for ${\Cal F}|_{M\setminus\Sigma}$ is a  
current $C\in {\Cal D}'(M\setminus\Sigma)$ (in the sense of de Rham \ref\dR, see also [56, Chap.\ IX]), 
homogeneous of dimension (and degree) equal to $1$, with the following property:
$$\imath_XC={\Cal L}_XC=0\,\,, \eqnu $$
for all smooth vector fields $X$ tangent to $\Cal F$, with compact support in $M\setminus\Sigma$.
\enddemo

The operators of contraction and Lie derivative, with respect to a smooth vector 
field $X$, denoted  $\imath_X$ and ${\Cal L}_X$ respectively, are extended to currents in the 
standard distributional (weak) sense [56, Chap.\ IX, \S 3]. The vector space of all {\it real}
basic currents for ${\Cal F}|_{M\setminus\Sigma}$ will be denoted by ${\Cal B}_{\Cal F}(M\setminus
\Sigma)$. 

 Let $q\in Q_{\kappa}$ be an orientable (holomorphic) quadratic differential on a (marked)
Riemann surface $M_q$. Let $\Sigma_q$ be the (finite) set of its zeroes. We introduce the following 
space $\Omega_q(M)$ of smooth test forms on $M$. Let $p\in \Sigma_q$ be a zero of (even) order $k\in 
{\Bbb N}$. There exists a canonical complex coordinate $z:{\Cal U}_p\to {\Bbb C}$ on a neighbourhood 
${\Cal U}_p$ such that $p\in {\Cal U}_p$, $z(p)=0$ and $q=z^k dz^2$ on ${\Cal U}_p$. Let $\pi_p:
{\Cal U}_p \to {\Bbb C}$ be the (local) covering map defined by $\pi_p(z):=z^{k/2+1}/(k/2+1)$. A form
$\alpha\in\Omega_q(M)$ if and only if, for all $p\in\Sigma_q$, there exists a smooth form $\lambda_p$
on a neighbourhood of $0\in {\Bbb C}$ such that $\alpha=\pi_p^{\ast}(\lambda_p)$ on ${\Cal U}_p'
\subset {\Cal U}_p$. The dual space ${\Cal S}'_q(M)$ of $\Omega_q(M)$ will be called the space of 
$q$-{\it tempered currents }on $M$. A {\it homogeneous }tempered current of dimension $d$ will be a 
continuous functional on the subset $\Omega^d_q(M)\subset\Omega_q(M)$ of homogeneous forms of degree 
$d\in {\Bbb N}$. 

Let ${\Cal V}_q(M)$ be the space of vector fields $X$ on 
$M\setminus\Sigma_q$ such that $\imath_X
\alpha$, ${\Cal L}_X\alpha\in \Omega_q(M)$ for all $\alpha\in \Omega_q(M)$. 

\demo{Definition $6.1'$} A current $C\in {\Cal S}'_q(M)$, of degree and dimension equal to $1$,
is basic for ${\Cal F}_{\pm q}$ if the identities $(6.1)$ holds in ${\Cal S}'_q(M)$ for all 
$X\in {\Cal V}_q(M)$, tangent to ${\Cal F}_{\pm q}$ on $M\setminus\Sigma_q$. The vector 
spaces of real ${\Cal F}_{\pm q}$-basic currents will be denoted, respectively, by ${\Cal B}_{\pm q}
(M)$. \enddemo

 Let ${\Cal Z}'(M\setminus\Sigma_q)\subset {\Cal D}'(M\setminus\Sigma_q)$ be the vector 
space of {\it closed }real homogeneous currents of dimension (and degree) equal to $1$. A current $C\in
{\Cal D}'(M\setminus\Sigma_q)$ is closed if the identity $dC=0$ holds in ${\Cal D}'(M\setminus
\Sigma_q)$ [11, Chap.\ IV, \S 18]. By the generalized de Rham theorem [11, \S 15, Th.\ 12], [56, 
Chap.\ IX, \S 3, Th.\ I], there exists a natural cohomology map 
$$j_q:{\Cal Z}'(M\setminus\Sigma_q)\to H^1(M\setminus\Sigma_q,{\Bbb R})\,\,.\eqnu $$
In fact, any closed current $C\in {\Cal Z}'(M\setminus\Sigma_q)$ defines a linear functional on the 
de Rham cohomology with compact supports $H^1_c(M\setminus\Sigma_q, {\Bbb R})$ and, by Poincar\'e
duality, $H^1_c(M\setminus\Sigma_q,{\Bbb R})^{\ast}\equiv H^1(M\setminus\Sigma_q,{\Bbb R})$. Let 
${\Cal Z}_q'(M)\subset {\Cal Z}'(M\setminus\Sigma_q)$ be the subspace of $q$-tempered currents $C$
of dimension (and degree) equal to $1$ which are {\it closed }in ${\Cal S}_q'(M)$, in the sense that 
$dC=0$ holds in ${\Cal S}_q'(M)$.

\proclaim {Lemma} The cohomology map $(6.2)$ has the property that
$$j_q: {\Cal Z}_q'(M)\to H^1(M,{\Bbb R})\subset H^1(M\setminus\Sigma_q,{\Bbb R})\,\,.\speqnu{6.2'}$$
\endproclaim

\demo{Proof} Let $C\in {\Cal Z}_q'(M)$. Then $C\in {\Cal Z}'(M\setminus\Sigma_q)$ and the cohomology
class $j_q(C)\in H^1(M\setminus\Sigma_q,{\Bbb R})$ is therefore well-defined. Let $c\in H^1_c(M
\setminus\Sigma_q,{\Bbb R})$ belong to the kernel of the (surjective) canonical map 
$$\epsilon: H^1_c(M\setminus\Sigma_q,{\Bbb R})\to H^1(M,{\Bbb R}) \,\,. \eqnu $$
Then $c=[dv]$, where $v$ is a smooth function on $M$ constant in a neighbourhood of each point $p\in 
\Sigma_q$. Since $v\in\Omega^0_q(M)$ and $C$ is closed in ${\Cal S}'_q(M)$, $C(dv)=dC(v)=0$. Hence 
$j_q(C)\in N(\epsilon)^{\perp}$, the annihilator of the kernel $N(\epsilon)$. The dual map 
$$\epsilon^{\ast}:H^1(M,{\Bbb R})^{\ast}\to H^1(M\setminus\Sigma_q,{\Bbb R})  \speqnu{6.3'}$$ 
coincides with the canonical injection $H^1(M,{\Bbb R})\to H^1(M\setminus\Sigma_q,{\Bbb R})$, under 
the identification $H^1(M,{\Bbb R})^{\ast}\equiv H^1(M,{\Bbb R})$ given by the intersection form on 
$H^1(M,{\Bbb R})$. Since $N(\epsilon)^{\perp}$ is equal to the range of the dual map $\epsilon^{\ast}$, 
$j_q(C)\in H^1(M,{\Bbb R})$.
\enddemo     

 Let ${\Cal B}_{\pm q}(M\setminus\Sigma_q):={\Cal B}_{\Cal F}(M\setminus\Sigma_q)$, ${\Cal F}
={\Cal F}_{\pm q}$.

\specialnumber{6.2'} 
\proclaim {Lemma} The following inclusions hold\/{\rm :}\/
\smallbreak
\item{\rm (1)} ${\Cal B}_{\pm q}(M)\subset {\Cal B}_{\pm q}(M\setminus\Sigma_q)${\rm ,} \smallbreak
\item{\rm (2)} ${\Cal B}_{\pm q}(M\setminus\Sigma)\subset {\Cal Z}'(M\setminus\Sigma_q)${\rm ,} \smallbreak
\item{\rm (3)} ${\Cal B}_{\pm q}(M)\subset {\Cal Z}_q'(M)${\rm .}
\endproclaim

\demo{Proof} $(1)$ follows directly from the definitions. In fact, since smooth forms with compact
support in $M\setminus\Sigma_q$ belong to the space $\Omega_q(M)$, the space of $q$-tempered currents ${\Cal
S}_q'(M)\subset {\Cal D}' (M\setminus\Sigma_q)$ and vector fields with compact support in $M\setminus\Sigma_q$ belong to
the  space ${\Cal V}_q(M)$. 

\smallbreak $(2)$ The identity
$$ {\Cal L}_X C \,=\, \imath_X dC + d( \imath_X C)  \eqnu $$
holds for currents [56, (IX,3;32)] in ${\Cal D}'(M\setminus\Sigma_q)$ if $C\in {\Cal D}'(M
\setminus\Sigma_q)$ and $X$ is a vector field with compact support in $M\setminus\Sigma_q$. Hence, 
if $C\in {\Cal B}_{\pm q}(M\setminus\Sigma_q)$, $\imath_X dC = 0$ in ${\Cal D}'(M\setminus\Sigma_q)$ 
for any vector field $X$, with compact support in $M\setminus\Sigma_q$, tangent to ${\Cal F}_{\pm q}$.
Since $C$ has dimension (and degree) equal to $1$, $dC$ has dimension $0$ (and degree $2$). It follows
that $dC=0$, hence $C$ is closed, in ${\Cal D}'(M\setminus\Sigma_q)$.
\smallbreak
$(3)$ The identity $(6.4)$ holds in ${\Cal S}'_q(M)$ if $C\in {\Cal S}'_q(M)$ and $X\in {\Cal
V}_q(M)$. If $C\in {\Cal B}_{\pm q}(M)$, it follows by $(6.4)$ that $\imath_X dC=0$ for all 
vector fields $X\in {\Cal V}_q(M)$, tangent to ${\Cal F}_{\pm q}$. Hence $dC=0$ in ${\Cal S}'_q(M)$.
\enddemo

\numbereddemo{ Definition} The (first) distributional basic cohomology (with real coefficients) 
of the measured foliations ${\Cal F}_{\pm q}|_{M\setminus\Sigma_q}$, ${\Cal F}_{\pm q}$, denoted
respectively by $H^1_{\pm q}(M\setminus\Sigma_q,{\Bbb R})$, $H^1_{\pm q}(M,{\Bbb R})$, are 
the images, under the cohomology map $(6.2)$, of the spaces of basic currents ${\Cal B}_{\pm q}(M
\setminus\Sigma_q)$, ${\Cal B}_{\pm q}(M)$: 
$$ \eqalign { H^1_{\pm q} (M\setminus\Sigma_q,{\Bbb R}) & :=j_q \bigl( {\Cal B}_{\pm q}(M\setminus
\Sigma_q)\bigr) \subset H^1(M\setminus\Sigma_q,{\Bbb R})\,\, , \cr
H^1_{\pm q}(M,{\Bbb R}) & :=j_q \bigl( {\Cal B}_{\pm q}(M) \bigr) \subset H^1(M,{\Bbb R})\,\,. \cr}
\eqnu $$  
\enddemo

 Basic currents in ${\Cal B}_q(M)$ [${\Cal B}_{-q}(M)$] are related to the invariant 
distributions for the vector field $S$ [$T$], constructed in [18, \S 4] as obstructions to the 
existence of smooth solutions to the cohomological equation $Su=f$ [$Tu=f$]. They therefore have a 
dynamical significance as they give obstructions to triviality of time changes for flows tangent to 
the foliation. In fact, time-change triviality is equivalent to the solvability of the cohomological 
equation [32, Def.\ 2.2.3 \& \S 2.9]. 

\numbereddemo{Definition} A quasi $S$-invariant [a quasi $T$-invariant] distribution\break ${\Cal D}\in
{\Cal D}'(M\setminus\Sigma)$ is a distributional solution of the equation $S{\Cal D}=0$ [$T{\Cal D}
=0$] in ${\Cal D}'(M\setminus\Sigma)$. An $S$-invariant [a $T$-invariant] distribution is a solution 
${\Cal D}\in{\Cal S}_q'(M)$ of the equation $S{\Cal D}=0$ [$T{\Cal D}=0$] in ${\Cal S}'_q(M)$. The 
vector space of quasi $S$-invariant [quasi $T$-invariant] distributions will be denoted by ${\Cal I}
_q(M\setminus\Sigma_q)$ [${\Cal I}_{-q}(M\setminus\Sigma_q)$]. The subspace of $S$-invariant 
[$T$-invariant] distributions will be denoted by ${\Cal I}_q(M)$ [${\Cal I}_{-q}(M)$].
\enddemo

A distribution $\Cal D$ will be considered as a current of degree $2$ (and dimension $0$). 
However, it is possible to identify distributions with currents of degree $0$ (and dimension $2$) by 
the isomorphism ${\Cal D}\to {\Cal D}^{\ast}$ determined by the area form $\omega_q$ on the 
(open) manifold $M\setminus\Sigma$ [56, Chap.\ IX, \S 2]. Given any distributions ${\Cal D}
\in {\Cal D}'(M\setminus\Sigma_q)$, let ${\Cal D}^{\ast}$ be the current of degree $0$ (and dimension
$2$) uniquely determined by the identity ${\Cal D}={\Cal D}^{\ast}\omega_q$. It is then possible to 
define the (exterior) products ${\Cal D}^{\ast}\eta_S$ or ${\Cal D}^{\ast}\eta_T$ as in [11, p.42], 
[56, Chap.\ IX, \S 3]. Such products are the currents, of dimension (and degree) equal to 
$1$, given by the formulas
$$\eqalign{ {\Cal D}^{\ast}\eta_S&:={\Cal D}^{\ast}\wedge\eta_S =\imath_S{\Cal D}\,\,, \cr 
{\Cal D}^{\ast}\eta_T&:={\Cal D}^{\ast}\wedge\eta_T=-\imath_T{\Cal D} \,\,.\cr} \eqnu $$

\proclaim {Lemma} The maps $C^{\pm}_q:{\Cal I}_{\pm q}(M\setminus\Sigma_q)\to {\Cal B}_{\pm q}
(M\setminus\Sigma_q)$ given by the formulas
$$\eqalign{ C^+_q({\Cal D})&:= {\Cal D}^{\ast}\eta_S\,\,,\,\,\,\,{\Cal D}\in {\Cal I}_q(M\setminus
                                \Sigma_q)\,\,, \cr
            C^-_q({\Cal D})&:= {\Cal D}^{\ast}\eta_T\,\,,\,\,\,\,{\Cal D}\in {\Cal I}_{-q}(M\setminus
                                \Sigma_q)\,\,, \cr } \eqnu $$
are bijective{\rm .} The inverse mappings  are given by the formulas
$$\eqalign { C^+ &\to -C^+\wedge\eta_T\in {\Cal I}_q(M\setminus\Sigma_q)\,\,, \cr
             C^- &\to C^-\wedge\eta_S\in {\Cal I}_{-q}(M\setminus\Sigma_q)\,\,.\cr} \speqnu{6.7'}$$
In addition $C^{\pm}_q({\Cal D})\in {\Cal B}_{\pm q}(M)\subset {\Cal B}_{\pm q}(M\setminus\Sigma_q)$
if and only if ${\Cal D}\in {\Cal I}_{\pm q}(M)\subset {\Cal I}_{\pm q}(M\setminus\Sigma_q)${\rm .} 
\endproclaim

\demo{Proof} Let $X$ be a smooth vector field tangent to ${\Cal F}_q$ [${\Cal F}_{-q}$].
Then, since ${\Cal F}_{\pm q}$ are $1$-dimensional, $X=WS$ [$X=WT$], where $W$ is a smooth function 
with compact support in $M\setminus\Sigma_q$ if $X$ has compact support and $W\in\Omega^0_q(M)$ if 
$X\in {\Cal V}_q(M)$. Since, by $(6.6)$, ${\Cal D}^{\ast}\eta_S=\imath_S{\Cal D}$ [${\Cal D}^{\ast}
\eta_T=-\imath_T{\Cal D}$] and $\imath_S^2=\imath_T^2=0$, 
$$\eqalign{ \imath_X({\Cal D}^{\ast}\eta_S)&=W\imath_S({\Cal D}^{\ast}\eta_S)= W\imath_S^2{\Cal D}
                                            =0\,\, \cr
[\imath_X({\Cal D}^{\ast}\eta_T)&=W\imath_T({\Cal D}^{\ast}\eta_T)=-W\imath_T^2{\Cal D}
                                            =0\,]\,,\cr} \eqnu $$
in ${\Cal D}'(M\setminus\Sigma)$ if ${\Cal D}\in {\Cal D}'(M\setminus\Sigma)$, in ${\Cal S}_q'(M)$ 
if $\Cal D\in {\Cal S}_q'(M)$.

 Since $\eta_S$ [$\eta_T$] is closed, if ${\Cal D}\in {\Cal I}_q(M\setminus\Sigma_q)$
[${\Cal D}\in {\Cal I}_{-q}(M\setminus\Sigma_q)$], we have 
$$\eqalign {d({\Cal D}^{\ast}\eta_S)&=d({\Cal D}^{\ast}\wedge\eta_S)=S{\Cal D}=0 \cr
         [d({\Cal D}^{\ast}\eta_T)&=d({\Cal D}^{\ast}\wedge\eta_T)=-T{\Cal D}=0] \cr} \speqnu{6.8'}$$
in ${\Cal D}'(M\setminus\Sigma)$. Hence, by $(6.4)$, ${\Cal D}^{\ast}\eta_S\in {\Cal B}_q(M\setminus
\Sigma_q)$ [${\Cal D}^{\ast}\eta_T\in {\Cal B}_{-q}(M\setminus\Sigma_q)$]. If ${\Cal D}\in {\Cal I}_q
(M)$ [${\Cal D}\in {\Cal I}_{-q}(M)$] the identities $(6.8')$ hold in ${\Cal S}_q'(M)$. Hence ${\Cal 
D}^{\ast}\eta_S\in {\Cal B}_q(M)$ [${\Cal D}^{\ast}\eta_T\in {\Cal B}_{-q}(M)$].

 Conversely, let $C\in{\Cal B}_q(M\setminus\Sigma_q)$ [$C\in {\Cal B}_{-q}(M\setminus
\Sigma_q)$]. Let $\Cal D\in {\Cal D}'(M\setminus\Sigma)$ be the distribution given by
$$\eqalign{ {\Cal D}& := -C\wedge \eta_T  \cr
             [{\Cal D}& := C\wedge \eta_S] .\cr} \eqnu $$

Let $\phi$ be any real-valued function with compact support in $M\setminus\Sigma_q$
and let $X:=\phi S$ [$X:=\phi T$]. Since $\imath_X\eta_T=\phi\imath_S\eta_T=\phi$ [$\imath_X
\eta_S=\phi\imath_T\eta_S=\phi$] and $\imath_X C=0$,
$$\eqalign{ \phi{\Cal D}^{\ast}\eta_S& =\imath_X{\Cal D}=-\imath_X(C\wedge \eta_T)=\phi C\, \cr
 [\phi{\Cal D}^{\ast}\eta_T& =-\imath_X{\Cal D}=-\imath_X (C\wedge \eta_S)=\phi C]\,, \cr} 
\speqnu{6.9'}$$
hence $C={\Cal D}^{\ast}\eta_S$ [$C={\Cal D}^{\ast}\eta_T$] in ${\Cal D}'(M\setminus\Sigma)$. In 
addition, ${\Cal D}\in {\Cal I}_q(M\setminus\Sigma_q)$ [${\Cal D}\in {\Cal I}_{-q}(M\setminus
\Sigma_q)$], since
$$\eqalign{ S{\Cal D} & =d({\Cal D}^{\ast}\eta_S)=dC=0 \cr
           [T{\Cal D}&=-d({\Cal D}^{\ast}\eta_T)=-dC=0] \cr} \eqnu $$
in ${\Cal D}'(M\setminus\Sigma_q)$. If $C\in{\Cal B}_q(M)$ [$C\in{\Cal B}_{-q}(M)$], since $S,\,T
\in {\Cal V}_q(M)$, the identities $(6.9')$, with $\phi\equiv 1$, and $(6.10)$ hold in ${\Cal S}'_q
(M)$. It follows that ${\Cal D}\in {\Cal I}_q(M)$ [${\Cal D}\in {\Cal I}_{-q}(M)$]. 
\enddemo 

 The space of all basic currents is filtered by a scale of subspaces of {\it finite 
order currents}, defined in terms of the weighted Sobolev spaces introduced in [18, \S 2].

\demo{{\rm 6.2.} Weighted Sobolev spaces of currents}
The space $L^2_q(M)$ is simply the space of square-summable complex-valued functions
with respect to the area-form $\omega_q$ considered in \S 2, i.e. $L^2_q(M):=L^2(M,\omega_q)$. 
The space $H^s_q(M)$, $s\in {\Bbb N}\setminus\{0\}$, was defined in [18, \S 2] 
as the completion of the space of smooth complex-valued functions $v$ on $M$ 
such that 
$$|v|_s:= \Bigl(\sum_{i+j\leq s} |S^iT^jv|^2_0\Bigr)^{1/2} \,<\,+\infty
\eqnu $$
with respect to the norm $|\cdot|_s$. Let $H^s(M)$ denote the standard Sobolev 
spaces on the compact manifold $M$ [70, Chap.\ IV, \S 1]. The following 
properties hold:
\medbreak
\item{(1)} $L^2(M)\subset L^2_q(M)$,
\smallbreak\item{(2)} $H^1(M)\equiv H^1_q(M)$,
\smallbreak\item{(3)} $H^s_q(M)\subset H^s(M)$, $s\geq 2$.
\medbreak

 The embedding $(1)$ is an immediate consequence of the definitions, since $\omega_q$ is a 
smooth $2$-form on $M$, vanishing at $\Sigma$. The isomorphism in $(2)$ can be proved by [18, 
Lemmas 2.1 \& 2.2]. The proof is based on Poincar\'e's inequality and on the following formulas (see [18, (2.7)]):
$$\eqalign { & S=|z|^{-k}\left\{\Re(z^{k/2}){\partial\over{\partial x}}\,-\, \Im(z^{k/2})
{\partial\over{\partial y}}\right\}\,, \cr &T=|z|^{-k}\left\{\Im(z^{k/2}){\partial\over{\partial x}}
\,+\,\Re(z^{k/2}) {\partial\over{\partial y}}\right\}\,,} \speqnu{6.11'}$$ 
where $z=x+iy$ is a canonical complex coordinate for the (orientable) quadratic
differential 
$q$ at a zero $p\in\Sigma_q$ of {\it even }order $k$. By the same formulas and 
$(2)$ one can also 
derive $(3)$. We also remark that the space of smooth functions with compact 
support in $M\setminus
\Sigma_q$ is {\it not }dense in $H^s_q(M)$ for $s\geq 2$, as a consequence of 
$(3)$ and of trace theorems for Sobolev spaces\break [1, Th.\ 5.4]. However, it can be
proved that $\Omega_q^0\subset H^s_q(M)$  is a dense subspace, for all 
$s\geq 0$. In fact, $C^{\infty}(M)\cap H^s_q(M)$ is dense in $H^s_q(M)$
by definition. The density of $\Omega_q^0$ in $C^{\infty}(M)\cap H^s_q(M)$
can be proved by Taylor expansion at each $p\in \Sigma_q$.

The {\it weighted Sobolev spaces }of $1$-forms, defined for all $s\geq 0$ by
$${\Cal H}^s_q(M):=\{\alpha\in {\Cal S}'_q(M)\,|\,(\imath_S\alpha,\imath_T\alpha) \in H^s_q(M)
\times H^s_q(M)\}\,\,, \eqnu $$ 
have a natural Banach space structure, isomorphic to the product $H^s_q(M)\times H^s_q(M)$. Since 
$\Omega_q^0\subset H^s_q(M)$ is dense, $\Omega_q^1\subset {\Cal H}^s_q(M)$ is also dense, for all
$s\geq 0$. The dual Sobolev spaces $H^{-s}_q(M)$, ${\Cal H}^{-s}_q(M)$ will be the distributional 
spaces defined, according to the standard Banach space theory, as the spaces of bounded functionals 
on $H^s_q(M)$, ${\Cal H}^s_q(M)$ respectively. Let ${\Cal I}_q^s(M)\subset {\Cal I}_q(M)$ [${\Cal 
I}_{-q}^s(M)\subset {\Cal I}_{- q}(M)$] be the space of $S$-invariant [$T$-invariant] distributions
of {\it order} $s\in {\Bbb N}$, that is the space of ${\Cal D}\in H^{-s}_q(M)$ such that $S{\Cal D}=0$ 
[$T{\Cal D}=0$] in $H^{-s-1}_q(M)$. Let ${\Cal B}^s_{\pm q}(M) \subset {\Cal B}_{\pm q}(M)\cap {\Cal H}^{-s}_q(M)$ be the space of ${\Cal F}_{\pm q}$-basic 
currents of {\it order} $s\in {\Bbb N}$, defined as follows:
$$\eqalign{ &{\Cal B}^s_q(M):=\{ C\, | \,\imath_S C=0 
\in H^{-s}_q(M),\,{\Cal L}_S C=0\in {\Cal H}^{-s-1}_q(M)\}, \cr
         &{\Cal B}^s_{-q}(M):=\{ C\,|\, \imath_T C=0 \in 
H^{-s}_q(M),\,
{\Cal L}_T C =0\in {\Cal H}^{-s-1}_q(M)\}. \cr} \speqnu{6.12'}$$  

 By Lemma 6.5, we have the following: 
\proclaim {Lemma} A current $C\in {\Cal B}^s_q(M)$ {\rm [}$C\in {\Cal B}^s_{-q}(M)$] if and only 
if the distribution $C\wedge\eta_T\in {\Cal I}^s_q(M)$ {\rm [}$C\wedge\eta_S\in {\Cal I}_{-q}^s(M)${\rm ].}
In addition{\rm ,} the map
$$ \eqalign{ C & \to -C\wedge \eta_T\,\,,\cr
             [C & \to C\wedge \eta_S] \,\,, \cr } \eqnu $$
is a bijection from ${\Cal B}^s_q(M)$ {\rm [}${\Cal B}^s_{-q}(M)${\rm ]} onto the vector space 
${\Cal I}^s_q(M)$ [${\Cal I}^s_{-q}(M)$] of $S$\/{\rm -}\/invariant {\rm [}$T$\/{\rm -}\/invariant{\rm ]} distributions of order
$s\in {\Bbb N}${\rm .}
\endproclaim

\numbereddemo {Definition} The $H^{-s}$ (first) basic cohomology (with real coefficients) of the 
measured foliations ${\Cal F}_{\pm q}$ is given by  
$$H^{1,s}_{\pm q}(M,{\Bbb R}):=j_q\bigl({\Cal B}^s_{\pm q}(M)\bigr)\subset H^1_{\pm q}
(M,{\Bbb R}) \subset H^1(M,{\Bbb R})\,\,.\eqnu $$
\enddemo

The $H^{-s}$ basic cohomology of the measured foliations ${\Cal F}_{\pm q}$ depends locally
only on the A. Katok's {\it fundamental class }of the foliation. The proof is in fact based on Katok's
{\it local classification theorem }for orientable measured foliations on surfaces, announced in [29, Th.\ 3], 
proved in [32, Th.\ 14.7.4] or [47, Th.\ 7.11.7].\pagebreak

\proclaim{Lemma} Let $q_t\in Q_{\kappa}${\rm ,} $t\in [0,1]${\rm ,} be a smooth family of quadratic 
differentials{\rm .} Let $\Sigma_t$ be the set of zeroes of $q_t${\rm .} With the natural identification of
the relative cohomology vector spaces $H^1(M,\Sigma_t,{\Bbb R})\equiv H^1(M,\Sigma_0,{\Bbb R})${\rm ,} 
for all $t\in [0,1]${\rm ,} we have\/{\rm :}
$$\eqalign {[\Im(q_t^{1/2})]\equiv\lambda^+_t[\Im(q_0^{1/2})]\in& H^1(M,\Sigma_0,{\Bbb R}), \hbox{ and } \cr
&  \hbox{all }\lambda^+_t>0\,\Rightarrow\,  H^{1,s}_{q_t}(M,{\Bbb R})\equiv H^{1,s}_{q_0}(M,{\Bbb R})\,, 
\cr
[\Re(q_t^{1/2})]\equiv\lambda^-_t[\Re(q_0^{1/2})]\in& H^1(M,\Sigma_0,{\Bbb R}), \hbox{ and } \cr
&  \hbox{all }\lambda^-_t>0\,\Rightarrow\,  H^{1,s}_{-q_t}(M,{\Bbb R})\equiv H^{1,s}_{-q_0}(M,{\Bbb R})\,. 
\cr} \eqnu $$            
\endproclaim

\demo{Proof} Since quadratic differentials in $Q_{\kappa}$ are determined up to isotopies and the 
$H^{-s}$ basic cohomologies $H^{1,s}_{\pm q}(M,{\Bbb R})$  are invariant under isotopies, we can 
assume that $\Sigma_t\equiv\Sigma_0$ and $q_t\equiv q_0$ on a neighbourhood of $\Sigma_0$. By a 
result of A. Katok [32, Th.\ 14.7.4] or [47, Th.\ 7.11.7], if $[\Im(q_t^{1/2})]=\lambda^+_t 
[\Im(q_0^{1/2})]\in H^1(M,\Sigma_0;{\Bbb R})$ [if $[\Re(q_t^{1/2})]=\lambda^-_t[\Re(q_0^{1/2})]\in 
H^1(M,\Sigma_0;{\Bbb R})]$, with $\lambda^+_t>0$\break  [$\lambda^-_t>0$], for all $t\in [0,1]$, then the
horizontal foliation ${\Cal F}_{q_t}$ is isotopic to ${\Cal F}_{q_0}$ [the vertical foliation
${\Cal F}_{-q_t}$ is isotopic to ${\Cal F}_{-q_0}$]. In both cases the Katok's classification 
result holds since $q_t\equiv q_0$ in a neighborhood of the set $\Sigma_t\equiv\Sigma_0$ of common
zeroes (hence the isotopy is equal to the identity on a neighbourhood of $\Sigma_0$). By the isotopy 
invariance of the $H^{-s}$ basic cohomology, the argument can be reduced to the case ${\Cal F}_{q_t}
={\Cal F}_{q_0}$ [${\Cal F}_{-q_t}={\Cal F}_{-q_0}$]. Then, since the property that $q_t\equiv q_0$
in a neighbourhood of $\Sigma_t\equiv \Sigma_0$ implies that the dual Sobolev spaces of currents
${\Cal H}^{-s}_{q_t}(M)\equiv{\Cal H}^{-s}_{q_0}(M)$, by the definition $(6.12')$ of basic currents 
of order $s\in {\Bbb N}$, it follows that
$$\eqalign { {\Cal B}^s_{q_t}(M)&={\Cal B}^s_{q_0}(M)  \cr
             [{\Cal B}^s_{-q_t}(M)&={\Cal B}^s_{-q_0}(M)] \,.\cr} \speqnu{6.15'}$$
In fact, let $\omega_t$, $\omega_0$ be the (degenerate) area-forms induced, respectively, by the 
quadratic differentials $q_t$, $q_0$. There exists a function $W_t:M\to {\Bbb R}^+$ such that $W_t
\equiv 1$ in a neighbourhood of $\Sigma_0$ and $\omega_t = W_t\omega_0$. Hence, ${\Cal F}_{q_t}=
{\Cal F}_{q_0}$ [${\Cal F}_{-q_t}={\Cal F}_{-q_0}$] implies $S_t=W_t^{-1}S_0$ [$T_t=W_t^{-1} T_0$]. 
Then $(6.15')$ follows by applying the definition $(6.12')$. By $(6.15')$ and Definition 6.7
of the $H^{-s}$ basic cohomologies, the statement is proved.
\enddemo

6.3. {\it A geometric estimate of the Poincar{\rm \'{\it e}} constant}.
The property of weighted Sobolev spaces which is the most relevant to the content of
this paper is a geometric estimate of the {\it Poincar{\rm \'{\it e}} constant }(the constant in the Poincar\'e
inequality for the Dirichlet form of the metric induced by a quadratic differential). Let $q\in 
Q_{\kappa}$. The Dirichlet form of the metric $R_q$, introduced in [18, (2.6)], is defined as 
the hermitian form on the Hilbert space $H^1_q(M)$ given by
$${\Cal Q}(u,v):= (Su,Sv)_q\, +\, (Tu,Tv)_q\,\,. \eqnu $$
The following is a stronger version of [18, Lemma 2.2], which adds a crucial explicit geometric
estimate of the constant.  

\proclaim {Lemma} There is a constant $K_{g,\sigma}>0$ such that the following holds{\rm .} Let 
$q\in Q^{(1)}_g$ be the square of a holomorphic differential{\rm ,} $\Sigma_q$ be the {\rm (}\/finite\/{\rm )} set
of its zeroes and let $\sigma:=\hbox{card}(\Sigma_q)${\rm .} Denote by $|\!|q|\!|$ the $R_q$\/{\rm -}\/length 
of the shortest geodesic segment with endpoints in $\Sigma_q${\rm .} Then{\rm ,} for any $v\in H^1_q(M)${\rm ,}
$$\left|v-\int_M v\,\omega_q\right|_0 \leq {{K_{g,\sigma}} \over {|\!|q|\!|}}\, {\Cal Q}(v,v)^{1/2}\,\,. 
\speqnu{6.16'}$$
\endproclaim  

\demo{Proof} It was shown in [18, \S 2] that many basic properties of the Laplace-Beltrami
operator for a (non-degenerate) smooth Riemannian metric $R$ on a compact surface hold for the
Dirichlet form $\Cal Q$. In particular, since the Max-Min principle holds [18, Lemma 2.4], 
we have the following standard relation between the constant in the Poincar\'e lemma and the first
non-zero eigenvalue $\lambda_1$ of $\Cal Q$. Let $Z^1_q\subset H^1_q(M)$ be the subspace of zero 
average functions. Then 
$$\min_{v\in Z^1_q} {{{\Cal Q}(v,v)}\over {|v|_0^2}}=\lambda_1 \,\,, \eqnu $$
where the minimum is achieved at $v=u_1$, where $u_1$ is any eigenfunction with eigenvalue 
$\lambda_1$. A standard lower estimate for the first non-zero eigenvalue is given, for smooth 
Riemannian metrics on compact manifolds, in terms of the Cheeger isoperimetric constant \ref\Cg. 
Cheeger's argument \ref\Cg, [4, Chap.\ III, D.4] carries over without modifications to 
the case of the degenerate Riemannian metric $R_q$. In fact, by [18, Th.\ 2.3 (iii)] all 
eigenfunctions of the Dirichlet form $\Cal Q$ are smooth on $M$, including at $\Sigma_q$. In 
addition, the degenerate metric $R_q$ induces a well-defined smooth conformal structure on 
$M$, hence the key step of Cheeger's proof, the co-area formula [8, Th.\ 6.3], holds. Let 
$h_q$ be the Cheeger isoperimetric constant for the metric $R_q$,
$$ h_q := \inf  {{L_q(\Gamma)}\over{\min (\hbox{vol}_q(M_1),\hbox{vol}_q(M_2))}}\,\,, 
\eqnu $$
where $L_q$ denotes the $R_q$-length and the infimum is taken over all compact $1$-dimensional 
submanifolds $\Gamma\subset M$, dividing $M$ into two submanifolds with boundary $M_1$, $M_2$ such 
that $\partial M_1=\partial M_2=\Gamma$ and $M_1 \cup M_2=M$. Then the Cheeger inequality reads
$$\lambda_1\geq {1\over 4}\,h_q^2 \,\,. \speqnu{6.18'}$$
\smallskip
\noindent The final step consists in a lower estimate of $h_q$. A $1$-dimensional submanifold
$\Gamma\subset M$ is a finite union of simple closed curves. If $\Gamma$ contains a closed curve 
$\gamma$ non-homotopic to zero in $M$, and since $\hbox{vol}_q(M)=1$, 
$${{L_q(\Gamma)}\over{\min (\hbox{vol}_q(M_1),\hbox{vol}_q(M_2))}}\geq L_q(\gamma)\geq |\!|q|\!|\,\,.
\eqnu $$
In fact, the shortest simple closed curve $\widehat{\gamma}$ freely homotopic to $\gamma$, 
which exists by [58, Th.\ 18.4], either contains singularities or it is a waist curve of an embedded flat
cylinder. The boundary of the cylinder is then the union of a finite number of geodesic segments
with endpoints in $\Sigma_q$. In either case $L_q(\gamma)\geq L_q(\widehat{\gamma})\geq |\!|q|\!|$.
If all simple closed curves in $\Gamma$ are homotopic to zero, since by an elementary lemma due
to S. T. Yau [73, (1)-(4)] the infimum in the definition $(6.18)$ can be equivalently taken on 
$\Gamma$ such that $M_1$, $M_2$ are connected, by the Jordan-Brouwer separation theorem, we can 
assume that $\Gamma$ is reduced to a single simple closed smooth curve bounding a simply connected
domain $\Omega\subset M$. Let $p_0\in \Omega$ be any given point and consider the holomorphic map
$\Phi_q:\Omega\to {\Bbb C}$ given by
$$\Phi_q(p):=\int_{p_0}^p q^{1/2} \,\,, \eqnu $$
where $q^{1/2}$ is a holomorphic square root of $q$ on $\Omega$. The map $\Phi_q$ is well-defined 
and holomorphic since $\Omega$ is simply connected and $q^{1/2}$ is a holomorphic differential. By 
the open mapping theorem for holomorphic functions, $\Phi_q(\Omega)\subset~{\Bbb C}$ is an open set 
and $\Phi_q$ is a branched covering with at most $2g-1$ sheets by the Riemann-Hurwitz relation 
[13, I.2.7]. Since, by definition, the metric $R_q=\Phi_q^{\ast}(R_e)$, where $R_e$ is the 
Euclidean metric on ${\Bbb C}$,
$${{L_q(\Gamma)}\over{\min (\hbox{vol}_q(\Omega),\hbox{vol}_q(M\setminus\Omega))}} \geq  
{{L_e(\Phi_q(\Gamma))}\over{\hbox{vol}_e(\Phi_q(\Omega))}}\geq 2\pi^{1/2}(2g-1)^{1/2}\,\,, \eqnu $$
where the last inequality follows from the classical isoperimetric inequality for the Euclidean 
metric in ${\Bbb R}^2$, that is $L_e^2(\partial D)\geq 4\pi\hbox{vol}_e(D)$ if $D\subset {\Bbb R}^2$ 
is any domain [8, \S 6.2], and from the estimate $\hbox{vol}_e(\Phi_q(\Omega))=(2g-1)^{-1}
\hbox{vol}_q(\Omega)\leq (2g-1)^{-1}$. 

 Let $d_q(M):=\max\{d_q(p,\Sigma_q)\,|\,p\in M\}$ be the maximum $R_q$-distance of any 
point of $M$ to a zero of $q$. Since $|\!|q|\!|/2\leq d_q(M)$ and, by [45, Cor. 5.6], there
exists a constant $K'_{g,\sigma}>0$ (depending only on $(g,\sigma)$) such that either $d_q(M)\leq
\sqrt{2/\pi}$ or $d_q(M)\leq K'_{g,\sigma}/|\!|q|\!|$, the estimate $(6.16')$ follows from $(6.17)$, 
$(6.18')$, $(6.19)$ and $(6.21)$. \enddemo\pagebreak

 \section{The structure of the space of basic currents of finite order}

The results we present below on the structure of the space of basic currents and on the basic cohomology of measured
foliations are derived from \ref\Ftwo and therefore do not depend on the previous sections of the paper (with the
exception of \S 6.1 and \S 6.2). 

\demo{{\rm 7.1.} Basic currents with non-vanishing cohomology class}
 Let $q\in Q_{\kappa}$, let $\Sigma:=\Sigma_q$ be finite set of
its zeroes and ${\Cal F}:=
{\Cal F}_q$ be its horizontal foliation. A current $C\in {\Cal D}'(M\setminus\Sigma)$ will be said 
to be of {\it finite order }$s\in {\Bbb N}$ if it can be extended to a
continuous linear functional 
on the Sobolev space of forms ${\Cal H}^s_{loc}(M\setminus\Sigma_q)$. The space of basic  
currents of order $s\in {\Bbb N}$ of the measured foliation ${\Cal F}|_{M\setminus\Sigma}$ will be 
denoted by ${\Cal B}^s_{\Cal F}(M\setminus\Sigma)$. The space ${\Cal B}^s_{\Cal F}(M):={\Cal B}^s
_q(M)$ of basic currents of order $s\in {\Bbb N}$ for ${\Cal F}={\Cal F}_q$ has been defined in Section 6.2. 
Let $j_{\Sigma}:{\Cal Z}'(M\setminus\Sigma)\to H^1(M\setminus\Sigma,{\Bbb R})$ be the cohomology
map $(6.2)$ and let
$$\eqalign{H^{1,s}_{\Cal F}(M\setminus\Sigma,{\Bbb R})&:=j_{\Sigma}\bigl({\Cal B}^s_{\Cal F}(M
\setminus\Sigma) \bigr)\subset H^1(M\setminus\Sigma,{\Bbb R})\,\,,\cr
H^{1,s}_{\Cal F}(M,{\Bbb R})&:=j_{\Sigma}\bigl({\Cal B}^s_{\Cal F}(M)\bigr)\subset H^1(M,{\Bbb R}) 
 \,\,,\cr} \eqnu $$
be the $H^{-s}$ basic cohomologies of, respectively, ${\Cal F}|_{M\setminus\Sigma}$ and $\Cal F$.
We will prove the following:
\enddemo

\proclaim {Theorem} There exists an integer $l>1$ such that the following holds{\rm .} Let $q_0\in
Q_{\kappa}${\rm .} For almost all $q\in {\rm SO}(2,{\Bbb R})q_0\subset Q_{\kappa}$ with respect
to the {\rm (}\/Haar\/{\rm )} Lebesgue measure{\rm ,} the basic cohomologies of order $s>l$
 of the horizontal measured 
foliation ${\Cal F}:={\Cal F}_q$ satisfy the following identities\/{\rm :}
$$\eqalign {& {\rm (i)}\, H^{1,s}_{\Cal F}(M,{\Bbb R})\equiv H_{\Cal F}:=\{ c\in H^1(M,{\Bbb R})\,|
           c\wedge [\eta_{\Cal F}]=0 \} \,\,, \cr 
          & {\rm (ii)}\, H^{1,s}_{\Cal F}(M\setminus\Sigma,{\Bbb R})\equiv H^1(M\setminus\Sigma,
                   {\Bbb R})\,\,.\cr} \speqnu{7.1'}$$
The form $\eta_{\Cal F}:=\Im(q^{1/2})$ is the closed $1$\/{\rm -}\/form determined by the measured foliation
$\Cal F${\rm ,} hence the cohomology class $[\eta_{\Cal F}]$ is the projection under the map $H^1(M,\Sigma;
{\Bbb R}) \to H^1(M,{\Bbb R})$ of the Katok\/{\rm '}\/s fundamental class of $\Cal F${\rm .}       
\endproclaim
 
\specialnumber{7.1'}
\demo{Remark} Let $\mu$ be a probability measure on the moduli space ${\Cal M}^{(1)}_{\kappa}$
with the property that the conditional measures induced by $\mu$ on the orbits of the circle group
${\rm SO}(2,{\Bbb R})$ are absolutely continuous. By Theorem 7.1 and Fubini's theorem, the identities $(7.1')$ 
hold for $\mu$-almost all $q\in {\Cal M}^{(1)}_{\kappa}$. The statement holds in particular if $\mu=
\mu^{(1)}_{\kappa}$ is the (unique) absolutely continuous $G_t$-invariant probability measure on 
${\Cal M}^{(1)}_{\kappa}$.
\enddemo

In [18, \S 4] $S$-invariant distributions were constructed by the following method.
Let $m^+\in {\Cal M}^+_q$ be a meromorphic function in $L^2_q(M)$ with poles at $\Sigma:=\Sigma_q$,
the set of zeroes of a quadratic differential $q\in Q_{\kappa}$. If $U\in H^{-s}_q(M)$ is a solution 
of the cohomological equation ${\rm SU}=m^+$ in $H^{-s-1}_q(M)$, then ${\Cal D}=\partial_q^+U\in H^{-s-1}_q
(M)$ is an $S$-invariant distribution on $M$. A refined version of this construction will enables us 
to detect the cohomology classes of basic currents corresponding to the $S$-invariant distributions 
constructed.

The spaces ${\Cal M}^{\pm}_{\Sigma}$ of all meromorphic, respectively anti-meromorphic, 
functions with poles at $\Sigma:=\Sigma_q$ is a subspace of ${\Cal S}_q'(M)$. In fact, a function 
$m^{\pm}\in {\Cal M}^{\pm}_{\Sigma}$ can be regarded as a distribution in $H^{-s}_q(M)$ if and only 
if at any point $p\in\Sigma$ the order of  a pole of $m^{\pm}$ does not exceed $k/2+ (k/2+1)s$, where 
$k$ is the (even) order of zero of the quadratic differential $q$ at $p$ [18, \S 4]. The 
distributions defined by $m^{\pm}\in {\Cal M}_{\Sigma}$ are given, for any function $v\in H^s_q(M)$, 
by the formula
$$(m^{\pm},v):=\lim_{\epsilon\to 0} \int_{M\setminus \Sigma_{\epsilon}}
m^{\pm}\,v\,\omega_q\,\,,\eqnu $$
where $\Sigma_{\epsilon}$ denotes the $\epsilon$-neighbourhood of the finite set $\Sigma$ with respect
to the geodesic distance. The distribution defined by $(7.2)$ is called the {\it principal value} of
$m^{\pm}$ (principal values of meromorphic and anti-meromorphic functions on the complex plane and 
the action of the Cauchy-Riemann operators on distributions are studied in [3, 3.6 \& 3.8]).
There exist isomorphisms $h^{\pm}_q:{\Cal M}^{\pm}_{\Sigma}\to {\Cal H}^{\pm}_{\Sigma}(M)$ onto the 
spaces of meromorphic, respectively anti-meromorphic, differentials with poles at $\Sigma$:
$$\eqalign { h^+_q(m^+)&:= m^+\,q^{1/2}\,\, , \cr
            h^-_q(m^-)&:=  m^-\, {\overline q}^{1/2}\,\,. \cr } \speqnu{7.2'}$$
In particular, $h_q^{\pm}$ maps the finite dimensional spaces ${\Cal M}^{\pm}_q$ of $L^2_q(M)$ 
meromorphic, respectively anti-meromorphic, functions onto the spaces ${\Cal H}^{\pm}_q(M)$ of 
holomorphic, respectively anti-holomorphic, differentials on the Riemann surface $M_q$.  

 Let $m^+\in {\Cal M}^+_{\Sigma}$ be a meromorphic function with poles at $\Sigma$. Let 
$U\in {\Cal D}'(M\setminus\Sigma)$ be a distributional solution of the equation ${\rm SU}=-\Re(m^+)$ in 
${\Cal D}'(M\setminus\Sigma)$. Let $U^{\ast}$ be the (unique) current of dimension $2$ and degree 
$0$ such that $U=U^{\ast}\omega_q$. The $1$-current $C\in {\Cal D}'(M\setminus\Sigma)$ uniquely 
determined by the identity 
$$ dU^{\ast}= -\Re h^+_q(m^+) \,+\,C \eqnu $$
is {\it basic }for ${\Cal F}_q|_{M\setminus\Sigma}$. In fact, by $(7.3)$, $C$ is closed (since
$h_q^+(m^+)$ is closed and $dU^{\ast}$ is exact) in ${\Cal D}'(M\setminus\Sigma)$ and $\imath_S C=0$,
since ${\rm SU}=-\Re (m^+)$. Then ${\Cal L}_SC=0$ by $(6.4)$. A short computation also shows that 
$C={\Cal D}^{\ast}\eta_S$, where $\Cal D$ is the quasi $S$-invariant distribution given by ${\Cal D}
= -C\wedge \eta_T=TU - \Im(m^+)$. If $h^+_q(m^+)$ is closed in ${\Cal S}'_q(M)$, which is the cases 
if $m^+\in L^2_q(M)$, and the equation ${\rm SU}=-\Re(m^+)$ holds in ${\Cal S}'_q(M)$, then $C\in {\Cal S}'
_q(M)$ and it is basic for ${\Cal F}_q$. By $(7.3)$, the cohomology class $[C]=[\Re h^+_q(m^+)]\in 
H^1(M\setminus\Sigma,{\Bbb R})$. 

 It is a standard fact in the theory of Riemann surfaces, related to Hodge theory, that 
each cohomology class in $H^1(M_q,{\Bbb R})$ can be represented by a harmonic differential $\Re(h^+)$, 
where $h^+$ is a holomorphic differential on the Riemann surface $M_q$ [13, III.2]. Moreover, 
it follows from the Riemann-Roch theorem [13, III.4] that all cohomology classes in $H^1(M_q
\setminus\Sigma;{\Bbb R})$ can be represented in the form $\Re(h^+)$ where $h^+$ is a {\it 
meromorphic }differential with at most {\it simple poles }at $\Sigma$. The strategy of our 
argument will therefore be based on the construction of solutions to the cohomological equation 
${\rm SU}=-\Re(m^+)$ in ${\Cal D}'(M\setminus\Sigma)$ (in ${\Cal S}_q'(M)$) for appropriate choices of 
a meromorphic function $m^+\in {\Cal M}^+_{\Sigma}$.
  The construction of basic currents depends therefore on the following:

\proclaimtitle{[18, Th. 4.1]}
\proclaim {Theorem}  There exists an integer $r>1$ such that the following 
holds{\rm .} Let $q_0\in Q_{\kappa}${\rm .} For almost all $q\in {\rm SO}(2,{\Bbb R})q_0\subset Q_{\kappa}$ 
with respect
to the {\rm (}\/Haar\/{\rm )} Lebesgue measure{\rm ,} the cohomological equation $Su=f$ has a distributional solution
$u\in  H^{-r}_q(M)$ for any given $f\in H^{r-1}_q(M)$ with $\int_M f\,\omega_q=0${\rm .}
\endproclaim

We would like to construct solutions of the equation ${\rm SU}=-\Re(m^+)$, as $h^+_q(m^+)$ varies 
over all possible meromorphic differentials on $M$ with at most simple poles at $\Sigma$. Theorem 7.2
cannot be directly applied since $\Re(m^+)\not \in H^r_q(M)$. In the case of a smooth differential, 
which corresponds to the case that $m^+\in L^2_q(M)$, the idea is to first solve the equation locally
around $\Sigma$.

 Let $p\in \Sigma$ be a zero of $q$ of {\it even }order $k$. With respect to a canonical 
holomorphic coordinate $z=x+iy$ at $p$, $q(z)=z^kdz^2$. By $(6.11')$ we have: 
$$\eqalign{ Sz=z^{-k/2}\,,&\,\,S{\overline z}={\overline z}^{-k/2}\,, \cr
   Tz={\imath} z^{-k/2}\,,&\,\,T{\overline z}=-{\imath}{\overline z}^{-k/2}\,. \cr} \eqnu $$
By the identities $(7.4)$, the equation ${\rm SU}=-\Re(m^+)$ can be solved locally by expanding the meromorphic
function $m^+\in L^2_q(M)$ in Laurent series around each point $p\in \Sigma$. In fact, if a $m^{\pm} 
\in {\Cal M}^{\pm}_q\subset L^2_q(M)$, then $m^{\pm}$ has a pole of order at most $k/2$ at a zero 
$p\in \Sigma$ of order $k$. Hence, by $(7.4)$, the local solution $u^{\pm}$ of the equation $Su^{\pm}
=m^{\pm}$, constructed formally by Laurent expansion, is respectively holomorphic, anti-holomorphic.
Let $u_0$ be the function obtained as an extension to $M$ of the local solution of the equation 
${\rm SU}=-\Re(m^+)$ by partition of unity. Then $u_0$ is smooth on $M$, hence $u_0\in H^1_q(M)$, although 
$u_0\not\in H^2_q(M)$, and $Su_0=-\Re(m^+)$ in a neighbourhood of $\Sigma$. It follows that the
function $f_1:=-\Re(m^+)-Su_0\in H^s_q(M)$, for all $s\in {\Bbb N}$. In fact, $f_1$ is smooth and
has compact support in $M\setminus\Sigma$. In order to apply Theorem 7.2 to the equation $Su_1=f_1$
we have to impose the zero average condition:
$$\int_M \Re(m^+ + Su_0)\,\omega_q=\int_M \Re h^+_q(m^+) \wedge \eta_S=0\,\,. \eqnu $$
Under condition $(7.5)$ Theorem 7.2 yields a solution $u_1\in H^{-r}_q(M)$. Hence the 
distribution $U=u_0+u_1 \in H^{-r}_q(M)$ is a solution of the equation ${\rm SU}=-\Re(m^+)$ in $H^{-r-1}_q
(M)$ and the identity $(7.3)$ determines a real ${\Cal F}_q$-basic current $C\in{\Cal H}^{-r-1}_q(M)$,  
cohomologous to the harmonic form $\Re h^+_q(m^+)$. By this method, we can construct basic currents
$C\in {\Cal B}^{r+1}_{\Cal F}(M)$ representing any given cohomology class in the subspace $\{c\in 
H^1(M,{\Bbb R})\,|\, c\wedge [\eta_{\Cal F}]=0\}$. In fact, such cohomology classes can be
represented by smooth harmonic differentials on $M$, $\Re h^+_q(m^+)$, $m^+\in L^2_q(M)$, satisfying
the zero average condition $(7.5)$.

 Conversely, it can be proved that if $C\in {\Cal B}_{\Cal F}(M)$, then the zero average 
condition $(7.5)$ holds. In fact, 
$$[C]\wedge[\eta_{\Cal F}]=C\wedge \eta_S=\imath_S C\wedge \omega_q=0\,\,. 
\speqnu{7.5'}$$ 
Theorem 7.1(i) is therefore proved.

In order to prove Theorem 7.1(ii) we need a regularization result proved in [18, Lemma 
4.3B]. The main idea is to express all distributions in the dual weighted Sobolev spaces, which
vanish on constant functions, as sums of iterated Cauchy-Riemann derivatives of smooth functions. In 
[18, Prop. 4.6A] we described the kernels of the Cauchy-Riemann operators $\partial^{\pm}_{-s}:
H^{-s+1}_q(M)\to H^{-s}_q(M)$ modulo distributions supported at $\Sigma$. From this perspective such 
kernels are given by {\it all }equivalence classes of, respectively, meromorphic and anti-meromorphic,
functions in $H^{-s+1}_q(M)$ modulo distributions supported at $\Sigma$. Such a conclusion follows 
immediately from the remark that $\partial^{\pm}$ are Cauchy-Riemann operators in the standard sense
on the open manifold $M\setminus\Sigma$; hence their distributional kernel is given by, respectively,
{\it all }holomorphic and anti-holomorphic functions. We prove below a refined version of such a 
result:

\proclaim {Lemma} Let $\partial^{\pm}_{-s}:H^{-s+1}_q(M)\to H^{-s}_q(M)${\rm ,} $s\geq 1${\rm ,} denote the 
Cauchy\/{\rm -}\/Riemann operators on dual weighted Sobolev spaces{\rm .}
 The range $R(\partial^{\pm}_{-s})\subset 
H^{-s}_q(M)$ of the Cauchy\/{\rm -}\/Riemann operators is given by all distributions vanishing on constant 
functions{\rm .} The kernel $N(\partial^{\pm}_{-s})\subset H^{-s+1}_q(M)$ consists respectively of 
meromorphic{\rm ,} anti\/{\rm -}\/meromorphic{\rm ,} functions in $H^{-s+1}_q(M)$ which have no terms of the form 
$z^{-\ell(k/2+1)}${\rm ,} ${\overline z}^{-\ell(k/2+1)}${\rm ,}  $\ell\in {\Bbb N}\setminus \{0\}${\rm ,} in the 
Laurent expansion with respect to the canonical holomorphic coordinate, at a zero of $q$ of {\rm (}\/even\/{\rm )} 
order $k\in {\Bbb N}\setminus\{0\}${\rm .}
\endproclaim

\demo {Proof} Let $\Phi\in H^{-s}_q(M)$. Then
$$(U,\partial^{\pm}v):= -(\Phi,v)\,\,,\,\,\,\,\hbox{ for all }
\,\,v\in H^s_q(M)\,\,, \eqnu $$
defines a linear bounded functional on the (closed) range of the linear operator $\partial^{\pm}_s:
H^s_q(M)\to H^{s-1}_q(M)$. In fact, the kernel of $\partial^{\pm}$ on $H^s_q(M)$ contains only 
constant functions [18, Prop. 4.3A] and by hypothesis $\Phi$ vanishes on constant functions. 
In addition, by Poincar\'e's inequality (Lemma 6.9),
$$|(U,\partial^{\pm}v)|=|(\Phi,v)|\leq C |v|_s \leq C_s |\partial^{\pm}v|_{s-1}\,\,.
\eqnu $$
\vglue6pt\noindent 
Finally, $U$ extends by the Hahn-Banach theorem to a linear bounded functional on $H^{s-1}_q(M)$,
i.e.\ to a distribution $U\in H^{-s+1}_q(M)$, which by construction satisfies the equation $\partial
^{\pm}_{-s}U=\Phi$ in $H^{-s}_q(M)$. Hence the stated description of the range $R(\partial^{\pm}_{-s})$ 
is proved (following [18, Prop. 4.6A]).

 Let $U\in H^{-s+1}_q(M)$ be a solution of $\partial^{\pm}U=0$. By the local properties of 
the Cauchy-Riemann operators on $M\setminus\Sigma$, the restriction of $U$ to $M\setminus\Sigma$ is 
a holomorphic, respectively an anti-holomorphic, function. Since $U\in H^{-s+1}_q(M)$, it is a 
meromorphic, respectively an anti-meromorphic, function with poles of order at most $k/2+(k/2+1)
(s-1)$ at any point $p\in \Sigma$ where $q$ has a zero of (even) order $k$. A local computation 
of $\partial^{\pm}U$ at points of $\Sigma$ concludes the proof. In fact, let $p\in \Sigma$ be a 
zero of order $k$ and let $v_{\alpha,\beta}$ be the Taylor coefficients of a smooth function $v$ 
supported in a neighbourhood of $p$, with respect to the complex coordinates $z$, ${\overline z}$, 
where $z$ is the canonical coordinate for $q$ at $p$ (determined by the conditions $q(z)=z^kdz^2$ and 
$z(p)=0$). Since $\partial^{\pm}=S\pm\,iT$ and the vector fields $S$, $T$ are given by $(7.4)$, with 
respect to the canonical coordinate, by applying Stokes theorem we get the following formulas:
\vglue-6pt

$$\eqalign{&\bigl(\partial^+ (z^{-n}),v\bigr)=-\lim_{\epsilon\to 0}
\int_{M\setminus\Sigma_{\epsilon}} z^{-n}\,\partial^+v\,\omega_q= 2\pi v_{n-k/2	-1,0}\,\,, \cr
\noalign{\vskip4pt}
&\bigl(\partial^- ({\overline z}^{-n}),v\bigr)=-\lim_{\epsilon\to 0}
\int_{M\setminus\Sigma_{\epsilon}} {\overline z}^{-n}\,\partial^-v\,\omega_q= 2\pi v_{0,n-k/2-1}\,\,.
\cr} \eqnu $$
\vglue6pt\noindent 
The reader can compare $(7.8)$ with similar formulas for the non-singular situation [3, \S\S 3.6 \& 
3.8]. Since $U\in H^{-s+1}_q(M)$ and it is meromorphic, respectively anti-meromorphic, it can 
have at $p$ a pole of order at most $k/2+(k/2+1)(s-1)$. On the other hand, a smooth function $v\in 
H^s_q(M)$ is not allowed to have an arbitrary Taylor series at $p$. In fact, $v\in H^s_q(M)$ if and 
only if:

\vglue-6pt

$$v_{\alpha,\beta}\not=0 \Longrightarrow \alpha+\beta\geq -k/2+s(k/2+1) \hbox{ or }  
\alpha,\beta\in {\Bbb Z}\cdot (k/2+1)\,.  \eqnu $$ 
\vglue6pt\noindent 
By $(7.8)$ and $(7.9)$, if $n-k/2-1=(\ell-1)(k/2+1)$, that is $n=\ell (k/2+1)$, $\ell\in 
{\Bbb N}\setminus\{0\}$,
\vglue-6pt

$$\eqalign{&\partial^+ (z^{-n})=2\pi {{(-1)^{\ell-1}}\over {(l-1)!}}(\partial^-)^{\ell-1}\delta_p
\not=0\,\,,\cr
\noalign{\vskip4pt}
 &\partial^-({\overline z}^{-n})=2\pi {{(-1)^{\ell-1}}\over {(l-1)!}}(\partial^+)
^{\ell-1} \delta_p\not=0\,\,,\cr} \eqnu $$  
where $\delta_p$ is the Dirac delta at $p\in \Sigma$. If $n\leq k/2+(k/2+1)(s-1)$, $z^{-n}$ and
${\overline z}^{-n}\in H^{-s+1}_q(M)$ locally and, since $\delta_p\in H^{-2}(M)\subset H^{-2}_q(M)$,
the identities $(7.10)$ hold in $H^{-s}_q(M)$. Such formulas can be compared with the analogous 
formulas for the derivatives of meromorphic and anti-meromorphic functions on the complex plane [3, 
Prop. 3.6.3].

 We have therefore proved that, if $n\leq k/2+(k/2+1)(s-1)$ and $n=\ell(k/2+1)$ for any 
$\ell \geq 1$, then meromorphic and anti-meromorphic functions, with a pole of order $n$ at a zero 
of (even) order $k$ of the orientable holomorphic quadratic differential $q$, do not belong to the 
kernel $N(\partial^{\pm}_{-s})\subset H^{-s+1}_q(M)$. On the other hand, if $n\not=\ell (k/2+1)$, 
then $\partial^+(z^{-n})=\partial^-({\overline z}^{-n})=0$ in $H^{-s}_q(M)$, by $(7.8)$ and 
$(7.9)$, since $n-k/2-1<-k/2+s(k/2+1)$.
\enddemo 

 By Lemma 7.3, if ${\Cal D}\in H^{-s}_q(M)$ vanishes on constant functions, then there exist
functions $u^{\pm}\in L^2_q(M)$ (which can be taken of zero average) such that $(\partial^{\pm})^s 
u^{\pm}={\Cal D}$. It is not possible in general to have a similar statement with $u^{\pm}$ smoother 
than $L^2_q$, since the operators $\partial^{\pm}:H^1_q(M)\to L^2_q(M)$ have a non-trivial cokernel. 
However, the following result holds:
\proclaim {Lemma} Let $s\geq 0$ and $r\geq 1$. Then{\rm ,} any function $v\in H^s_q(M)$ such that 
$\int_M v\,\omega_q=0$ can be (non-uniquely) represented as follows:
$$v=\sum_{i+j=r}({\partial}^+)^i({\partial}^-)^j\,v_{ij}\,\,,\,\,\,\,
\hbox{ where } \,\,v_{ij}\in H^{s+r}_q(M)\,. \eqnu $$
\endproclaim

The proof of this lemma is given in [18, Prop.\ 4.3A \& Lemma 4.3B] and will be 
omitted.

 We can now prove a stronger version of Theorem 7.2:

\specialnumber{7.5}
\proclaim {Corollary} There exists an integer $r>1$ {\rm (}\/given by Theorem {\rm 7.2)} such that the following
holds{\rm .} Let $q_0\in Q_{\kappa}${\rm .} For almost all $q\in {\rm SO}(2,{\Bbb R})q_0\subset Q_{\kappa}$\break
 with respect
to the {\rm (}\/Haar\/{\rm )} Lebesgue measure{\rm ,} the cohomological equation\break
 ${\rm SU}=\Phi$ has a
distributional solution 
$U\in H^{-s-2r+1}_q(M)$ for any given $\Phi\in H^{-s}_q(M)$ vanishing on constant functions{\rm .}
\endproclaim

\demo {Proof} Since $\Phi$ vanishes on constant functions, by iterated application of Lemma 7.3, 
there exists a zero average function $f\in L^2_q(M)$ such that
$$(\partial^+)^s f=\Phi\,\,. \eqnu $$
Then, by applying Lemma 7.4, with $s=0$, to the function $f$, there exist zero average functions 
$f_{ij}\in H^{r-1}_q(M)$ such that 
$$\Phi=(\partial^+)^s f=\sum_{i+j=r-1}({\partial}^+)^{i+s}({\partial}^-)^j\,f_{ij}
\,\,. \eqnu $$
By Theorem 7.2, it is   possible to find distributional solutions $U_{ij}\in H^{-r}_q(M)$ of the 
equations ${\rm SU}_{ij}=f_{ij}$. Then
$$U:=\sum_{i+j=r-1}({\partial}^+)^{i+s}({\partial}^-)^jU_{ij}\,
\in\, H^{-2r-s+1}_q(M) \eqnu $$ 
is a distributional solution of the equation ${\rm SU}=\Phi$ in $H^{-2r-s}_q(M)$.
\enddemo 
\smallskip
\noindent All cohomology classes $c\in H^1(M\setminus\Sigma,{\Bbb R})$ can be represented by
the real part of a meromorphic differential $h^+$ with at most simple poles at $\Sigma$. In
fact, by the Riemann-Roch theorem [13, III.4.8] the complex dimension of the space of 
meromorphic differentials with at most simple poles at $\Sigma$ is equal to $2g+\sigma-1$,
where $\sigma:=\hbox{card}(\Sigma)$, hence it is equal to the dimension of the cohomology
$H^1(M\setminus\Sigma,{\Bbb C})$. If $h^+$ has at most simple poles at $\Sigma$, then $h^+=
h^+_q(m^+)$, where $m^+\in {\Cal M}_{\Sigma}$ has a pole of order at most $k/2+1$ at any zero
of $q$ of order $k$. We would like to construct a solution in ${\Cal D}'(M\setminus\Sigma)$
of the equation ${\rm SU}=-\Re (m^+)$. If $\Re (m^+)$ does not vanish on constant functions, then
no solution of finite order $U\in H^{-s}_q(M)$ exists. In fact, if $U\in H^{-s}_q(M)$, then
${\rm SU}\in H^{-s-1}_q(M)$ vanishes on constant functions. In order to bypass this difficulty, we 
construct a distributional solution of the modified equation ${\rm SU}=-\Re(m^+)+\delta$, where 
$\delta$ is a (Dirac) distribution supported on $\Sigma$ with the property that the distribution
$\Phi:=-\Re(m^+)+\delta\in H^{-2}_q(M)$ vanishes on constant functions. The existence of a
solution $U\in H^{-2r-1}_q(M)$ is given by Corollary 7.5. Since $U$ is a solution of the
modified equation ${\rm SU}=-\Re(m^+)+\delta$ in $H^{-2r-2}_q(M)$ and $\delta$ is supported on 
$\Sigma$, the distribution $U$   solves the equation ${\rm SU}=-\Re(m^+)$ in ${\Cal D}'(M\setminus
\Sigma)$. Let $C\in {\Cal H}^{-2r-2}_q(M)$ be the currents obtained from $U$ by the identity 
$(7.3)$. By construction of $U$ and $(7.3)$, $C\in {\Cal B}^s_{\Cal F}(M\setminus\Sigma)$, $s\geq 2r+2$, and 
the cohomology class $[C]=[\Re(h^+)]\in H^1(M\setminus\Sigma,{\Bbb R})$. Theorem 7.1(ii) is 
therefore proved and the proof of Theorem 7.1 is completed. 
\enddemo

\demo{{\rm 7.2.} Basic currents with vanishing cohomology class}
Let $\Cal D$ be a quasi\break $S$-invariant distribution and let ${\Cal D}^{\ast}$ be the 
current of dimension $2$ (and degree~$0$) corresponding to $\Cal D$ under the isomorphism induced 
by $\omega_q$. Then $C:=d{\Cal D}^{\ast}\in {\Cal B}_{\Cal F}(M\setminus\Sigma)$. In fact, 
$$\eqalign{ \imath_S\, d\,{\Cal D}^{\ast}&={\Cal L}_S{\Cal D}^\ast=0\,\,,\cr
  {\Cal L}_Sd\,{\Cal D}^{\ast}&=d\,\imath_S\, d\,{\Cal D}^{\ast}=0 \cr}\eqnu $$
in ${\Cal D}'(M\setminus\Sigma)$. Moreover, by its definition as a differential, $C$ has vanishing 
cohomology class. If $\Cal D$ is $S$-invariant, then $(7.15)$ holds in ${\Cal S}_q'(M)$ and $C\in 
{\Cal B}_{\Cal F}(M)$. It follows that the dimension of the vector space of basic currents with 
vanishing cohomology class is (at least) countable. In fact, it is possible to generate a sequence 
of linearly independent basic currents by iteration of the construction just explained. In other
terms, if $\Cal D$ is any (quasi) $S$-invariant distribution, then $\{T^{\ell}{\Cal D}\}$, ${\ell
\in \Bbb N}$, is a sequence of linearly independent (quasi) $S$-invariant distributions. We prove 
below that this procedure generates all basic currents for ${\Cal F}_{M\setminus\Sigma}$ [for 
$\Cal F$] (hence all quasi $S$-invariant [$S$-invariant] distributions) of 
{\it finite order}.
\enddemo

\specialnumber{7.6}
\proclaim {Lemma} Let $C$ be any real current of order $s\in {\Bbb N}${\rm ,} dimension {\rm (}\/and degree\/{\rm )}\/ equal to
{\rm 1,} closed in
${\Cal D}' (M\setminus\Sigma)${\rm .} Then there exist a current $U^{\ast}$ of order $s-1${\rm ,} dimension $2$ {\rm
(}\/and  degree $0${\rm )} and a meromorphic differential $h^+${\rm ,} with at most simple poles at $\Sigma${\rm },
such that 
$$dU^{\ast}=-\Re(h^+) + C \eqnu $$
in ${\Cal D}'(M\setminus\Sigma)${\rm .} If $C\in {\Cal H}^{-s}_q(M)$ is closed in ${\Cal S}_q'(M)${\rm ,} then 
there exist a current  $U^{\ast}\in H^{-s+1}_q(M)$ and a differential $h^+${\rm ,} holomorphic on $M_q${\rm ,} 
such that $(7.16)$ holds in ${\Cal H}^{-s}_q(M)${\rm .}
\endproclaim

\demo {Proof} If $C$ is closed in ${\Cal D}'(M\setminus\Sigma)$, since, by the Riemann-Roch theorem
[13, III.4], any cohomology class in $H^1(M\setminus\Sigma;{\Bbb R})$ can be represented as 
$\Re(h^+)$, where $h^+$ is a meromorphic differential with at most simple poles at $\Sigma$, it 
follows from [11, \S 15, Th. 12] that there exists a current $U^{\ast}\in {\Cal D}'(M\setminus
\Sigma)$ of dimension $2$ (and degree $0$) such that $(7.16)$ holds in ${\Cal D}'(M\setminus\Sigma)$.
In addition, if $C$ is of order $s\in {\Bbb N}$, then $U^{\ast}$ is of order $s-1$ by $(7.16)$.

 If $C\in {\Cal H}^{-s}_q(M)$ is closed in ${\Cal S}'_q(M)$, we cannot immediately deduce 
our statement from the standard de Rham theory. Instead, we argue as follows. Let $z$ be the canonical
coordinate at a zero $p\in \Sigma$ of (even) order $k$ and let $\pi:{\Cal U}_p\to D_p$ be the local 
covering map defined by $\pi_p(z):=z^{k/2+1}/(k/2+1)$ on a neighbourhood of $p\in M$ onto a 
neighbourhood $D_p\subset {\Bbb C}$ of the origin. The push-forward $(\pi_p)_{\ast}(C)\in {\Cal D}'
(D_p)$ is well-defined (see [11, Chap. III, \S 11] or [56, Chap. IX, \S 4]) and closed, hence 
exact. There exists therefore a current $V_p\in {\Cal D}'(D_p)$ such that $dV_p=(\pi_p)
_{\ast}(C)$. By the definition of the distributional space ${\Cal S}'_q(M)$, the pull-back $\pi_p
^{\ast}(V_p)$ is also (locally) defined as a current in ${\Cal S}_q'(M)$ over ${\Cal U}_p$.
In fact, a form $\alpha\in \Omega^2_q(M)$, locally over ${\Cal U}_p$, if and only if $\alpha=\pi_p
^{\ast}(\lambda_p)$, where $\lambda_p$ is a (unique) smooth $2$-form on $D_p$. Then, by definition, 
$\pi_p^{\ast}(V_p)(\alpha):=V_p(\lambda_p)$. A computation shows that $d\pi_p^{\ast}(V _p)=C$ 
in ${\Cal S}_q'(M)$, locally over ${\Cal U}_p$. Let $\{{\Cal U}_p\}$, $p\in \Sigma$, be a finite family
of disjoint open sets and, for each $p\in\Sigma$, let ${\Cal U}'_p\subset\subset{\Cal U}_p$ be a 
relatively compact open subset such that $p\in {\Cal U}'_p$. Let $U^{\ast}_p\in {\Cal S}_q'(M)$ be a 
current such that $U^{\ast}_p\equiv \pi_p^{\ast}(V_p)$ in ${\Cal U}'_p$ and $U^{\ast}_p\equiv 
0$ in $M\setminus {\overline {\Cal U}}_p$. Then the current of dimension $2$ (and degree $0$)
$$U^{\ast}_0 := \sum_{p\in\Sigma}\,U^{\ast}_p\,\, \in \,\, {\Cal S}_q'(M)  \eqnu $$
is such that $C-dU^{\ast}_0\in {\Cal S}_q'(M)$ has compact support in $M\setminus\Sigma$ and can 
therefore be extended to a current in ${\Cal D}'(M)$. In addition, since $C$ is closed and $dU^{\ast}
_0$ is exact in ${\Cal S}'_q(M)$, by Lemma 6.2 the cohomology class $[C-dU^{\ast}_0]=[C]\in H^1(M,
{\Bbb R})$. By [11, \S 15, Th. 12] and by the standard representation theorem of cohomology 
classes on a Riemann surface [13, III.2], there exists a current $U^{\ast}_1 \in {\Cal D}'(M)$ 
and a differential $h^+$, holomorphic on $M_q$, such that
$$dU^{\ast}_1= -\Re (h^+) + C - dU^{\ast}_0 \speqnu{7.17'}$$
in ${\Cal D}'(M)$. Let $U^{\ast}:= U^{\ast}_0+U^{\ast}_1\in {\Cal S}_q'(M)$. By $(7.17')$, the
identity $(7.16)$ holds in ${\Cal S}_q'(M)$. Since $C\in {\Cal H}^{-s}_q(M)$, $\Re(h^+)$
is smooth on $M$ and $\Omega^1_q$ is dense in ${\Cal H}^s_q(M)$, it follows that $dU^{\ast}\in 
{\Cal H}^{-s}_q(M)$. Hence $U^{\ast}\in H^{-s+1}_q(M)$, and $(7.16)$ holds in ${\Cal H}^{-s}_q(M)$.  
\enddemo

 Let ${\Cal F}:={\Cal F}_q$ be the horizontal foliation and $\Sigma:=\Sigma_q$ be the set 
of zeroes of a quadratic differential $q\in Q_{\kappa}$. We have:

\specialnumber{7.7}
\proclaim {Theorem} Let ${\Cal B}^s_{\Cal F}(M\setminus\Sigma)$ {\rm [}${\Cal B}^s_{\Cal F}(M)${\rm ]} 
be the vector space of all real basic currents of order $s\in {\Bbb N}$ for ${\Cal F}|_{M\setminus
\Sigma}$ {\rm [}for $\Cal F${\rm ].} There exist exact sequences
$$\eqalign{ &0\to {\Bbb R} \to {\Cal B}^{s-1}_{\Cal F}(M\setminus\Sigma) \to
 {\Cal B}^s_{\Cal F}(M\setminus\Sigma) \to H^{1,s}_{\Cal F}(M\setminus\Sigma;{\Bbb R})\,\,, \cr
            &0\to {\Bbb R} \to {\Cal B}_{\Cal F}^{s-1}(M) \to
{\Cal B}_{\Cal F}^s(M) \to H^{1,s}_{\Cal F}(M,{\Bbb R})\,\,.\cr} \eqnu $$
\endproclaim

\demo {Proof} Let $\delta_s:{\Cal B}^{s-1}_{\Cal F}(M\setminus\Sigma)\to {\Cal B}^s_{\Cal F}(M
\setminus\Sigma)$ be the map defined as follows: 
$$\delta_s (C):= d\,U_C^{\ast}\,\,,\,\,\,\,U_C:=- C\wedge\eta_T\,\,.\speqnu{7.18'}$$
The map $(7.18')$ is well-defined. In fact, if $C\in {\Cal B}^{s-1}_{\Cal F}(M\setminus\Sigma)$
[$C\in {\Cal B}^{s-1}_{\Cal F}(M)$], then $\delta_s(C)$ is by definition closed in ${\Cal D}'
(M\setminus\Sigma)$ [in ${\Cal H}^{-s}_q(M)$). Moreover, by Lemma 6.5, $U_C=-C\wedge\eta_T$ is a quasi 
$S$-invariant [an $S$-invariant] distribution of order $s-1$, hence $\imath_S dU_C^{\ast}={\Cal L}_SU^\ast_C=0$. 
It follows by $(6.4)$ that $\delta_s (C)\in {\Cal B}^s_{\Cal F}(M\setminus\Sigma)$ [$\delta_s(C)
\in {\Cal B}^s_{\Cal F}(M)$].

 The image under $\delta_s$ of ${\Cal B}^{s-1}_{\Cal F}(M\setminus\Sigma)$ [of ${\Cal B}^{s-1}
_{\Cal F}(M)$] consists of all basic currents $C^s\in{\Cal B}^s_{\Cal F}(M\setminus\Sigma)$ [$C^s\in 
{\Cal B}^s_{\Cal F}(M)$] with vanishing cohomology class. In fact, by Lemma 7.6, if $C^s$ has
vanishing cohomology class, then there exists a distribution $U_C$ of order $s-1$ such that $dU_C
^{\ast}=C^s$ in ${\Cal D}'(M\setminus\Sigma)$ [in ${\Cal H}^{-s}_q(M)$]. Since $C^s$ is a basic 
current,
$${\Cal L}_S U^\ast_C = \imath_S (dU_C^{\ast})=\imath_S C^s=0 \eqnu $$
in ${\Cal D}'(M\setminus\Sigma)$ [in $H^{-s}_q(M)$]. It follows that $U_C$ is a quasi $S$-invariant 
distribution [an $S$-invariant distribution] of order $s-1$, hence there exists a current $C^{s-1}
\in {\Cal B}^{s-1}_{\Cal F}(M\setminus\Sigma)$ [$C^{s-1}\in {\Cal B}^{s-1}_{\Cal F}(M)$] such that
$U_C=-C^{s-1}\wedge\eta_T$. By $(7.18')$, $C^s=\delta_s(C^{s-1})$.   

 The kernel of $\delta_s$ is given by the one-dimensional space ${\Bbb R}\cdot\eta_S$, i.e. 
by all basic currents which are scalar multiples of the closed $1$-form $\eta_S$. In fact, if $\delta_s
(C)=0$, then the distribution $U_C:=-C\wedge\eta_T\in {\Bbb R}\cdot\omega_q$. Since $C\in {\Cal B}^{s-1}
_{\Cal F}(M\setminus\Sigma)$ [$C\in {\Cal B}^{s-1}_{\Cal F}(M)$], it follows that $C= \imath_S U_C$ in ${\Cal D}'
(M\setminus\Sigma)$\break [in ${\Cal H}^{-s+1}_q(M)$], and hence $C\in {\Bbb R}\cdot\eta_S$. Let 
$\imath_{\Cal F}:{\Bbb R}\to {\Cal B}_{\Cal F}(M\setminus\Sigma)$ be the linear map given by $\imath_{\Cal
F}(1):=\eta_S$. We have proved that the following sequences are {\it exact}:
$$\eqalign{&0\to {\Bbb R}_{\overrightarrow{\,\,\,\,\imath_{\Cal F}\,\,\,\,}} 
{\Cal B}^{s-1}_{\Cal F}(M\setminus\Sigma)_{\overrightarrow{\,\,\,\, \delta_s\,\,\,\,}} 
{\Cal B}^s_{\Cal F}(M\setminus\Sigma)_{\overrightarrow{\,\,\,\, j_{\Sigma}\,\,\,\,}} 
H^{1,s}_{\Cal F}(M\setminus\Sigma;{\Bbb R})\,\,,\cr
&0\to {\Bbb R} _{\overrightarrow{\,\,\,\,\imath_{\Cal F}\,\,\,\,}}  
{\Cal B}^{s-1}_{\Cal F}(M)_{\overrightarrow{\,\,\,\,\delta_s\,\,\,\,}} 
{\Cal B}^s_{\Cal F}(M)  _{\overrightarrow{\,\,\,\,j_{\Sigma}\,\,\,\,}}  
H^{1,s}_{\Cal F}(M,{\Bbb R})\,\,.\cr} \eqnu $$
The cohomology maps 
$$\eqalign{ {\Cal B}^s_{\Cal F}(M\setminus\Sigma)&_{\overrightarrow{\,\,\,\,j_{\Sigma}\,\,\,\,}} 
H^{1,s}_{\Cal F}(M\setminus\Sigma)\,\,,\cr
{\Cal B}^s_{\Cal F}(M) & _{\overrightarrow{\,\,\,\,j_{\Sigma}\,\,\,\,}} 
H^{1,s}_{\Cal F}(M) \,\,,\cr} \speqnu{7.20'}$$
are surjective by definition. \enddemo

 \section{The non-uniform hyperbolicity of the Kontsevich-Zorich cocycle\\ 
and
an application to currents}

In Section 8.1 we will describe the invariant unstable and stable sub-bundles of the 
Kontsevich-Zorich cocycle in terms of the basic cohomology of the horizontal and vertical foliations 
of orientable quadratic differentials. This result is independent of the proof of the non-uniform 
hyperbolicity of the cocycle, concluded in Section 8.2. We then derive in Section 8.3 an Oseledec-type theorem 
for the cocycle $G^c_t$ on the bundle ${\Cal Z}^1_{\kappa}(M)$ of closed currents of order $1$, 
introduced in Section 8.1.

\demo{{\rm 8.1.} The stable and unstable sub\/{\rm -}\/bundles of the Kontsevich\/{\rm -}\/Zorich cocycle as bundles 
of basic currents}
 In Section~7, we have obtained results on the $H^{-s}$ basic cohomology of measured foliations 
for $s\gg 1$. The stable and unstable sub-bundles of the Kontsevich-Zorich cocycle are related to the 
basic cohomology in the case $s=1$.

 Let $H^{-1}_{\kappa}(M):=Q^{(1)}_{\kappa}\times H^{-1}(M)/\Gamma_g$, where the action of 
the mapping class group $\Gamma_g$ on the product is defined by {\it pull-back }on both factors. The
mapping class group is a quotient of the group of orientation-preserving diffeomorphisms and the 
action of a diffeomorphism on currents by pull-back is described in [11, Chap. III, \S 11] or [56, 
Chap. IX, \S 5]. The bundle $H^{-1}_{\kappa}(M)$ is the bundle of all distributions of order 
$1$ on $M$ over the stratum ${\Cal M}^{(1)}_{\kappa}$ of the moduli space. We introduce a cocycle 
$G^0_t$,
over the Teichm\"uller geodesic flow $G_t$, defined on the {\it Hilbert }bundle $H^{-1}_{\kappa}(M)$:
$$G^0_t:= (G_t\times \hbox{id}) /\Gamma_g \,\,\hbox{ on } \,\,H^{-1}_{\kappa}(M):= 
\bigl(Q^{(1)}_{\kappa}\times H^{-1}(M)\bigr)/\Gamma_g\,\,. \eqnu $$
We also introduce a cocycle $G^c_t$ over the Teichm\"uller flow, defined on the Hilbert bundle ${\Cal 
H}_{\kappa}^{-1}(M)$ of $H^{-1}$ currents over the stratum ${\Cal M}^{(1)}_{\kappa}$. The fiber of
the pull-back of the bundle ${\Cal H}_{\kappa}^{-1}(M)$ to the stratum $Q^{(1)}_{\kappa}$ of the 
Teichm\"uller space will be isomorphic by definition to the vector space ${\Cal H}_q^{-1}(M)$ at any 
$q\in Q^{(1)}_{\kappa}$. Since the  weighted Sobolev space of $1$-forms ${\Cal H}_q^1(M)$ can be identified, by its 
definition $(6.12)$, with the tensor product ${\Bbb R}^2\otimes H^1(M)$, the dual space ${\Cal H}_q
^{-1}(M)$ can be identified with the tensor product ${\Bbb R}^2\otimes H^{-1}(M)$. Let
$$G^c_t:= \hbox{diag}(e^{-t},e^t)\otimes G^0_t\,\,\hbox{ on } {\Cal H}_{\kappa}^{-1}(M):=
{\Bbb R}^2\otimes H^{-1}_{\kappa}(M)\,\,. \eqnu $$ 
The cocycle $G^c_t$ can be described as the cocycle obtained by parallel transport of currents, 
with respect to the trivial connection, along the orbits of the Teichm\"uller flow. In fact, by
the formulas $(2.7)$ describing the Teichm\"uller flow $G_t$ on $Q^{(1)}_{\kappa}$, we have:
$$G^c_t | {\Cal H}^{-1}_q(M) \equiv \hbox{id}: {\Cal H}^{-1}_q(M) \to {\Cal H}^{-1}_{G_t(q)}(M)\,\,.
\speqnu{8.2'}$$   
Let ${\Cal Z}^1_{\kappa}(M)\subset {\Cal H}^{-1}_{\kappa}(M)$ be the infinite-dimensional sub-bundle 
over ${\Cal M}^{(1)}_{\kappa}$ with fiber at $q\in {\Cal M}^{(1)}_{\kappa}$ given by the vector space of {\it closed currents of order }$1$: 
$${\Cal Z}^1_q(M):={\Cal Z}_q'(M) \cap  {\Cal H}^{-1}_q(M)\,\,. \eqnu $$ 
The sub-bundle ${\Cal Z}^1_{\kappa}(M)$ is $G^c_t$-invariant and the cocycle induced by $G^c_t$
on the $H^{-1}$ de Rham cohomology bundle is isomorphic to the Kontsevich-Zorich cocycle. The latter
is the essential motivation for the formulas $(8.1)$ and $(8.2)$ which define, respectively, the 
cocycles $G^0_t$ and $G^c_t$. Let
$$j_{\kappa}:{\Cal Z}^1_{\kappa}(M)\to {\Cal H}^1_{\kappa}(M,{\Bbb R}) \speqnu {8.3'}$$
be the natural bundle map, given fiber-wise by $(6.2')$, onto the cohomology bundle ${\Cal H}^1
_{\kappa}(M,{\Bbb R})$ (defined in \S 1). 

 Let ${\Cal B}^1_{\kappa,\pm}(M)\subset {\Cal Z}^1_{\kappa}(M)$ be the 
sub-bundles with fiber at $q\in {\Cal M}^{(1)}_{\kappa}$ given by the vector 
spaces ${\Cal B}^1_{\pm q}(M)$ of ${\Cal F}_{\pm q}$-basic currents of order 
$1$, defined by $(6.12')$. Let 
$${\Cal B}^1_{\kappa}(M):={\Cal B}^1_{\kappa,+}(M) + {\Cal B}^1_{\kappa,-}(M)
\,\,. \eqnu $$ 
\enddemo

\proclaim {Lemma} Let $\mu$ be any Borel probability ergodic invariant measure 
for the 
Teichm{\rm \"{\it u}}ller flow on a stratum ${\Cal M}^{(1)}_{\kappa}${\rm ,} with the property 
that the conditional 
measures induced by $\mu$ on the orbits of the action of the circle group 
${\rm S}(2,{\Bbb R})$ on 
${\Cal M}^{(1)}_{\kappa}$ are absolutely continuous with respect to the {\rm (}\/Haar\/{\rm )} 
Lebesgue measure{\rm .}     
\medbreak
\item{\rm (1)} The identity $j_{\kappa}\circ G^c_t=G^{KZ}_t \circ j_{\kappa}$ 
holds on 
${\Cal Z}^1_{\kappa}(M)${\rm .}
\smallbreak
\item{\rm (2)} The sub\/{\rm -}\/bundles ${\Cal B}^1_{\kappa,\pm}(M)\subset {\Cal Z}'_{\kappa}(M)$ are 
$G^c_t$\/{\rm -}\/invariant{\rm ,} $\mu$\/{\rm -}\/measurable and have $\mu$\/{\rm -}\/almost everywhere constant  {\rm (}\/finite\/{\rm
)} rank\/{\rm .}\/
\smallbreak
\item{\rm (3)} The cocycle $G^c_t$ has $\mu$\/{\rm -}\/almost everywhere strictly
positive {\rm [}\/strictly negative\/{\rm ]} 
Lyapunov spectrum on the invariant sub\/{\rm -}\/bundle ${\Cal B}^1_{\kappa,+}(M)$ {\rm [}${\Cal B}^1_{\kappa,-}(M)${\rm ].} 
\smallbreak
\item{\rm (4)} The splitting $(8.4)$ is direct and the restriction of the map 
$j_{\kappa}$ to the sub\/{\rm -}\/bundle ${\Cal B}^1_{\kappa}(M)$ is $\mu$\/{\rm -}\/almost 
everywhere injective\/{\rm .}
  
\endproclaim

\demo{Proof} $(1)$ It can be directly derived from the definition of the cocycles $G_t^{KZ}$, $G^c_t$ 
and of the cohomology map $j_{\kappa}$.
\medbreak $(2)$ The $G_t^c$-invariance of the sub-bundles ${\Cal B}^1_{\kappa,\pm}(M)$ follows 
immediately from the definitions. We will write down in detail the proof of the measurability of 
the bundle ${\Cal B}^1_{\kappa,+}(M)$. In fact, since ${\Cal B}^1_{-q}(M)={\Cal B}^1_{r(q)}(M)$, 
where the diffeomorphism $r:{\Cal M}^{(1)}_{\kappa}\to {\Cal M}^{(1)}_{\kappa} $ is given by 
the action of the counterclockwise rotation of a right angle, the measurability of ${\Cal B}^1
_{\kappa,-}(M)$ then follows. By Lemma 6.6, the measurability of ${\Cal B}^1_{\kappa,+}(M)$ can 
be reduced to the measurability of the bundle ${\Cal I}^1_{\kappa}(M)$ with fiber at $q\in {\Cal M}
^{(1)}_{\kappa}$ given by the vector space ${\Cal I}^1_q(M)$ of $S$-invariant distributions 
of order $1$. Let 
$$R^1_q(S):= \{Sv\,|\,v\in H^2_q(M)\}\subset H^1 (M) \,\,. \eqnu $$
We claim that any $q_0\in {\Cal M}^{(1)}_{\kappa}$ has a neighborhood 
${\Cal U}_0$ such that the following holds. There exists a basis ${\Cal V}(q)
:=\{v_k(q)\}$, $k\in{\Bbb N}$, of the vector
space $H^2_q(M)$ with the property that, for all $k\in {\Bbb N}$, the 
mappings $q\to Sv_k(q) \in H^1(M)$ are continuous on ${\Cal U}_0$.
In fact, there exists a neighbourhood $\widehat{\Cal U}_0$ of $q_0$ in the 
space of quadratic differentials with the prescribed pattern of zeroes, endowed
with the smooth topology, such that the following holds. Let $\Sigma_0$ be
the set of zeroes of $q_0$.
There exists on $\widehat{\Cal U}_0$ a continuous mapping $q\to\phi_q\in
\hbox{Diff}^+_0(M)$ such that the set of zeroes $\Sigma_{\hat q}$ 
of the quadratic differential ${\hat q}:=\phi_q^{\ast}(q)$ coincides 
with $\Sigma_0$ and ${\hat q}\equiv q_0$ in a neighbourhood of 
$\Sigma_0$. By definition, the quadratic differentials $q$ and 
${\hat q}$ yield the same point in   Teichm\"uller space. The weighted 
Sobolev spaces $H^s_{\hat q}(M)\equiv H^s_{q_0}(M)$ (although they 
do not coincide as Hilbert spaces). 
Let ${\Cal V}_0:=\{v_k\}$, $k\in {\Bbb N}$, be any basis of $H^2_{q_0}(M)$ 
and let ${\Cal V}(q):=\{v_k(q)\}$, where $v_k(q):=v_k\circ \phi_q^{-1}$,
for all $k\in {\Bbb N}$. Then ${\Cal V}(q)$ is a basis of $H^2_q(M)$ and
the mappings $q\to Sv_k(q)\in H^1(M)$ are continuous on $\widehat
{\Cal U}_0$. Let ${\Cal U}_0:= {\widehat{\Cal U}_0}/\hbox{Diff}^+(M)\subset
{\Cal M}^{(1)}_{\kappa}$. The claim is therefore proved.    

Let $v_0\equiv 1$ in the choice of the basis $\nu_0$.  Let then 
${\Cal W}(q)=\{w_k(q)\}$ be the orthonormal system in the Hilbert space 
$H_q^1(M)$   
constructed by applying the Gram-Schmidt orthonormalization procedure to the system $\{Sv_k\,|\, k\in {\Bbb N}\setminus\{0\}\}\subset
R^1_q(S)$. The system
${\Cal W}(q)$ is  well-defined at all $q\in {\Cal U}_0$ such that ${\Cal F}_q$ is {\it quasi-minimal}, since in that case the system
$\{Sv_k\,|\, k\in {\Bbb N}\setminus\{0\}\}$ is linearly independent. The subset of quadratic differentials with quasi-minimal horizontal
foliation has full measure with respect to any Borel measure $\mu$ with purely continuous conditional measures on $\mu$-almost every
orbit of the circle group. In fact, for any given
$q$, the subset of the orbit ${\rm SO}(2,{\Bbb R})\cdot q$ of the circle qroup with non-quasi-minimal horizontal foliation is countable as a
consequence of the structure theorem for measured foliations (or area-preserving vector fields) on closed orientable surfaces
\ref\Ktone, \ref\ZK, [47, Th.\ 3.1.7]. It follows that ${\Cal W}(q)$ is well-defined and the maps $q\to w_k(q)\in H^1(M)$ are continuous
$\mu$-almost everywhere. The orthogonal projection $\pi_q:H^1_q(M)\to R^1_q(S)^{\perp}$, can be written as
$$ \pi_q(v)= v-\sum_{k=1}^{\infty} \bigl(v,w_k(q)\bigr)_1\,w_k(q)\,\,,\,\,\,\,v\in H^1(M)\,,
\eqnu $$
where $(\cdot,\cdot)_1$ is the inner product in $H^1_q(M)$. The subspace 
$R^1_q(S)^{\perp}$ depends measurably on $q\in {\Cal U}_0$, since by $(8.6)$ 
the orthogonal projection $\pi_q$ does. The space 
of continuous linear functionals on $H^1(M)$ vanishing on the closure of the 
subspace $R^1_q(S)$ is 
by definition the space ${\Cal I}^1_q(M)$ of $S$-invariant distributions of 
order $1$. Since, by 
Lemma 6.4, the vector space ${\Cal B}^1_q(M)$ is canonically identified with 
${\Cal I}^1_q(M)$, the map $q\to {\Cal B}^1_q(M)$ is $\mu$-measurable. 

The almost everywhere finite dimensionality of ${\Cal B}^1_{\pm q}(M)$ is a consequence
of the injectivity of the cohomology map $j_q:{\Cal B}^1_{\pm q}(M)\to H^1(M,{\Bbb R})$. By
Theorem~7.7, if $C^+\in {\Cal B}^1_q(M)$ [$C^-\in {\Cal B}^1_{-q}(M)$] and $j_q(C^{\pm})=0$, then 
there exists $U^{\pm}\in L^2_q(M)$ such that ${\rm SU}^+=0$ [$TU^-=0$] and $dU^{\pm}=C^{\pm}$. Since, 
by [35, Th.~2], ${\Cal F}_q$ [${\Cal F}_{-q}$] is ergodic for $\mu$-almost all $q\in 
{\Cal M}^{(1)}_{\kappa}$, it follows that $j_q(C^{\pm})=0$ implies $C^{\pm}=0$ $\mu$-almost 
everywhere. Moreover, the dimension of ${\Cal B}^1_{\pm q}(M)$ is $\mu$-almost everywhere
constant by its invariance under the Teichm\"uller flow $G_t$ and the ergodicity of the system
$(G_t,\mu)$.  
\medbreak $(3)$ The argument is based on the construction of appropriate infinitesimal Lyapunov 
functions in the sense of \ref\Ktthree. Let $q\in {\Cal M}^{(1)}_{\kappa}$, $C^{\pm}\in {\Cal B}^1
_{\pm q}(M)$. Let ${\Cal D}^{\pm}\in H^{-1}_q(M)$, $U^{\pm}\in L^2_q(M)$ and $A^{\pm}\in {\Bbb R}$ 
be defined as follows:
$$\eqalign{&C^+\wedge\eta_T:= {\Cal D}^+=\partial^+_q U^+ + A^+\,\,,\,\,\,\,
                                              U^+\perp {\Cal M}^+_q\,\,,\cr
           &C^-\wedge\eta_S:= {\Cal D}^-=\partial^-_q U^- +  A^-\,\,,\,\,\,\,
                                             U^-\perp {\Cal M}^-_q\,\,,\cr} \eqnu $$
where $A^{\pm}={\Cal D}^{\pm}(1)$ and $U^{\pm}\in L^2_q(M)$ are uniquely determined by the
orthogonality conditions to the kernels ${\Cal M}^{\pm}_q\subset L^2_q(M)$ of the Cauchy-Riemann
operators (consisting of meromorphic, respectively anti-meromorphic, functions). The existence
of the functions $U^{\pm}\in L^2_q(M)$ is given by the existence of distributional solutions of the 
equations $\partial^{\pm} U={\Cal D}$ whenever ${\Cal D}\in H^{-1}_q(M)$ and ${\Cal D}(1)=0$ (Lemma
7.3 or [18, Prop.\ 4.6A]). Let ${\Cal L}^{\pm}: {\Cal B}^1_{\kappa, \pm}\to {\Bbb R}$ be the 
non-negative functions defined by 
$${\Cal L}^{\pm}_q(C^{\pm}):= |U^{\pm}|_0^2 + (A^{\pm})^2\,\,. \speqnu{8.7'}$$
Since the right-hand side of $(8.7')$ is invariant under the action of the mapping class group 
$\Gamma_g$, it defines a function on the measurable bundle ${\Cal B}^1_{\kappa,\pm}(M)$. The 
functions ${\Cal L}^{\pm}$ are measurable, since they extend to
continuous functions on the smooth 
bundle ${\Cal H}^{-1}_{\kappa}(M)$ with fiber ${\Cal H}^{-1}_q(M)$ at all $q\in{\Cal M}^{(1)}_
{\kappa}$ and ${\Cal B}^1_{\kappa,\pm}(M)$ are measurable sub-bundles by $(2)$. In addition, 
${\Cal L}^{\pm}_q(\cdot)$ are continuous functions, homogeneous of degree $2$ on each fiber and 
strictly positive on all non-zero $C^{\pm}\in {\Cal B}^1_{\pm q}(M)$. We claim that the function 
${\Cal L}^+$ [${\Cal L}^-$] is strictly increasing [decreasing] along the forward orbits of the 
cocycle $G^c_t$ if the horizontal foliation ${\Cal F}_q$ [the vertical foliation ${\Cal F}_{-q}$] 
is {\it quasi\/{\rm -}\/minimal}.

 By the formulas $(2.7)$, $(2.7')$, which describe the evolution of a 
quadratic differential $q$ and of the corresponding Cauchy-Riemann operators 
$\partial^{\pm}_q$ under the Teichm\"uller flow $G_t$ on $Q^{(1)}_{\kappa}$, 
and by the definition of $G^c_t$, we have:
$$\eqalign{{\Cal D}^+_t&= C^+ \wedge\eta_T(t)= e^t\,C^+\wedge \eta_T= e^t\,
{\Cal D}^+_0\,\,,\cr
 {\Cal D}^-_t&= C^-\wedge\eta_S(t)=e^{-t}\,C^-\wedge\eta_S= e^{-t}\,
{\Cal D}^-_0\,\,.\cr}\eqnu $$
Let $q_t:=G_t(q)$ and ${\Cal M}^{\pm}_t:=N(\partial^{\pm}_t)\subset L^2_q(M)$, 
where $\partial^{\pm}_t=\partial^{\pm}_q$ at $q=q_t$. Let $V^{\pm}_t\in 
L^2_q(M)$ be the unique solutions of the equations
$$ {\Cal D}^{\pm}_0= \partial ^{\pm}_t V^{\pm}_t + A^{\pm}_0\,\,,\,\,\,\, 
V^{\pm}_t \perp{\Cal M}^{\pm}_t \,\,. \eqnu $$
We have
$$ {\Cal D}^{\pm}_t= e^{\pm t} {\Cal D}^{\pm}_0= e^{\pm t}\{ \partial^{\pm}_t 
V^{\pm}_t + 
A^{\pm}_0\}=\partial^{\pm}_t U^{\pm}_t + A^{\pm}_t\,\,, \eqnu $$
where $U^{\pm}_t:= e^{\pm t} V^{\pm}_t$ and $A^{\pm}_t:={\Cal D}^{\pm}_t(1)$. 
It follows finally from $(8.7)$--$(8.10)$ that  
$${\Cal L}^{\pm}\circ G^c_t (C^{\pm})=e^{\pm 2t}\{|V^{\pm}_t|_0^2 + 
(A^{\pm}_0)^2\} \,\,. \eqnu $$ 

 Let $\pi^{\pm}_t:L^2_q(M)\to {\Cal M}^{\pm}_t$ be the orthogonal projections. We claim that 
$V^{\pm}_t$ is a  solution  of the following O.D.E.  in $L^2_q(M)$:
$$\eqalign{& u^{\pm}_t=\partial^{\pm}_t v^{\pm}_t + \pi^{\mp}_t(u^{\pm}_t) \,\,, \cr
         & {d\over{dt}} u^{\pm}_t =\partial^{\mp}_t v^{\pm}_t - M^{\pm}_t(u^{\pm}_t)\,\,,\cr} \eqnu $$
where $M^{\pm}_t:=M^{\pm}_q$, at $q=q_t$, and $M^{\pm}_q:L^2_q(M)\to{\Cal M}^{\pm}_q$ is the linear 
operator defined as follows. Let $\{m^{\pm}_1, \ldots ,m^{\pm}_g\}$ be an orthonormal basis of ${\Cal M}
^{\pm}_q\subset L^2_q(M)$. Let $m^{\pm}_i=\partial^{\pm}v^{\pm}_i + \pi^{\mp}_q(m^{\pm}_i)$ be the
decompositions of the functions $m^{\pm}_i$ according to the splitting $(2.3)$. Then
$$M^{\pm}_q(u):= \sum_{i=1}^g (u,\partial^{\pm}v^{\pm}_i)_q\,m^{\pm}_i\,\,. \speqnu{8.12'}$$
The equation  $(8.12)$ yields, as $(2.12)$, an O.D.E. of the form $u'=F(u,t)$ in the Banach space $L^2_q(M)$. In fact, its first line is just the
Hilbert space orthogonal decomposition according to the splitting $(2.3)$. Hence $F(t,u)$ \pagebreak is smooth and linear in the second
variable. By the standard theorems on O.D.E.'s in Banach spaces [40, Chap. IV] the Cauchy problem for $(8.12)$ has (locally) a unique solution
for any initial data.

 By the variational formulas $(2.10)$ for the Cauchy-Riemann operators, any solution 
$u^{\pm}_t$ of $(8.12)$ satisfies the identity ${d\over{dt}}(\partial^{\pm}_t u^{\pm}_t)\equiv 0$ 
(no matter how the operator $M^{\pm}_q$ is defined) and, by the choice $(8.12')$ of $M^\pm_q$, 
$u^{\pm}_t\perp {\Cal M}^{\pm}_t$. It follows that $V^{\pm}_t$ is the unique solution $u^{\pm}_t$
of $(8.12)$ satisfying the initial condition $u^{\pm}_0=V^{\pm}_0$.

 Let $u^{\pm}_t$ be any solution of $(8.12)$. Since $u^{\pm}_t\perp {\Cal M}^{\pm}_t$,
by $(2.3)$ there exists a unique zero average function $w^{\pm}_t\in H^1(M)$ such that $u^{\pm}_t=
\partial^{\mp}_tw^{\pm}_t$. We claim that, if $\partial^{\pm}_q u^{\pm}_0$ is a real-valued distribution, then $w^{\pm}_t$ is real-valued for all $t\in {\Bbb R}$. In fact, $\partial^{\pm}_t
u^{\pm}_t\equiv \partial^{\pm}_q u^{\pm}_0$, which is real-valued, and $\partial^{\pm}_t u^{\pm}_t
=\triangle_t w^{\pm}_t$, where $\triangle_t:=\partial^+_t \partial^-_t$ is the Laplace operator of 
the metric $R_t$ induced by the quadratic differential $q_t=G_t(q)$. Since the Laplacian is a {\it 
real }operator, it suffices to show that, for any $q\in Q^{(1)}_{\kappa}$, the unique zero average 
solution $w_0\in H^1_q(M)$ of the equation $\triangle_q w_0=0$ is $w_0\equiv 0$. In fact, $|\partial
^+_q w_0|_0^2=-(\triangle_q w_0,{\overline {w_0}})=0$, hence $w_0$ is a meromorphic function in 
$H^1_q(M)$, that is, by [18, Prop.\ 3.2], a constant function. The zero average condition then 
implies $w_0\equiv 0$. We can then compute:
$$\eqalign{ {d\over{dt}} |u^{\pm}_t|_0^2=&2\Re\,(u^{\pm}_t,{d\over{dt}}u^{\pm}_t)_q=
2\Re\,(\partial^{\mp}_t w^{\pm}_t,\partial^{\mp}_tv^{\pm}_t)_q \cr
=&2\Re\,(\partial^{\pm}_t w^{\pm}_t,\partial^{\pm}_tv^{\pm}_t)_q=2\Re\,(\partial^{\pm}_tw^{\pm}_t,
\partial^{\mp}_tw^{\pm}_t)_q  \cr
=&2\Re\,\{|S_tw^{\pm}_t|^2_0-|T_tw^{\pm}_t|^2_0 \pm 2i(S_tw^{\pm}_t,T_tw^{\pm}_t)_q\} \cr
=&2\{2|S_t w^{\pm}_t|_0^2-|u^{\pm}_t|_0^2\}=2\{|u^{\pm}_t|_0^2-2|T_t w^{\pm}_t|_0^2\}\,\,.\cr}
\eqnu $$
We have exploited the fact that, since $w^{\pm}_t$ is real-valued, $(S_tw^{\pm}_t,
T_tw^{\pm}_t)_q\in~{\Bbb R}$. In addition, if $\partial^{\pm}_q u_0\not=0$ and ${\Cal F}_q$ [${\Cal F}_{-q}$] is quasi-minimal; then $S_t
w^{\pm}_t
\not=0$\break [$T_tw^{\pm}_t\not=0$]. In fact, otherwise, 
by the quasi-minimality property, $w^{\pm}_t=0$ hence $u_t=0$ and $\partial^{\pm}u^{\pm}_0\equiv\partial
^{\pm}_t u^{\pm}_t=0$.

Since the functions $V^{\pm}_t$ satisfy $(8.12)$, we have by $(8.11)$ and $(8.13)$:
$$\eqalign{ {d\over{dt}}{\Cal L}^+\circ G_t^c(C^+)&=e^{2t}\{4|S_t w^+_t|_0^2 + 2(A^+_0)^2\}>0\,\,,\cr 
         {d\over{dt}}{\Cal L}^-\circ G_t^c(C^-)&=-e^{-2t}\{4|T_t w^-_t|_0^2 + 2(A^-_0)^2\}<0\,\,. \cr}
\speqnu{8.13'}$$
We have therefore constructed an infinitesimal  Lyapunov function 
${\Cal L}^+$ [${\Cal L}^-$] on the finite-dimensional measurable bundle ${\Cal B}^1_{\kappa,+}(M)$
[${\Cal B}^1_{\kappa,-}(M)$],  strictly increasing [decreasing] along the Teichm\"uller orbit $G_t(q)$, for all $q\in {\Cal M}^{(1)}_{\kappa}$
with quasi-minimal horizontal  [vertical] foliation. The set of $q\in {\Cal M}^{(1)}_{\kappa}$ with quasi-minimal
horizontal and vertical foliations has full $\mu$-measure. In fact, by the
structure theorem for measured foliations (or area-preserving vector fields)
on closed orientable surfaces, for any orbit of the circle group on ${\Cal M}^{(1)}_{\kappa}$ the subset of quadratic differentials $q$ with 
non-quasi-minimal horizontal or vertical foliation is countable \ref\Ktone, \ref\ZK, [47, Th.\ 3.1.7]. Since, by the
definition $(8.7')$, ${\Cal L}^{\pm}(C^{\pm})>0$ if $C^{\pm}\in {\Cal B}^1_{\kappa,\pm}\setminus
\{0\}$, by [31, Th.\ 2.1], the Lyapunov exponents of the cocycle $G^c_t$ restricted to the 
invariant bundle ${\Cal B}^1_{\kappa,+}(M)$ [${\Cal B}^1_{\kappa,-}(M)$] are all strictly positive 
[strictly negative] $\mu$-almost everywhere. 
\medbreak $(4)$ Let $q\in {\Cal M}^{(1)}_{\kappa}$ be a quadratic differential with {\it quasi-minimal} 
horizontal and vertical foliations. Let $C^{\pm}\in {\Cal B}^1_{\pm q}(M)$ be such that $j_q(C^++C^-)
=0$. Then the cohomology class $c=j_q(C^+)=-j_q(C^-)\in H^1(M_q,{\Bbb R})$ belongs by $(1)$ and $(3)$ to 
both the stable and the unstable invariant bundle of the Kontsevich-Zorich cocycle, hence $c=0$. 
By the injectivity of the cohomology map $j_{\kappa}$ on ${\Cal B}^1_{\kappa,\pm}(M)$ (proved 
above), it follows that $C^+=C^-=0$. We have therefore proved that ${\Cal B}^1_q(M) \cap {\Cal B}^1
_{-q}(M)=\{0\}$ and that $j_{\kappa}$ is injective on ${\Cal B}^1_q(M)\oplus {\Cal B}^1_{-q}
(M)$, for $\mu$-almost all $q\in {\Cal M}^{(1)}_{\kappa}$. 
\enddemo 

By Lemma 8.1 every cohomology class $c\in H^1(M_q,{\Bbb R})$, which can be represented
by a current of order $1$, basic for the horizontal [vertical] foliation of a `generic' quadratic 
differential $q\in {\Cal M}^{(1)}_{\kappa}$, belongs to the unstable [stable] invariant sub-bundle 
of the Kontsevich-Zorich cocycle. In the remaining part of this section we will prove the converse
statement. The argument is based on the geometric estimate of the Poincar\'e constant obtained in
Lemma 6.9 and on the {\it logarithmic law for geodesics }in the moduli space due to H. Masur \ref\Msthree.   

 Let $E^+_q [E^-_q] \subset H^1(M_q,{\Bbb R})$ be the (Lagrangian) subspaces 
corresponding to the subset of strictly positive [strictly negative] Lyapunov exponents in the 
Oseledec's decomposition of the Kontsevich-Zorich cocycle at an Oseledec {\it regular }point $q\in 
{\Cal M}^{(1)}_{\kappa}$. By Oseledec's theorem \ref\Os, [32, Th.\ S.2.9] the set of regular 
points has full measure with respect to any ergodic $G_t$-invariant measure $\mu$ on ${\Cal M}^{(1)}
_{\kappa}$.    

\proclaim{Lemma} Let $\mu$ be any ergodic $G_t$\/{\rm -}\/invariant measure with the 
property that the conditional measures induced by $\mu$ on the orbits of the 
circle action are absolutely continuous{\rm .} For $\mu$\/{\rm -}\/almost all $q\in {\Cal M}
^{(1)}_{\kappa}${\rm ,} every cohomology class\break $c^{\pm}\in E^{\pm}_q$ can be 
represented by a current $C^{\pm}\in{\Cal B}^1_{\pm q}(M)${\rm .} In addition{\rm ,} the 
following estimates holds{\rm .} There exists a $\mu$\/{\rm -}\/measurable function $K:{\Cal M}
^{(1)}_{\kappa}\to {\Bbb R}^+$ such that{\rm ,} if $c\in E^+_q\oplus E^-_q$ has the 
harmonic representation $c=[\Re\,(m^+ q^{1/2})]${\rm ,} $m^+\in  {\Cal M}^+_q${\rm ,} 
then for a unique current $C\in {\Cal B}^1_{+q}(M)\oplus {\Cal B}^1_{-q}(M)$ 
and a unique zero average function $U\in L^2_q(M)$ we have\/{\rm :}\/
$$\eqalign{& dU=C \,-\,\Re\,(m^+ q^{1/2}) \,\,, \cr 
           & |U|_0\leq K(q) |\!|c|\!|_q\,\,, \cr }\eqnu $$
where $|\!|c|\!|_q$ is the Hodge norm of $c\in H^1(M_q,{\Bbb R})${\rm .} 
\endproclaim  

\demo{Proof} Let $q\in {\Cal M}^{(1)}_{\kappa}$ be a regular point of the 
Kontsevich-Zorich cocycle, 
$c=c^+_q(m^+)\in E^+_q [E^-_q]$, where $m^+\in {\Cal M}^+_q$ is a meromorphic function and\break
$c^+_q:{\Cal M}^+_q\to H^1(M,{\Bbb R})$ is the isomorphism given by $(2.4)$. By $(2.5)$ the Hodge
norm $|\!|c|\!|_q=|m^+|_0$. Let $c_t:=G_t^{KZ}(c)=[\Re(m^+_t q_t^{1/2})]$, where $m^+_t\in{\Cal M}^+_t$,
the space of meromorphic functions with respect to the complex structure determined by the quadratic 
differential $q_t=G_t(q)$. Since the $L^2_q$ norm is invariant under the Teichm\"uller flow, there 
exist a $\mu$-measurable function $K_1(q)>0$ and an exponent $0<\lambda<1$ such that, if $c\in 
E^+_q$ [$c\in E^-_q$], 
$$|\!|c|\!|_{q_t}=|m^+_t|_0\leq K_1(q) |m^+|_0\,\exp (-\lambda |t|)\,\,,\,\,\,\,t\leq 0\,\,[t\geq 0] 
\,\,. \eqnu $$
Let $U_t\in L^2_q(M)$ be the unique function with zero average satisfying
$$ dU_t = c_t-c=\Re\,(m^+_t q_t^{1/2})\,-\,\Re\,(m^+q^{1/2})\,\,. \eqnu $$
By $(2.12)$, $(2.13)$ there exists a smooth one-parameter family of functions\break
$v_t\in H^1(M)$ such that the function $U_t$ satisfies the following Cauchy problem in $L^2_q(M)$:
$$\eqalign { {d\over{dt}}U_t &= 2\Re\,(v_t)\,\,, \cr
                         U_0 &= 0\,\,. \cr } \speqnu{8.16'}$$
If $v_t$ is chosen with zero average, since by definition the Dirichlet form ${\Cal Q}(v,v)=|\partial^{\pm}_q v|_0^2$ for any orientable
quadratic differential, by the Poincar\'e inequality (Lemma 6.9) and by the orthogonality  of the decomposition in the first line of $(2.12)$,  
$$|v_t|_0\leq K_{g,\sigma} |\!|q_t|\!|^{-1} |\partial^+_t v_t|_0 \leq K_{g,\sigma} |\!|q_t|\!|^{-1}  
|m^+_t|_0 \,\,, \eqnu $$
where $|\!|q_t|\!|$ denotes the length of the shortest segment with endpoints in the set of zeroes 
of the quadratic differential $q_t$ with respect to the induced metric. By 
$(8.16')$ and $(8.17)$,
there exists a measurable function $K_2(q)>0$ such that, if $c\in E^+_q$ [$c\in E^-_q$], 
$$\left|{d\over{dt}} U_t\right|_0\leq 2|v_t|_0\leq K_2(q)\,|m^+|_0\,|\!|q_t|\!|^{-1}\exp (- \lambda |t|) \,\,, \,\,\,\, t\leq 0\,\,[t\geq 0]
\,\,.\eqnu
$$ Since $U_0=0$, by Minkowski's integral inequality, $(8.18)$ implies the 
following estimate:
$$ |U_t|_0 \leq K_2(q) |m^+|_0\,\left|\int_0^t e^{-\lambda |s|}|\!|q_s|\!|^{-1}\, 
ds\,\right|\,\,,\,\,\,\,t\leq 0 \,[t\geq 0] \,\,. \speqnu{8.18'}$$
We claim that the integral in $(8.18')$ converges as $t \to -\infty$ [$t\to +\infty$]. In \ref\Msthree 
H. Masur proved a {\it logarithmic law }for the Teichm\"uller geodesic flow on the moduli space. The 
following estimate is the simplest step of Masur's argument [44, Prop.\ 1.2]:
$$ \limsup_{t\to \pm \infty} {{-\log|\!|q_t|\!|}\over {\,\log |t|}}\,\leq \, {1\over 2}\,\,,  
\eqnu $$
The statement holds for $\mu$-almost all quadratic differentials $q\in {\Cal M}^{(1)}_{\kappa}$.
It follows that our claim concerning $(8.18')$ holds $\mu$-almost everywhere. An alternative 
proof of the $\mu$-almost everywhere convergence of the integral in $(8.18')$ can be derived 
from a recent result of A. Eskin and H. Masur [12, Lemma 5.4], which
strengthens a previous result by 
W. Veech [68, Th.\ 0.2], [69, Corollary 2.8], by applying integration by parts and Birkhoff's ergodic 
theorem. By $(8.18')$ and $(8.19)$, there exists a $\mu$-measurable function $K_3(q)>0$ such that, 
if $c\in E^+_q$ [$c\in E^-_q$], 
$$|U_t|_0 \leq K_3(q) |\!|c|\!|_q\,\,,\,\,\,\,t\leq 0 \,[t\geq 0] \,\,. \eqnu $$   
\smallskip
\noindent Let $U\in L^2_q(M)$ be any weak limit of $U_t$ as $t\to -\infty$ [$t\to +\infty$], which 
exists by $(8.20)$ since all bounded subsets of the Hilbert space $L^2_q(M)$ 
are sequentially weakly compact. By contraction in $(8.16)$,
$$\eqalign{ {\rm SU}_t &= -\Re\,(m^+)\,+\, e^t\,\Re\,(m^+_t)\,\,,\,\,\,\,t\leq 0\,,\cr \noalign{\vskip5pt}
           [TU_t &= \Im\,(m^+) \,-\, e^{-t}\,\Im\,(m^+_t)\,\,,\,\,\,\, t\geq 0] \,,\cr}\eqnu $$
then by taking the limit as $t\to -\infty$ [$t\to +\infty$],
$$ \eqalign{ {\rm SU} &= -\Re\,(m^+)\,\,,\cr  \noalign{\vskip5pt}
             [TU &= \Im\,(m^+)]\,\,.\cr} \speqnu{8.21'}$$
The identities $(8.21')$ hold in the distributional space $H^{-1}_q(M)$. It follows that the current 
$C\in {\Cal H}^{-1}_q(M)$ determined by the identity
$$dU=C\,-\,\Re\,(m^+ q^{1/2}) \eqnu $$
is a basic current of order $1$ for the horizontal foliation ${\Cal F}_q$ [for the vertical foliation
${\Cal F}_{-q}$], representing the cohomology class $c\in E^+_q$ [$c\in E^-_q$]. In fact, 
by $(8.21')$ and its definition $(8.22)$, the current $C$ is closed and $\imath_S C=0$ [$\imath_T C=
0$]. Hence, by $(6.4)$, ${\Cal L}_SC=0$ [${\Cal L}_TC=0$]. If the horizontal foliation ${\Cal F}_q$ 
[the vertical foliation ${\Cal F}_{-q}$] is ergodic (a property which holds for $\mu$-almost all $q \in {\Cal M}^{(1)}_{\kappa}$ [35, Th.\ 2]) the limit function $U\in L^2_q(M)$ is uniquely determined by $(8.21')$
and by the zero average condition. In addition, by $(8.20)$, $U$ satisfies the 
estimate in $(8.14)$.    \pagebreak

 Since any cohomology class $c\in E^+_q\oplus E^-_q$ can be split as $c=c^++c^-$, 
where $c^{\pm}\in E^{\pm}_q$ and $|\!|c^{\pm}|\!|_q\leq\delta(q) |\!|c|\!|_q$, where $\delta(q)
>0$ is the {\it distorsion }(in the sense of [41, p.\ 269]) of the splitting $E^+_q\oplus 
E^-_q$, the cohomology map
$$j_q: {\Cal B}^1_q(M)\oplus {\Cal B}^1_{-q}(M) \,\to\,E^+_q\oplus E^-_q 
\eqnu $$
is surjective and the estimate $(8.14)$ holds, for $\mu$-almost all $q\in {\Cal M}^{(1)}_{\kappa}$. 
The map $(8.23)$ was proved to be $\mu$-almost everywhere injective in Lemma 
8.1,$\,(4)$.
\enddemo
 
\specialnumber{8.2'}
\demo{Remark} Lemma 8.2 can be viewed as a {\it regularity }result. In fact, by Theorem 
7.1$(i)$, for almost all $q\in {\Cal M}^{(1)}_{\kappa}$, every cohomology class $c\in H_q:=\{c\in 
H^1(M_q,{\Bbb R})\,|\,c\wedge [\Im(q^{1/2})]=0\}$, can be represented by a basic current $C\in 
{\Cal B}_q(M)$ of finite order $l>1$. The subspace $H_q$ has codimension $1$, hence it has dimension 
$2g-1$. By Lemma 8.2 and Theorem 8.5 below, if $c$ belongs to a Lagrangian subspace $E^+_q\subset 
H_q$ (which has dimension equal to $g\geq 2$), then $C$  actually has order $1$.
\enddemo   

 We can finally derive from Lemma 8.1 and Lemma 8.2 the following characterization 
of the invariant unstable [stable] bundle $E^+_{\kappa}$ [$E^-_{\kappa}$] of 
the Kontsevich-Zorich cocycle on any stratum ${\Cal M}^{(1)}_{\kappa}$.  

\proclaim {Theorem} Let $\mu$ be any ergodic $G_t$\/{\rm -}\/invariant Borel probability measure on
${\Cal M}^{(1)}_{\kappa}$ with the property that the conditional measures induced on the orbits
of the circle action are absolutely continuous{\rm .} The invariant unstable {\rm [}\/stable\/{\rm ]} sub\/{\rm -}\/bundle of the 
Kontsevich\/{\rm -}\/Zorich cocycle coincides $\mu$\/{\rm -}\/almost everywhere with the bundle of $H^{-1}$ basic 
cohomologies for the horizontal {\rm [}\/vertical\/{\rm ]} foliation of quadratic differentials\/{\rm :}
$$\eqalign{ E_{\kappa}^+&=j_{\kappa}\bigl({\Cal B}^1_{\kappa,+}(M)\bigr):=
{\Cal H}^{1,1}_{\kappa,+}(M,{\Bbb R})\,\,, \cr
E_{\kappa}^-&=j_{\kappa}\bigl({\Cal B}^1_{\kappa,-}(M)\bigr):=
{\Cal H}^{1,1}_{\kappa,-}(M,{\Bbb R})\,\,. \cr} \eqnu $$
\endproclaim

 Theorem 8.3 implies that the unstable subspace $E^+_q$ [the stable subspace 
$E^-_q$] depends locally only on the fundamental class of the horizontal measured foliation 
${\Cal F}_q$ [of the vertical measured foliation ${\Cal F}_{-q}$], a property conjectured in 
[38, Th.\ p.\ 8, (3)]. In fact, by Lemma 6.8 and Theorem 8.3, we have:

\specialnumber{8.3'} 
\proclaim{Corollary} For almost all {\rm (}\/regular\/{\rm )} points $q_0\in {\Cal M}^{(1)}_{\kappa}$ of the 
Kontsevich\/{\rm -}\/Zorich coycle{\rm ,} there exists a neighbourhood ${\Cal U}_0\subset{\Cal M}^{(1)}_{\kappa}$
with the property that $q_0\in{\Cal U}_0$ and{\rm ,} for almost all {\rm (}\/regular\/{\rm )} points $q\in {\Cal U}_0${\rm ,}
$$\eqalign{ [\Im(q^{1/2})]=[\Im(q_0^{1/2})]\in H^1(M,\Sigma_{\kappa},{\Bbb R}) &\Longrightarrow 
                                                E^+_q=E^+_0\,\,, \cr      
             [\Re(q^{1/2})]=[\Re(q_0^{1/2})]\in H^1(M,\Sigma_{\kappa},{\Bbb R}) &\Longrightarrow 
                                                E^-_q=E^-_0 \,\, ,\cr}
\eqnu $$
where $E^{\pm}_q$ {\rm [}$E^{\pm}_0${\rm ]} denotes the invariant unstable{\rm ,} respectively stable{\rm ,} subspace of the 
Kontsevich\/{\rm -}\/Zorich cocycle at $q$ {\rm [}$q_0${\rm ].}
\endproclaim 

  8.2. {\it The Kontsevich\/{\rm -}\/Zorich cocycle is non\/{\rm -}\/uniformly hyperbolic}.
\vglue-8pt
  \proclaim{Lemma} The set of quadratic differentials $q_0$ having the following properties is dense
in ${\Cal M}^{(1)}_{\kappa}$\/{\rm :}\/
  \item{\rm (1)} The horizontal foliation ${\Cal F}_{q_0}$ is Lagrangian {\rm (}\/Definition {\rm 4.3);}
 \item{\rm (2)} Any neighbourhood $\,{\Cal U}^{(1)}_{\kappa}\subset {\Cal M}^{(1)}
_{\kappa}$ of 
$q_0$ contains a compact positive measure set ${\Cal P}^{(1)}_{\kappa}$ such that 
\itemitem{\rm (a)}  all $q\in {\Cal P}^{(1)}_{\kappa}$ are Birkhoff generic points for the Teichm{\rm \"{\it u}}ller flow and 
Oseledec regular points for the Kontsevich\/{\rm -}\/Zorich cocycle{\rm ;}
\itemitem{\rm (b)} the stable and the unstable subspace $E^{\pm}_q$ of the Kontsevich\/{\rm -}\/Zorich cocycle 
depend continuously on $q\in {\Cal P}^{(1)}_{\kappa}${\rm ;}
  \smallbreak \itemitem{\rm (c)} the Poincar{\rm \'{\it e}} dual $P({\Cal L}_0)$ of the Lagrangian subspace ${\Cal L}_0:={\Cal L}({\Cal F}
_{q_0})${\rm ,} generated by the regular trajectories of ${\Cal F}_{q_0}${\rm ,} is transverse to $E^-_q$ for
all $q\in {\Cal P}^{(1)}_{\kappa}${\rm .}

\endproclaim

{\it Proof}. Let ${\Cal F}_{\kappa}(M)$ be the set of all isotopy equivalence classes of orientable 
measured foliations ${\Cal F}$ on $M$ with a finite set of canonical saddle-like singularities of 
multiplicities $\kappa:=(k_1, \ldots ,k_{\sigma})$. Let ${\Cal F}'_{\kappa}(M)\subset {\Cal F}_{\kappa}
(M)$ be the subset given by the measured foliations ${\Cal F}$ such that there exists $q\in 
Q_{\kappa}$ with ${\Cal F}_q={\Cal F}$ (or ${\Cal F}_{-q}={\Cal F}$). The space $Q_{\kappa}$ has
locally a product structure modeled on ${\Cal F}'_{\kappa}(M)\times {\Cal F}'_{\kappa}(M)$. Let
${\Cal U}_{\kappa}={\Cal U}^+_{\kappa}\times {\Cal U}^-_{\kappa}$, ${\Cal U}^{\pm}_{\kappa}\subset
{\Cal F}'_{\kappa}(M)$, be {\it any }open subset of $Q_{\kappa}$  having a product structure. By 
Luzin's theorem, there exists a compact set of positive measure ${\Cal P}^-_{\kappa} \subset {\Cal U}
^-_{\kappa}$ such that the basic cohomology map ${\Cal F}\to H^{1,1}_{\Cal F}(M,{\Bbb R})$ is 
continuous on ${\Cal P}^-_\kappa$. Let ${\Cal F}_0^-$ be a density point of ${\Cal P}^-_{\kappa}$ and let 
$E^-_0=H^{1,1}_{{\Cal F}_0^-}(M,{\Bbb R})$. By Lemma 4.4, there exists a Lagrangian foliation 
${\Cal F}^+_0\in {\Cal U}^+_{\kappa}$ such that the Poincar\'e dual $P({\Cal L}_0)$ of the 
Lagrangian subspace ${\Cal L}_0:={\Cal L}({\Cal F}^+_0)$ is transverse to $E^-_0$. Let $q_0\in 
{\Bbb R}^+\cdot {\Cal U}_{\kappa}$ be the quadratic differential uniquely defined by the condition  $({\Cal F}
_{q_{\lower2pt\hbox{$\scriptscriptstyle 0$}}},{\Cal F}_{-q_{\lower2pt\hbox{$\scriptscriptstyle 0$}}})\in {\Bbb R}^+\cdot ({\Cal
F}^{+}_0,{\Cal F}^-_0)$ and the condition that the total area
$A(q_0)=1$.

Let ${\Cal U}^{(1)}_{\kappa}\subset {\Cal M}^{(1)}_{\kappa}$ be a neighbourhood of $q_0$. By 
construction and Theorem~8.3, the map $q\to E^-_q$ is defined and continuous on a positive measure 
subset of ${\Cal U}^{(1)}_{\kappa}$ having density $1$ at $q_0$ and the subspace $P({\Cal L}_0)$ is 
transverse to $E^-_{q_0}$. Hence, by Birkhoff's ergodic theorem, Oseledec's theorem and Luzin's theorem,
there exists a compact set of positive measure ${\Cal P}^{(1)}_{\kappa}\subset {\Cal U}^{(1)}_{\kappa}$,
consisting~of quadratic differentials that are Birkhoff generic points for the Teichm\"uller flow and 
Oseledec regular points for the Kontsevich-Zorich coycle, such that the maps $q\to E^{\pm}_q$ are 
continuous and $P({\Cal L}_0)$ is transverse to $E^-_q$, for all $q\in {\Cal P}^{(1)}_{\kappa}$.           
\phantom{phun} \hfill\qed 
  
 \vglue-8pt
\proclaim {Theorem} The Kontsevich\/{\rm -}\/Zorich cocycle is non\/{\rm -}\/uniformly hyperbolic on every connected 
component  ${\Cal C}^{(1)}_{\kappa}$ of a stratum ${\Cal M}^{(1)}_{\kappa}${\rm .} In fact{\rm ,} the Lyapunov 
exponents of the normalized restriction of the absolutely continuous invariant measure $\mu^{(1)}
_{\kappa}$ to ${\Cal C}^{(1)}_{\kappa}$ form a symmetric subset of the real line and satisfy the 
following inequalities\/{\rm :}
$$ \eqalign{ 1=\lambda_1>\lambda_2\geq&\cdots\geq\lambda_g>0>\lambda_{g+1}=-\lambda_g\geq\cdots \cr
              &\cdots \geq \lambda_{2g-1}=-\lambda_2 > \lambda_{2g}=-\lambda_1=-1 \,\,.\cr} \eqnu $$
\endproclaim

\demo{Proof} The argument will proceed by contradiction. If $(8.26)$ does not 
hold on a connected
component ${\Cal C}^{(1)}_{\kappa}\subset{\Cal M}^{(1)}_{\kappa}$, there exists a natural number 
$k<g$ such that $\lambda_k\not=0$ and $\lambda_{k+1}=0$. Let $E^+_k \subset {\Cal H}^1_{\kappa}
(M,{\Bbb R})$ be the $k$-dimensional unstable invariant sub-bundle of the Kontsevich-Zorich cocycle, 
which corresponds to the set of strictly positive Lyapunov exponents $\{\lambda_1, \ldots ,\lambda_k\}$ 
in the Oseledec's decomposition. By Corollary 5.3 and 5.5, the function $\Phi_k$  defined by
$(5.7)$ satisfies, for $\mu^{(1)}_{\kappa}$-almost all $q\in {\Cal C}^{(1)}_{\kappa}$, the
following identity:
$$ \Phi_k\bigl(q, E^+_k(q)\bigr) \equiv \Lambda_1(q) +\cdots +\Lambda_g(q)\,\,. \eqnu $$
Let $\{c_1, \ldots ,c_g\}\subset H^1(M_q,{\Bbb R})$ be any orthonormal system, with respect to the
Hodge norm, such that $\{c_1, \ldots ,c_k\}$ is a basis of $E^+_k(q)$ and let $c_i=[\Re (m^+_i 
q^{1/2})]$ be the harmonic representation $(2.4)$. By formulas $(2)$ and $(3)$ 
in Lemma $5.2'$, $(8.27)$ implies that
$$ B_q(m^+_i,m^+_j)=0 \,\,,\,\,\,\, k+1\leq i,j\leq g\,\,. \speqnu{8.27'}$$
The argument will therefore be concluded by the construction of a   subset of 
the connected component ${\Cal C}^{(1)}_{\kappa}$ with positive measure (with respect to the measure $\mu^{(1)}_{\kappa}$) 
where $(8.27')$ does not hold.

 We claim that there exists a compact  set of positive measure ${\Cal P}^{(1)}_{\kappa} 
\subset {\Cal C}^{(1)}_{\kappa}$  with the following properties:
\medbreak
\item{(1)} All $q\in {\Cal P}^{(1)}_{\kappa}$ are Oseledec regular points, hence the unstable 
subspace $E^+_k(q)\subset H^1(M,{\Bbb R})$ is defined and it is a $k$-dimensional (isotropic) 
subspace. The mapping $q\to E^+_k(q)\subset H^1(M,{\Bbb R})$ is continuous on ${\Cal P}^{(1)}
_{\kappa}$.
\smallbreak \item{(2)} For any $\varepsilon>0$, the flat structure induced by $q\in {\Cal P}^{(1)}_{\kappa} $, 
that is the flat Riemannian surface with singularities $(M,R_q)$, contains $g\geq 2$ embedded 
{\it metric cylinders }$A_1^{\varepsilon}, \ldots ,A_g^{\varepsilon}$ (in the sense of the {Definition }in 
[45, p.\ 457]) with non-horizontal waist curves $\gamma_1^{\varepsilon}, \ldots ,\gamma_g^{\varepsilon}$ such 
that $M\setminus\cup\{\gamma_1^{\varepsilon}, \ldots ,\gamma_g^{\varepsilon}\}$ is
homeomorphic to a sphere 
minus $2g$ disjoint disks and

\vglue-15pt
$$d_q\bigl(E^+_k(q),P(\Gamma^{\varepsilon})\bigr) < \varepsilon\,\,, \eqnu $$
\item{}
where $P(\Gamma^{\varepsilon})\subset H^1(M,{\Bbb R})$ is the Poincar\'e dual of the (Lagrangian)
subspace $\Gamma^{\varepsilon}\subset H_1(M,{\Bbb R})$ generated by the homology classes of the waist
curves $\gamma_1^{\varepsilon}, \ldots ,\gamma_g^{\varepsilon}$ and $d_q$ denotes the metric induced on the 
Grassmannian by the Hodge inner product on $H^1(M_q,{\Bbb R})$.   
 
 Let $q_0\in {\Cal C}^{(1)}_{\kappa}$ a quadratic differential with the properties listed
in the statement of Lemma 8.4. The existence of such a differential in every connected component
of ${\Cal M}^{(1)}_{\kappa}$ follows immediately from the density property asserted by the lemma. 
Since the horizontal foliation ${\Cal F}_{q_0}$ is Lagrangian, the flat structure induced by $q_0$
contains at least $g\geq 2$ flat horizontal cylinders with waist curves 
$\gamma_1(q_0), \ldots ,
\gamma_g(q_0)$ which generate a Lagrangian subspace ${\Cal L}_0\subset 
H_1(M_{q_0},{\Bbb R})$. There 
exists a neighbourhood ${\Cal U}^{(1)}_{\kappa}\subset {\Cal C}^{(1)}_{\kappa}$ of $q_0$ such that 
any $q\in {\Cal U}^{(1)}_{\kappa}$ has $g\geq 2$ distinct closed (possibly 
non-horizontal) regular 
trajectories $\gamma_1(q), \ldots ,\gamma_g(q)$, isotopic respectively to $\gamma_1
(q_0), \ldots ,\gamma_g
(q_0)$. Let ${\Cal P}^{(1)}_{\kappa}\subset {\Cal U}^{(1)}_{\kappa}$ be a 
compact set of positive
measure with the properties listed in Lemma 8.4.  

 Let $q\in {\Cal P}^{(1)}_{\kappa}$ and let $T<0$ be a backward return time of the 
Teichm\"uller orbit $G_t(q)$ to the compact  set  of positive measure ${\Cal P}^{(1)}_{\kappa}$. The 
quadratic differential $q_T:=G_T(q)$ has therefore $g\geq 2$ closed non-horizontal regular 
trajectories $\gamma_1(q_T), \ldots ,\gamma_g(q_T)$ which generate in the homology vector space the 
Lagrangian plane ${\Cal L}_0$ with Poincar\'e dual $P({\Cal L}_0)$ transverse 
to $E^-(q_T)$.
Let $\Gamma^T \subset H_1(M_q,{\Bbb Z})$ be the Lagrangian subspace generated by the homology 
classes of the closed non-horizontal regular $q$-trajectories $\gamma_1^T(q), \ldots ,\gamma_g^T(q)$, 
obtained from $\gamma_1(q_T), \ldots ,\gamma_g(q_T)$ by pull-back under the Teichm\"uller map $f_T:(M,q)
\to (M,q_T)$. Let $P(\Gamma_T)\subset H^1(M_q,{\Bbb R})$ be the Poincar\'e dual of the subspace
$\Gamma_T$. By definition of the Kontsevich-Zorich cocycle, $P(\Gamma_T)=
G^{KZ}_{-T}
\bigl(P({\Cal L}_0)\bigr)$. Since $P({\Cal L}_0)$ is transverse to $E^-(q_T)$, by the Oseledec's 
theorem, for any $\varepsilon >0$ there exists $T_{\varepsilon}<0$ such that, if the return time $T<
T_{\varepsilon}$, then the Lagrangian subspace $P(\Gamma^{\varepsilon}):= P(\Gamma_T)$ satisfies 
the estimate $(8.28)$. We therefore choose a return time $T<T_{\varepsilon}$ and let $\{\gamma_1
^{\varepsilon}, \ldots ,\gamma_g^{\varepsilon}\}:=\{\gamma_1^T(q), \ldots ,\gamma_g^T(q)\}$. The embedded metric 
cylinders $A_1^{\varepsilon}, \ldots ,A_g^{\varepsilon}$ are then determined by the closed regular 
non-horizontal $q$-trajectories $\gamma_1^{\varepsilon}, \ldots ,\gamma_g^{\varepsilon}$. We have thus 
constructed a compact  set of positive measure  ${\Cal P}^{(1)}_{\kappa}\subset {\Cal C}^{(1)}_{\kappa}$ 
with the required properties.

 Let $\{\gamma_1, \ldots ,\gamma_g\}$ be a set of disjoint simple closed curves on $M$ such that 
$M\setminus\{\gamma_1, \ldots ,\gamma_g\}$ is homeomorphic to a sphere minus $2g$ disjoint (paired) disks. Let 
${\Cal U}_{\kappa}(\gamma_1, \ldots ,\gamma_g)$ be the open set of quadratic differentials $q\in Q_{\kappa}$ 
having closed regular non-horizontal trajectories $\gamma_1(q), \ldots ,\gamma_g(q)$ isotopic respectively 
to $\gamma_1, \ldots ,\gamma_g$. We will construct a smooth deformation 
$$\Psi: (0,1]^g\times {\Cal U}_{\kappa}(\gamma_1, \ldots ,\gamma_g) \to {\Cal U}_{\kappa}(\gamma_1, \ldots ,
\gamma_g)  \eqnu $$
having the properties listed below. Let $q_t:=\Psi_t(q)$, $t=(t_1, \ldots ,t_g)\in (0,1]^g$, then   
\medbreak
\item{(1)} if $t_i=1$ for all $i=1, \ldots ,g$, then $q_t=q$. 
\smallbreak \item{(2)}  if $t_i\to 0$ for all $i=1, \ldots ,g$, then the Riemann surface carrying $q_t$ converges
to a Riemann surface with nodes, pinched along the curves $\gamma_1, \ldots ,\gamma_g$, hence $q_t$
converges in the moduli space ${\Cal M}_g$ to a quadratic differential $q_0 \in {\Cal S}_g$, the 
boundary component of the moduli space ${\Cal M}_g$ consisting of regular quadratic differentials 
on the Riemann sphere with $2g$ paired punctures, having poles of order $2$ at all punctures 
(considered in \S 4).
\smallbreak \item{(3)} the zero set of $q_t$ coincides with the zero set $\Sigma_q$ of $q$ for all $t\in 
(0,1]^g$ and the cohomology class $[\Im(q_t^{1/2})]=[\Im(q^{1/2})] \in H^1(M,\Sigma_q,{\Bbb R})$.
\medbreak
 We will construct $\Psi_t(q)$ as a composition of deformations $\Psi^{(i)}_t$, $i=\break1, \ldots ,g$,
pinching along $\gamma_i$. Let $q\in {\Cal U}_{\kappa}(\gamma_1, \ldots ,\gamma_g)$ and let $\theta_i
\in (0,\pi)$ be the angle between the oriented horizontal foliation ${\Cal F}_q$ and the closed regular
oriented trajectory $\gamma_i(q)$, isotopic to $\gamma_i$. Let $H^+_s$, $s\in {\Bbb R}$, be the
horocycle flow on $Q_g$, given by the action of the one-parameter subgroup of ${\rm SL}(2,{\Bbb R})$
$$ H^+_s:=  \bmatrix 1\,\,s\cr 0\,\,1 \endbmatrix  \,\,. \eqnu  $$
Let $s_i:=\cot\theta_i\in {\Bbb R}$. Then the horizontal foliation of the quadratic differential 
$H^+_{-s_i}(q)\in Q_{\kappa}$ coincides with ${\Cal F}_q$ and $\gamma_i(q)$ is a closed regular 
{\it vertical }trajectory of $H^+_{-s_i}(q)$. In addition, since the set of closed regular 
trajectories is invariant under the horocycle action, $H^+_{-s_i}(q)\in {\Cal U}_{\kappa}(\gamma_1,
\ldots,\gamma_g)$. Let ${\Cal V}_{\kappa}(\gamma_i)\subset {\Cal U}_{\kappa}(\gamma_1, \ldots ,\gamma_g)$ be 
the set of quadratic differentials having a closed regular vertical trajectory $\gamma_i(q)$ isotopic 
to $\gamma_i$ defined in the proof of Lemma~$4.4'$ and let $\Phi^{(i)}: (0,1] \times {\Cal V}_{\kappa}
(\gamma_i)\to {\Cal V}_{\kappa}(\gamma_i)$ be the deformation defined by the identity $(a)$ in
$(4.31)$. We then define: 
$$\eqalign {&(a)\,\,\Psi^{(i)}_{t_i}(q):= H^+_{s_i}\circ \Phi^{(i)}_{t_i} \circ H^+_{-s_i}(q)\,\,,
                    \,\,\,\,(t_i,q)\in (0,1]\times {\Cal U}_{\kappa}(\gamma_1, \ldots ,\gamma_g)\,, \cr
& (b)\,\,\Psi_t(q):=\Psi^{(1)}_{t_1}\circ \cdots \circ\Psi^{(g)}_{t_g}(q)\,\,,\,\,\,\,\,\,\,\,\,\,\,
\,\,\, (t,q)\in (0,1]^g\times {\Cal U}_{\kappa}(\gamma_1, \ldots ,\gamma_g)\,.\cr} \eqnu $$
It can be checked that the pinching deformation $\Psi:(0,1]^g\times {\Cal U}_{\kappa}(\gamma_1, \ldots ,
\gamma_g) \to {\Cal U}_{\kappa}(\gamma_1, \ldots ,\gamma_g)$ is well-defined and has the required 
properties. In fact, by construction, if $q\in {\Cal U}_{\kappa}(\gamma_1, \ldots ,\gamma_g)$ and 
$\theta_j\in (0,\pi)$ is the oriented angle of the oriented horizontal foliation ${\Cal F}_q$ 
and the closed regular non-horizontal oriented\break $q$-trajectory $\gamma_j(q)$, then $\theta_j$
is preserved by the deformation $\Psi^{(i)}_{t_i}(q)$, $i,j=1, \ldots ,g$. Hence, the composition
in $(b)$ is well-defined. The properties $(1)$ and $(2)$ hold by the definition $(8.31)$ and by 
the corresponding properties of the deformations $\Phi^{(i)}_{t_i}$. By its construction, in
the proof of Lemma $4.4'$, $\Phi^{(i)}_{t_i}(q)=q$, if $t_i=1$, and $\Phi^{(i)}_{t_i}(q)$ converges 
in the compactified moduli space  as $t_i\to 0$, to a regular quadratic differential with $2$ paired 
poles of order $2$ and real strictly positive residues on a Riemann surface of genus $g-1$ with $2$ 
punctures. Hence $q_t:=\Psi_t(q)$ converges, as $t_i\to 0$ for all $i=1, \ldots ,g$, to a regular 
quadratic differential $q_0\in {\Cal S}_g$ (on a punctured Riemann sphere \pagebreak with $2g$ paired punctures)
and real  strictly positive residues at all punctures. Finally, $(3)$ follows from the corresponding 
property of $\Phi^{(i)}_{t_i}$, established in the proof of Lemma $4.4'$, since the horocycle flow 
$H^+_s$ preserves the horizontal measured foliation.

Let $\{a_1, \ldots ,a_g,b_1, \ldots ,b_g\}$ be a canonical homology basis on a Riemann surface $M$ 
and $\{\theta_1, \ldots ,\theta_g\}$ be the dual basis of holomorphic differentials, defined by the 
standard condition $\theta_i(a_j)=\delta_{ij}$. Let ${\Cal A}\subset H_1(M,{\Bbb R})$ be the\break
Lagrangian subspace generated by the homology classes of the curves\break $\{a_1, \ldots ,a_g\}$. By the 
definition of the Poincar\'e duality, the Poincar\'e dual $P({\Cal A})\subset H^1(M,{\Bbb R)}$ 
is generated by the cohomology classes of the harmonic forms $\{\Im(\theta_1), \ldots ,\Im(\theta_g)\}$. 
Let $q\in Q_{\kappa}$ and let $\phi^+_i:=\theta_i/q^{1/2}\in L^2_q(M)$ be the corresponding 
meromorphic functions. Since $\Im(\theta_i)=-\Re(\imath\,\theta_i)$, the Lagrangian subspace 
$P({\Cal A})$ is represented, under the ${\Bbb R}$-linear isomorphism defined in $(2.4)$, by the 
$g$-dimensional {\it real }subspace $V^+_q({\Cal A})\subset{\Cal M}^+_q$ generated by the system 
$\{\imath\,\phi^+_1, \ldots ,\imath\,\phi^+_g\}$. Let $\{m^+_1, \ldots ,m^+_g\}$ be any orthonormal basis of 
$V^+_q({\Cal A})$ over $\Bbb R$, with respect to the (symmetric) inner product induced on $V^+_q
({\Cal A})$ by the hermitian structure of $L^2_q(M)$. Let $B$ be the symmetric complex matrix defined
by $B_{ij}:=B_q(m^+_i,m^+_j)$. The matrix $B$ depends on the quadratic differential and on the choice
of the orthonormal basis $\{m^+_1, \ldots ,m^+_g\}$ of $V^+_q({\Cal A})$. By $(4.4)$ and $(4.5)$, we have:
$$ B= C^{-1}\,{{d\Pi}\over {d\mu_q}}\,(C^t)^{-1} \,\,, \eqnu $$
where $C$ is the matrix of base change defined in $(4.5)$ and $\Pi$ is the period matrix with respect 
to the canonical basis $\{a_1, \ldots ,a_g,b_1, \ldots ,b_g\}$. Since in 
the present case the matrix $C$ is {\it purely imaginary}, by $(4.5)$ we have that $CC^t=-CC^{\ast}=
-\Im(\Pi)$. Let $q_0\in {\Cal S}_g$ be a regular quadratic differential on the $2g$-punctured Riemann
sphere with positive residues at all punctures. We claim that, if $q\to q_0$ as $M_q$ is pinched 
along the homotopically non-trivial curves $a_1, \ldots ,a_g$, then
$$\lim _{q\to q_0} B_q(m^+_i,m^+_j) \,=\,\delta_{ij}\,\,, \eqnu $$
uniformly with respect to the orthonormal basis $\{m^+_1, \ldots ,m^+_g\}$ of $V^+_q({\Cal A})$. In fact, 
by $(4.5')$ and $(4.33)$,
$$(\phi^+_i,\phi^+_j)_q \,+\, {1\over{2\pi}} \delta_{ij} \log  |t_i|  \eqnu $$
is bounded as $t\to 0$, where $t=(t_1, \ldots ,t_g)$ are the pinching parameters, described in \S 4,
corresponding respectively to the curves $a_1, \ldots ,a_g$. It follows that the purely imaginary matrix $C$ of base change converges, up to
left multiplication by a diagonal purely imaginary matrix, to the {\it  compact }subgroup ${\rm O}(g,{\Bbb R})\subset {\rm GL}(g,{\Bbb R})$ of
real
$g
\times g$ orthogonal matrices.  In fact, let $\Delta$ be the $g\times g$ diagonal matrix with entries
$$
\Delta_{ij} := i \left| {\log |t_i|\over 2\pi}\right|^{1/2}\delta_{ij}\,\,. \speqnu{8.34'}$$
By the boundedness of (8.34), as $t\to 0$,
$$
\Delta^{-1}C\to O(g,\Bbb R)\,\, . \eqnu
$$
On the other hand, by Lemma 4.2 and Lemma 4.2$'$, as $t\to 0$,
$$
\Delta^{-1} {d\Pi\over d\mu_q} \Delta^{-1}\to I_{g\times g}\,\,. \speqnu{8.35'}
$$
The claim $(8.33)$ follows then from $(8.32)$, $(8.35)$ and $(8.35')$.

 Let $q^{(0)}\in {\Cal C}^{(1)}_{\kappa}$ be a density point of the compact 
set of positive measure ${\Cal P}^{(1)}_{\kappa}$ constructed above. Let $\varepsilon>0$ and let $\{\gamma_1^{\varepsilon},
\ldots,\gamma_g^{\varepsilon}\}$ be the system of disjoint closed curves with the properties stated in the 
description of ${\Cal P}^{(1)}_{\kappa}$ at  (2) (see p.~75). Let $q^{(0)}_t=\Psi_t(q^{(0)})$ be the pinching 
deformation along the curves $\gamma_1^{\varepsilon}, \ldots ,\gamma_g^{\varepsilon}$ constructed above. Since
$q^{(0)}_t$ converges as $t\to 0$ to a regular quadratic differential $q^{(0)}_0\in {\Cal S}_g$ with 
real strictly positive residues and there exists a canonical homology basis $\{a_1, \ldots ,a_g,b_1, \ldots ,
b_g\}$ such that $a_i=\gamma_i^{\varepsilon}$, for all $i\in\{1, \ldots ,g\}$, by 
$(8.33)$ there exists 
$\tau>0$ such that the following holds. Let $t\in (0,1]^g$ with $|t|<\tau$ and $q:=q^{(0)}_t$. Let 
$\{c_1, \ldots ,c_g\}\subset H^1(M_q,{\Bbb R})$ be {\it any }orthonormal basis, with respect to the Hodge 
norm, of the Poincar\'e dual $P(\Gamma^{\varepsilon})$ of the Lagrangian subspace $\Gamma^{\varepsilon} 
\subset H_1(M_q,{\Bbb R})$ generated by the homology classes of the curves $\gamma_1^{\varepsilon}, \ldots ,
\gamma_g^{\varepsilon}$ and let $c_i=[\Re (m^+_i q^{1/2})]\in H^1(M_q,{\Bbb R})$ be the harmonic 
representation $(2.4)$. Then
$$ |B_q(m^+_i,m^+_i)|>{1\over 2}\,\,,\,\,\,\, i=1, \ldots ,g\,. \eqnu $$
Let us choose $t^{\ast}\in (0,1]^g$ such that $|t^{\ast}|<\tau$ and let $\Psi_{\ast}:=\Psi_t$ for 
$t=t^{\ast}$. There exists an open neighbourhood ${\Cal U}^{(0)}_{\kappa}\subset {\Cal U}_{\kappa}
(a_1, \ldots ,a_g)$ such that $q^{(0)}\in {\Cal U}^{(0)}_{\kappa}$, $(8.28)$ holds for all $q\in 
{\Cal U}^{(0)}_{\kappa}\cap {\Cal P}^{(1)}_{\kappa}$ and $(8.36)$ holds for all $q\in \Psi_{\ast}
({\Cal U}^{(0)}_{\kappa})$. In fact, the first assertion follows from the fact that, by construction,
$E^+_k(q)$ depends continuously on $q\in {\Cal P}^{(1)}_{\kappa}$, the latter from the continuity
of the bilinear form $B_q$. Finally, since the inequality $(8.28)$ is preserved at any Oseledec 
regular point of the form $\Psi_t(q)$, $q\in {\Cal P}^{(1)}_{\kappa}$, by property $(3)$ of $\Psi_t$,
Lemma 6.8 and Theorem 8.3, if $\varepsilon >0$ is chosen sufficiently small, the subset of regular 
points of the set $\Psi_{\ast}({\Cal U}^{(0)}_{\kappa}\cap{\Cal P}^{(1)}_{\kappa})$ yields a  subset of ${\Cal C}^{(1)}_{\kappa}$ of
positive 
measure (with
respect to $\mu^{(1)}_{\kappa}$) where the  identities $(8.27')$ fail for all $i=j\in \{k+1, \ldots ,g\}$.   
\enddemo

\demo{{\rm 8.3.} The Oseledec's theorem for the bundle of closed currents of order $1$}
By Theorem 8.3, which describes the invariant sub-bundles of the Kontsevich-Zorich cocycle 
in terms of basic currents of order $1$, and Theorem 8.5, which proves the (non-uniform) hyperbolicity
of the cocycle, we can derive an Oseledec's theorem for the restriction of the cocycle $G^c_t$ to the
bundle of all {\it closed }currents of order $1$. This is a crucial step in order to obtain results 
on the deviation of ergodic averages for measured foliations. In fact, by the standard trace theorem 
for Sobolev spaces [1, Th.\ 5.4 (5)], the (return) leaves of a\break $1$-dimensional foliation on a 
$2$-dimensional manifold can be regarded as currents of order $1$.
\enddemo

\proclaim {Lemma} The forward and backward Lyapunov spectra of the cocycle $G^0_t$ on the Hilbert
bundle $H^{-1}_{\kappa}(M)${\rm ,} defined by $(8.1)${\rm ,} are contained in the closed interval $[-1,1]${\rm ,} in 
the sense that{\rm ,} for all $q\in {\Cal M}^{(1)}_{\kappa}$ and any ${\Cal D}\in H^{-1}_q(M)${\rm ,}
$$-1\leq \liminf_{t\to \pm\infty}{1\over |t|}\,\log |G^0_t({\Cal D})|_{-1}\leq 
\limsup_{t\to \pm \infty}{1\over |t|}\,\log |G^0_t({\Cal D})|_{-1}\leq 1\,\,.
\eqnu $$
The forward and backward Lyapunov spectra of the cocycle $G^c_t$ on the Hilbert bundle ${\Cal H}^{-1}
_{\kappa}(M)${\rm ,} defined by $(8.2)${\rm ,} $(8.2')${\rm ,} are contained in the closed interval $[-2,2]${\rm ,} in the 
sense that{\rm ,} for all $q\in {\Cal M}^{(1)}_{\kappa}$ and any $C\in {\Cal H}^{-1}_q(M)${\rm ,}
$$-2\leq \liminf_{t\to \pm \infty}{1\over |t|}\,\log |G^c_t(C)|_{-1}\leq\limsup_{t\to \pm\infty}
{1\over |t|}\,\log |G^c_t(C)|_{-1}\leq 2\,\,.\speqnu{8.37'}$$
There is a continuous {\rm (}\/orthogonal\/{\rm )} $G^c_t${\rm -}invariant splitting 
$${\Cal H}^{-1}_{\kappa}(M):= I^-_{\kappa}(M)\oplus I^+_{\kappa}(M) \eqnu $$
into the closed infinite\/{\rm -}\/dimensional sub{\rm -}bundles $I^{\pm}_{\kappa}(M)$ with fiber at $q\in 
{\Cal M}^{(1)}_{\kappa}$ given by
$$\eqalign{ I^+_q(M)&:= \{ C\in {\Cal H}^{-1}_{\kappa}(M)\,|\, \imath_S C =0\}\,\,, \cr
         I^-_q(M)&:= \{ C\in {\Cal H}^{-1}_{\kappa}(M)\,|\, \imath_T C =0\}\,\,. \cr}\speqnu{8.38'}$$
The restriction of the cocycle $G^c_t$ to $I^+_{\kappa}(M)$ [$I^-_{\kappa}(M)$] has non{\rm -}negative 
{\rm [}\/non{\rm -}positive\/{\rm ]} forward Lyapunov spectrum{\rm .} A corresponding statement{\rm ,} reversing the role of 
$I^+_{\kappa}(M)$ and $I^-_{\kappa}(M)${\rm ,}  holds for the backward spectrum{\rm .}
\endproclaim

\demo{Proof} Since the $L^2_q$ norm is invariant under the Teichm\"uller flow $G_t$, by $(2.7')$ and the definition $(8.1)$ of $G^0_t$, we
have:
$$e^{-|t|}\, |{\Cal D}|_{-1}\leq |G^0_t({\Cal D})|_{-1}\leq e^{|t|}\,|{\Cal D}|_{-1}\,,\eqnu $$
where $|\cdot|_{-1}$ denotes the norm induced on the fiber at any $q\in {\Cal M}^{(1)}_{\kappa}$ by 
the Hilbert structure of $H^{-1}_q(M)$. The bounds $(8.37)$ follow immediately 
from~$(8.39)$. \pagebreak

The sub-bundles $I^{\pm}_{\kappa}(M)$ are $G^c_t$-invariant by definition and, by the
definition $(8.2)$ of $G^c_t$, the Lyapunov spectra of $G^c_t$ on $I^{\pm}_{\kappa}(M)$ can be 
obtained from the Lyapunov spectrum of $G^0_t$ by translation of $\pm 1$ respectively. The 
orthogonal splitting $(8.38)$ holds since $M$ has dimension $2$ and the bounds 
$(8.37')$ follow from the bounds $(8.37)$ and the splitting $(8.38)$. 
\enddemo

\proclaim{Theorem} The infinite dimensional closed sub\/{\rm -}\/bundle ${\Cal Z}_{\kappa}^1(M)\subset 
{\Cal H}_{\kappa}^{-1}(M)$ of closed currents of order $1$ is $G^c_t$\/{\rm -}\/invariant and has a measurable 
$G^c_t$\/{\rm -}\/invariant splitting\/{\rm :}\/ 
$${\Cal Z}_{\kappa}^1(M)={\Cal B}^1_{\kappa,+}(M)\oplus {\Cal B}^1_{\kappa,-}(M)\oplus {\Cal E}
_{\kappa}^1(M) \,\,. \eqnu $$

\medbreak \item{\rm (1)} The measurable sub\/{\rm -}\/bundles ${\Cal B}^1_{\kappa,\pm}(M)$ have fibers
equal to ${\Cal B}^1_{\pm q}(M)${\rm ,} the finite dimensional vector space of basic 
currents of order $1$ for ${\Cal F}_{\pm q}${\rm ,} at $\mu^{(1)}_{\kappa}$\/{\rm -}\/almost 
all quadratic differentials $q\in {\Cal M}^{(1)}_{\kappa}${\rm .} The sub\/{\rm -}\/bundle 
${\Cal E}_{\kappa}^1(M)\subset {\Cal Z}_{\kappa}^1(M)$ is the closed infinite 
dimensional continuous bundle of exact currents{\rm ,} which has everywhere defined 
fibers ${\Cal E}_q^1(M):=\{dU\,|\,U\in  L^2_q(M)\}${\rm . }
\smallbreak \item{\rm (2)} The restriction of the cocycle $G^c_t$ to the bundle $${\Cal B}^1
_{\kappa}(M)={\Cal B}^1_{\kappa,+}(M)\oplus {\Cal B}^1_{\kappa,-}(M)$$ is 
measurably isomorphic to the Kontsevich\/{\rm -}\/Zorich cocycle $G_t^{KZ}$ on the real 
cohomology bundle ${\Cal H}^1_{\kappa}(M,{\Bbb R})${\rm ,} hence it has 
the Lyapunov spectrum $(8.26)${\rm .} The invariant sub\/{\rm -}\/bundles ${\Cal B}^1_{\kappa,
\pm}(M)$ correspond to the  {\rm (}\/strictly\/{\rm )} positive{\rm ,} respectively {\rm (}strictly\/{\rm )} 
negative{\rm ,} Lyapunov exponents{\rm .} The Lyapunov spectrum of $G^c_t$ on 
${\Cal E}_{\kappa}^1(M)$ is reduced to the single Lyapunov exponent $0${\rm .} 

\endproclaim

\demo{Proof} The bundle ${\Cal Z}_{\kappa}^1(M)$ and the splitting $(8.40)$ are $G^c_t$-invariant
by the definition $(8.2)$, $(8.2')$ of the cocycle $G^c_t$. By Lemma 8.1, Theorem 8.3 and Theorem
8.5, the restriction of $G^c_t$ to the sub-bundle ${\Cal B}^1_{\kappa}(M)$ is measurably isomorphic
to the Kontsevich-Zorich cocycle, hence it has Lyapunov spectrum given by 
$(8.26)$. In fact, the
restriction of the cohomology map $j_{\kappa}$, defined by $(8.3')$, to 
${\Cal B}^1_{\kappa}(M)$ is 
injective by Lemma 8.1$(4)$, and surjective by Theorem 8.3 and Theorem 8.5. 
By Lemma 8.1(1),
it follows that $j_{\kappa}$ is an isomorphism of the cocycle $G^c_t$ onto the Kontsevich-Zorich
cocycle $G^{KZ}_t$. By Lemma 8.1(3), the Lyapunov exponents of $G^c_t$ on ${\Cal B}^1_{\kappa,+}
(M)$ [on ${\Cal B}^1_{\kappa,-}(M)$] are strictly positive [strictly negative]. 

 The forward and backward Lyapunov spectra of the restriction of the cocycle $G^c_t$ to the 
sub-bundle ${\Cal E}^1_{\kappa}(M)$ of exact currents are reduced to the Lyapunov exponent $0$. 
In fact, since $|dU|_{-1}\leq |U|_0$ and the $L^2_q$ norm is invariant under the Teichm\"uller flow, 
$$\limsup_{t\to \pm\infty} {1\over |t|} \log |G^c_t(dU)|_{-1}\leq \limsup_{t\to\pm\infty} {1\over |t|} 
\log |U|_0 = 0\,\,. \eqnu $$
On the other hand, if the horizontal and vertical foliations ${\Cal F}_{\pm q}$ are ergodic, $dU\not
=0$ implies that the currents $({\rm SU})^{\ast}\eta_T \in I^-_q(M)$ and  $(TU)^{\ast}\eta_S\in  I^+_q(M)$,
which are the projections of $dU$ given by the splitting $(8.38)$ are non-zero. 
Hence, by Lemma 8.6, 
$$\eqalignno{ &\liminf_{t\to +\infty}{1\over t}\log |G^c_t(dU)|_{-1} \geq \liminf_{t\to+\infty}{1\over 
t} |G^c_t\bigl((TU)^{\ast}\eta_S\bigr)|_{-1} \geq 0 \,\,, \cr
\noalign{\vskip-6pt}
&&\speqnu{8.41'}\cr
\noalign{\vskip-6pt}
            & \liminf_{t\to-\infty} {1\over {|t|}}\log |G^c_t(dU)|_{-1} \geq 
\liminf_{t\to-\infty}
{1\over {|t|}} |G^c_t\bigl(({\rm SU})^{\ast}\eta_T)|_{-1}  \geq 0\,\,. \cr
\noalign{\vskip-24pt}} 
$$
\enddemo

\specialnumber{8.7'}
\demo{Remark} There is no exhaustive theory of Lyapunov exponents for cocycles on Hilbert 
or Banach bundles. In particular, the available general results (\ref\Sl and references therein) do 
not seem to apply to the cocycles $G^0_t$ and $G^c_t$. It is possible to deduce from Theorem 8.7 some
additional information concerning the Lyapunov spectrum of $G^0_t$, namely that its forward Lyapunov
spectrum contains the Lyapunov exponents
$$  -1+\lambda_g<\cdots <-1+\lambda_2<0<1-\lambda_2<\cdots <1-\lambda_g\,\,, \eqnu $$
but no Oseledec's decomposition of the bundle $H^{-1}_{\kappa}(M)$ seems to exist. The Lyapunov 
spectrum of $G^c_t$ can be obtained from that of $G^0_t$ as the union of the translate by $+1$
and of the translate by $-1$. 
\enddemo

 \section{The deviation of ergodic averages and open questions}

In this conclusive section of the paper we apply the results we have obtained on the
Kontsevich-Zorich cocycle and on the cocycle $G^c_t$ on the bundle of closed currents of order $1$
to prove results concerning the deviation of the ergodic averages of measured foliations and of
area-preserving smooth vector fields on higher genus surfaces. By Theorem 8.7 
and by Sobolev estimates on the first return orbits as $1$-dimensional currents of order $1$, we derive Sobolev estimates valid for a sequence of `best' return times. Then by a standard technical argument we 
deduce estimates for all times.

\demo{{\rm 9.1.} The $L^2$ mean deviation and the cohomological equation}
 By the Gottschalk-Hedlund theorem, there is a relation between the deviation of ergodic 
averages and the existence of solutions of the cohomological equation. The topological version of
the theorem, which states that the cohomological equation $Xu=f$, for a minimal flow generated 
by a vector field $X$, has a continuous  solution $u$ if and only if the ergodic integrals of $f$ are 
uniformly bounded, was first proved in [22, Th.\ 14.11].  A proof can also be found in [32, Th.\ 2.9.4].  
A similar theorem can be proved for $L^2$ solutions, although we were not able to find a reference.  
By the $L^2$ version of the Gottschalk-Hedlund theorem, Theorem 8.3 yields a negative answer to 
a question, left open in \ref\Ftwo, concerning the existence of $L^2$ solutions of the cohomological equation. We therefore prove that the statement of [17, Th.\ A] is incorrect and should be replaced 
by [18, Th.\ A].    
\enddemo

\proclaim{Theorem} For $\mu^{(1)}_{\kappa}$\/{\rm -}\/almost all quadratic differentials $q\in {\Cal M}
^{(1)}_{\kappa}$ there exists {\rm (}\/at least\/{\rm )} a function $f\in C^{\infty}_0(M\setminus\Sigma_q)${\rm ,}
with  zero average{\rm ,} such that the cohomological equation $Su=f$ has no solution $u\in L^2_q(M)${\rm .} Let 
$\Phi_q^{\Cal T}$ be the {\rm (}\/almost everywhere defined\/{\rm )} flow of the vector field $S${\rm ,} by the
Gottschalk\/{\rm -}\/Hedlund theorem{\rm ,}
$$\limsup_{{\Cal T}\to +\infty}\,\,|\int_0^{\Cal T} f(\Phi_q(p,\tau))\,d\tau\,|_0\,=\, +\infty\,\,. 
\eqnu$$   
\endproclaim

\demo{Proof} Assume that the cohomological equation $Su=f$ has a solution $u\in L^2_q(M)$ for any 
$f\in C^{\infty}_0(M\setminus\Sigma_q)$. Then by the methods of Section 7.1 it would be possible to 
construct a $2g-1$ dimensional subspace of basic currents of order $1$ for the horizontal foliation 
${\Cal F}_q$. In fact, given any meromorphic function $m^+\in {\Cal M}^+_q$ with zero average, the
cohomological equation $Su=-\Re(m^+)$ would have a solution $u\in L^2_q(M)$. The solution could be 
constructed as follows. Let $u_0$ be a smooth local solution of the equation $Su_0=-\Re(m^+)$ in
a neighbourhood of the zero set $\Sigma_q$, which can be obtained by Laurent expansion at the
points of $\Sigma_q$ and by the formulas $(7.4)$. Let then $u_1\in L^2_q(M)$ be a solution of the 
equation $Su_1=f_1$, where $f_1:=-\Re(m^+)-Su_0\in C_0^{\infty} (M\setminus\Sigma_q)$. The function
$u:=u_0+u_1\in L^2_q(M)$ is therefore a solution of the cohomological
equation $Su=-\Re(m^+)$. By
formula $(7.3)$, we could then construct a space of basic currents of order $1$ for ${\Cal F}_q$ 
of dimension $2g-1$. In fact, we would obtain that the basic cohomology $H^{1,1}_q(M,{\Bbb R})$ has
codimension $1$ in $H^1(M,{\Bbb R})$, hence dimension $2g-1$. On the other hand, by Theorem 8.3,
for almost all $q\in {\Cal M}^{(1)}_{\kappa}$, $H^{1,1}_q(M,{\Bbb R})$ coincides with the unstable
subspace of the Kontsevich-Zorich cocycle, hence it has dimension (at most) equal to the genus 
$g\geq 2$. Since $2g-1>g\geq 2$, we have obtained a contradiction. Hence, we have proved that
there exists $f\in C_0^{\infty}(M\setminus\Sigma_q)$ with zero average such that the cohomological
equation $Su=f$ has no solution $u\in L^2_q(M)$.

By the Gottschalk-Hedlund theorem, if $(9.1)$ were false, the cohomological equation $Su=f$ 
would have a solution $u\in L^2_q(M)$. In fact, in that case the one-parameter family
$$u_{\Cal T}(p):= -{1\over {\Cal T}} \int_0^{\Cal T} \int_0^{\tau} f(\Phi_q(p,s))\,ds\,d\tau\,\, 
\eqnu $$
would be bounded in $L^2_q(M)$. In addition, a computation shows that, if ${\Cal F}_q$ is ergodic, 
$$Su_{\Cal T}(p)= f(p)-{1\over {\Cal T}}\int_0^{\Cal T} f(\Phi_q(p,\tau))\,d\tau \,\to\,f \,
\hbox{ in }\,L^2_q(M)\,\,. \speqnu{9.2'}$$
Any weak limit $u\in L^2_q(M)$ of the family $(9.2)$ would therefore be a solution of the equation 
$Su=f$. Since we proved above that such an equation has no solution $u\in L^2_q(M)$, it follows
that $(9.1)$ holds. \enddemo 

The non-existence of $L^2_q(M)$ solutions for a cohomological equation\break $Su=f$, with $f$ 
smooth and of zero average is therefore related to the existence of {\it unbounded deviation }in the 
ergodic integrals of $f$. In the case of the torus $T^2$ ($g=1$) this phenomenon does not occur for
a full measure set of flows (given by a Diophantine condition). In the case of higher genus, as 
it was conjectured in \ref\KZone, the deviation of the ergodic integrals from the linear growth 
behaviour predicted by the ergodic theorem obeys, for a typical function, a {\it power law }with 
exponent $\lambda_2>0$. 

\demo{{\rm 9.2.} Sobolev estimates for the {\rm (}\/first\/{\rm )} return orbits}
 The key idea of our argument consists in studying the dynamics of the flow $G^c_t$ on 
the infinite-dimensional non-closed sub-bundles $\Gamma^{\pm}_{\kappa}\subset {\Cal H}^{-1}_{\kappa}
(M)$ generated by segments of leaves of the horizontal and vertical foliations. Let $\Cal T>0$. We 
denote by $\gamma_{\pm q}^{\Cal T}$ a positively oriented segment of length $\Cal T>0$ of a leaf of 
the measured foliation ${\Cal F}_{\pm q}$ respectively. By trace theorems for Sobolev spaces, the 
vector spaces $\Gamma^{\pm}_q$ generated by the segments $\gamma_{\pm q}^{\Cal T}$ are subspaces of 
the space ${\Cal H}^{-1}_q(M)$ of\break $1$-dimensional currents of order $1$. In addition, they are 
invariant under the action of the cocycle $G^c_t$. Let $\Phi_{\pm q}^{\Cal T}$, ${\Cal T}\in{\Bbb R}$,
 be the (almost everywhere defined) flows of the vector fields $S$, $T$ respectively. The foliations 
${\Cal F}_{\pm q}$ are almost everywhere the orbit foliations of the flows $\Phi_{\pm q}^{\Cal T}$. 
Let $\alpha:=f^+\eta_T +f^-\eta_S \in {\Cal H}^1_q(M)$, then
$$\gamma_{\pm q}^{\Cal T}(\alpha)=\int_0^{\Cal T}f^{\pm}(\Phi_{\pm q}(p^{\pm},\tau))\,d\tau\,\,,
\eqnu $$
where $p^{\pm}\in M$ are the starting points of the oriented segments $\gamma_{\pm q}^{\Cal T}$. The 
ergodic averages of the functions $f^{\pm}\in H^1_q(M)$ can therefore be understood by studying the 
dynamics of the `renormalization' cocycle $G^c_t$ on $\Gamma^{\pm}_{\kappa}$.

The first step consists in an estimate of the norm of $\gamma_{\pm q}^{\Cal T}$ in 
${\Cal H}^{-1}_q(M)$. Let $|\!|q|\!|$ denote the $R_q$-length of the shortest saddle connection of 
the quadratic differential $q\in {\Cal M}^{(1)}_{\kappa}$.
\enddemo

\proclaim{Lemma} There exists a constant $K>0$ such that{\rm ,} for all quadratic differentials
$q\in {\Cal M}^{(1)}_{\kappa}${\rm ,}
$$ |\gamma_{\pm q}^{\Cal T}|_{{\Cal H}^{-1}_q(M)} \,\leq\, K\left(1+ {{\Cal T}\over {|\!|q|\!|}}\right)\,\,. 
\eqnu $$
\endproclaim

\demo{Proof} Let $d_q(M)$ be the maximum distance of any point in $M$ to a zero. By 
[45, Cor. 5.6] there is a constant $K_1>1$ such that $d_q(M)\leq K_1/|\!|q|\!|$ or $d_q(M)
\leq \sqrt{2/\pi}$. 
Since $|\!|q|\!|/2\leq d_q(M)$, it follows that there is $K_2>1$ such that $|\!|q|\!| \leq 
K_2$. Let $\delta:=|\!|q|\!|/3$. Let $\gamma_+:=\gamma_{+q}^{\Cal T}$ be a {\it regular }segment of 
length $\Cal T>0$ of a leaf of the horizontal foliation ${\Cal F}_q$. The case of a regular segment
$\gamma_-:=\gamma_{-q}^{\Cal T}$ of a leaf of the vertical foliation is similar. For all but a finite
number of points $p\in\gamma_+$, the vertical segment $\gamma_-^{\delta}(p)$ of length $2\delta$, 
centered at $p$, is well-defined. At the exceptional points $p_1,\ldots ,p_N\in \gamma_+$ (ordered in 
the positive direction along $\gamma_+$) the vertical leaf through $p_i$ meets a singularity at a 
distance from $p_i$ less than $\delta$. Let $\gamma_+(p_i,p_{i+1})$ be the open horizontal segment 
with endpoints $p_i$, $p_{i+1}$. By the definition of $\delta>0$, the $R_q$-length $L_q(\gamma_+
(p_i,p_{i+1}))\geq \delta$. Otherwise, there would be two singularities at a distance strictly 
less than $3\delta=|\!|q|\!|$, in contradiction with the definition of $|\!|q|\!|$. By possibly 
introducing additional points, we can assume that the following holds:
$$\eqalign{\delta\leq & L_q\bigl(\gamma_+(p_i,p_{i+1})\bigr)\leq 2\delta\,\,,\,\,\,\,1\leq i 
                                                                                  \leq N-1\,\,, \cr 
          &L_q\bigl(\gamma_+(p_0,p_1)\bigr)\,,\,\,L_q\bigl(\gamma_+(p_N,p_{N+1})\bigr)\leq 
                                                                   2\delta\,\,, \cr} \eqnu $$
where $p_0$, $p_{N+1}$ are the endpoints (ordered in the positive direction) of $\gamma_+$. Let $R_i
\subset {\Bbb R}^2$, $0\leq i\leq N$, be the open rectangles defined by 
$$R_i:=\{(x,y)\in {\Bbb R}^2\,|\, 0<x<L_q\bigl(\gamma_+(p_i,p_{i+1})\bigr)\,,\,\, -\delta\leq y\leq 
\delta\}\,\,.  \eqnu $$ 
We claim that there are isometric embeddings of the rectangles $R_i$ into $(M,R_q)$ such that the 
image of the horizontal segments $\gamma_i=R_i\cap \{y=0\}$ coincide with $\gamma_+(p_i,p_{i+1})$. 
In fact, let $p\not=p'\in \gamma_+(p_i,p_{i+1})$. If $\gamma_-^{\delta}(p)\cap\gamma_-^{\delta}(p')
\not=\emptyset$, then, by the Pythagorean theorem, there would be a homotopically non-trivial regular
closed geodesic of length at most $\sqrt{4\delta^2 +4\delta^2}<3\delta=|\!|q|\!|$. Hence, there would
be an embedded cylinder with waist of length $<|\!|q|\!|$ [58, Chap.\ IV, \S 9.3]. The boundary 
of such a cylinder should then contain at least one saddle connection of smaller or equal length, in 
contradiction with the definition of $|\!|q|\!|$.

Let $R_{a,b}:=\{(x,y)\in {\Bbb R}^2\,|\,0<x<a\,,\,\,-b<y<b$. By a rescaling argument, that 
is, by reducing to the case of $a=b=1$ by an affine change of coordinate and by the Sobolev trace 
theorem [1, Th.\ 5.4 (5)], there exists a constant $K_3>0$, such that
$$\left|\int_0^a f(x,0)dx\right|\,\leq K_3\,\left({a\over b}\right)^{1/2} \max\{a,b,1\}\,|f|_{ H^1(R_{a,b})}\,\,.
\eqnu $$
By applying the inequality $(9.7)$ to each of the rectangles $R_i\subset M_q$, $0\leq i\leq N$, 
by $(9.3)$ we get 
$$\left|\gamma_i(\alpha)\right|=\left|\int_{\gamma_i}f^+\eta_T\right| \leq K_4 |f^+|_1\leq K_4|\alpha|_1\,\,, 
\speqnu{9.7'}$$
where $K_4:=2^{1/2}K_2 K_3$. The estimate $(9.4)$ is then a consequence of $(9.7')$ and of the 
inequality $N-1\leq {\Cal T}/\delta$, which follows from the first inequality of $(9.5)$. 
\enddemo 

 A point $p\in M$ is {\it regular }with respect to a measured foliation $\Cal F$ if it does
not belong to a singular leaf of $\Cal F$. Let $q\in {\Cal M}^{(1)}_{\kappa}$. A point $p\in M$ will 
be said to be $q$-{\it regular }if it is regular with respect to the horizontal and vertical foliations
${\Cal F}_{\pm q}$. The set of $q$-regular points is of full measure and it is equivariant under the
action of the mapping class group and of the Teichm\"uller flow. Let $p\in M$ be $q$-regular and let
$I_{\pm q}(p)$ be the vertical [horizontal] segment of length $|\!|q|\!|/2$ centered at $p$. A {\it
forward horizontal {\rm [}\/vertical\/{\rm ]} return time }of the point $p\in M$ is defined to be any real
number
${\Cal T}_{\pm q}(p)>0$ such that $\Phi_{\pm q}\bigl(p,{\Cal T}_{\pm q}(p)\bigr)\in I_{\mp q}(p)$. 
If ${\Cal T}>0$ is any horizontal [vertical] return time of a $q$-regular point $p\in M$, the 
horizontal [vertical] forward segment $\gamma_{\pm q}^{\Cal T}(p)$, with initial point $p$, will 
be called a {\it forward horizontal {\rm [}\/vertical\/{\rm ]} return trajectory }at $p$. There is a natural map 
from the set of horizontal [vertical] return trajectories into the  set of homotopically non-trivial 
closed curves. The {\it closing }of any horizontal [vertical] trajectory segment $\gamma_{\pm q}
^{\Cal T}(p)$ is defined by
$$\widehat{\gamma}_{\pm q}^{\Cal T}(p):=\gamma_{\pm q}^{\Cal T}(p)\cup\gamma\bigl(p,\Phi_{\pm q}
(p,{\Cal T})\bigr)\,\,, \eqnu $$
where $\gamma\bigl(p,\Phi_{\pm q}(p,{\Cal T})\bigr)$ is the shortest geodesic segment joining the 
endpoints $p$, $\Phi_{\pm q}(p,{\Cal T})$ of the trajectory segment (which exists by [58, Th.\ 18.2.2]).   

Let ${\Cal T}_{\pm q}^{(1)}(p)$ be the {\it forward horizontal {\rm [}\/vertical\/{\rm ]} first return 
time }of the $q$-regular point $p\in M$, defined to be the real number
$${\Cal T}_{\pm q}^{(1)}(p):=\min\{{\Cal T}>0\,|\, \Phi_{\pm q}(p,{\Cal T})\in I_{\mp q}(p)\} \,\,. 
\eqnu $$
The corresponding forward horizontal [vertical] trajectory $\gamma_{\pm q}^{(1)}(p)$ with initial 
point $p$ will be called the {\it forward horizontal [vertical] first return trajectory }at~$p$. We 
have the following estimate for first return times:

\specialnumber{9.2'}
\proclaim{Lemma} There exists a measurable function $K_r:{\Cal M}_{\kappa}\to {\Bbb R}^+$ such 
that{\rm ,} if\/ ${\Cal T}_{\pm q}^{(1)}(p)$ is the forward horizontal {\rm [}\/vertical\/{\rm ]} first return
time of a\break
$q$\/{\rm -}\/regular point $p\in M${\rm ,} then
$$ |\!|q|\!|/2\leq {\Cal T}_{\pm q}^{(1)}(p) \leq K_r(q)\,\,. \eqnu $$
\endproclaim

\pagebreak

\demo{Proof} Since $\widehat{\gamma}_{\pm q}^{\Cal T}(p)$ is a closed curve not homotopic to zero, 
there exists a shortest geodesic  in its homotopy class [58, Th.\ 18.4.]. If such a geodesic  
contains singular points, then its length is $\geq |\!|q|\!|$. Otherwise, it is the waist curve of 
an embedded flat cylinder. Since the boundary of any flat cylinder is the union of saddle connections,
the lower bound $|\!|q|\!|$ is valid also in this case. Since $I_{\pm q}(p)$ has length $|\!|q|\!|/2$,
 ${\Cal T}_{\pm q}^{(1)}(p)\geq |\!|q|\!|/2$.

For all $q\in {\Cal M}_{\kappa}$, the functions $p\to {\Cal T}_{\pm q}^{(1)}(p)$ are
almost everywhere defined and simple, in the sense that they are locally constant and take only a 
finite set of values. In fact, if $I^{\mp}$ is an interval transverse to ${\Cal F}_{\pm q}$ of
length $\delta>0$, the Poincar\'e return map on $I^{\mp}$ is an orientation-preserving isometry 
with (at most) a finite number of discontinuities, hence it is an interval exchange transformation
with (at most) a finite number of sub-intervals. On each sub-interval ${\Cal T}_{\pm q}^{(1)}(p)$
is constant. The function 
$$K_r(q):= \max \{{\Cal T}_{\pm q}^{(1)}(p)\,|\, p\in M\,\hbox{ is }\,q\hbox{-regular}\,\} \,\,
\speqnu{9.10'}$$
is therefore well-defined and everywhere finite. It can be proved that $K_r$ is lower semicontinuous,
hence measurable, and $(9.10)$ holds by the definition $(9.10')$. 
\enddemo 

9.3. {\it Special sequences of {\rm `}\/close\/{\rm '} return times}.
 Special sequences of return times for the horizontal [vertical] foliation of a generic
quadratic differential can be constructed by considering return times of the Teichm\"uller geodesic 
flow. This will clarify in what sense the Teichm\"uller flow $G_t$ and the cocycles $G^{KZ}_t$, 
$G^c_t$ play the role of a `renormalization dynamics' for orientable measured foliations. Let $q\in 
{\Cal M}^{(1)}_{\kappa}$ be a Birkhoff generic point of the Teichm\"uller flow $G_t$ and let 
${\Cal S}_{\kappa}(q)\subset {\Cal M}^{(1)}_{\kappa}$ be a smooth compact hypersurface of codimension
$1$, containing $q$ and transverse to the Teichm\"uller flow. Let $(t_k)$ be the sequence of return 
times of the orbit $\{G_t(q)\,|\,t\in {\Bbb R}\}$ to ${\Cal S}_{\kappa}(q)$. Since, by definition,  
$$({\Cal F}_{q_t},{\Cal F}_{-q_t})=(e^{-t}{\Cal F}_q, e^t {\Cal F}_{-q}) \,\,, \eqnu $$
if $t=t_k<0$ is a {\it backward }return time of $G_t(q)$, any forward first return trajectory of 
the horizontal foliation ${\Cal F}_{q_t}$ is a forward return trajectory of the horizontal foliation 
${\Cal F}_q$, provided $|t_k|$ is sufficiently large. In a similar way, if $t=t_k>0$ is a {\it 
forward }return time of $G_t(q)$, any forward first return trajectory of the vertical foliation 
${\Cal F}_{-q_t}$ is a forward return trajectory of the vertical foliation ${\Cal F}_{-q}$. In all 
the arguments which follow we will consider the case of the horizontal foliation, the case of the 
vertical foliation being similar.   

 By {\it closing }the return trajectories of the horizontal [vertical] foliation, as in 
$(9.8)$, we obtain {\it closed }currents of order $1$. The evolution of such currents under the 
action of the Teichm\"uller flow is therefore   described by the cocycle $G^c_t$ on the bundle 
${\Cal Z}_{\kappa}^1(M)$, studied in Section 8.3 (in particular Theorem 8.7). We 
describe below a relevant additional property of the Oseledec's splitting~$(8.40)$. 

In Theorem 8.7 we have proved that the Lyapunov spectrum of the 
cocycle $G^c_t$ on the 
bundle ${\Cal E}^1_{\kappa}(M)$ of exact currents of order $1$ is reduced to the single Lyapunov 
exponent $0$. In fact, there exists a {\it Lyapunov norm }on ${\Cal E}^1_{\kappa}(M)$, invariant 
under $G^c_t$. Let $C$ be an exact current of order $1$, then $C=dU$, where the function $U\in 
L^2_q(M)$ has zero average.

\proclaim{Lemma} The norm $|\cdot|_q${\rm ,} defined fiber\/{\rm -}\/wise on the bundle 
${\Cal E}^1_{\kappa}(M)$ of exact currents of order $1$ as
$$|dU|_q:=|U|_{L^2_q(M)}\,\,,\,\,\,\,\int_M U\,\omega_q=0\,\, , \eqnu $$
is invariant under $G^c_t${\rm .} There exists a continuous function $K:{\Cal M}^{(1)}_{\kappa}\to
{\Bbb R}^+$ such that{\rm ,} for all $q\in{\Cal M}^{(1)}_{\kappa}$ and $dU\in
{\Cal E}^1_q(M)${\rm ,} we have\/{\rm :}
$$ K(q)\,|dU|_q \leq |dU|_{{\Cal H}^{-1}_q(M)} \leq |dU|_q \,\,. \speqnu{9.12'}$$
\endproclaim

\demo{Proof} The norm $|dU|_q$ is invariant under the flow $G^c_t$, since $G^c_t(dU)=dU$ and the 
inner product of the Hilbert space $L^2_q(M)$ is invariant under the action of the Teichm\"uller 
flow $G_t$ on ${\Cal M}_{\kappa}$. The second inequality in $(9.12')$ follows immediately from the
definition $(9.12)$.

Let $q\in {\Cal M}_{\kappa}$ be a quadratic differential. By [18, Prop. 4.3A], the 
following splitting holds:
$$H^0_q:=\{U\in L^2_q(M)\,|\,\int_M U\omega_q=0\}= R(\partial^+_q)\oplus R(\partial^-_q)\,\,, 
\eqnu $$
where $R(\partial^{\pm}_q)\subset L^2_q(M)$ is the range of the Cauchy-Riemann operator $\partial_q
^{\pm}$ on the weighted Sobolev space $H^1_q(M)$. By [18, Prop. 3.2], the subspace $R(\partial
^+_q)$ [$R(\partial^-_q)$] is the orthogonal complement in $L^2_q(M)$ of the finite-dimensional 
subspace ${\Cal M}_q^-$ [${\Cal M}_q^+$] of anti-meromorphic [meromorphic] functions. Since the
Cauchy-Riemann operators $\partial^{\pm}_q$ depend continuously on the quadratic differential $q\in
{\Cal M}_\kappa$, there  exists a continuous function $K_1:{\Cal M}_{\kappa}\to {\Bbb R}^+$ such that the
following holds. If 
$U\in H^0_q$, then there exist functions $v^{\pm}\in H^1_q(M)$ such that 
$$ \eqalign{ U& = \partial^+_q v^+  + \partial^-_q v^-\,\,,\cr
                |v^{\pm}|_1&\leq K_1(q)\, |U|_{L^2_q(M)} \,\,. \cr}\eqnu $$  
Hence
$$\eqalignno{  |U|^2_{0}=&(U,\partial^+v^+ + \partial^-v^-)_q   \cr
\noalign{\vskip-6pt}
&&\speqnu{9.14'}\cr
\noalign{\vskip-12pt}
             \leq &\,2|dU|_{-1}(|v^+|_1 + |v^-|_1) \leq 4K_1(q)|U|_0|dU|_{-1}\,\,.\cr
\noalign{\vskip-18pt}} 
$$
\enddemo 

\bigbreak   Let $q\in{\Cal M}^{(1)}_{\kappa}$ be a regular point for the Kontsevich-Zorich cocycle (in 
the sense of the Oseledec's theorem). Let 
$$\lambda'_1=1>\lambda'_2>\cdots \lambda'_s>0>-\lambda'_s>\cdots >-\lambda'_1=-1 \eqnu $$
be the {\it distinct }Lyapunov exponents at $q$ of the Kontsevich-Zorich cocycle. By Corollary 2.2 
and Theorem 8.5, there exists $s\in \{2,\ldots ,g\}$ such that $(9.15)$ holds. It has been conjectured
in \ref\KZone, on the basis on numerical evidence, that the Lyapunov spectrum of the Kontsevich-Zorich
cocycle is simple, hence $s=g$. By Theorem 8.7, there exists an Oseledec's splitting
$${\Cal Z}^1_q(M)=E^+_1(q)\oplus\cdots \oplus E^+_s(q)\,\bigoplus\, E^-_1(q)\oplus\cdots \oplus E^-_s(q)\,
\bigoplus\, {\Cal E}^1_q(M)\,\,,\speqnu{9.15'}$$
where $E^{\pm}_i(q)$, $i\in \{1,\ldots s\}$, is the eigenspace of the cocycle $G^c_t$ corresponding to 
the Lyapunov exponent $\pm\lambda'_i$ and ${\Cal E}^1_q(M)$ is the closed infinite-dimensional 
subspace of exact currents of degree $1$. The subspaces $E^{\pm}_i(q)$ are finite dimensional with
dimension equal to the {\it multiplicity }of the Lyapunov exponent $\lambda'_i$, while the subspace
${\Cal E}^1_q(M)$ is infinite-dimensional. 

Let $\Pi^{\pm i}_q:{\Cal Z}^1_q(M)\to E^{\pm}_i(q)$, $i\in \{1,\ldots ,s\}$, $\Pi^{\Cal E}_q:{\Cal Z}^1_q(M)\to
{\Cal E}^1_q(M)$ be the projections determined by the splitting $(9.15')$.  Such projections can be extended to
the (Hilbert) space ${\Cal H}^{-1}_q(M)$ by composition with the  orthogonal projection onto the closed
subspace ${\Cal Z}^1_q(M)$. Let $\delta_{\kappa}:{\Cal M} ^{(1)}_{\kappa}\to{\Bbb R}^+$ be the {\it distorsion
}(in the sense of [41, Chap.\ IV, \S 11, p. 269])  of the splitting $(9.15')$ with respect to the ${\Cal H}^{-1}_q$
norm on ${\Cal B}^1_q(M)\oplus {\Cal M}^1_{-q}(M)$ and the  Lyapunov norm $|\cdot|_q$ (introduced in Lemma
9.3) on
${\Cal E}^1_q(M)$:
$$\delta_{\kappa}(q):=\sup_C  {{\sum_{i=1}^g |\Pi^{+i}_q(C)|_{-1}^2\,+\, 
\sum_{i=1}^g |\Pi^{-i}_q(C)|_{-1}^2\,+\, |\Pi^{\Cal E}_q(C)|_q^2} \over {|C|_{-1}^2}}\,\,,
\speqnu{9.15''}$$
where the supremum is taken over all currents $C\in {\Cal Z}^1_q(M)$.

Let ${\Cal P}^{(1)}_{\kappa}\subset{\Cal M}^{(1)}_{\kappa}$ be a compact set satisfying the following
conditions:
\medbreak
\item{(1)} All $q\in {\Cal P}^{(1)}_{\kappa}$ are Birkhoff generic points for the Teichm\"uller flow $G_t$ and 
Oseledec regular points for the cocycle $G^c_t$;
\medbreak
\item{(2)} ${\Cal P}^{(1)}_{\kappa}$ is transverse to the Teichm\"uller flow and has positive 
transverse measure;
\medbreak
\item{(3)} the distorsion $\delta_{\kappa}$, defined by $(9.15'')$, is continuous (hence bounded) on 
${\Cal P}^{(1)}_{\kappa}$;
\medbreak
\item{(4)} the function $K_r:{\Cal M}^{(1)}_{\kappa}\to{\Bbb R}^+$, that is the upper bound $(9.10)$
on the first return times of the horizontal and vertical foliations, is continuous (hence bounded) on 
${\Cal P}^{(1)}_{\kappa}$.
\medbreak
 It follows from the ergodicity of the Teichm\"uller flow (Theorem 1.1), from Theorem 8.7, 
Lemma $9.2'$ and Luzin's theorem that the union of all sets ${\Cal P}^{(1)}_{\kappa}$ with the 
properties 
 (1)--(4)  is a full measure subset of ${\Cal M}^{(1)}_{\kappa}$.

Let $q\in {\Cal P}^{(1)}_{\kappa}$ and $(t_k)_{k\in {\Bbb N}}$ be the sequence of return 
times of the {\it backward }orbit $\{G_t(q)\,|\,t\leq t_1=0\}$ to ${\Cal P}^{(1)}_{\kappa}$. Let $p
\in M$ be any $q$-regular point. A {\it principal sequence }of forward return times at $p$ for the 
horizontal foliation ${\Cal F}_q$ is the sequence
$${\Cal T}^{(k)}_q(p):= {\Cal T}^{(1)}_{G_t(q)}(p) \exp{|t|}\,\,,\,\,\,\,t=t_k \,\,. \eqnu $$  
The horizontal forward  trajectory at $p$ corresponding to a principal return time ${\Cal T}
^{(k)}_q(p)$ will be called a {\it horizontal principal return trajectory }at $p$ and denoted by 
$\gamma_q^{(k)}(p)$. We remark that a horizontal principal return trajectory $\gamma_q^{(k)}(p)$ 
coincides with the horizontal {\it first }return trajectory at $p$ of the quadratic differential 
$G_t(q)$, $t=t_k<0$. A similar construction holds for the vertical foliation ${\Cal F}_{-q}$ by 
considering forward return times of the Teichm\"uller flow.

\proclaim {Lemma} Under the conditions  {\rm (1)--(4)}  there exists a constant $K_{\Cal P}>0$ such 
that the following holds{\rm .} Let $q\in {\Cal P}^{(1)}_{\kappa}$ and $\gamma^{\Cal T}_q(p)${\rm ,}
 ${\Cal T}>0${\rm ,} 
be a forward   trajectory with initial point a $q$\/{\rm -}\/regular point $p\in M${\rm .}
 There exists a finite 
sequence of points $(p^{(k)}_j)\subset\gamma^{\Cal T}_q(p)${\rm ,} $1\leq k\leq n${\rm ,} $1\leq j\leq m_k${\rm ,} 
such that the principal return trajectories $\gamma_q^{(k)}(p^{(k)}_j)\subset \gamma^{\Cal T}_q(p)$ 
do not overlap and{\rm ,} in addition{\rm ,}
$$\eqalign{ & \gamma^{\Cal T}_q(p)=\sum_{k=1}^n\sum_{j=1}^{m_k}\gamma_q^{(k)}(p^{(k)}_j)\,\,
              +\,\,b^{\Cal T}_q(p)\,\,,\cr 
            & m_k\leq K_{\Cal P}\exp(|t_{k+1}|-|t_k|)\,\,, \cr \cr
            & L_q\bigl(b^{\Cal T}_q(p)\bigr)\leq K_{\Cal P}\,\,.\cr} 
\speqnu {9.16'}$$
\endproclaim

\demo{Proof} The proof is based on the following estimate on principal return times. By $(9.16)$, 
Lemma $9.2'$ and condition $(4)$, there exists a constant $K_{pr}>0$ such that, for all $q\in {\Cal P}^{(1)}
_{\kappa}$, all $q$-regular points $p\in M$ and all $k\in {\Bbb N}$,
$$K_{pr}^{-1}\exp{|t_k|}\leq {\Cal T}^{(k)}_q(p) \leq K_{pr}\exp{|t_k|}\,\,.\eqnu $$
Let $n=\max\{k\in {\Bbb N}\,|\, {\Cal T}^{(k)}_q(p)\leq {\Cal T}\,\}$. The maximum exists (i.e., is finite) by
$(9.17)$. Let $p^{(n)}_1:=p$. The sequence $(p^{(k)}_j)$ with the properties stated in $(9.16')$ can
be constructed by a finite iteration of the following procedure. Let $p^{(k)}_j$ be the last point 
in the sequence already determined and let 
$$p^{(k)}_{j+}:=\Phi_q\bigl(p^{(k)}_j,{\Cal T}_q^{(k)}(p^{(k)}_j)\bigr)\in \gamma^{\Cal T}_q(p)\,\,.
\eqnu $$  
Let then $k'\in \{1,\ldots ,k\}$ be the largest integer such that 
$$\Phi_q\bigl(p^{(k)}_{j+},{\Cal T}_q^{(k')}(p^{(k)}_{j+})\bigr)\in 
\gamma^{\Cal T}_q(p)\,\,. \speqnu{9.18'}$$ 
If $k'=k$, let $p^{(k)}_{j+1}:=p^{(k)}_{j+}$. If $k'<k$, let $m_k:=j$, $m_h=0$ (no points) for all 
$k'<h<k$ and $p^{(k')}_1:=p^{(k)}_{j+}$. The iteration step is concluded. By $(9.17)$ we have:
$$K_{pr}^{-1}\exp{|t_k|}\,m_k \leq {\Cal T}^{(k+1)}_q(p^{(k)}_1)\leq  K_{pr}\exp{|t_{k+1}|} \,\,. 
\eqnu $$
The length of the remainder $b^{\Cal T}_q(p)$ has to be less than any upper bound for the length of
first return times. Hence, by $(9.19)$, Lemma $9.2'$ and condition $(4)$, the estimates in $(9.16')$ are proved 
and the argument is concluded.
\enddemo 

9.4. {\it The main theorem on the deviation of ergodic averages}.
 
\proclaim{Theorem} For almost all quadratic differentials $q\in {\Cal M}^{(1)}_{\kappa}$ and 
almost all points $p\in M_q$ the following holds{\rm .} Let $\gamma_q^{\Cal T}(p)$ be the horizontal 
forward trajectory segment of length ${\Cal T}>0$ with initial point $p${\rm .} Then{\rm ,} for all $i\in
\{1,\ldots s\}${\rm ,} 
$$\eqalign{&\limsup_{{\Cal T}\to +\infty}{{\log |\Pi^{+ i}_q(\gamma_q^{\Cal T}(p))|_{-1}}\over 
              {\log {\Cal T}}}\,=\, \lambda'_i\,\,, \cr
           &\limsup_{{\Cal T}\to +\infty}{{\log |\Pi^{-i}_q(\gamma_q^{\Cal T}(p))|_{-1}}\over
              {\log {\Cal T}}}\,=\, 0\,\,, \cr
           &\limsup_{{\Cal T}\to +\infty}{{\log |\Pi^{\Cal E}_q(\gamma_q^{\Cal T}(p))|_{-1}} \over
              {\log {\Cal T}}}\,=\,0\,\,.\cr}\eqnu $$
\endproclaim            
\demo{Proof} The argument   consists of two parts. We   first prove the upper bounds implicit 
in $(9.20)$, then the lower bound.
\medbreak
 {\it Upper bound}. By Lemma 9.4 the estimate can be reduced to the case of a principal 
return trajectory $\gamma^{(k)}_q(p)$. Let ${\Cal P}^{(1)}_{\kappa}$ be a compact set satisfying 
conditions  (1)--(4)  listed in Section 9.3. By Lemmas 9.2, 9.2$'$ and the conditions $(3)$, $(4)$, there 
exists a constant $K_1>0$ such that the following holds. Let $q\in {\Cal P}^{(1)}_{\kappa}$ and 
$p\in M$ be a $q$-regular point. Let $\gamma^{(1)}_q(p)$ be the first return horizontal trajectory 
at $p$. Then
$$\eqalign{ |\Pi^{\pm i}_q\bigl({\widehat\gamma}^{(1)}_q(p)\bigr)|_{-1} &\leq K_1\,\,, \cr
            |\Pi^{\Cal E}_q\bigl({\widehat\gamma}^{(1)}_q(p)\bigr)|_q &\leq K_1\,\,, \cr}
\eqnu $$
where the closing operation $\gamma\to {\widehat\gamma}$ of a return trajectory is described 
by $(9.8)$
and $|\cdot|_q$ is the Lyapunov norm on ${\Cal E}^1_q(M)$ defined in Lemma 9.3. The currents
of order $1$ given by the principal return trajectories can be obtained by the action of the cocycle
$G^c_t$ on the first return trajectories:
$$\gamma^{(k)}_q(p)=G^c_{-t}\bigl(\gamma^{(1)}_{G_t(q)}(p)\bigr)\,\,,\,\,\, t=t_k\leq 0\,\,, 
\eqnu $$
where $(t_k)$ is the sequence of backward return times of the Teichm\"uller orbit $G_t(q)$   ${\Cal
P}^{(1)}_{\kappa}$ . Hence, by Theorem 8.7, Lemma 9.3 and the invariance of the Lyapunov  norm $|\cdot|_q$
under the cocycle $G^c_t$,  if $\lambda^+_i>
\lambda'_i>0$ and $-\lambda'_i<\lambda^-_i<0$, there exists a constant $K_2>0$, such that
$$\eqalign{ &|\Pi^{\pm i}_q\bigl({\widehat\gamma}^{(k)}_q(p)\bigr)|_{-1} \leq K_2\exp(\lambda_i^{\pm}
|t_k|)\,\,,\cr
            &|\Pi^{\Cal E}_q\bigl({\widehat\gamma}^{(k)}_q(p)\bigr)|_{-1} \leq K_2\,\,. \cr} 
\speqnu{9.22'}$$
 By 
Lemma~9.4 the   trajectory $\gamma_q^{{\Cal T} }(p)$ can be split as a union of {\it principal 
return trajectories }and a {\it uniformly bounded }remainder. Hence, by Lemma 9.2, $(9.11)$ and 
$(9.22')$, there exists a constant $K_3>0$
such that
$$\eqalign{&|\Pi^{\pm i}_q\bigl(\gamma_q^{\Cal T}(p)\bigr)|_{-1}\leq K_3 \sum_{1\leq k\leq n} m_k 
          \exp(\lambda_i^{\pm} |t_k|)\,\,, \cr 
    &|\Pi^{\Cal E}_q\bigl(\gamma_q^{\Cal T}(p)\bigr)|_{-1}\leq K_3\sum_{1\leq k\leq n} m_k\,\,.\cr}
\speqnu{9.22''}$$
Since the transverse set ${\Cal P}^{(1)}_{\kappa}$ has positive transverse measure, by the ergodicity
of the Teichm\"uller flow (Theorem 1.1) and Birkhoff ergodic theorem, the backward return times 
$(t_k)$ have the following property:  
$$ \lim_{k\to +\infty} {{|t_k|}\over k}\,=\, {1\over{\mu}}\,\, ,\eqnu $$
where $\mu:=\mu({\Cal P}^{(1)}_{\kappa})>0$ is the transverse measure of ${\Cal P}^{(1)}_{\kappa}$. 
Let $0<\delta\ll 1/\mu$. By $(9.23)$ there exists $N_{\delta}\in {\Bbb N}$ such that, for all $k\geq 
N_{\delta}$,  
$$ ({1\over {\mu}} -\delta)k \leq |t_k| \leq ({1\over {\mu}}+\delta)k \,\,, \speqnu{9.23'}$$
hence, by the estimates proved in Lemma 9.4 and $(9.23')$, there exist constants $K^{(1)}_{\mu,
\delta}$, $K^{(2)}_{\mu,\delta}>0$ such that, for all $n>N_{\delta}$,
$$\eqalign{&|\Pi^{\pm i}_q\bigl(\gamma_q^{\Cal T}(p)\bigr)|_{-1}\leq K^{(1)}_{\mu,\delta}\,+\,
K^{(2)}_{\mu,\delta}\,|\exp\{({{\lambda_i^{\pm}}\over{\mu}}+(\lambda_i^{\pm}+2)\delta)n\}\,-
\,1|\,\,, \cr
         &|\Pi^{\Cal E}_q\bigl(\gamma_q^{\Cal T}(p)\bigr)|_{-1}\leq K^{(1)}_{\mu,\delta}\,+\,
           K^{(2)}_{\mu,\delta}\exp (2\delta n)\,\,.\cr}
\eqnu $$
On the other hand, by the choice of $n\in {\Bbb N}$ in Lemma 9.4, $$n=\max\{k\in {\Bbb N}\,|\,
{\Cal T}^{(k)}_q(p)\leq {\Cal T}_r\},$$ hence, by $(9.16)$ and Lemma 9.2$'$, there exists a constant 
$K_4>0$ such that
$$K_4\exp\{({1\over{\mu}}-\delta)n\}\leq K_4\exp(|t_n|)\leq {\Cal T}   \,\,. 
\speqnu{9.24'}$$
By $(9.24)$ and $(9.24')$,
$$\eqalign{&\limsup_{{\Cal T}\to +\infty}{{\log |\Pi^{+ i}_q\bigl(\gamma_q^{\Cal T}(p)\bigr)|_{-1}}
\over {\log{\Cal T}}}\,\leq\, {{\lambda_i^++(\lambda_i^++2)\mu\delta}\over {1-\mu\delta}}
\,,\cr
&\limsup_{{\Cal T}\to +\infty}{{\log |\Pi^{-i}_q\bigl(\gamma_q^{\Cal T}(p)\bigr)|_{-1}}\over
{\log {\Cal T}}}\,\leq\, {{\max\{0,\lambda_i^-+(\lambda_i^-+2)\mu\delta\}}\over {1-\mu\delta}}\,,\cr
&\limsup_{{\Cal T}\to +\infty}{{\log |\Pi^{\Cal E}_q\bigl(\gamma_q^{\Cal T}(p)\bigr)|_{-1}} \over
              {\log {\Cal T}}}\,\leq \,{{2\mu\delta}\over {1-\mu\delta}}\,.\cr}\eqnu $$
Since $0<\delta\ll 1/\mu$ and $\lambda_i^+>\lambda'_i>0$, $\lambda^-_i<-\lambda'_i<0$ are arbitrarily
chosen, the upper bound in $(9.20)$ is proved.      
\medbreak
 {\it Lower bound}. We claim the following statement holds. For $\mu^{(1)}_{\kappa}$-almost 
all quadratic differentials $q\in {\Cal M}^{(1)}_{\kappa}$, there exists a full measure set $A_q
\subset M_q$ (consisting of $q$-regular points) with the following property. For each $p\in A_q$, 
there exists a diverging sequence ${\Cal T}_k:={\Cal T}_k(p)$ of positive times such that  
$$\lim_{k\to +\infty}{{\log |\Pi^{+ i}_q\bigl(\gamma_q^{{\Cal T}_k}(p)\bigr)|_{-1}}\over
              {\log {\Cal T}_k}}\,=\, \lambda'_i\,\,. \eqnu $$
We remark that it is sufficient to find a {\it positive }measure set ${\Cal P}^{(1)}_{\kappa}$ in 
every connected component of the moduli space ${\Cal M}^{(1)}_{\kappa}$ with the property that 
$(9.26)$ holds on a set $A_q\subset M_q$ of {\it positive }measure. In fact, since almost all 
quadratic differentials have ergodic horizontal foliation \ref\Mstwo, if $(9.26)$ holds on a positive
measure subset of $M_q$, it holds on a full measure subset. Then, by the ergodicity of the 
Teichm\"uller flow on each connected component of the moduli space (Theorem 1.1), if $(9.26)$ holds 
for almost all $p\in M_q$ on a positive measure subset ${\Cal P}^{(1)}_{\kappa}$ of a connected 
component ${\Cal C}^{(1)}_{\kappa}\subset {\Cal M}^{(1)}_{\kappa}$, then it holds on a full measure 
subset of ${\Cal C}^{(1)}_{\kappa}$.

 Let $q_0 \in {\Cal C}^{(1)}_{\kappa}$ be a quadratic differential with Lagrangian 
horizontal foliation satisfying the properties listed in Lemma 8.4. The flat structure induced
by $q_0$ contains at least $g\geq 2$ embedded flat cylinders $A^{(0)}_j$ with waist curves
$\gamma_j^{(0)}$, $j\in \{1,\ldots ,g\}$, such that the subspace ${\Cal L}_0\subset H_1(M_{q_0},
{\Bbb R})$ generated by $\{\gamma_1^{(0)},\ldots ,\gamma_g^{(0)}\}$ is Lagrangian. Let $h_j>0$ be
the height of the cylinder $A^{(0)}_j$. Then $A^{(0)}_j\equiv \gamma^{(0)}_j\times (-h_j/2,
h_j/2)$. There exists a compact neighbourhood ${\Cal U}^{(1)}_{\kappa}$ of $q_0$ in ${\Cal C}^{(1)}
_{\kappa}$ such that all $q\in {\Cal U}^{(1)}_{\kappa}$ have embedded flat cylinders $A_j(q)$ with 
waist curves isotopic to $\gamma^{(0)}_j$ and, in addition, the following two properties hold:
\medbreak
\item{(1)} the closing ${\widehat \gamma}^{(1)}_q(p)$ of the $q$-horizontal first return trajectory
 of any point $p\in A_j(q)$ is isotopic to a waist curve of $A_j(q)$ in $M_q\setminus \Sigma_q$; 
\smallbreak \item{(2)} there is a strictly positive lower bound for the area of all cylinders $A_j(q)\subset 
M_q$, $q\in {\Cal U}^{(1)}_{\kappa}$. 
\medbreak

Let ${\Cal P}^{(1)}_{\kappa}\subset {\Cal U}^{(1)}_{\kappa}$ be a compact positive measure
set satisfying the conditions listed in Lemma 8.4 and, in addition, the conditions $(1)-(4)$ listed
in Section 9.3. Since, by  (2c)  in Lemma 8.4, the Poincar\'e dual $P({\Cal L}_0)$ of the Lagrangian 
subspace ${\Cal L}_0$, generated by the waist curves $\{\gamma_1^{(0)},\ldots ,\gamma_g^{(0)}\}$, is 
transverse to the stable subspace $E^-(q)$, the condition $(1)$ above implies the following. Let 
$q\in {\Cal P}^{(1)}_{\kappa}$ and $\Pi^{\pm i}_q$, $i\in \{1,\ldots ,s\}$, be one of the projections 
determined by the Oseledec splitting $(9.15')$. There exists a cylinder $A_{j(i)}(q)$ such that 
$$\Pi^i_q\bigl({\widehat\gamma}^{(1)}_q(p)\bigr)\not=0\,\,,\,\,\,\, \hbox{ for all }\,
p\in A_{j(i)}(q)\,\,. \eqnu $$
\noindent By Luzin's theorem, conditions $(3)$, $(4)$ in Section 9.3, Lemmas 9.2, $9.2'$ and $(9.27)$, 
there exist a compact positive measure set ${\Cal P}'_{\kappa}\subset{\Cal P}^{(1)}_{\kappa}$, on 
which the splitting $(9.15')$ depends continuously on the quadratic differential, and a constant 
$K_5>0$ such that, if $q\in {\Cal P}'_{\kappa}$, we have: 
$$ K_5^{-1} \leq |\Pi^i_q\bigl({\widehat\gamma}^{(1)}_q(p)\bigr)|_{-1}\leq K_5\,\,,\,\,\,\,
p\in A_{j(i)}(q)\,\,. \eqnu $$
Let $q\in {\Cal P}'_{\kappa}$ and $(t_n)$ be the sequence of {\it backward }return times of the orbit 
$G_t(q)$ to the positive measure set ${\Cal P}'_{\kappa}$. Let $A_{j(i)}^{(n)}:=A_{j(i)}(q_n)$, 
$q_n:=G_{t_n}(q)$. Let $p\in A_{j(i)}^{(n)}$ and ${\widehat\gamma}^{(n)}_q(p)$ 
be the closing of the first 
return trajectory at $p$ of the horizontal foliation of $G_{t_n}(q)$. Then, by 
Theorem 8.7, for any
$\lambda_i^-<\lambda_i<\lambda^+$, there is a constant $K_6>0$ such that
$$K_6^{-1}\exp(\lambda_i^-|t_n|) \leq |\Pi^i_q\bigl({\widehat\gamma}^{(n)}_q(p)\bigr)|_{-1} \leq K_6
\exp(\lambda_i^+|t_n|) \,\,. \speqnu{9.28'}$$
Let ${\Cal T}^{(n)}_q(p)$ be the return time of the return trajectory $\gamma^{(n)}_q(p)$, $p\in 
A_{j(i)}^{(n)}$. By $(9.17)$ there exists a constant $K_{pr}>0$ such that, for all $q$-regular $p\in
M_q$, 
$$K_{pr}^{-1}\exp{|t_n|}\leq {\Cal T}^{(n)}_q(p) \leq K_{pr}\exp{|t_n|} \,\,. \speqnu{9.28''}$$
All cylinders $A_{j(i)}^{(n)}$ have positive area in $M_q$, uniformly bounded away from zero 
(as a consequence of the corresponding property of the cylinders $A_{j(i)}(q)$ for all $q\in {\Cal U}
^{(1)}_{\kappa}$ and of the invariance of the area form $\omega_q$ under the Teichm\"uller flow). 
Hence the set 
$$A_q:= \bigcap_{N\in {\Bbb N}} \bigcup_{n\geq N} A_{j(i)}^{(n)}  \eqnu $$
such that $p\in A_q$ if and only if $p\in A_{j(i)}^{(n)}$ for infinitely many $n\in {\Bbb N}$, has 
{\it positive measure }in $M_q$. In fact, the $R_q$-area of $A_q$ is greater than any lower bound for
the $R_q$-area of the sets $A_{j(i)}^{(n)}$. This is a simple instance of the Kochen-Stone inequality
[19, \S 6.2, Lemma 4]. By $(9.28')$ and $(9.28'')$, for each $p\in A_q$ there exists therefore 
a diverging sequence $(n_k)$ such that
$$\lim_{k\to +\infty}{{\log |\Pi^{+ i}_q\bigl(\gamma_q^{(n_k)}(p)\bigr)|_{-1}}\over
              {\log {\Cal T}^{(n_k)}}}\,=\, \lambda'_i\,\,. \eqnu $$
Let ${\Cal T}_k(p):={\Cal T}^{(n_k)}$. We have proved $(9.26)$ for a positive measure set ${\Cal P}'
_{\kappa}$ of quadratic differentials in each connected component of the stratum ${\Cal M}^{(1)}_
{\kappa}$ and, given any $q\in {\Cal P}'_{\kappa}$, for a positive measure set $A_q\subset M_q$. By 
the above remark the argument is concluded. \enddemo

\demo{ Remark $9.5'$} In proving Theorem 9.5, we have in fact proved the following statement. 
For almost all $q\in {\Cal M}^{(1)}_{\kappa}$ and any given $\varepsilon>0$ there exists ${\Cal T}
_{\varepsilon}:={\Cal T}_{\varepsilon}(q)>0$ such that, for all $i\in \{1,\ldots ,s\}$ and all $q$-regular 
$p\in M_q$, 
$$\eqalign { & \sup _{{\Cal T}\geq {\Cal T}_{\varepsilon}} {{\log|\Pi^{+i}_q\bigl(\gamma_q^{\Cal 
T}(p)\bigr)|_{-1}} \over {\log {\Cal T}}}\,\leq\, \lambda'_i+\varepsilon\,\,, \cr
& \sup _{{\Cal T}\geq {\Cal T}_{\varepsilon}} {{\log|\Pi^{-i}_q\bigl(\gamma_q^{\Cal 
T}(p)\bigr)|_{-1}} \over {\log {\Cal T}}}\,\leq\, \varepsilon\,\,, \cr
& \sup _{{\Cal T}\geq {\Cal T}_{\varepsilon}} {{\log|\Pi^{\Cal E}_q\bigl(\gamma_q^{\Cal 
T}(p)\bigr)|_{-1}} \over {\log {\Cal T}}}\,\leq\, \varepsilon\,\,. \cr} \eqnu $$
In addition, for each $i\in \{1,\ldots ,s\}$, there exist a cylinder $A^{(i)}_{\varepsilon}(q)\subset M_q$, 
with $\hbox{vol}_q\bigl(A^{(i)}_{\varepsilon}(q)\bigr)$ uniformly bounded away from zero (independently 
of $\varepsilon>0$), and a real number ${\Cal T}^{(i)}_{\varepsilon}:={\Cal T}^{(i)}_{\varepsilon}(q)>0$ such 
that, for all $p\in A^{(i)}_{\varepsilon}(q)$, 
$$ {{\log |\Pi^{+ i}_q\bigl(\gamma_q^{{\Cal T}^{(i)}_{\varepsilon}}(p)\bigr)|_{-1}}
\over {\log {\Cal T}^{(i)}_{\varepsilon}}}\,\geq\, \lambda'_i-\varepsilon  \eqnu$$
and there exists a constant $K^{(i)}_q>0$ such that, for all $p_1$, $p_2\in A^{(i)}_{\varepsilon}(q)$,
$$\left|\Pi^{+ i}_q\bigl(\gamma_q^{{\Cal T}^{(i)}_{\varepsilon}}(p_1)\bigr)\,-\,\Pi^{+ i}_q
\bigl(\gamma_q^{{\Cal T}^{(i)}_{\varepsilon}}(p_2)\bigr)\right|_{-1} \,\leq \, K^{(i)}_q\,\,. \speqnu{9.32'}$$
The bound in $(9.32')$ follows from the remark that all return trajectories ${\widehat\gamma}_q
^{{\Cal T}^{(i)}_{\varepsilon}}(p)$, $p\in A^{(i)}_{\varepsilon}(q)$, are isotopic in $M_q\setminus \Sigma_q$ to
the waist curve of  a flat cylinder (of uniformly bounded area).

The argument also proves that the cohomology class of the closed current ${\widehat\gamma}
^{\Cal T}_q(p)$ stays within a bounded distance of the unstable subspace $E^+_q$ of the 
Kontsevich-Zorich cocycle in the cohomology group $H^1(M,{\Bbb R})$. In fact, by choosing $\delta>0$ 
sufficiently small in $(9.24)$, we obtain the following estimate. There is a constant $K_q>0$ such 
that, for all $i\in \{1,\ldots ,s\}$ and all $q$-regular points $p\in M_q$,
$$\sup _{{\Cal T}>0} |\Pi^{-i}_q\bigl(\hat{\gamma}_q^{\Cal T}(p)\bigr)|_{-1}\leq K_q\,\,. \eqnu $$
\enddemo
 
9.5. {\it The deviation for area-preserving vector fields and open questions}.
 Let $\omega$ be a smooth area form on a compact orientable surface $M$ of genus $g\geq 2$
with at most finite order degeneracies at a finite set $\Sigma:=\{p_1,\ldots ,p_{\sigma}\}\subset M$. Let 
$\imath:=(\imath_1,\ldots ,\imath_{\sigma})$ satisfy the following properties: each $\imath_k$ is a 
negative integer and $\sum \imath_k =2-2g$. Let ${\Cal E}^{\imath}_{\omega}(M,\Sigma)$ the set of 
vector fields $X$ such that the closed $1$-form $\eta_X:=\imath_X\omega$ is smooth on $M$ and the 
orbit foliation ${\Cal F}_X:=\{\eta_X=0\}$ is a measured foliation in the sense of \ref\Th with a 
saddle-type singularity of index $\imath_k$ at each point $p_k\in \Sigma$, $k=1,\ldots ,\sigma$. The 
measure class on ${\Cal E}^{\imath}_{\omega}(M,\Sigma)$ is defined by pull-back of Lebesgue measure 
class under the map ${\Cal E}^{\imath}_{\omega}(M,\Sigma) \to H^1(M,\Sigma;{\Bbb R})$ given by the 
Katok's {\it fundamental class }\ref\Ktone:
$$ X\to [\eta_X]\in  H^1(M,\Sigma;{\Bbb R})\,\,. \eqnu $$    
For almost all $X\in {\Cal E}^{\imath}_{\omega}(M,\Sigma)$, the orbit foliation ${\Cal F}_X$ can be 
realized as the horizontal foliation of an orientable quadratic differential $q\in Q^{(1)}_{\kappa}$
[27, Chap.~II] with $\kappa=-2\imath$. Let ${\Cal C}^{(1)}_X$ be the connected component of the
stratum ${\Cal M}^{(1)}_{\kappa}$ of the moduli space such that $q\in {\Cal C}^{(1)}_X$ and let
$$\{\bigl(\lambda'_i(X),m_i(X)\bigr)\,|\,i=1,\ldots ,s\} \subset {\Bbb R}^+\times {\Bbb N}\setminus\{0\}
\eqnu $$
be the {\it positive }Lyapunov spectrum (with multiplicities) of the Kontsevich-Zorich cocycle on the
connected component ${\Cal C}^{(1)}_X$. The distinct Lyapunov exponents
$\lambda_1'(X)=1>\lambda_2'(X)>\cdots >\lambda_s'(X)>0$ and the corresponding multiplicities
 $m_1(X)=1,\ldots ,m_s(X)$ (such that $\sum m_i(X)=g$)  are defined for almost all $X\in {\Cal
 E}^{\imath}_{\omega}(M,\Sigma)$ and only depend on the connected 
component ${\Cal C}^{(1)}_X$. 
The Kontsevich-Zorich conjecture $(1.2)$ states that $m_i(X)=1$ for all $i\in 
\{1,\ldots ,s\}$ (hence $s=g$). We have proved only the weaker statement that all 
the Lyapunov exponents are non-zero (Theorem 8.5).  

 If $X\in {\Cal E}^{\imath}_{\omega}(M,\Sigma)$ and ${\Cal F}_X={\Cal 
F}_q$, there exists a 
strictly positive smooth function $W=W_X$, with zeroes of finite order at 
$\Sigma$, such that $\omega
=W^{-1}\omega_q$, hence $X=WS$, where $S$ is the normalized horizontal vector 
field for the metric 
$R_q$, defined in Section 2, and $\omega_q$ is the area form of $R_q$. Let
$$H^1_W(M):=\{f\,|\,W^{-1}f \in H^1(M)\} \eqnu $$ 
be the related weighted Sobolev space of once weakly differentiable functions.
\pagebreak

\proclaim {Theorem} For almost all vector fields $X\in {\Cal E}^{\imath}_{\omega}(M,\Sigma)${\rm ,} 
the vector space ${\Cal I}_X^1(M)$ of $X$\/{\rm -}\/invariant distributions which belong to the dual weighted 
Sobolev space $H^{-1}_W(M)$ {\rm (}\/i.e.\ solutions ${\Cal D}^X\in H^{-1}_W(M)$ of the equation $S{\Cal 
 D}=0${\rm )}
has precisely dimension $g\geq 2${\rm .} There exists a splitting
$${\Cal I}_X^1(M)={\Cal I}_X^1(\lambda'_1)\oplus {\Cal I}_X^1(\lambda'_2)\oplus \cdots \oplus
{\Cal I}_X^1(\lambda'_s)\,\,, \eqnu $$
where the dimension of the subspace ${\Cal I}_X^1(\lambda'_i)$ is equal to the multiplicity $m_i(X)$ 
of the $i^{{\rm th}}$ Lyapunov exponent $\lambda'_i(X)>0$ of the Kontsevich\/{\rm -}\/Zorich 
cocycle on the 
connected component ${\Cal C}^{(1)}_X${\rm ,} such that the following holds{\rm .} Let $\Phi_X^{\Cal T}$ denote 
the flow of the vector field $X${\rm .} Let $f\in H^1_W(M)$ be any function such that 
$${\Cal D}^X(f)=0\,\,,\,\,\,\,\hbox{ for all }\,\, {\Cal D}^X\in {\Cal I}_X^1(\lambda'_1)\oplus 
{\Cal I}_X^1(\lambda'_2)\oplus \cdots \oplus {\Cal I}_X^1(\lambda'_i)\,\,,\eqnu $$
then{\rm ,} if $0\leq i<s${\rm ,} for all $p\in M$ with regular forward trajectory{\rm ,}
$$\limsup_{{\Cal T}\to +\infty}{{\log |\int_0^{\Cal T} f(\Phi_X(p,\tau))\,d\tau|}\over
              {\log {\Cal T}}}\,\leq \, \lambda'_{i+1}(X)\,\,, \eqnu $$
and there exists 
${\Cal D}^X_{i+1}\in {\Cal I}_X^1(\lambda'_{i+1})\setminus\{0\}$ such that{\rm ,} if ${\Cal D}^X_{i+1}(f)\not= 0${\rm ,}
then equality holds in $(9.39)$ for almost all $p\in M${\rm .}

 If $i=s$ in $(9.38)${\rm ,} then{\rm ,} for all $p\in M$ with regular forward trajectory{\rm ,}    
$$\limsup_{{\Cal T}\to +\infty}{{\log |\int_0^{\Cal T} f(\Phi_X(p,\tau))\,d\tau|}\over
              {\log {\Cal T}}}\,=\, 0\,\,. \speqnu{9.39'}$$
\endproclaim

\demo{Proof} The $X$-invariant distributions ${\Cal D}^X\in H^{-1}_W(M)$ are related to\break $S$-invariant
distributions ${\Cal D}^S\in H^{-1}_q(M)$ by the following formula:
$${\Cal D}^X(f)={\Cal D}^S(W^{-1}f)\,\,,\,\,\,\,f\in  H^{1}_W(M)\,.\eqnu $$
Let $E^+_i(q)\subset {\Cal B}^1_q(M)$, $i\in\{1,\ldots ,s\}$, be the eigenspace, corresponding to the 
Lyapunov exponent $\lambda'_i(X)>0$, of the Oseledec's splitting $(9.15')$ for the restriction of 
the cocycle $G^c_t$ on closed currents of order $1$ to the 
connected component ${\Cal C}^{(1)}_X$.
Let $C^+_{i,j}\in E^+_i(q)$, $j=1,\ldots ,m_i(X)$, be a basis for the space of basic currents of order 
$1$ of the horizontal foliation ${\Cal F}_q={\Cal F}_X$. By Lemma~6.6, the distributions 
$${\Cal D}_{i,j}^S:=  C^+_{i,j}\wedge \eta_T \in H^{-1}_q(M) \eqnu $$ 
are $S$-invariant. Let ${\Cal D}^X_{i,j}$ be the corresponding $X$-invariant distributions given by 
$(9.40)$ and let ${\Cal I}_X^1(\lambda'_i)$, $i\in \{1,\ldots ,s\}$, be the vector space generated by
$\{{\Cal D}_{i,j}^X\,|\, 1\leq j\leq m_i(X)\}$. Let $\alpha_f:=f\eta_T$, $f\in H^1(M)$. By $(9.23)$, 
$$\gamma_q^{\Cal T}(\alpha_f)=\int_0^{\Cal T} f(\Phi_S(p,\tau))\,d\tau\,\,. \eqnu $$
Since, by $(9.41)$, $C^+_{i,j}(\alpha_f)={\Cal D}^S_{i,j}(f)$, it follows from Theorem 9.5 and $(9.31)$
that $(9.39)$ and $(9.39')$ hold in the case $X=S$ for all functions in $H^1(M)$. The $S$-invariant
distributions ${\Cal D}^S_{i+1} \in {\Cal I}_S^1(\lambda'_{i+1})$, $0\leq i <s$, such that 
${\Cal D}^S_{i+1}(f)\not= 0$ implies equality in $(9.39)$, are constructed as follows. By Theorem 
9.5 and Remark $9.5'$, for each $i\in \{1,\ldots ,s\}$, there exist a sequence of flat cylinders $A^{(i)}_k
\subset M_q$, with area uniformly bounded away from zero, and a sequence ${\Cal T}^{(i)}_k>0$ such 
that, for any sequence of $q$-regular points $p^{(i)}_k\in A^{(i)}_k$,  
$$ \lim_{k\to +\infty}{{\log |\Pi^{+ i}_q\bigl(\gamma_q^{{\Cal T}^{(i)}_k}(p^{(i)}_k)\bigr)|_{-1}}
\over    {\log {\Cal T}^{(i)}_k } }\,=\, \lambda'_i \,\,. \eqnu $$
Since $\Pi^{+i}_q$ projects onto a finite-dimensional subspace, by passing if necessary to a
subsequence, we can define
$${\Cal D}^S_i := \lim_{k\to +\infty}  { {\Pi^{+ i}_q\bigl(\gamma_q^{{\Cal T}^{(i)}_k}(p^{(i)}_k)\bigr)} \over 
{\,\,\,\,|\Pi^{+ i}_q\bigl(\gamma_q^{{\Cal T}^{(i)}_k}(p^{(i)}_k)\bigr)|_{-1}} } \,\in {\Cal I}_S^1(\lambda'_i) 
\,\,. \speqnu{9.43'}$$
By $(9.32')$ the limit in $(9.43')$ does not depend on the choice of the sequence of points 
$p^{(i)}_k \in A^{(i)}_k$. As in $(9.29)$ we define
$$ A^{(i)}_q:= \cap _{N\in {\Bbb N}} \cup _{k\geq N} A^{(i)}_k\,\,.  \eqnu$$
Since the area of the sets $A^{(i)}_k$ is uniformly bounded away from zero, $A^{(i)}_q$ has positive
measure. By Theorem 9.5, $(9.32')$, $(9.42)$, $(9.43)$ and $(9.43')$, for all $q$-regular $p\in 
A^{(i)}_q$, there exists a sequence ${\Cal T}_k:={\Cal T}_k(p)>0$ such that, if $f\in H^1(M)$ and 
${\Cal D}^S(f)=0$ for all ${\Cal D}^S\in {\Cal I}_S^1(\lambda'_1)\oplus \cdots 
\oplus {\Cal I}_S^1
(\lambda'_{i-1})$,
$$ {1\over \left| \Pi^{+1}_q\left( \gamma^{{\Cal T}_k}_q(p)\right)\right|_{-1}} \,\int_0^{{\Cal T}_k}
f(\Phi_S(p,\tau))
\,d\tau \,\,
\to \,\, {\Cal D}^S_i(f)\,\,. \eqnu $$  
Since $A^{(i)}_q$ has positive measure and, for almost all $q\in {\Cal M}^{(1)}_{\kappa}$, the
flow $\Phi_S$ is ergodic (see \ref\KMS, \ref\Mstwo),   $(9.45)$ holds for almost all $p\in M_q$
along a sequence ${\Cal T}_k:={\Cal T}_k(p)$. We have therefore completed the proof of Theorem
9.6 in the case $X=S$.

 Let $f\in H^1_W(M)$, then by definition $W^{-1}f\in H^1(M)$. Let $\gamma^{\Cal T}_X(p)$ be 
the forward $X$-orbit of the (generic) point $p\in M$ up to time ${\Cal T}>0$. We have
$$\eqalign{&\int_0^{\Cal T}f(\Phi_X(p,\tau))\,d\tau= \int_{\gamma^{\Cal T}_X(p)} W^{-1}f\,\eta_T\,\,, \cr
&\int_0^{\Cal T}W(\Phi_X(p,\tau))\,d\tau=L_q\bigl(\gamma^{\Cal T}_X(p)\bigr)\,\,.\cr}\eqnu $$
In addition, by the Birkhoff ergodic theorem, for almost all $p\in M$, 
$$\lim_{{\Cal T}\to+\infty}{1\over{\Cal T}}\int_0^{\Cal T}W(\Phi_X(p,\tau))\,d\tau\,=\,
\int_M W\omega\,=\, \int_M\omega_q\,=1\,\,, \speqnu{9.46'}$$
Hence Theorem 9.6 for a vector field $X=WS$ follows from the case $X=S$ proved above.
\enddemo  

 We conclude the paper by listing some related questions we have left unanswered:
 
\numbereddemo{Question} Does the Kontsevich-Zorich cocycle have a {\it simple }Lyapunov 
spectrum, with respect to the canonical absolutely continuous invariant measure? The affirmative 
answer, conjectured by M. Kontsevich and A. Zorich in the joint paper \ref\KZone, is strongly supported
by the numerical evidence reported in \ref\Zrone, \ref\Zrtwo, \ref\Zrfour, \ref\KZone. Our paper 
proves the conjecture in the case of $g=2$. It is unclear to the author whether similar methods can 
yield a proof in the case of genus $g\geq 3$.
\enddemo

\numbereddemo{Question} Is the Kontsevich-Zorich cocycle non-uniformly 
hyperbolic (or has it simple Lyapunov spectrum) with respect to {\it other 
invariant ergodic measures }of the Teichm\"uller flow? In particular, it would be interesting, as\break 
J. Smillie suggested to the author, to know the answer for {\it periodic Teichm\"uller disks }such as 
those constructed by W. Veech in \ref\Vcthree, which correspond to a class of rational triangular 
billiards. By the methods developed in this paper, a preliminary step would be to answer the question 
whether the support of the invariant measure $\mu$ under consideration is contained in the determinant
locus. In case the answer to the latter question turns out to be negative, it would then follow at 
least that $\lambda^{\mu}_2>0$ (Theorem 3.3). It seems likely that {\it in the case of Veech 
surfaces}, other known methods would be more effective. The methods described in \ref\GM can 
probably be applied to prove that the Lyapunov spectrum is {\it simple}. We owe this remark to A. 
Eskin. 
\enddemo

\numbereddemo{Question} Do there exist basis $\{C^{\pm}_1,\ldots ,C^{\pm}_g\}\subset 
{\Cal B}^1_{\pm q}(M)$ of basic currents of order $1$ for the horizontal, respectively the 
vertical, measured foliation ${\Cal F}_{\pm q}$ at an Oseledec regular point $q\in 
{\Cal M}^{(1)}_{\kappa}$ with the property that in fact $ C^{\pm}_i\in {\Cal H}^{\lambda_i-1}_q(M)$, 
the dual Sobolev space of currents with negative exponent $\lambda_i-1$ ? This question was
formulated independently by M. Kontsevich when we told him about the results of Section 8.1. According 
to Lemma 8.2, this question can be asked at any regular $q\in {\Cal M}^{(1)}_{\kappa}$ at which 
the Lyapunov exponents $\lambda_1,\ldots ,\lambda_g$ are non-zero. A partial form of the question 
can also be asked in the case that, for any $k\in\{1,\ldots ,g\}$, the exponents $\lambda_1,\ldots,
\lambda_k$ are non-zero and $\lambda_{k+1}=0$. In the case $k=1$, \pagebreak we remark that we can choose 
$C^+_1:=\eta_S$ and $C^-_1:=\eta_T$, which are smooth currents. Hence we indeed have that 
$C^{\pm}_1\in {\Cal H}^{\lambda_1-1}_q(M)={\Cal H}^0_q(M)$. An affirmative answer in the case 
$k>1$ could possibly be obtained by strengthening the proof of Lemma 8.2. 
\enddemo   

\numbereddemo{Question} What is the dynamical significance, if any, of basic currents
which are of order $l>1$? We have proved in Section 7.1 that all cohomology classes in a codimension
$1$ subspace of $H^1(M,{\Bbb R})$ can be represented by basic currents of order higher than $1$.
A natural speculation is that they are related to lower order (polynomial, 
sub-polynomial?) deviations of ergodic averages for sufficiently smooth 
functions. 
\enddemo


\references

Ad 
\name{R.\ Adams}, {\it Sobolev Spaces}. {\it Pure and Applied Mathematics} {\bf 65}, Academic Press,  
New York-London, 1975.

BE 
\name{R.\ Barre} and \name{A.\ El Kacimi Alaoui}, Foliations, in 
{\it Handbook of Differential Geometry}\break (F.\ Dillen and L.\ 
Verstraelen, eds.),
North-Holland, Amsterdam, to appear.

BG
\name{C.\ Berenstein} and \name{R.\ Gay}, {\it Complex Variables. An 
Introduction}, {\it Grad.\ Texts in Math.\/}  {\bf 125}, 
Springer-Verlag, New York, 1991.

BGM 
\name{M.\ Berger}, \name{P.\ Gauduchon}, and \name{E.\ Mazet}, 
{\it Le Spectre d'une
Vari\'et\'e Riemannienne}, {\it Lecture Notes in Math\/}.\  {\bf 194}, 
Springer-Verlag, New York, 1971.

Beone
\name{L.\ Bers}, Spaces of degenerating Riemann surfaces, in 
{\it Discontinuous Groups and Riemann Surfaces}   
(Univ.\ 
Maryland, College Park, Md., 1973), {\it Ann.\  of Math.\  Studies\/} 
{\bf 79}, 43--55, Princeton Univ.\  Press, Princeton, NJ, 1974.

Betwo 
\bibline, Finite-dimensional Teichm\"uller spaces and generalizations,
{\it Bull.\ A.\ M.\ S\/}.\  {\bf 5} (1981),  131--172.

Bu
\name{M.\ Burger}, Horocycle flow on geometrically finite surfaces, 
{\it Duke Math.\ J\/}.\
{\bf 61} (1990),  779--803.

Ch 
\name{I.\ Chavel}, {\it Riemannian Geometry - A Modern
Introduction}, {\it Cambridge Tracts in Math\/}.\ {\bf 108}, 
Cambridge Univ.\ Press, Cambridge, 
1993.

Cg 
\name{J.\ Cheeger}, A lower bound for the smallest
eigenvalue of the Laplacian, in {\it Problems in Analysis\/} 
(Papers dedicated to  Salomon Bochner, 1969), 195--199, 
 Princeton Univ.\  Press, Princeton, NJ, 1970.

Dk
\name{M.\ Denker}, The central limit theorem for dynamical systems,
in {\it Dynamical Systems and Ergodic Theory } (Warsaw, 1986), 
{\it Banach Center Publ\/}.\  {\bf 23},  33--62, PWN, Warsaw, 1989.

dR 
\name{G.\ de Rham}, {\it Variet\'es Differentiables}.  {\it Formes{\rm ,} courants{\rm ,}
 formes harmoniques} (French).  Troisime \'edition revue et augment\'ee.  {\it Publications de l\/{\rm '}\/Institut
de Mathematique de l\/{\rm '}\/Universit{\rm \'{\it e}} de Nanacago}, III.  {\it Actualit\'es Scientifiques et Industrielles},
{\bf 1222b}, Hermann, Paris, 1973.

EM 
\name{A.\ Eskin} and \name{H.\ Masur}, Asymptotic formulas on flat surfaces, {\it Ergodic Theory Dynam.\
Systems} {\bf 21} (2001), 443--478.

FK 
\name{H.\ M.\ Farkas} and \name{I.\ Kra}, {\it Riemann Surfaces} 
(second edition), {\it Grad.\ Texts in Math\/}.\ {\bf 71}, 
Springer-Verlag, New York, 1992.

FLP
\name{A.\ Fathi}, \name{F.\ Laudenbach}, and \name{V.\ Po\'enaru}, {\it Travaux de
Thurston sur les surfaces}  (French).  S\'eminaire Orsay.  Preprint of Travaux de Thurston sure les
surfaces,   Soc.\ Math.\ France, Paris, 1979. {\it Ast{\rm\'{\it e}}risque\/} {\bf 66--67} (1991). 
Soci\'et\'e   Math\'ematique  de France, Paris, 1991.

Fa 
\name{J.\ D.\ Fay}, {\it Theta Functions on Riemann Surfaces}, {\it Lecture
Notes in Math\/}.\ {\bf 352}, Springer-Verlag, Berlin-New York, 1973.

Fd
\name{B.\ R.\ Fayad}, {\it Reparam\'etrage de flots irrationnels sur le 
tore} (English),  Th\`eses de Doctorat,  \'Ecole Polytechnique, Paris, Juin 2000.

Fone 
\name{G.\ Forni}, The cohomological equation for
area-preserving flows on compact surfaces, {\it Electron.\  Res.\ 
Announc.\  of the A.\ M.\ S.} {\bf 1} (1995),  114--123.

Ftwo 
\bibline, Solutions of the cohomological equation
for area-preserving flows on compact surfaces of higher genus,
 {\it Ann.\ of Math\/}.\ 
{\bf 146} (1997),  295--344.

FG 
\name{B.\ Fristedt} and \name{L.\ Gray}, {\it A Modern Approach to Probability 
Theory}, {\it Probab.\ Appl\/}., Birkh\"auser, Boston, 1997.

Ga 
\name{F.\ Gardiner}, {\it Teichm\"uller Theory and Quadratic
Differentials}, John Wiley \& Sons, Inc., New York, 1987.

GM 
\name{I.\ Ya.\ Gol'dshe\u{\i}d} and \name{G.\ A.\ Margulis}, Lyapunov 
exponents of a product
of random matrices (Russian), {\it  Uspekhi Mat.\ Nauk} {\bf 44} (1989), no.\ 5(269), 13--60 (English
translation: {\it Russian Math.\ Surveys} {\bf 44} (1989), no.\ 5,  11--71).

GH
\name{W. H. Gottschalk} and \name{G. A. Hedlund}, Topological dynamics, {\it A. M. S.
Colloq.\ Publ}.\  {\bf 36}, A. M. S.\   
 Providence RI, 1955.

Gr 
\name{P.\ A.\ Griffiths}, Periods of integrals on algebraic manifolds: summary
of main results and discussion of open problems, {\it Bull.\ A.M.S}.\ 
{\bf 76} (1970), 228--296.

Hl 
\name{S.\ Helgason}, {\it Groups and Geometric Analysis.  Integral Geometry, 
Invariant Differential Operators, and Spherical Functions}, {\it Pure and Appl.\ Math\/}.\ {\bf 113},
Academic Press, 
Inc., Orlando, Florida, 1984.

Hrone
\name{M.-R.\ Herman}, Sur la conjugaison 
diff\'erentiable des diff\'eomorphismes 
du cercle \`a des rotations (French), {\it IHES Publ.\ Math\/}.\  {\bf 49} (1979),  5--233.

Hrtwo 
\bibline, Une m\'ethode pour minorer les exposants de Lyapunov et quelques 
exemples montrant le caract\`ere local d'un th\'eor\`eme d'Arnol'd et de Moser 
sur le tore de dimension $2$ (French), {\it  Comment.\ Math.\ Helv\/}.\  {\bf 58} (1983), 453--502.

HM 
\name{J.\ Hubbard} and \name{H.\ Masur}, Quadratic differentials
and foliations, {\it Acta Math\/}.\  {\bf 142} (1979),  221--274.

IT 
\name{Y.\ Imayoshi} and \name{M.\ Taniguchi}, {\it An Introduction to Teichm\"uller
Spaces}, Springer-Verlag, Tokyo, 1992.

Ktone 
\name{A.\ Katok}, Invariant measures of flows on orientable surfaces (Russian), 
{\it Dokl.\ Akad.\ Nauk SSSR\/} {\bf 211} (1973), 775--778 (English translation: 
{\it Sov.\ Math.\ Dokl\/}.\  {\bf 14} (1973), 1104--1108).

Kttwo
\bibline, Interval exchange transformations 
and some special flows are not mixing, 
{\it Israel J.\ Math\/}.\  {\bf 35} (1980),  301--310.

Ktthree
\bibline, Infinitesimal Lyapunov 
functions, invariant cone families and stochastic properties of smooth 
dynamical systems (with the collaboration of K.\ Burns), 
{\it Ergod.\ Theory  Dynam.\ Systems} {\bf 14} (1994), 757--785.

KHone 
\name{A.\ Katok} and \name{B.\ Hasselblatt}, {\it Introduction
to the Modern Theory of Dynamical Systems}, Cambridge Univ.\
Press, Cambridge, 1995.

KHtwo
\bibline, {\it Handbook of Dynamical Systems, Vol.\ 1, Survey I, Principal 
Structures}, to appear.

KR
\name{M.\ S.\ Keane} and \name{G.\ Rauzy}, Stricte ergodicit\'e des 
\'echanges d'intervalles (French), {\it Math.\ Z\/}.\  {\bf 174} (1980),  203--212.

KMS 
\name{S.\ Kerckhoff}, \name{H.\ Masur}, and \name{J.\ Smillie}, Ergodicity of
Billiard Flows and Quadratic Differentials, {\it Ann.\ of Math\/}.\ 
{\bf 124} (1986), 293--311.

Kh 
\name{A.\ Ya.\ Khinchin}, {\it Continued Fractions}, The Univ.\  of Chicago 
Press, Chicago, Ill.-London, 1964.

Kc
\name{A.\ V.\ Ko\v cergin}, Mixing in special flows over a rearrangement of segments
and in smooth flows on surfaces (Russian), {\it Mat. Sb.\ {\rm (}\/N.S}.)  {\bf 96}({\bf 138})  (1975), 
 471--502, 504 (English translation: {\it Math.\ USSR-Sb\/}.\  {\bf 25} (1975), 
441--469).

KZone 
\name{M.\ Kontsevich}, Lyapunov exponents and
Hodge theory, in  {\it The Mathematical Beauty of Physics} (Saclay, 1996),
{\it Adv.\ Ser.\ Math.\ Phys\/}.\  {\bf 24}, 318--332, World Sci.\ Publ.,
River Edge, NJ, 1997.

KZtwo
\name{M. Kontsevich} and \name{A. Zorich}, Connected components of the moduli space of abelian 
differentials with prescribed singularities, preprint.

Ln 
\name{S.\ Lang}, {\it Differentiable Manifolds}, Springer-Verlag, New York, 1985.

Mn 
\name{R.\ Ma\~n\'e}, {\it Ergodic Theory and Differentiable
Dynamics}, {\it Ergeb.\ Math.\ Grenzgeb.\/} {\bf 8},\break Springer-Verlag, Berlin,
 1987.

Msone 
\name{H.\ Masur}, Extension of the Weyl-Petersson
metric to the boundary of the Teichm\"uller space, {\it  Duke Math.\ J\/}.\ 
{\bf 43} (1976),  623--635.

Mstwo 
\bibline, Interval exchange transformations
and measured foliations, {\it Ann.\ of Math\/}.\  {\bf 115}  (1982),
169--200.

Msthree 
\bibline, Logarithmic law for geodesics in moduli space, in {\it Mapping 
Class Groups and Moduli Spaces of Riemann Surfaces}  
(G\"ottingen, 1991/Seattle, WA, 1991),   {\it Contemp.\ 
Math\/}.\ {\bf 150}, 229--245, A.\ M.\ S., Providence, RI, 1993.

MS 
\name{H.\ Masur} and \name{J.\ Smillie}, Hausdorff dimension of 
sets of nonergodic measured foliations, {\it Ann.\ of Math\/}.\  {\bf 134}
(1991), 455--543.

Ng 
\name{S.\ Nag}, {\it The Complex Analytic Theory of Teichm\"uller Spaces},
John Wiley \& Sons, Inc., New York, 1988.

NZ 
\name{I.\ Nikolaev} and \name{E.\ Zhuzhoma}, {\it Flows on $2$-dimensional
Manifolds.  An Overview}, {\it Lectures Notes in Math\/}.\ {\bf 1705}, Springer-Verlag, Berlin,  
1999.

No
\name{S.\ P.\ Novikov}, The Hamiltonian formalism and a multivalued analogue
of Morse theory (Russian), {\it Uspekhi Mat.\ Nauk\/} {\bf 37} (1982),   3--49 (English 
translation: {\it Russian Math.\ Surveys\/} {\bf 37} (1982), 1--56).

Os 
\name{V.\ I.\ Oseledec},  A multiplicative ergodic theorem. Characteristic
Ljapunov, exponents of dynamical systems (Russian), {\it Trudy Moskov Mat.\
Ob\v{s}\v{c}\/}.\  {\bf 19} (1968),   179--210 (English translation: {\it Trans.\ Moscow Math.\ 
Soc\/}.\  {\bf 19} (1968),  197--231).

Pe
\name{K.\ Petersen}, {\it Ergodic Theory}, {\it Cambridge Studies in
Adv.\ Math\/}.\ {\bf 2}, Cambridge Univ.\
Press, Cambridge, 1989.

Rtone
\name{M.\ Ratner}, Rigidity of time changes for horocycle flows, {\it Acta Math\/}.\  
{\bf 156} (1986),  1--32.

Rttwo
\bibline, The rate of mixing for geodesic and horocycle flows, 
{\it Ergod.\ Theory Dynamic.\ Systems\/}  
{\bf 7} (1987),   267--288.

Reone
\name{B.\ L.\ Reinhart}, Harmonic integrals on foliated manifolds,
{\it Amer.\ J.\ Math}.\ {\bf 81} (1959), 529--536.

Retwo 
\bibline, {\it Differential Geometry of Foliations},
Springer-Verlag, New York, 1983.

Sl 
\name{K.-U.\ Schauml\"offel}, Multiplicative ergodic
theorems in infinite dimensions, in {\it Lyapunov Exponents} 
(Oberwolfach, 1990), {\it Lecture
Notes in Math}.\ {\bf 1486}, 187--195, Springer-Verlag, Berlin, 1991.

Sc 
\name{L.\ Schwartz}, {\it Th\'eorie des Distributions} (French),
 {\it Publ.\ de l\/{\rm '}\/Institut de Math.\ de\break
l\/{\rm '}\/Universit{\rm \'{\it e}} de Strasbourg}, No. IX-X,  Nouvelle \'edition, enti\'erement
corrig\'ee, refondue et augment\'ee,  Hermann, Paris, 1966.

Sn
\name{S.\ Schwartzman}, Asymptotic cycles, {\it  Ann.\ of Math\/}.\  {\bf 66} (1957), 
270--284.

St 
\name{K.\ Strebel}, {\it Quadratic Differentials}, {\it Ergeb.\ Math.\ Grenzgeb.}\ {\bf 5},  
Springer-Verlag, Berlin,  1984.

Su 
\name{D.\ Sullivan}, The Dirichlet problem at infinity for a negatively
curved manifold, {\it J.\ Differential Geom\/}.\  {\bf 18} (1983), 723--732.

Taone
\name{M.\ Taniguchi}, A note on the second variational formulas of functionals on Riemann surfaces,
{\it  Kodai Math.\ J\/}.\ 
{\bf 12} (1989), 283--295.

Tatwo
\name{M.\ Taniguchi}, The behavior of the extremal length function on
arbitrary Riemann surface, in {\it Prospects in Complex
Geometry} (Katata and Kyoto, 1989), {\it
Lecture Notes in Math\/}.\  {\bf 1468}, 160--169,
Springer-Verlag, Berlin, 1991.

Th 
\name{W.\ Thurston}, On the geometry and dynamics of
diffeomorphisms of surfaces, {\it Bull.\ A.\ M.\ S\/}.\ 
{\bf 19} (1988), 417--431.

Td 
\name{P.\ Tondeur}, {\it Geometry of Foliations}, 
{\it Monographs in Math\/}.\ {\bf 90}, Birkh\"auser Verlag,
Basel, 1997.

Tr
\name{A.\ Tromba}, {\it Teichm\"uller Theory in Riemannian Geometry}, 
{\it Lectures Math.\ ETH Z\"urich\/}, Birkha\"user Verlag, Basel, 1992.

Vcone 
\name{W.\ A.\ Veech}, Gauss measures for transformations on the
space of interval exchange maps, {\it Ann.\ of Math\/}.\  {\bf 115} (1982), 201--242.

Vctwo
\bibline, The Teichm\"uller geodesic flow, {\it Ann.\ of Math\/}.\   {\bf 124} (1986),  
 441--530.

Vcthree 
\bibline, Teichm\"uller curves in moduli space, Eisenstein series
and an application to triangular billiards, {\it Invent.\ Math\/}.\ {\bf 97} (1989), 
 553--583.

Vcfour 
\bibline, Moduli spaces of quadratic differentials,
{\it J.\ Analyse Math\/}.\  {\bf 55} (1990),  117--171.

Vcfive 
\bibline, Siegel measures, {\it Ann.\ of Math\/}.\  {\bf 148} (1998),
 895-944.

We 
\name{R.\ O.\ Wells}, {\it Differential Analysis on Complex
Manifolds} (second edition), {\it Grad.\ Texts in Math\/}.\ {\bf 65}, 
Springer-Verlag, New York-Berlin, 1980.

YaA 
\name{A.\ Yamada}, Precise variational formulas for abelian differentials,
{\it Kodai Math.\ J\/}.\  {\bf 3} (1980),  114--143.

YaS
\name{S.\ Yamada}, Weil-Petersson convexity of the energy functional on classical and universal
 Teichm\"uller
spaces,
{\it J.\ Differential Geom\/}.\  {\bf 51} (1999),  35--96.

Yu 
\name{S.\ T.\ Yau}, Isoperimetric constants and the first
eigenvalue of a compact Riemannian manifold, 
{\it Ann.\ Sci.\ \'Ecole Norm.\ Sup}.\
{\bf 8} (1975), 487--507.

Ygone     
\name{L.-S.\ Young}, Entropy of continuous flows on compact $2$-manifolds, {\it Topology\/} 
{\bf 16} (1977),  469--471.

Ygtwo
\bibline, Statistical properties of dynamical systems
with some hyperbolicity, {\it Ann.\ of Math\/}.\ {\bf 147} (1998),   585--650.

Ygthree
\bibline, Developments in chaotic dynamics, {\it Notices  A.\ M.\ S\/}.\ 
 {\bf 45}, no.\ 10 (1998),   1318--1328.

ZK
\name{A.\ Zemljakov} and \name{A.\ Katok}, Topological transitivity of
billiards in polygons (Russian), {\it Mat.\ Zametki\/} {\bf 18} (1975),   291--300 
(English translation: {\it Math.\ Notes\/} {\bf 18} (1975), 760--764).

Zrone 
\name{A.\ Zorich}, Asymptotic flag of an orientable measured foliation
on a surface, in  {\it Geometric Study of Foliations} (Tokyo, 1993), 479--498, 
World Sci.\  Publ., River Edge, NJ,  	
1994.

Zrtwo 
\bibline, Finite Gauss measure on the space of interval exchange 
transformations. Lyapunov exponents, {\it Ann.\ Inst.\ Fourier {\rm (}\/Grenoble}\/) {\bf 46} 
(1996),  325--370.

Zrthree 
\bibline, Deviation for interval exchange transformations,
{\it Ergod.\ Theory Dynam.\ Systems\/} {\bf 17} (1997), 1477--1499.

Zrfour 
\name{A.\ Zorich}, On hyperplane sections of periodic surfaces,
in {\it Solitons{\rm ,} Geometry and Topology\/{\rm :}\/ on the Crossroad}, 
 {\it A.\ M.\ S.\ Transl\/}.\  {\bf 179}, 173--189,  A.\ M.\ S.\ 
Providence, RI, 1997.

Zrfive
\bibline, How do the leaves of a closed $1$-form wind around a surface?, in
{\it Pseudoperiodic Topology}, 
{\it  A.\ M.\ S.\ Transl\/}.\ {\bf 197}, 135--178,  A.\ M.\ S., Providence, RI, 1999.

\endreferences

\bye